\documentclass[reqno]{amsart}
\usepackage{amsmath,amssymb}
\usepackage{amsmath,amssymb,cite}
\usepackage{graphics,epsfig,color}
\usepackage[mathscr]{euscript}
\usepackage{xcolor}
\usepackage{a4wide}
\usepackage{tikz}
\usepackage{stmaryrd}

\theoremstyle{plain}
\newtheorem{theorem}{Theorem}[section]
\newtheorem{proposition}[theorem]{Proposition}
\newtheorem{lemma}[theorem]{Lemma}
\newtheorem{corollary}[theorem]{Corollary}

\newtheorem{conclusion}[theorem]{Conclusion}
\theoremstyle{definition}

\theoremstyle{remark}
\newtheorem{remark}[theorem]{Remark}
\newtheorem{example}[theorem]{Example}

\makeatletter
\@addtoreset{equation}{section}
\makeatother

\DeclareMathOperator*{\esssup}{ess\,sup}

\newcommand{\mc}[1]{{\mathcal #1}}
\newcommand{\mb}[1]{{\mathbf #1}}
\newcommand{\mf}[1]{{\mathfrak #1}}
\newcommand{\bs}[1]{{\boldsymbol #1}}
\newcommand{\bb}[1]{{\mathbb #1}}
\newcommand{\ms}[1]{{\mathscr #1}}
\newcommand{\mss}[1]{{\textsf #1}}
\newcommand{\mtt}[1]{{\mathtt #1}}

\newcounter{as}[section]

\newcommand{\<}{\langle}
\renewcommand{\>}{\rangle}

\newcommand{\cb}[1]{{\color{blue} #1}}

\title{Non-equilibrium fluctuations of gradient exclusion processes in
dimension $d\le 3$}

\author{ Claudio Landim and Sunder Sethuraman}

\address{Claudio Landim
  \hfill\break\indent IMPA \hfill\break\indent Estrada Dona Castorina
  110, \hfill\break\indent
J. Botanico, 22460 Rio de Janeiro, Brazil\hfill\break\indent
{\normalfont and} \hfill\break\indent
Univ. Rouen Normandie, \hfill\break\indent
CNRS Normandie Univ,  LMRS UMR 6085, \hfill\break\indent
F-76000 Rouen, France.} 
\email{landim@impa.br}

\address{Sunder Sethuraman \hfill\break\indent Department of
Mathematics \hfill\break\indent University of Arizona
\hfill\break\indent 621 N. Santa Rita Ave, \hfill\break\indent Tucson,
AZ, 85721 USA.}  \email{sethuram@arizona.edu}

\begin{document}

\begin{abstract}
We consider a gradient, speed-change, symmetric exclusion process on
the discrete torus $\bb T^d_n$, $d\le 3$, whose empirical measure
evolves, on the diffusive scale, according to a non-linear parabolic
equation.  We prove that, starting from a sequence of initial states
whose relative entropy with respect to the inhomogeneous product
measure associated with the initial profile is $o(n^{d/2})$, and whose
density fluctuation field converges, the density fluctuation field
converges, in $L^1(0,T; \mc H_{-\mtt r})$, to the solution of a linear
stochastic partial differential equation with time-dependent
coefficients driven by a conservative space-time white noise.  
\end{abstract}

\noindent


\maketitle
\thispagestyle{empty}

\section{Introduction}
\label{s:0}

Consider a mass conservative interacting particle system evolving on
the $d$-dimensional discrete torus $\bb T^d_n$.  Let
$\pi^n(\eta^n(t)) = n^{-d}\sum_{x\in \bb
T^d_n}\eta_x^n(t)\delta_{x/n}$ be the empirical measure in diffusive
scale, where $\eta^n(t)=\{\eta^n_x(t)\}$ denotes the configuration of
(unlabeled) particles at time $n^2t$. When $\pi^n$ converges to the
solution $u$ of a partial differential equation, such a limit
constitutes the `hydrodynamic' limit for the model \cite{kl}.

The next order term in this description is the behavior of the 
density fluctuation field,
\begin{equation*}
X^n_t(F) \;=\; \frac{1}{n^{d/2}} \sum_{x\in \bb T^d_n} F(x/n)\,
\big[\, \eta_x^n(t) \,-\, u(t,x/n)\, \big] \;.
\end{equation*}
One expects $X^n_t$ to converge to the solution of a linear stochastic
partial differential equation, obtained by linearizing the non-linear
parabolic hydrodynamic equation around $u$ and adding a conservative
noise whose intensity is prescribed by a certain `compressibility' of
the model.

When the process starts from an invariant state this program has a
long history.  The equilibrium fluctuations of the density field were
first obtained by Brox and Rost \cite{BR84} for zero-range and
exclusion dynamics, and later extended to gradient and non-gradient
reversible systems by Chang \cite{c94} and Lu \cite{Lu1994};
see Chapter 11 of \cite{kl} for an account of this
theory.

Out of equilibrium, however, the picture is far less complete, and the
problem has remained one of the main open questions in the theory of
hydrodynamic limits.  Of the few results available, Chang and Yau
\cite{cy92}, DeMasi, Presutti, and Scacciatelli \cite{DeMasi},
Ferrari, Presutti, and Vares \cite{Ferrari}, Ravishankar
\cite{Ravishankar} are restricted to one dimension or rely on special
features of the dynamics, such as `duality'.  More recently, Jara and
Menezes \cite{jm1}, \cite{jm2} have derived the non-equilibrium
fluctuation limit in exclusion-type systems in $d\leq 3$.  Also, with
respect to asymmetric exclusion processes in $d=1$, fluctuation limits
for `height' functions have been shown in `KPZ' scales; see Matetski,
Quastel, and Remenik \cite{MQR} and references therein.  These last
works also take advantage of special properties of the exclusion
system or integrability of certain quantities.

The reason for the difficulty in the general setting is as follows.
Writing the semi-martingale decomposition of $X^n_t(F)$, the drift is
not a function of the density field itself but a space-time average of
local functions.  Turning this average into a function of $X^n_t$ is
the content of the `Boltzmann--Gibbs principle', introduced in
\cite{BR84}.  Namely, a local function $\tau_x g$ may be replaced, at
first order and after averaging in space and time, by its projection
on the conserved quantity, here the number of particles, and therefore
written in terms of the fluctuation field.  In equilibrium, this
replacement follows from the classical Kipnis--Varadhan bound on the
variance of an additive functional of a reversible Markov process, see
\cite[Proposition A1.6.1]{kl}.  Out of equilibrium, however, this tool
is not available, and one has to control otherwise $\mu_n(t)$, the
state of the process at macroscopic time $t$.

As in the Guo, Papanicolaou, and Varadhan \cite{gpv} `entropy method' for
the hydrodynamic limit, Chang and Yau \cite{cy92}, for a $d=1$
Ginzburg-Landau model, impose the bound
$H(\mu^n\, |\, \nu^n_{\bar u}) = O(n)$ of the relative entropy of
the initial state $\mu^n=\mu^n(0)$ with respect to an invariant
reference measure, and derive the non-equilibrium fluctuation limit.
Such a condition, where control of $n^dX^n_t(F)$ is required in the
argument, seems too strong in higher dimensions.  However, the state
$\mu^n(t)$ may be compared to a more natural reference measure, the
`local equilibrium' given by the inhomogeneous product measure
$\nu^n_{u(t,\cdot)}$ associated with the solution of the hydrodynamic
equation.  Though such a measure is not invariant, the relative
entropy $H_n(\mu^n(t) \,|\, \nu^n_{u(t,\cdot)})$ may be estimated.

Yau's `relative entropy method' \cite{Yau1991} provides an $o(n^d)$
bound for this entropy production, sufficient to derive the
hydrodynamic limit.  But, such a bound is not suitable for the
homogenizations required in the fluctuation limit.  What is needed is
that $H_n(\mu^n(t) \,|\, \nu^n_{u(t,\cdot)})=o(n^{d/2})$, since then
the entropy inequality would resolve fluctuations of size $n^{d/2}$.
An improvement, due to Jara and Menezes \cite{jm1, jm2}, bounds the
relative entropy of order $g_d(n)n^{d-2}$ where $g_d(n)=n$, $\log(n)$,
and $1$ in $d=1$, $d=2$, and $d\geq 3$ respectively, via a `path
decomposition' analysis amenable in certain exclusion-type systems,
but less so more generally, say for systems allowing multiple
occupancy on sites.  However, such bounds permitted the derivation of
the non-equilibrium fluctuations of a reaction-diffusion model and a
form of weakly-asymmetric exclusion in $d\leq 3$.  These bounds were
subsequently used in Jara and Landim \cite{jl23} to prove that the
density field of a stirring dynamics perturbed by a voter model
converges to the stochastic heat equation, and also in Funaki, Landim,
and Sethuraman \cite{FLS} to identify certain non-equilibrium
fluctuations for Glauber+Kawasaki dynamics with respect to continuum
mean-curvature flow derivations, both in dimensions $d\leq 3$.  In
these works, the perturbation is small, in the sense that the
exclusion part of the generator is accelerated faster than the other
one, so that the reference measures may be taken corresponding to the
exclusion dynamics.

The purpose of this article is to derive the non-equilibrium
fluctuation limit via a robust `multi-scale' homogenization argument
suitable for a general class of gradient models in $d\leq 3$.  The
strength of the approach is that it relies only on large deviations
estimates, a log-Sobolev inequality and equivalence of ensembles
expansions, and not on more special features of the process.  More
specifically, we focus in this article on a `speed-changed' symmetric
exclusion system, as a representative in the class, for which a
log-Sobolev inequality is already known, and for which useful forms of
equivalence of ensembles estimates can be shown, as we derive later in
the paper.

The use of multi-scale schemes in problems of hydrodynamics go back at
least to the seminal paper of Guo, Papanicolaou, and Varadhan
\cite{gpv}, and can be seen in other works in different contexts such
as Landim, Olla, and Yau \cite{LOY}, Grunewald, Otto, Villani, and
Westdickenberg \cite{GOVW}, Dizdar, Menz, Otto, and Wu \cite{DMOW},
among others.  However, the formulation of multi-scale methods that we
provide here to address systematically the question of non-equilibrium
fluctuations is novel.

Fix strictly positive local functions $c_{0,e_j}$, $1\le j\le d$,
which do not depend on the occupation variables $\eta_0$,
$\eta_{e_j}$, and satisfy the gradient condition
$c_{0,e_j}(\eta)\, [\eta_{e_j} - \eta_0] = \tau_{e_j} h_j - h_j$ for
some local functions $h_j$.  In diffusive scale, a particle at $x$
jumps to $x\pm e_j$ at rate $n^2 c_{x,x\pm e_j}(\eta)$.  The Bernoulli
product measures are reversible for this dynamics, and its empirical
measure converges to the solution $u$ of the non-linear diffusion
equation
$\partial_t u = \sum_j \partial_{x_j} ( D_j(u)\, \partial_{x_j} u)$,
where $D_j(\alpha)$ is the expectation of $c_{0,e_j}$ with respect to
the Bernoulli measure of density $\alpha$.

Our main result, Theorem \ref{t1}, asserts that in dimension $d\le 3$
the density fluctuation field converges to the solution of the
stochastic partial differential equation
\begin{equation}
\label{eq:intro1}
\partial_t Y_t \;=\; \mc A^*_t \, Y_t \;+\; \sum_{j=1}^d
\partial_{x_j} \Big(\, \sqrt{2\, \chi(u(t))\, D_j(u(t))}\;
\xi^j_t \,\Big)\;, \qquad Y_0 \,=\, X_0 \;,
\end{equation}
where $\mc A_t$ is the elliptic operator obtained by linearizing the
hydrodynamic equation at $u(t,\cdot)$,
$\chi(\alpha)D_j(\alpha)=\alpha(1-\alpha)D_j(\alpha)$ is the static
compressibility of the speed-changed exclusion dynamics, and $\xi$ is
a $d$-dimensional space-time white noise.  Both the drift and the
intensity of the noise are time-dependent, involving the solution of
the hydrodynamic equation.  In equilibrium, the equation
\eqref{eq:intro1} reduces to a linearization around a constant
solution $u(t)\equiv \bar u$, with constant drift and noise intensity
depending on $\bar u$.

We discuss now three aspects of the statement of Theorem \ref{t1}.

\smallskip\noindent{\it The hypotheses on the initial state.}  We
require only that the relative entropy of $\mu_n$ with respect to
$\nu^n_{u_0(\cdot)}$ be $o(n^{d/2})$, and that the initial fluctuation
field $X^n_0$ converges weakly in $\mc H_{-\mtt r}$ (a Sobolev space of
negative norm $-\mtt r$) to some random element $X_0$.  No further
property of $\mu_n$ is assumed.  in particular $\mu_n$ is not required
to be a product measure, and the limit $X_0$ is not required to be
Gaussian, nor centered.  As
$H_n(\nu^n_{\rho(\cdot)} \,|\, \nu^n_{u_0(\cdot)})$ is of order
$n^d \Vert \rho - u_0\Vert^2_{L^2}$, the entropy condition allows the
initial profile to differ from $u_0$ by $o(n^{-d/4})$, a perturbation
much larger than the scale $n^{-d/2}$ of the fluctuations themselves.
The limit is the mild solution $Y_t = P^*_{0,t} X_0 + \mss M_t$, where
$(P_{s,t})$ is the semigroup generated by $\mc A_t$ and $\mss M$ is a
Gaussian field independent of $X_0$.  Here, the initial data is
transported deterministically, the term $\mss M$ reflects the noise
from the dynamics, and covariances of the sum
$P^*_{0,t}X_0 + \mss M_t$ correspond to the linear SPDE
\eqref{eq:intro1}.

\smallskip\noindent{\it The topology.}  We prove the convergence in
$L^1(0,T;\mc H_{-\mtt r})$ and not, as is customary, in
$D([0,T],\mc H_{-\mtt r})$.  This is a consequence of the estimates at
our disposal, and not a technical convenience.  The bound we obtain
for the space-time average of local functions, Corollary \ref{s94},
contains a term which does not vanish as the interval of time
integration tends to zero.  The estimate therefore yields no modulus
of continuity in time, and the compactness criteria in the Skorohod
space do not apply.  In $L^1(0,T;\mc H_{-\mtt r})$, on the other hand,
compactness only requires a bound on the mean $L^1$-modulus of
continuity, by Simon's criterion \cite{s86}, which is what the
Boltzmann--Gibbs principle (Theorem \ref{t02}) provides.  A similar
difficulty was circumvented in \cite{jl23}, \cite{FLS} by considering
tightness of the time-integral of the fields instead of the fields
themselves.  We refer to Section \ref{sec-c3}, where the properties of
$L^1(0,T;\mc H_{-\mtt r})$ needed here are collected, including the
fact that a probability measure on this space is determined by the
laws of the integrals of the field over time-intervals.

\smallskip\noindent{\it The restriction to $d\le 3$.}  On the one
hand, fundamentally, the bound on the entropy production which we are
able to prove is not strong enough in dimension $d\ge 4$ to yield the
Boltzmann--Gibbs principle.  On the other hand, the deterministic
remainder produced by the replacement of the discrete Laplacian by the
continuous one, and by the replacement of
$E_{\nu^n_{u(s,\cdot)}}[\tau_x h_j]$ by $\mtt p_{h_j}(u(s,x/n))$, is
of order $n^{(d/2)-2}$, and therefore vanishes only for $d\le 3$; see
\eqref{240} and the remark which follows \eqref{145}.

\smallskip 
\noindent {\it Discussion of results and outline.} The first main part
of the article is the entropy estimate, Theorem \ref{t01}: for every
$T>0$, there exists a finite constant $C_1$ such that
\begin{equation*}
H_n(\mu_n(t) \,|\, \nu^n_{u(t, \cdot)}) \,+\,
\int_0^t   n^2\, I_n  (\mu_n(s);  \nu^n_{u(s, \cdot)} ) \; ds
\,\le\,  \big\{\, H_n(\mu_n \,|\, \nu^n_{u_0(\cdot)}) \, +\, C_1\,
\varkappa_n^{(d)}\,\big\} \, e^{C_1\,  t}
\end{equation*}
for all $0\le t\le T$, where $\varkappa^{(1)}_n = n^{(1/3)+b}$ for
some $b<1/6$, $\varkappa^{(2)}_n = n^{4/5}$ and
$\varkappa^{(3)}_n = n^{3/2} \kappa^{-3}$, the constant $C_1$ being
independent of $\kappa$.  Under the hypotheses of Theorem \ref{t1},
this yields that
$H_n(\mu_n(t) \,|\, \nu^n_{u(t, \cdot)}) = o(n^{d/2})$ uniformly on
$[0,T]$, and, in addition, a bound on the time-integral of the
Dirichlet form which measures the smoothness of the state of the
process.  Such a combination is reminiscent of the approach in
\cite{gpv} for the hydrodynamic limit.  The entropy alone controls
only the deviations of the density field at a fixed scale, whereas the
multi-scale argument formulated 
makes use of the Dirichlet form bound to help bound certain
local functions in the different scales.

The proof of Theorem \ref{t01} occupies Sections \ref{sec2} to
\ref{sec1}.  It starts, in Section \ref{sec2}, from the differential
inequality for the entropy with respect to the family of
time-dependent reference measures $\nu^n_{u(t,\cdot)}$.  With this
choice of reference measures, allowing to cancel zeroth order terms,
the remaining expression is a space-time average of the microscopic
projection $\Pi_{u(s,\cdot),g}$ for certain functions $g$ onto the
conserved particle numbers.  The later Boltzmann--Gibbs principle is
formulated to control precisely such quantities in a general context.

Sections \ref{sec4} and \ref{sec9} contain the `first' and `recursive'
steps of the multi-scale analysis, referring to preliminary estimates
given in Section \ref{sec3}.  Local functions are replaced, at the
expense of the Dirichlet form, by their average over a box of a larger
scale, and the procedure is iterated.  Each step improves the
estimate, while the accumulated cost remains of order
$\varkappa^{(d)}_n$, which is dimension dependent.  In particular,
with only the initial scale in $d=1$, one further iteration in $d=2$,
and two iterations in $d=3$, the final scaling is able to reach box
widths of order $n^{2/3 -b}$ for a small $b>0$ in $d=1$, $n^{4/5}$ in
$d=2$, and $\kappa^{-3} n^{1/2}$ for a $\kappa>0$.  The equivalence of
ensembles for inhomogeneous product measures, needed at each scale to
estimate canonical conditional expectations, is established in
Appendices \ref{sec5} and \ref{sec12}.  A log-Sobolev estimate is also
used in parts of the argument to bound expressions in terms of the
Dirichlet form.

The second main part of the article, the Boltzmann--Gibbs principle,
Theorem \ref{t02}, is proved in Section \ref{sec7} in a quantitative
form, with its own multi-scale scheme.  The argument here differs from
that of Theorem \ref{t02} in that we do not have an apriori Dirichlet
form bound to use.  However, since the relative entropy estimates in
Theorem \ref{t02} are now available, a successful multi-scale scheme
leveraging these estimates with the same scales used in Theorem
\ref{t02} can be fashioned.

Finally, the proof of Theorem \ref{t1} is completed in Section
\ref{sec11}, where the tightness of the sequence in
$L^1(0,T;\mc H_{-\mtt r})$, as well as the characterization of its
limit points is examined.

\smallskip {\it Future directions.}  Several questions remain open.
The most natural one is the extension to $d\ge 4$, which requires both
a sharper entropy estimate and a finer description of the
deterministic corrections.  The non-gradient case, in which the
diffusion coefficient is given by a variational formula rather than by
a closed expression, is also a natural direction. Although there has
been success in deriving stochastic Burgers and KPZ equations in $d=1$
(see Bertini and Giacomin \cite{Bert-Gia}, Amir, Corwin, and Quastel
\cite{ACQ}, Gon\c calves and Jara \cite{GJ}, Gon\c calves, Jara, and
Sethuraman \cite{GJS}, Yang \cite{Yang}, Bernardin, Funaki, and Sethuraman \cite{BFS},
Cannizzaro, Gon\c calves,
Misturini, and Occelli \cite{cgmo}, Huang, Matetski, and Weber \cite{HMW}
among others), the derivation more generally of non-linear stochastic
partial differential equations from microscopic dynamics is another
major challenge in the field.

\section{Notation and Results}
\label{s:1}

Let $\cb{\bb T^d_n} : = (\bb Z / n \bb Z)^d$,
$n \in \bb N =\{1, 2, \dots \}$, be the $d$-dimensional discrete torus
with $n^d$ points.  Let $\color{blue} \Omega_n = \{0,1\}^{\bb T^d_n}$
be the state space. Elements of $\Omega_n$ are represented by the
Greek letters $\eta = (\eta_x : x\in \bb T^d_n)$, $\xi$. Hence,
$\eta_x =1$ if site $x$ is occupied for the configuration $\eta$, and
$0$ otherwise.

Denote by $\color{blue} \{e_1,\dots,e_d\}$ the canonical basis of
$\bb R^d$.  Let $c_{0,e_j}\colon \{0,1\}^{\bb Z^d} \to \bb R$,
$1\le j\le d$, be strictly positive, cylinder functions `local'
functions which depend only on a finite number of variables
$\eta_z$]. Assume that there exists a constant $\mf c_0$ such that
\begin{equation}
\label{14}
c_{0,e_j}(\eta) \;\ge\; \mf c_0 \;>\; 0
\end{equation}
for all $\eta\in \{0,1\}^{\bb Z^d}$, $1\le j\le d$.  Assume that $c_{0,e_j}$
does not depend on the variables $\eta_0$ and $\eta_{e_j}$ and that
the following gradient conditions are in force. For each $j$, there
exist cylinder functions $h_{j,k}$, $1\le k\le d$, such that
\begin{equation}
\label{11}
c_{0,e_j} (\eta) \, [\, \eta_{e_j} \,-\, \eta_{0}\,] \,=\,
\sum_{k=1}^d 
\big\{\, (\tau_{e_k} \, h_{j,k})(\eta) \;-\; h_{j,k}(\eta) \,\big\}  \;.
\end{equation}
In this formula, $\cb{\{\tau_z : z\in \bb Z^d\}}$ represents the group of
translations acting on the configurations of $\{0,1\}^{\bb Z^d}$ or
$\{0,1\}^{\bb T_n^d}$:
\begin{equation}
\label{12}
(\tau_x \eta)_z \;=\;\eta_{x+z}\;, \quad x\,,\, z\,\in\, \bb A^d\;,
\;\; \eta\,\in\, \{0,1\}^{\bb A^d} \;, \quad \bb A = \bb T_n
\;\;\text{or}\;\; \bb Z\,.
\end{equation}
We use the same notation for translations acting on $\Omega_n$ and on
$\{0,1\}^{\bb Z^d}$, but the context will clarify to which one we are
referring.  When $\bb A=\bb T_n$, the summation has to be understood
modulo $n$.

Denote by $L_n$ the generator of the speed-change, symmetric
exclusion process given by
\begin{equation}
\label{10}
\cb{ (L_nf)\, (\eta) } \; :=\; n^2\, \sum_{x\in\bb T_n^d} \sum_{j=1}^d 
c_{x,x+e_j} (\eta)  \, \{f(T_{x,x+e_j}\eta)-f(\eta)\}\;.
\end{equation}
In this formula, on the bond $\{x, x+e_j\}$, the rate
$\cb{c_{x,x+e_j} (\eta) = c_{0,e_j} (\tau_x \eta)}$, and
$T_{ x,y}\eta$ represents the configuration of particles obtained from
$\eta$ by exchanging the values of $\eta_x$ and $\eta_y$:
\begin{equation*}
(T_{x,y}\eta)_{z} \;=\;
\begin{cases}
\eta_{z} & z\not =x \,,\, y\;,\\
\eta_{x} & z=y\;, \\
\eta_{y} & z=x\;.
\end{cases}
\end{equation*}

In the special case where $c_{0,e_j}(\eta) = 1$ for all $1\le j\le d$,
we recover the symmetric simple exclusion process on $\bb T^d_n$.  Let
$\color{blue} (\eta^n(t); t\geq 0)$ be the $\Omega_n$-valued,
continuous-time Markov chain whose generator is $L_n$.

Denote by $\color{blue} \nu^n_\alpha$, $0\le \alpha\le 1$, the Bernoulli
product measure on $\Omega_n$ with density $\alpha$. This is the product
measure whose marginals are Bernoulli distributions with parameter
$\alpha$. A straightforward computation shows that these measures
satisfy the detailed balance conditions for the speed-change exclusion
process because the cylinder functions $c_{0,e_j}$ are assumed not to depend
on $\eta_0$, $\eta_{e_j}$. In particular, they are stationary for this
dynamics.

For a cylinder function $h\colon \{0,1\}^{\bb Z^d} \to \bb R$, denote
by $\mtt p_h \colon [0,1] \to \bb R$ the function which maps $\alpha$
to the expectation of $h$ with respect to the Bernoulli product
measure $\nu_\alpha$ on $\{0,1\}^{\bb Z^d}$ with density $\alpha$:
\begin{equation}
\label{05}
\cb{\mtt p_h (\alpha) } \,:=\, E_{\nu_\alpha} [\,h\, ]\, , \quad \alpha\in[0,1]\,.
\end{equation}
Clearly, $\mtt p_h$ is a polynomial for each cylinder function
$h$.

\begin{example}
Here are two examples of jump rates $c_{0,e_1}(\cdot)$ satisfying
\eqref{11}:
\begin{equation*}
c_{0,e_1}(\eta) \,=\, a \,+\, [\,\eta_{-e_1} + \eta_{2e_1}\,]\,,
\quad
c_{0,e_1}(\eta) \,=\, a \,+\, [\,\eta_{e_2}  + \eta_{e_1+e_2}  +
\eta_{-e_2}  + \eta_{e_1-e_2} \,]
\end{equation*}
where $a>0$ to guarantee that the jump rates are strictly positive. In
the first case,
\begin{equation*}
h_{1,1}(\eta) \,=\, a \, \eta_0 \,-\, \eta_{-e_1}\eta_{e_1} \,+\, 
\eta_{-e_1} \eta_0 \,+\, \eta_0 \eta_{e_1}\,, \quad
h_{1,k}(\eta) \,=\, 0\,, \;\; 2\le k\le d\,.
\end{equation*}
In the second one, 
\begin{equation*}
h_{1,1}(\eta) \,=\, a \, \eta_0 \,+\, \eta_{-e_2}\eta_{0} \,+\,
\eta_{0}\eta_{e_2}\,, \;\;
h_{1,2}(\eta) \,=\, \eta_{e_1-e_2}\eta_{0} \,-\,
\eta_{e_1}\eta_{-e_2} \,, \quad
h_{1,k}(\eta) \,=\, 0\,, \;\; 3\le k\le d\,.
\end{equation*}
Mind that $\mtt p_{h_{j,k}} (\alpha) = 0$ for all $0\le \alpha\le 1$,
$j\neq k$. This seems to be a common feature of all gradient models.
\end{example}

From now on, instead of \eqref{11} we assume the stronger condition
that for each $j$, there exists a cylinder function $h_{j}$,
$1\le j\le d$, such that
\begin{equation}
\label{11b}
c_{0,e_j} (\eta) \, [\, \eta_{e_j} \,-\, \eta_{0}\,] \,=\, 
(\tau_{e_j} \, h_{j})(\eta) \;-\; h_{j}(\eta) \;.
\end{equation}
By  the proof of \cite[Proposition 5.7]{jlt21}, 
\begin{equation}
\label{125}
\cb{ D_j (\alpha) } \,:=\, \mtt p_{h_j}' (\alpha)
\,=\, E_{\nu_\alpha} [
\, c_{0,e_j} \,] \,, \quad\text{that is}\,, \quad
D_j (\alpha)  \,=\, \mtt p_{c_{0, e_j}} (\alpha)  \,, \quad
0\le \alpha \le 1 \,.
\end{equation}

Denote by $\color{blue} D([0,T], \Omega_n)$, $T>0$, the set of
right-continuous trajectories $\mf e\colon  [0,T] \to \Omega_n$ with
left-limits, endowed with the Skorohod topology.  For a probability
measure $\mu_n$ on $\Omega_n$, denote by $\color{blue} \bb P_{\mu_n}$
the measure on $D([0,T], \Omega_n)$ induced by the Markov chain
$\eta^n(t)$ and the initial distribution $\mu_n$.

Let $\bb T^d$ be the $d$-dimensional continuous torus. 
Denote by $\color{blue} \ms M (\bb T^d)$ the space of positive
measures on $\bb T^d$ with total mass bounded by $1$, equipped with
the weak topology.  Introduce the empirical measure
$\pi^n(\eta) \in \ms M(\bb T^d)$ associated with an element
$\eta \in \Omega_n$ as the measure on $\bb T^d$ defined by
\begin{equation*}
{\color{blue} \pi^n} \,=\, \pi^n(\eta) \,:=\, 
\frac{1}{n^d}\sum_{x\in \bb T^d_n}\eta_x \, \delta_{x/n}\;,
\end{equation*}
where $\color{blue}\delta_{x}$ stands for the Dirac mass at point $x$.  

Fix a density profile $u_0 \colon \bb T^d \to [0,1]$. Consider a
sequence of initial configurations $\zeta^n\in \Omega_n$, and suppose
that the associated empirical measure $\pi^n(\zeta^n)$ converges to
$u_0(x)\, dx$ weakly (and therefore in probability).  Then, for all
$t\ge 0$, the random measure $\pi^n(\eta^n(t))$, where $\eta^n(\cdot)$
is the process started from $\zeta^n$, converges in probability to the
absolutely continuous measure $u(t,x)\, dx$ whose density $u$ solves
the parabolic equation
\begin{equation}
\label{127}
\left\{
\begin{aligned}
& {\displaystyle \partial_t u = 
\sum_{j=1}^d \partial^2_{x_j} \mtt p_{h_j} (u)
\,=\,
\sum_{j=1}^d \partial_{x_j} ( D_{j} (u)
\, \partial_{x_j} u) \;,}
\\
& u(0, x) = u_0 (x) \;,
\end{aligned}
\right.
\end{equation}
where $D_j (\cdot)$ has been introduced in \eqref{125}.

For $T>0$ and a non-integer, positive number $\beta$, denote by
$\cb{H^\beta} = H^\beta(\bb T^d)$ the H\"older space of exponent
$\beta$, and by $\cb{ \| F \|_{(\beta)} }$ the H\"older-norm of an
element $F$ in $H^\beta(\bb T^d)$, see \cite[Eq.  (1.9), Chapter
I]{lsu68}.  Similarly, denote by
$\cb{H^{\beta/2,\beta}} = H^{\beta/2,\beta} ([0,T]\times \bb T^d)$ the
H\"older space of exponent $(\beta/2, \beta)$, and by
$\cb{ \| \mf u \|_{(\beta)} }$ the H\"older-norm of an element
$\mf u \colon [0,T]\times \bb T^d \to \bb R$ in $H^{\beta/2,\beta}$,
see \cite[Eq. (1.10), Chapter I]{lsu68}.  Whenever the underlying space
needs to be specified, we represent $H^{\beta/2,\beta}$ by
$\cb{H^{\beta/2,\beta}([0,T]\times \bb T^d)}$.

For the existence and uniqueness of solutions to equation \eqref{127}
we rely on \cite[Theorem V.6.1]{lsu68}, which states that this
equation has a unique solution in $H^{1+\beta/2, 2+\beta}$ if the
initial condition $u_0$ belongs to $H^{2+\beta}$. To verify condition
(e) of this theorem, consider the linear equation
$\partial_t v = \sum_{1\le j\le d} \partial_{x_j} [\, D_{j} (u_0) \,
\partial_{x_j} v \,]$. By \cite[Theorem IV.5.1]{lsu68}, this linear
equation with initial condition $u_0$ has a solution in
$H^{1+\beta/2, 2+\beta}$.

Assume that there exists a constant $0 < \mf r <1/2$ such that
\begin{equation}
\label{255a}
\mf r \,\le\, u_0(x) \,\le\, 1-\mf r\qquad
\text{for all} \;\;  x \in \bb T^d\,.
\end{equation}
By the maximum principle, \cite[Theorem 3.5]{pw12},
\begin{equation}
\label{255b}
\mf r \,\le\, u(t,x) \,\le\, 1-\mf r\qquad
\text{for all} \;\;  (t,x) \in \bb R_+\times \bb T^d\,.
\end{equation}

Our aim is to study the density fluctuations of this model, when the
process starts from a measure close to the product Bernoulli measure
associated with the density profile $u_0(\cdot)$. For a continuous
function $\rho\colon \bb T^d \to [0,1]$, denote by
$\nu^n_{\rho(\cdot)}$ the product measure on $\Omega_n$ given by
\begin{equation*}
\cb{\nu^n_{\rho(\cdot)} \{\eta_x=1\}}
\,=\, \rho(x/n)\, , \quad x\in \bb T^d_n\,.
\end{equation*}

\subsection*{The density fluctuation field}

Denote by $\color{blue} C(\bb T^d)$ the space of continuous,
real-valued functions on $\bb T^d$ and by $\color{blue} C^k(\bb T^d)$,
$k\ge 0$, the space of real-valued functions on $\bb T^d$
whose partial derivatives up to order $k$ exist and are
continuous. Elements of $C(\bb T^d)$ are represented by the letters
$F$, $G$. For $F \in C^k(\bb T^d)$, denote by $\Vert F\Vert_{C^k}$ the
norm of $F$:
\begin{equation*}
\cb{\Vert F \Vert_{C^k} } \,:=\,
\sum_{(j_1, \dots , j_d)} \Vert  \partial^{j_1}_{x_1}
\partial^{j_2}_{x_2} \cdots \partial^{j_d}_{x_d} F \Vert_{L^\infty (\bb T^d)}\,,
\end{equation*}
where the sum is carried out over all vectors $(j_1, \dots, j_d)$ such
that $j_a\ge 0$, $j_1 + j_2 + \cdots + j_d \le k$. In this formula
$\cb{\partial^{j_a}_{x_a} F}$ represents the $j_a$-th partial
derivative of $F$ with respect to $x_a$. For $T>0$, $p\ge 0$, and a
function $F\colon [0,T]\times \bb T^d \to \bb R$ such that
$F(t, \cdot) \in C^p(\bb T^d)$ for all $0\le t\le T$ let
\begin{equation}
\label{256}
\cb{\Vert F\Vert_{T,p}} \, := \, \sup_{0\le t\le T} \Vert F(t, \cdot)
\Vert_{C^p}\,.
\end{equation}

Denote by $\color{blue} L^2(\bb T^d)$ the space of complex-valued,
square-integrable, measurable functions on $\bb T^d$ endowed with the
usual scalar product, represented by $\<\,\cdot\,,\, \cdot\,\>$.  Let
$\color{blue} \mc H_{\mtt r}$, $\mtt r>0$, be the completion of
$C^\infty(\bb T^d)$ with respect to the scalar product
$\<\,\cdot\,,\, \cdot\,\>_{\mc H_{\mtt r}}$ defined by
\begin{equation*}
\< F\,,\, G\>_{\mc H_{\mtt r}} \;=\; \sum_{m\in \bb Z^d}   \,
\gamma_m^{\mtt r} \,
\mb F_m \, \overline{\mb G_m} \;,
\end{equation*}
where $\color{blue} \gamma_m \;:=\; 1\,+\, \Vert m\Vert^2$ and
$\mb F_m = \< F\,,\, \phi_m\>$,
$\color{blue} \phi_m(x) = \exp\{ 2\pi i x\cdot m\}$ for
$m\in \bb Z^d$. The sum converges because $F$, $G$ belong to
$C^\infty(\bb T^d)$. It is well known that $\mc H_{\mtt r}$ is
compactly embedded in $\mc H_{\mtt s}$ for $\mtt r> \mtt s$.

Denote by $\color{blue} \mc H_{- \mtt r}$, $\mtt r>0$, the dual space of
$\mc H_{\mtt r}$. Elements of $\mc H_{- \mtt r}$ are represented by
the letters $X$, $Y$. For $X$ in $\mc H_{- \mtt r}$ and $F$ in
$\mc H_{\mtt r}$, $X(F)$ can be represented as
\begin{equation*}
X(F) \;=\; \sum_{m\in \bb Z^d} X(\phi_m) \, \mb F_m \;.
\end{equation*}

Denote by $X^n_t$ the random element of $\mc H_{- \mtt r}$,
$\mtt r>d/2$, defined by
\begin{equation*}
X_t^n(F) \;=\; \frac{1}{n^{d/2}}
\sum_{x\in \bb T^d_n} F(x/n)\, [\, \eta_x^n(t) \,-\, u(t,x/n) \, ]\;, \quad
F\,\in\, \mc H_{\mtt r} \;.
\end{equation*}
This formula defines an $\mc H_{- \mtt r}$-valued process
$\{X_t^n; t \geq 0\}$. We call this process the {\em density
fluctuation field}.

Fix $T>0$, $\mtt r>d/2$, and denote by $L^1( 0, T ;\mc H_{-\mtt r})$
the space of $\mc H_{- \mtt r}$-valued functions in $L^1([0,T])$:
$\int_0^T \Vert F(t)\Vert_{\mc H_{- \mtt r}} \, dt <\infty$ (cf. \cite{s86}).  Let
$\mu_n$ be a sequence of probability measures on $\Omega_n$. Denote by
$\cb{\bb Q_{\mu_n}}$ the probability measure on
$L^1( 0, T ;\mc H_{- \mtt r})$ induced by the process
$\{X_t^n; t \geq 0\}$ and the initial distribution $\mu_n$.

For probability measures $\mu$, $\nu$ on $\Omega_n$ such that $\nu$
has full support, denote by $H_n(\mu|\nu)$ the relative entropy of
$\mu$ with respect to $\nu$:
\begin{equation*}
\cb{H_n(\mu|\nu)} \,:=\,
\sum_{\eta\in\Omega_n} \mu(\eta)\, \ln \frac{\mu(\eta)}{\nu(\eta)}
\,\cdot
\end{equation*}
In many places we write $\cb{H_n(f|\nu)}$ instead of $H_n(\mu|\nu)$,
where $f$ is the Radon--Nikodym derivative of $\mu$ with respect to
$\nu$.

Throughout the article, for a sequence of non-negative real numbers
$s_n$, we adopt the notation $\cb{s_n = o(n^a)}$, $a>0$, to mean that
$\lim_{n\to\infty} s_n/n^a =0$.

\begin{theorem}
\label{t1}
Let $d\le 3$. Fix $T>0$, a density profile
$u_0\colon \bb T^d\to (0,1)$, and a sequence of probability measures
$\mu_n$ on $\Omega_n$. Assume that

\begin{itemize}
\item [(a)] $H_n(\mu_n \,|\, \nu^n_{u_0(\cdot)}) = o(n^{d/2})$;

\item [(b)] there exist $\mtt r_0 > d/2$ and a random
$\mc H_{-\mtt r_0}$-valued density field $X_0$ such that $X^n_0$
converges weakly  to $X_0$ in $\mc H_{- \mtt r_0}$;

\item [(c)] the density profile $u_0$ belongs to
$H^{p+\beta}(\bb T^d)$, for some integer $p\ge \max\{4, \mtt r_0 + d/2
+ 2\}$, $0<\beta<1$, and satisfies \eqref{255a}. 

\end{itemize}
Let \begin{equation}
\label{242b}
\mtt r \;>\; \cb{\mtt r^\star } \,:=\, 
\max\Big\{\, 6 \,+\, (3d+1)/2\,,\;\;
\mtt r_0 \,+\, d\,+\,4 \,\Big\}\;.
\end{equation}
Then, the sequence of probability measures $\bb Q_{\mu_n}$ converges
weakly in $L^1( 0, T ;\mc H_{- \mtt r})$ to the measure $\bb Q$ induced by
the solution to  the equation
\begin{equation}
\label{58}
\left\{
\begin{aligned}
& \partial_t Y_t \;=\; \mc A^*_t \, Y_t   \;+\;  \sum_{j=1}^d
\partial_{x_j} \Big(\, \sqrt{2\, \chi(u(t))\, D_j(u(t)) }
\, \xi^j_t \,\Big)\;, \\
& Y_0\;=\; X_0\;.
\end{aligned}
\right.
\end{equation}
In this formula, $\cb{\chi(\rho) := \rho(1-\rho)}$ is the static
compressibility of the exclusion dynamics,
$\xi = (\xi^1, \dots, \xi^d)$ is a $d$-dimensional standard space-time
white noise, $\mc A_t$, $t\ge 0$, is the second-order linear
elliptic operator defined by
\begin{equation}
\label{41}
\cb{(\mc A_t \, F)(x)} \; :=\; \sum_{j=1}^d  D_{j} (u(t,x))\,
\partial^2_{x_j} F (x) \;,\quad x\in \bb T^d\,, 
\end{equation}
for functions $F$ in $C^2(\bb T^d)$, and $\mc A^*_t$ stands for its
adjoint, so that \eqref{58} has to be understood in the
weak sense:
\begin{equation*}
d Y_t (F) \;=\; Y_t(\mc A_t F)   \, dt \;+\;  \sum_{j=1}^d
\int_{\bb T^d} \sqrt{2\, \chi(u(t))\, D_j(u(t)) }\, 
(\partial_{x_j} F) (x) \, \xi^j (dt, dx) \;,
\end{equation*}
$F\in \mc H_{\mtt r}$.
\end{theorem}

\begin{remark}
It follows from the assumption
$\lim_{n\to \infty} n^{-d/2} H_n(\mu_n \,|\, \nu^n_{u_0(\cdot)}) =0$
and the entropy inequality 
that under $\mu_n$ the empirical measure $\pi^n$ converges to
$u_0(x)\, dx$ in $\ms M(\bb T^d)$.

One may view the process $Y_t$ as a time-dependent form of a
generalized Ornstein-Uhlenbeck process.
\end{remark}

\section{Sketch of the proof of Theorem \ref{t1}}
\label{sec13}

In this section, we provide a brief sketch of the proof which
highlights the importance of the Boltzmann--Gibbs principle.  

Fix a test function $F\colon \bb T^d \to \bb R$ and write the
semi-martingale decomposition of the process $X^n_t(F)$:
\begin{equation}
\label{30}
X^n_t(F)  \, =\;  X^n_0(F) \; + \; \int_0^t
(\partial_s + L_n)  X^n_s(F)\, ds\, +\, M^n_t(F)  \;.
\end{equation}
By \cite[Lemma A1.5.1]{kl}, the process $t \mapsto M^n_t(F)$ is a
martingale.  On the other hand, by \eqref{11b} and a summation by
parts,
\begin{equation*}
(\partial_s + L_n)  X^n_s(F) \,=\, 
\frac{1}{n^{d/2}} \, \sum_{j=1}^d \sum_{x\in\bb T^d_n}
(\Delta^n_j F) (x/n) \, \tau_x h_{j} (\eta^n(s)) 
\,-\, \frac{1}{n^{d/2}} \, \sum_{x\in\bb T^d_n} F(x/n) \, (\partial_s
u)(s,x/n) \; ,
\end{equation*} 
where
\begin{equation}
\label{17}
\cb { (\Delta^n_j F) (x/n)}  \,:=\,
n^2\, \big\{ \, F((x+e_j)/n) \,+\, F((x - e_j)/n)
\,-\, 2\, F(x/n) \, \big\} \,.
\end{equation}

For a cylinder function $g$ and a continuous density profile
$\rho\colon \bb T^d \to (0,1)$, let $\Pi_{\rho (\cdot), g} (x, \eta)$
be the cylinder function given by
\begin{equation}
\label{179}
\cb{\Pi_{\rho (\cdot), g}  (x,\eta) } \, =\, (\tau_x g) (\eta)
- E_{\nu^n_{\rho(\cdot)} } [ \tau_x g] \, -\, \mtt p_g' (\rho (x/n))
\, [\eta_x - \rho(x/n) ]\,,
\end{equation}
where $\mtt p_g(\cdot)$ is the polynomial introduced in \eqref{05}.
Note that $\Pi_{\rho(\cdot) , g} (x, \eta)$ has zero mean with respect
to $\nu^n_{\rho(\cdot)} $. It represents the part of the observable
$\tau_x g$ which is orthogonal to the observable $\eta_x$ associated
to the unique conserved quantity, the mass $\sum_{x\in\bb T^d_n} \eta_x$.

With this notation, rewrite $(\partial_s + L_n)  X^n_s(F)$ as
\begin{equation}
\label{07}
\begin{aligned}
& \frac{1}{n^{d/2}} \, \sum_{j=1}^d \sum_{x\in\bb T^d_n}
(\Delta^n_j F) (x/n) \, \Pi_{u(s, \cdot), h_j}  (x,\eta^n(s))
\\
& \,+\, \frac{1}{n^{d/2}} \, \sum_{j=1}^d \sum_{x\in\bb T^d_n}
(\Delta^n_j F) (x/n) \, \mtt p'_{h_j} (u(s,x/n))  \,[\, \eta_x^n(s) - u(s,x/n)\,]
\\
& \,+\, \frac{1}{n^{d/2}} \, \sum_{j=1}^d \, \sum_{x\in\bb T^d_n}
\Big\{\, (\Delta^n_j F) (x/n) \, E_{\nu^n_{u(s,\cdot)}} \, [\, \tau_x h_{j} \,]
\, -\, F(x/n) \,  (\partial_s u)(s,x/n) \,\Big\}\,.
\end{aligned}
\end{equation}

The last term is deterministic.  Since $u(s,x)$ is the
solution of the hydrodynamic equation \eqref{127} and
$\partial_{x_j}  ( D_j (u) \, \partial_{x_j} u)$ can be written as
$\partial^2_{x_j} \, \mtt p_{h_j} (u)$, this last term is equal to
\begin{equation*}
\frac{1}{n^{d/2}} \, \sum_{j=1}^d \, \sum_{x\in\bb T^d_n}
\Big\{\,   (\Delta^n_j \,  F) (x/n) \,  E_{\nu^n_{u(s,\cdot)}} \, [\, \tau_x h_{j} \,]
\, -\, F (x/n)\, 
\partial^2_{x_j} \, \mtt p_{h_j} (u(s,x/n)) \,\Big\} \;.
\end{equation*}
Let $G$ be a function in $C^4(\bb T^d)$. A fourth-order Taylor
expansion yields that there exists a finite constant $C_0$ such that
$\sup_{x\in \bb T^d_n} |\, (\Delta^n_j G) (x/n) - (\partial^2_{x_j} G)(x/n)
\,|\,\le\, C_0\, \Vert G\Vert_{C^4}/n^2$. Therefore, replacing
$\partial^2_{x_j} \mtt p_{h_j} (u(s,x/n)) $ by
$\Delta^n_j \mtt p_{h_j} (u(s,x/n))$, and summing by parts yields that
the previous expression is equal to
$\cb{\mtt D^{(n)}_s(F) } := \mtt D^{1,n}_s(F) +\mtt D^{2,n}_s(F)$, where
\begin{gather*}
\cb{\mtt D^{1,n}_s(F) } \,:=\, \frac{1}{n^{d/2}} \, \sum_{j=1}^d \, \sum_{x\in\bb T^d_n}
(\Delta^n_j \,  F) (x/n) \, \big\{\,    E_{\nu^n_{u(s,\cdot)}} \, [\, \tau_x h_{j} \,]
\, -\, E_{\nu^n_{u(s,x/n)}} \, [\, h_{j} \,]  \,\big\} \,,
\\
\cb{\mtt D^{2,n}_s(F)} \,:=\, 
\frac{1}{n^{d/2}} \, \sum_{j=1}^d \, \sum_{x\in\bb T^d_n}
F (x/n)\,  \big\{\,
[\, \Delta^n_{j} -\partial^2_{x_j} \,] \, \mtt p_{h_j} (u(s,x/n)) \,\big\}
\;.
\end{gather*}
Recall from \eqref{256} the definition of $\Vert F\Vert_{T,p}$, and
note that
$|\, \mtt D^{2,n}_s(F)\,| \le C_0 \, \Vert u \Vert_{T,4}\, \Vert
F\Vert_\infty\, n^{(d/2)-2}$ for some finite constant $C_0$ depending
only on the cylinder functions $h_j$, $1\le j\le d$. By Corollary
\ref{s24}, $|\, \mtt D^{1,n}_s(F)\,|$ is bounded by
$\mf c_1 \, \Vert F\Vert_{C^2} \, n^{(d/2)-2}$ for some finite
constant $\mf c_1$ which depends only on $\mf r$,
$\Vert u \Vert_{T,2}$, and $h_j$, $1\le j\le d$. In conclusion, since
$\mtt p'_{h_j} (\alpha) = D_j(\alpha)$, replacing
$(\Delta^n_j F) (x/n) $ by $(\partial^2_{x_j} F) (x/n) $, in view of
the definition \eqref{41} of $\mc A_s F$,
\begin{align}
\label{145}
(\partial_s + L_n)  X^n_s(F) \, =\, X^n_s(\mc A_s F) 
\,+\, 
\frac{1}{n^{d/2}} \, \sum_{j=1}^d \sum_{x\in\bb T^d_n}
(\Delta^n_j F) (x/n) \, \Pi_{u(s, \cdot), h_j}  (x,\eta^n(s))
+ \mtt D^{(n)}_s(F) \,,
\end{align}
where, from now on, $\mtt D^{(n)}_s(F)$ stands for the sum of
$\mtt D^{1,n}_s(F)$, $\mtt D^{2,n}_s(F)$ and of the error committed in
replacing $(\Delta^n_j F) (x/n) $ by $(\partial^2_{x_j} F) (x/n) $ in
the second line of \eqref{07}. The error $\mtt D^{(n)}_s(F)$ is bounded
uniformly in the configuration: 
\begin{align}
\label{240}
|\, \mtt D^{(n)}_s(F) \,| \,\le\,  \mf c'_1 \Vert F\Vert_{C^4}
n^{(d/2)-2}
\end{align}
for all functions $F$ in $C^4 (\bb T^d)$. In this formula,
$\mf c'_1$ is a finite constant which depends only on $\mf r$,
$\Vert u \Vert_{T,4}$, and $h_j$, $1\le j\le d$.

In view of the previous computation, \eqref{30} becomes
\begin{equation}
\label{145b}
X^n_t(F)  \, =\;  X^n_0(F) \; + \; \int_0^t
X^n_s( \mc A_s F)\, ds\, +\, M^n_t(F)  \,+\, \int_0^t \mtt B^{(n)}_s(F) \,
ds \,+\, \int_0^t \mtt D^{(n)}_s(F)  \, ds \;,
\end{equation}
where
\begin{align}
\label{241}
\cb{\mtt B^{(n)}_s(F)} \,:=\, 
\frac{1}{n^{d/2}} \, \sum_{j=1}^d \sum_{x\in\bb T^d_n}
(\Delta^n_j F) (x/n) \, \Pi_{u(s, \cdot), h_j}  (x,\eta^n(s))\,.
\end{align}
Denote by $(P_{s,t} : 0\le s\le t)$ the semigroup
associated to the time-dependent linear operator $\mc A_s$. Fix $t>0$,
a function $F$ in $C^4(\bb T^d)$ and let $F_s = P_{s, t} F$,
$0\le s\le t$, so that $F_0 = P_{0,t} F$, $F_t =F$,
$\partial_s F_s = - \mc A_s F_s$. The previous calculations, replacing
the time-independent test function $F$ by $F_s$, yield that
\begin{equation*}
X^n_t(F)  \, =\;  X^n_0(P_{0,t}F) 
\, +\, M^n_t( P_{\cdot, t} F)  \,+\, \int_0^t \mtt B^{(n)}_s( P_{s, t} F) \,
ds \,+\, \int_0^t \mtt D^{(n)}_s( P_{s, t} F)  \, ds \;,
\end{equation*}
where $M^n_t( P_{\cdot, t} F) $ is the Dynkin martingale associated to
the time-dependent test function $F_s = P_{s,t} F$, $0\le s\le t$.

Therefore, since the remainder $\mtt D^{(n)}_s( P_{s, t} F)$ vanishes
as $n\to\infty$ uniformly in time, the proof of Theorem \ref{t1} is
reduced to (a) the tightness of the process $X^n_t$, (b) the
characterization of the limit of the sum of the first two terms on the
right-hand side, and (c) the proof that the time integral of the term
$\mtt B^{(n)}_s( P_{s, t} F) $ vanishes.

The tightness of the process is examined in Section \ref{sec11}, where
we also show that the processes $X^n_0(P_{0,t}F)$ and
$M^n_t( P_{\cdot, t} F) $ are asymptotically independent. The
hypothesis that the density field converges at time $0$ yields that
$X^n_0(P_{0,t}F)$ converges.  It will be shown in Section \ref{sec11}
that $M^n_t( P_{\cdot, t} F) $ converges as well taking advantage of
classical theorems on the convergence of martingales. The main
difficulty in the proof is (c). This is the content of Theorem
\ref{t02} below, the so-called Boltzmann--Gibbs principle, whose proof
is given in Section \ref{sec7}. The proof requires a bound on the
entropy of the state of the process with respect to a reference
measure, whose proof runs from Sections \ref{sec2} to \ref{sec1}.

\smallskip We complete the sketch by stating the Boltzmann--Gibbs
principle, derived for the first time by Brox and Rost \cite{BR84} in
the context of equilibrium fluctuations. It asserts that the local
observables 
\begin{equation*}
X^{n,f}_t(G) \,:=\,
\frac{1}{n^ {d/2}} \sum_{x\in \bb T^d_n} G (x/n)\, \big\{ \,
f(\tau_x \eta^n(t)) -  E_{\nu^n_{u(t,\cdot)}} [\tau_x f] \,\big\}
\end{equation*}
may be replaced by their projection on the conserved quantity, in the
present context the number of particles.  It reads as follows.  Let
$\color{blue} C^{j,k} ([0,t]\times \bb T^d)$, $j$, $k\ge 0$, be the
space of real-valued functions on $[0,t]\times \bb T^d$ whose time,
space partial derivatives up to order $j$, $k$ exist and are
continuous, respectively.

\begin{theorem}[Boltzmann--Gibbs principle]
\label{t02}
Fix $d\le 3$. Assume that the initial condition $u_0$ satisfies
\eqref{255a}, and belongs to $H^{4+\beta}$ for some $0<\beta<1$.
Let $u(t,x)$ be the solution of the hydrodynamic equation \eqref{127}, and
$\mu_n$ be a sequence of probability measures on $\Omega_n$ such
that $H_n(\mu_n \,|\, \nu^n_{u_0(\cdot)} ) = o(n^{d/2})$. Then,
\begin{equation*}
\lim_{n\to\infty} \bb E_{\mu_n} \Big[\,\Big| \int_0^t
\frac{1}{n^{d/2}} \sum_{x\in \bb T^d_n} G (s, x/n)\,
\Pi_{u(s, \cdot), h}  (x,\eta^n(s)) \, ds \, \Big| \, \Big] \;=\; 0
\end{equation*}
for all $t>0$, functions $G$ in $C^{0,2}([0,t] \times \bb T^d)$, and
cylinder functions $h\colon \{0,1\}^{\bb Z^d} \to \bb R$.
\end{theorem}

\begin{remark}
The bound of order $n^{(d/2)-2}$ for $\mtt D^{(n)}_s(F) $ in
\eqref{145} is bad in dimension $d\ge 4$.  In dimension $d\ge 4$, it
might be necessary to take the density profile $u(t, \cdot )$ as the
solution of a semi-discrete in space PDE, and to consider a higher
order expansion of $E_{\nu^n_{\rho(\cdot)}} \, [\, h_j \,]$.
\end{remark}

\subsection*{Entropy production}

The proofs of Theorems \ref{t1} and \ref{t02} rely on the entropy
bound stated in the next result. The proof of this estimate extends
across a good portion of the article. Denote by $(S^n_t: t\ge 0)$ the
semigroup associated to the generator $L_n$, and fix an initial
measure $\mu_n$ on $\Omega_n$.  Hereafter, we represent
$H_n( \mu_n S^n_t \,|\, \nu^n_{u(t, \cdot) } )$ by
$\cb{H_n( f_t \,|\, \nu^n_{u(t, \cdot) } )}$, where
$\cb{f^n_t} := d(\mu_n S^n_t)/d \nu^n_{u(t, \cdot) } $, and we omit
the index $n$ of $f^n_t$. 

For a probability measure $\mu$ on
$\Omega_n$, denote by $I_n(h; \mu)$ the carr\'e du champ of a
function $h\colon \Omega_n \to \bb R$ averaged with respect to $\mu$:
\begin{equation}
\label{27}
\begin{gathered}
{\color{blue} I_n(h; \mu)}
\;:=\; \frac{1}{2}\, \sum_{j=1}^d \sum_{x\in\bb T^d_n} 
\int  c_{x,x+e_j} (\eta)\,  \Big[\,  \sqrt{h (T_{x,x+e_j} \eta )} \,-\,
\sqrt{h(\eta)}\, \Big]^2 \; \mu (d\eta) \;.
\end{gathered}
\end{equation}

\begin{theorem}
\label{t01}
Fix $T>0$. Assume that the initial condition $u_0$ satisfies
\eqref{255a}, and belongs to $H^{4+\beta}$ for some $0<\beta<1$.  Let
$u(t,x)$ be the solution of the hydrodynamic equation
\eqref{127}. Then, in dimension $d=1$, $2$, there exists a finite
positive constant $C_1$, depending only on $\mf r$, $h_j$,
$c_{0, e_j}$, $1\le j\le d$, $\Vert u \Vert_{T,4}$, and $T$, such that
\begin{equation}
\label{133}
H_n(f_t | \nu^n_{u(t, \cdot) }) \,+\,
\int_0^t   n^2\, I_n  (f_s;  \nu^n_{u(s, \cdot) } ) \; ds
\,\le\,  \big\{ H_n(f_0 | \nu^n_{u_0(\cdot) }  ) \, +\, C_1\,
\varkappa_n^{(d)}\,\big\} \, e^{C_1\,  t}
\end{equation}
for all $0\le t\le T$, $n\ge 1$.  In dimension $1$, a constant
$0<b<1/6$ is fixed, $\varkappa_n^{(1)} = n^{(1/3) +b}$, and the
constant $C_1$ may depend on $b$. In dimension $2$,
$\varkappa_n^{(2)} = n^{4/5}$.

In dimension $3$, fix $\kappa>1$.  Then, there exists a finite
positive constant $C_1$, depending only on $\mf r$, $h_j$,
$c_{0, e_j}$, $1\le j\le d$, $\Vert u \Vert_{T,4}$, and $T$, and an
integer $n_0$ depending on the previous variables and on $\kappa$ for
which \eqref{133} holds for all $n\ge n_0$, $0\le t\le T$, where
$\varkappa_n^{(3)} = n^{3/2}/\kappa^3$.
\end{theorem}

\begin{remark}
The proof of this result is based on a block lemma which allows one to
replace space averages of a local function $h$ by a function of the
particles' empirical density over mesoscopic cubes.  The method
presented in the next sections allows one to reach cubes of width
\begin{equation*}
n^{(2/3)-b} \;\;\text{in dimension $d=1$,}\quad
n^{3/5} \;\;\text{for $d=2$, }\quad
\kappa\, n^{1/2} \;\;\text{for $d=3$}\,,
\end{equation*}
where $\kappa\to\infty$ after $n$.
\end{remark}

\begin{corollary}
\label{c-t01}
Assume that $d=3$ and that
$H_n(f_0 | \nu^n_{u_0(\cdot) } ) = o(n^{d/2})$. Then, under the
hypotheses of Theorem \ref{t01}, there exists an increasing sequence
$(\alpha_n)_{n\ge 1}$ such that $\alpha_n\to\infty$, and
\begin{equation*}
\limsup_{n\to \infty} 
\frac{\alpha_n }{n^{d/2}}
\,\Big\{\, \sup_{0\le t\le T} H_n (f_t \,|\, \nu^n_{u (t, \cdot)}  )
\,+\, \int_0^T   n^2\, I_n  (f_s;  \nu^n_{u(s, \cdot) } ) \; ds\, \Big\}
\,=\, 0\,.
\end{equation*}
\end{corollary}

\begin{proof}
By Theorem \ref{t01}, the result holds if we set $\alpha_n=1$ for all
$n\ge 1$, letting $\kappa\to\infty$ after $n\to\infty$.  If
$(a_n)_{n\ge 1}$ is a sequence of real numbers converging to $0$, it
is easy to show that there exists a sequence $(\alpha_n)_{n\ge 1}$
satisfying the conditions of the corollary and such that
$\alpha_n\, a_n \to 0$.
\end{proof}

\section{Entropy production}
\label{sec2}

In this section, we introduce reference measures and present a bound
on the time derivative of the entropy of the distribution of the
system with respect to the reference measure. Fix an initial condition
$u_0\colon \bb T^d\to (0,1)$ satisfying the hypotheses of Theorem
\ref{t1}. 

We start with some notation.  Let
$\cb{\psi^{\rm ref}} \colon \bb R_+ \times \Omega_n \to \bb R_+$ be a
sequence of smooth densities with respect to the measure
$\nu^{n}_{1/2}$. This means that $\psi^{\rm ref} (t, \eta) > 0$ for
all $(t,\eta) \in \bb R_+ \times \Omega_n$,
$E_{\nu^{n}_{1/2}} [\psi^{\rm ref} (t, \cdot) ] = 1$ for all $t\ge 0$,
and the function $t\mapsto \psi^{\rm ref} (t,\eta)$ belongs to
$C^1(\bb R_+)$ for each $\eta \in \Omega_n$. We often represent
$\psi^{\rm ref} (t, \eta)$ by $\cb{\psi_t^{\rm ref} (\eta)}$.  Denote
by $\nu^{n,t}_{\rm ref}$, $t\ge 0$, the family of reference
probability measures on $\Omega_n$ given by
\begin{equation}
\label{18}
\nu^{n,t}_{\rm ref}(\eta) \,=\, \psi^{\rm ref}_t (\eta)
\, \nu^{n}_{1/2}(\eta)
\, , \quad \eta\in \Omega_n\,.
\end{equation}
The next result is the so-called Jara--Yau lemma for time-dependent
reference measures. We refer to \cite[Lemma A.1]{jm2} for a
proof. 

\begin{lemma}
\label{l03}
For all $t\ge 0$,
\begin{equation}
\label{2-22}
\frac{d}{dt}\, H_n(f_t | \nu^{n,t}_{\rm ref}) \;\leq\; -\,2\,  
n^2\, I_n(f_t ; \nu^{n,t}_{\rm ref})
\;+\;  \int
\big\{ L^{{\rm ref}, *}_{n,t} \, \bs 1 \, -\, \partial_t   \ln \psi^{\rm ref}_t \,\big\}
\; f_t \, d \nu^{n,t}_{\rm ref} \;,
\end{equation}
where, recall, $\bs 1 \colon \Omega_n \to \bb R$ is the constant
function equal to $1$, and $\cb{L^{{\rm ref}, *}_{n,t}}$ represents the
adjoint of $L_n$ in $L^2(\nu^{n,t}_{\rm ref})$.
\end{lemma}

Fix $T>0$, and denote by $u$ the solution of the hydrodynamic equation
\eqref{127}. Since $u_0$ belongs to $H^{p+\beta}$ for some $p\ge 4$
and satisfies \eqref{255a}, by the results presented in the paragraph
following \eqref{127}, $u$ belongs to $C^{1,4}([0,T] \times \bb T^d)$,
and fulfills \eqref{255b}. Throughout this section,
$\cb{\mf c_{1,q}}$, $q\ge 0$, represents a finite constant, whose
value may change from line to line, and which depends only on
$\mf r$ and $\sup_{0\le t\le T} \Vert u(t, \cdot)\Vert_{C^q}$, and
the model (the dimension and the local functions $c_{0,e_j}$,
$h_j$). Since $\Vert u(t, \cdot )\Vert_\infty \le 1$, the constant
$\mf c_{1,0}$ depends only on $\mf r$ and the model.

Denote by $\color{blue} L^*_{n,t}$ the adjoint of $L_n$ in
$L^2(\nu^{n}_{u(t, \cdot)})$, and by $\cb{\psi^n_t}$ the function
introduced in \eqref{18} with $\nu^{n}_{u (t, \cdot)} $ taken as
reference measure: $\nu^{n,t}_{\rm ref} = \nu^{n}_{u(t, \cdot)}$.

\subsection*{Computation of $L^*_{n,t}  \bs 1$ and $\partial_t   \ln
\psi^n_t$}

An elementary computation yields that
\begin{equation}
\label{19}
\partial_t   \ln \psi^{n}_t  (\eta) \,=\,
\sum_{x\in\bb T^d_n} \frac{\partial_t u(t,x/n)}{\chi(u(t,x/n))}\,
[\, \eta_x - u(t,x/n)\,]\,,
\end{equation}
and that for every function $F\colon \Omega_n \to \bb R$,
\begin{equation*}
L_{n,t}^{*} F(\eta)
\,=\, n^2 \sum_{j=1}^d \sum_{x\in\bb T^d_n}
c_{x,x+e_j}(\eta)\, 
\Big[F(T_{x,x+e_j} \eta) \,
\frac{\nu^{n}_{u(t, \cdot)} (T_{x,x+e_j} \eta)}
{\nu^{n}_{u(t, \cdot)} (\eta)} - F(\eta)\Big] \,.
\end{equation*}
Replacing $F$ by $1$ yields that
\begin{align*}
L_{n,t}^{*} \, \bs 1 \, &=\,
n^2\, \sum_{j=1}^d \sum_{x\in\bb T^d_n} c_{x,x+e_j}(\eta)\,   
\Big\{ \frac{\nu^n_{u(t,\cdot)} (T_{x,x+e_j}\eta)}{\nu^n_{u(t, \cdot)}
(\eta)} \,-\, 1\,\Big\}
\\
&=\,
n^2\, \sum_{j=1}^d \sum_{x\in\bb T^d_n} c_{x,x+e_j}(\eta)\,    
\Big\{ e^{ - [U(t,[x+e_j]/n) - U(t,x/n)] [\eta_{x+e_j} - \eta_x]} \,-\,
1\,\Big\}\,, 
\end{align*}
where $\cb{ U(t,x) = \ln \{ u(t,x)/[1-u(t,x)]\}}$.

Let $g_j$, $1\le j\le d$, be the cylinder function given by
\begin{align*}
\cb{g_j (\eta) } \,:=\, c_{0,e_j}(\eta)\, [\eta_{e_j} - \eta_0]^2\,.
\end{align*}
By \eqref{125},
\begin{equation}
\label{126}
\mtt p_{g_j} (\alpha) \,=\, 2\, \chi(\alpha)\, D_j(\alpha) \,,\quad
0<\alpha<1\,.
\end{equation}

We claim that
\begin{equation}
\label{106}
\begin{aligned}
L_{n,t}^{*} \, \bs 1 \,  &=\,
\sum_{j=1}^d \sum_{x\in\bb T^d_n}
(\partial^2_{x_j} U)  (t,x/n) \, \{\, \tau_x h_j -
E_{\nu^{n}_{u(t, \cdot)}}[\tau_x h_j ]  \,\}
\\
& +\, \frac{1}{2}\,
\sum_{j=1}^d \sum_{x\in\bb T^d_n}
\big[ (\partial_{x_j} U)  (t,x/n) \big]^2\,
\{\,  \tau_x g_j - E_{\nu^{n}_{u(t, \cdot)}}[\tau_x g_j]  \,\}
\,+\, I_n(t) \,+\,  \mf R^{(4)}_{t,n}\, n^{d-2}\,,
\end{aligned}
\end{equation}
where
\begin{equation}
\label{106c}
\cb{I_n(t)} \,:=\, \frac{1}{2n}\,
\sum_{j=1}^d \sum_{x\in\bb T^d_n}
(\partial_{x_j} U)  (t,x/n) \,
(\partial^2_{x_j} U)  (t,x/n)  \,
\{\,  \tau_x g_j - E_{\nu^{n}_{u(t, \cdot)}}[\tau_x g_j]  \,\}\,.
\end{equation}
Here and throughout this section, $\cb{\mf R^{(p)}_{t,n}}$, $p\ge 1$,
is a remainder such that
\begin{equation}
\label{106d}
\sup_{n\ge 1}  \, \sup_{0\le t\le T} \,
\sup_{\eta\in\Omega_n} |\mf R^{(p)}_{t,n}| \,\le\, \mf c_{1,p}\,.
\end{equation}

To prove \eqref{106}, observe that
$| \, U(t,[x+e_j]/n) - U(t,x/n) \,| \le \mf r^{-2} \, \Vert \nabla
u\Vert_\infty /n= \mf c_{1,1}/n$. Perform a Taylor expansion up to the
third order and note that
$[\eta_{x+e_j} - \eta_x]^3 = \eta_{x+e_j} - \eta_x$ to obtain that
\begin{align}
\label{105}
L_{n,t}^{*} \, \bs 1 \, =\,
& -\, n^2\, \sum_{j=1}^d \sum_{x\in\bb T^d_n}
[U(t,[x+e_j]/n) - U(t,x/n)] \, 
c_{x,x+e_j}(\eta)\,  [\eta_{x+e_j} - \eta_x] 
\\
\nonumber
&\quad +\, n^2\,  \frac{1}{2}\,
\sum_{j=1}^d \sum_{x\in\bb T^d_n}
[U(t,[x+e_j]/n) - U(t,x/n)]^2
c_{x,x+e_j}(\eta)\,   [\eta_{x+e_j} - \eta_x]^2
\\
&\quad -\, n^2\,  \frac{1}{6}\,
\sum_{j=1}^d \sum_{x\in\bb T^d_n}
[U(t,[x+e_j]/n) - U(t,x/n)]^3
c_{x,x+e_j}(\eta)\,  [\eta_{x+e_j} - \eta_x] \,+\, \mf R^{(1)}_{t,n} \, n^{d-2}\,.
\nonumber
\end{align}
By \eqref{11b} and a summation by parts, the first line is equal to
\begin{align*}
\sum_{j=1}^d \sum_{x\in\bb T^d_n}
(\Delta^n_{j} U) (t,x/n) \, \tau_x h_j   \;,
\end{align*}
where the discrete differential operator $\Delta^n_j$ has been
introduced in \eqref{17}. Since the function $U(t, \cdot)$ is of class
$C^4$, a Taylor expansion yields that this expression is equal to
\begin{align*}
\sum_{j=1}^d \sum_{x\in\bb T^d_n}
(\partial^2_{x_j} U)  (t,x/n) \, \tau_x h_j    \,+\, \mf
R^{(4)}_{t,n} \, n^{d-2} \;.
\end{align*}

Similarly, a Taylor expansion yields that the second line in
\eqref{105} is equal to
\begin{equation*}
\frac{1}{2}\,
\sum_{j=1}^d \sum_{x\in\bb T^d_n}
\big[ (\partial_{x_j} U)  (t,x/n) \big]^2\, \tau_x g_j
\,+\,
\frac{1}{2n}\,
\sum_{j=1}^d \sum_{x\in\bb T^d_n}
(\partial_{x_j} U)  (t,x/n) \,
(\partial^2_{x_j} U)  (t,x/n)  \, \tau_x g_j
\,+\, \mf R^{(3)}_{t,n}\, n^{d-2} \;,
\end{equation*}
where the cylinder function $g_j$ has been introduced just before
equation \eqref{126}.  We turn to the third term in \eqref{105}. By
\eqref{11b}, we may perform a summation by parts and conclude that it
is bounded by $ \mf R^{(2)}_{t,n} \, n^{d-2}$.

Up to this point, we have proved that
\begin{equation}
\label{106b}
\begin{aligned}
L_{n,t}^{*} \, \bs 1 \,  &=\,
\sum_{j=1}^d \sum_{x\in\bb T^d_n}
(\partial^2_{x_j} U)  (t,x/n) \, \tau_x h_j 
\, +\, \frac{1}{2}\,
\sum_{j=1}^d \sum_{x\in\bb T^d_n}
\big[ (\partial_{x_j} U)  (t,x/n) \big]^2\, \tau_x g_j
\\ & +\,
\frac{1}{2n}\,
\sum_{j=1}^d \sum_{x\in\bb T^d_n}
(\partial_{x_j} U)  (t,x/n) \,
(\partial^2_{x_j} U)  (t,x/n)  \, \tau_x g_j  \,+\, 
\mf R^{(4)}_{t,n} \, n^{d-2}\;.
\end{aligned}
\end{equation}

As $L^*_{n,t}$ is the adjoint of $L_n$ in $L^2(\nu^{n}_{u(t, \cdot)})$,
if we represent by $\< \cdot \,, \, \cdot \>_{\nu^{n}_{u(t, \cdot)}}$
the scalar product in $L^2(\nu^{n}_{u(t, \cdot)})$,
$E_{\nu^{n}_{u(t, \cdot)}}[L_{n,t}^{*} \, \bs 1] = \< \bs 1\,,
L_{n,t}^{*} \, \bs 1\>_{\nu^{n}_{u(t, \cdot)}} = \< L_n \bs 1\,, \,
\bs 1\>_{\nu^{n}_{u(t, \cdot)}} =0$. Therefore, in \eqref{106b} we may
subtract $E_{\nu^{n}_{u(t, \cdot)}}[L_{n,t}^{*} \, \bs 1]$, that is,
the expectation with respect to $\nu^{n}_{u(t, \cdot)}$ of the
cylinder functions.  This completes the proof of the claim
\eqref{106}.

A last remark. Recall from \eqref{106c} the definition of $I_n(t)$.
By Lemma \ref{s47} with $\gamma=1$,
\begin{equation}
\label{108}
\int |\, I_n(t)\,|  \, f \, d\nu^n_{u (t, \cdot)}  \,\le\,
\, H_n(f \,|\, \nu^n_{u (t,\cdot)}  )
\,+\, \ln (2) \,+\, \mf c_{1,2} \, n^{d-2}
\end{equation}
for all densities $f$ with respect to the measure
$\nu^n_{u (t,\cdot)} $ and all $0\le t\le T$.  

As $u$ is the solution to equation \eqref{127}, by \eqref{19} and
\eqref{127},
\begin{equation}
\label{107}
\begin{aligned}
\partial_t   \ln \psi^n_t \, &=\,
\sum_{j=1}^d \sum_{x\in\bb T^d_n}
\frac{  (\partial^2_{x_j} u)(t,x/n) }    {\chi(u(t,x/n))}\,
\mtt p_{h_j}'  (u (t,x/n))  
\big[\, \eta_x - u(t,x/n)\,\big]
\\
\, & +\,\sum_{j=1}^d \sum_{x\in\bb T^d_n}
\Big(\frac{ (\partial_{x_j} u)(t,x/n)}    {\chi(u(t,x/n))} \Big)^2\,
\mtt p_{h_j}''  (u (t,x/n))  \chi(u(t,x/n))
\big[\, \eta_x - u(t,x/n)\,\big]
\;.
\end{aligned}
\end{equation}

Recall the definition of $\Pi_{\rho, g} (x,\eta)$ introduced in
\eqref{179} and that $\cb{\mf c_{1,p}}$, $p\ge 0$, represents a finite
constant, whose value may change from line to line, and which depends
only on $\mf r$, $\sup_{0\le t\le T} \Vert u(t, \cdot)\Vert_{C^p}$ and
the model. When the constant also depends on another variable, say a
cylinder function $h$, we stress the dependence explicitly by writing
$\mf c_{1,p} = \mf c_{1,p} (h)$.

\begin{lemma}
\label{s41}
Recall that $u$ is the solution of the parabolic equation
\eqref{127}. Then,
\begin{equation*}
L_{n,t}^{*} \, \bs 1  \,-\, \partial_t   \ln \psi^n_t  \, =\,
\mf W_t \, +\, I_n (t)\,+\,  \mf R^{(4)}_{t,n} \, n^{d-2} \;,
\end{equation*}
where
\begin{equation*}
\begin{aligned}
\mf W_t \, & =\,
\sum_{j=1}^d \sum_{x\in\bb T^d_n}
\frac{(\partial^2_{x_j} u)  (t,x/n)}
{\chi(u  (t,x/n)) }\,  \Pi_{u  (t,\cdot) , h_j} (x, \eta)
\\
& +\, \frac{1}{2}\, 
\sum_{j=1}^d \sum_{x\in\bb T^d_n}
\Big(\, \frac{(\partial_{x_j} u)  (t,x/n)}
{\chi(u  (t,x/n)) }\Big)^2
\, \Pi_{u  (t,\cdot) , g_j} (x, \eta) 
\\
& -\, \sum_{j=1}^d \sum_{x\in\bb T^d_n}
\Big(\, \frac{(\partial_{x_j} u)  (t,x/n)}
{\chi(u  (t,x/n)) }\Big)^2 \, [1-2 u  (t,x/n)] \,
\, \Pi_{u  (t,\cdot) , h_j} (x, \eta)  \,,
\end{aligned}
\end{equation*}
$I_n(t)$ is given by \eqref{106c} and satisfies the
bound \eqref{108}, and $\mf R^{(4)}_{t,n}$ is a remainder
satisfying the bound \eqref{106d}.
\end{lemma}

\begin{proof}
Recall equation \eqref{106} for $L_{n,t}^{*} \, \bs 1$, and the
estimate \eqref{108} for $I_n(t)$. We turn to the first two sums in
\eqref{106}. Since
$\partial^2_{x_j} U = \chi(u)^{-1}\, \partial^2_{x_j} u -
(1-2u)\, \chi(u)^{-2}\, [\partial_{x_j} u]^2$, and
$\partial_{x_j} U = \chi(u)^{-1}\, \partial_{x_j} u$, equation
\eqref{106} can be rewritten as 
\begin{align}
\label{129}
L_{n,t}^{*} \, \bs 1 \,  &=\,
\sum_{j=1}^d \sum_{x\in\bb T^d_n}
\frac{(\partial^2_{x_j} u)  (t,x/n)}
{\chi(u  (t,x/n)) }\,
\, \{\, \tau_x h_j -
E_{\nu^{n}_{u(t, \cdot)}}[\tau_x h_j ]  \,\}
\\
& +\, \frac{1}{2}\, 
\sum_{j=1}^d \sum_{x\in\bb T^d_n}
\Big(\, \frac{(\partial_{x_j} u)  (t,x/n)}
{\chi(u  (t,x/n)) }\Big)^2
\{\,  \tau_x g_j - E_{\nu^{n}_{u(t, \cdot)}}[\tau_x g_j]
\,\}
\nonumber
\\
& -\, \sum_{j=1}^d \sum_{x\in\bb T^d_n}
\Big(\, \frac{(\partial_{x_j} u)  (t,x/n)}
{\chi(u  (t,x/n)) }\Big)^2 \, [1-2 u  (t,x/n)] \,
 \{\, \tau_x h_j - E_{\nu^{n}_{u(t, \cdot)}}[\tau_x h_j ] \,\}
\, +\, I_n (t)\,+\,  \mf R^{(4)}_{t,n} \, n^{d-2}\;,
\nonumber
\end{align}
where $\mf R^{(4)}_{t,n}$ satisfies the bound stated in the lemma.

The first line of the previous expression combined with the first line
of \eqref{107} yields $\Pi_{u(t, \cdot), h_j} (x, \eta)$. By
\eqref{125}, $\mtt p_{h_j}'(\alpha) = D_j(\alpha)$, so that
$\mtt p_{h_j}''(\alpha) = D'_j(\alpha)$. By \eqref{126},
$(1/2)\, \mtt p_{g_j}' (\alpha) - (1-2\alpha) D_j(\alpha) =
\chi(\alpha) D_j' (\alpha)$. Therefore, the difference of the second
and third lines of the previous formula combined with the second line
of \eqref{107} yields the second and third lines of the statement of
the lemma. This completes the proof of the lemma.
\end{proof}

\begin{remark}
The expression $L_{n,t}^{*} \, \bs 1 \,-\, \partial_t \ln \psi^n_t$
reveals the equation that $u$ has to satisfy. Indeed, choose a product
reference measure associated with a time-dependent density profile
$u(t,x)$, and compute $L_{n,t}^{*} \, \bs 1$ as in equation
\eqref{106}, expressing it in terms of the local functions $h_j$,
$g_j$ and the spatial derivatives of $u$. Next, compute
$\partial_t \ln \psi^n_t$, which is expressed in terms of the time
derivative of $u$ and the linear function $\eta_x - u(t,x/n)$. It then
remains to determine the PDE that $u$ must satisfy for the difference
$L_{n,t}^{*} \, \bs 1 \,-\, \partial_t \ln \psi^n_t$ to become a
function of $\Pi_{u(t, \cdot), h} (x,\eta)$, $h$ a cylinder function.
\end{remark}

\section{Entropy and large deviations estimates}
\label{sec3}

In this section, we present some estimates needed throughout the
article. Here, $\rho\colon \bb T^d \to [0,1]$ is a smooth density
profile bounded away from $0$ and $1$: there exists
$0< \cb{\mf r } <1/2$ such that
\begin{equation}
\label{83}
\mf r  \, \le\,  \rho(\mtt x ) \,\le\,  1-\mf r\quad
\text{for all $\mtt x \in \bb T^d$}\,.
\end{equation}

Throughout the article, for $r\ge 1$, $\Sigma^+_{r} $, $\Sigma_{r} $
represent the subsets of $\bb Z^d$ given by
\begin{equation}
\label{190}
\cb{ \Sigma^+_{r} } \,:=\, \{0, \dots, r-1\}^d\,, \quad
\cb{\Sigma_{r} } \,:=\, \{-r, \dots, r\}^d\,.
\end{equation}
Denote the support of a local function
$h\colon \{0,1\}^{\bb Z^d} \to \bb R$ by
$\cb{\Lambda_h} \subset \bb Z^d$. This means that $h$ depends on the
configuration $\eta$ only through the variables $\eta_x$,
$x\in \Lambda_h$.  Denote by $\cb{\ell_h}\ge 0$ the smallest non-negative
integer $\ell$ such that $\Lambda_h \subset \Sigma_\ell$.

\subsection*{Entropy estimates}

In this subsection we estimate expressions of the form
\begin{align*}
\int  \Big  | \, \sum_{x \in \bb T^d_n} \mf a (x) \,
\big( \, \tau_x  h - E_{\nu^n_{\rho (\cdot)}} [\tau_x h] \, \big)\,\Big |
\, f\; d\nu^n_{\rho (\cdot)} \,,
\end{align*}
when the function $\mf a\colon \bb T^d_n \to\bb R$ is bounded.  The
proofs rely entirely on the entropy inequality and the correlation
decay of the reference product measure $\nu^n_{\rho (\cdot)} $.

Below, sums are performed modulo $n$.  Hence, if $h$ is a cylinder
function and $x\in \bb T^d_n$, the support of $\tau_x h$ may not be
contained in $\{0, \dots, n-1\}^d$, but it is contained in $\bb T^d_n$
due to the periodic boundary conditions assumed on $\bb T^d_n$.

\begin{lemma}
\label{s47}
Fix a local function $h\colon \{0,1\}^{\bb Z^d} \to \bb R$, a
function $\mf a \colon \bb T^d_n \to \bb R$ and let
$\Vert \mf a \Vert_\infty := \max_{x\in \bb T^d_n} |\mf a(x)|$.
Then,
\begin{equation*}
\int  \Big  | \, \sum_{x \in \bb T^d_n} \mf a (x) \,
\big( \, \tau_x  h - E_{\nu^n_{\rho (\cdot)}} [\tau_x h] \, \big)\,\Big |
\, f\; d\nu^n_{\rho (\cdot)} 
\,\le \, \frac{1}{\gamma}\, \Big\{\, H_n(f \,|\, \nu^n_{\rho (\cdot)}
)  \,+\, \ln(2) \, \Big\} \,+ \, c_1\, \gamma\, e^{c_2 \, \gamma}\, n^d
\end{equation*}
for all $\gamma>0$, and density $f$ with respect to
$\nu^n_{\rho (\cdot)} $. In this formula,
\begin{equation*}
c_1 \,=\,  2\, (2\ell_h + 1)^d \, \Vert \mf a \Vert_\infty^2
\, \Vert h \Vert_\infty^2 \,, \quad
c_2\,=\,  2 \, (2\ell_h + 1)^d \, \Vert \mf a \Vert_\infty
\, \Vert h \Vert_\infty \, .
\end{equation*}
\end{lemma}

\begin{proof}
By the entropy inequality, the expression appearing on the left-hand
side of the statement of the lemma is bounded above by
\begin{equation*}
\frac{1}{\gamma}\, H_n(f \,|\, \nu^n_{\rho (\cdot)}  )  \,+\,
\frac{1}{\gamma}\, 
\ln \int  \exp\Big\{ \, \gamma\, \Big  | \, \sum_{x \in \bb T^d_n} \mf a (x) \,
\big( \, \tau_x  h   - E_{\nu^n_{\rho (\cdot)}} [\tau_x  h  ] \, \big)\,\Big
|\, \Big\}  \; d\nu^n_{\rho (\cdot)} 
\end{equation*}
for all $\gamma>0$.  Since $e^{|y|}\leq e^{y} + e^{-y}$, for any
bounded random variable $X$,
\begin{equation}
\label{154}
\ln \bb E \big [ \,e^{|X|} \, \big] \,\leq \,
\ln \bb E \big[ \, e^{X} + e^{-X} \big ] \,\leq \,
\ln \Big(\,  2 \, \max_{b=\pm 1} \,
\bb E \big[ \, e^{b X} \, \big ] \,\Big ) \,=\, 
\ln(2) \,+\,
\max_{b=\pm 1} \, \ln
\bb E \big[ \, e^{b X} \, \big ] \,.
\end{equation}
Therefore, the expression in the penultimate displayed equation is
bounded by
\begin{equation}
\label{155}
\frac{1}{\gamma}\, \Big\{\, H_n(f \,|\, \nu^n_{\rho (\cdot)}
)  \,+\, \ln(2) \, \Big\}  \,+\,
\frac{1}{\gamma}\, 
\max_{b=\pm 1} \,  \ln \int
\exp \Big\{  b\, \gamma\, \sum_{x \in \bb T^d_n}
\mf a (x) \,
\big( \, \tau_x  h   - E_{\nu^n_{\rho (\cdot)}} [\tau_x  h  ] \, \big)
\Big\}  \; d\nu^n_{\rho (\cdot)} \,.
\end{equation}

Recall the definition of $\ell_h$, introduced just before the
statement of the lemma, and write
\begin{equation*}
\sum_{x \in \bb T^d_n}
\mf a (x) \,
\big( \, \tau_x  h   - E_{\nu^n_{\rho (\cdot)}} [\tau_x  h  ] \, \big)
\,=\,
\sum_{y\in {\Sigma_{\ell_h}}} \sum_{z\in B_{n,y}}
\mf a (x(y,z)) \,
\big( \, \tau_{x(y,z)}  h   - E_{\nu^n_{\rho (\cdot)}} [\tau_{x(y,z)}
h  ] \, \big)\,. 
\end{equation*}
In this equation, $x(y,z) = y + (2\ell_h + 1) z$, and the second sum
is carried over all $z\in \bb Z^d$ such that
$y + (2\ell_h + 1) z \in \Sigma^+_n$.  By H\"older inequality, the
last term in the penultimate displayed equation is bounded by
\begin{equation*}
\max_{b=\pm 1} \,  \frac{1}{\gamma \, (2\ell_h + 1)^d}
\sum_{y\in {\Sigma_{\ell_h}}}  \ln \int
\exp \Big\{  b\, \gamma\, (2\ell_h + 1)^d\, \sum_{z\in B_{n,y}} 
\mf a (x(y,z)) \,
\big( \, \tau_{x(y,z)}  h   - E_{\nu^n_{\rho (\cdot)}} [\tau_{x(y,z)} h  ] \, \big)
\Big\}  \; d\nu^n_{\rho (\cdot)} \,.
\end{equation*}
Clearly, if $z\neq z'$, then $x(y,z)-x(y,z')$ does not belong to
$\Sigma_{2 \ell_h}$, which implies that
$\{\tau_{x(y,z)} h: z\in B_{n,y}\}$ are independent random variables
with respect to any product measure. The previous expression is thus
equal to
\begin{equation*}
\max_{b=\pm 1} \,  \frac{1}{\gamma \, (2\ell_h + 1)^d}
\sum_{x\in \bb T^d_n}  \ln \int
\exp \Big\{  b\, \gamma\, (2\ell_h + 1)^d\, 
\mf a (x) \,
\big( \, \tau_x  h   - E_{\nu^n_{\rho (\cdot)}} [\tau_x  h  ] \, \big)
\Big\}  \; d\nu^n_{\rho (\cdot)} \,.
\end{equation*}
Since $e^r \le 1 + r + (1/2) \, r^2 \, e^{|r|} $ for all $r\in \bb R$,
$\tau_x h - E_{\nu^n_{\rho (\cdot)}} [\tau_x h ] $ has mean zero with
respect to $\nu^n_{\rho (\cdot)}$, and $\ln (1+c) \le c$ for $c>0$,
the previous expression is bounded by
\begin{equation*}
\frac{\gamma}{2}\, e^{2\, \gamma\, (2\ell_h + 1)^d
\Vert \mf a \Vert_\infty \, \Vert h \Vert_\infty }\,
\sum_{x \in \bb T^d_n} (2\ell_h + 1)^d \, \mf a (x)^2\, 
\int \big( \, \tau_x  h   - E_{\nu^n_{\rho (\cdot)}} [\tau_x  h  ] \, \big)^2
\; d\nu^n_{\rho (\cdot)} \,.
\end{equation*}
To complete the proof of the lemma, it remains to estimate
$|\, \tau_x h - E_{\nu^n_{\rho (\cdot)}} [\tau_x h ] \,|$ by
$ 2 \Vert h\Vert_\infty$, and to recollect all previous estimates.
\end{proof}

Assume that $n = r_1\, r_2$ for some positive integers $r_1\ge 2$,
$r_2\ge 2$. Note that $r_1$ and $r_2$ depend on $n$, though it does
not appear in the notation. Assume, furthermore, that $r_1$ is odd, so
that $\cb{r_1 = 2b_1+1}$. With this hypothesis we have the convenient
relation $\Sigma_{b_1} = -\, (b_1, \dots, b_1) + \Sigma^+_{r_1}$.
Divide the torus $\bb T^d_n$ into $r_2^d$ disjoint cubes of side length
$r_1$, represented by $\cb{B_{\mtt k}}$,
$\mtt k =(k_1, \dots, k_d) \in \Sigma^+_{r_2}$.  Thus,
$B_{\mtt 0} = \{0, \dots, r_1-1\}^d$,
$B_{\mtt k} = r_1\, \mtt k + B_{\mtt 0} $. Consider the sets
$B_{\mtt k} $ as subsets of $\bb T^d_n$, and write
$\bb T^d_n = \cup_{\mtt k} B_{\mtt k}$.

Denote by $J(B)$ the average over a cube $B$ of a continuous function
$J\colon \bb T^d\to \bb R$, and by $\mss m (B)$ the density of
particles in $B$:
\begin{equation}
\label{156}
\cb{J (B) }  \,:=\, \frac{1}{ |B| }\,
\sum_{x\in B} J(x/n)\,,
\quad
\cb{\mss M (B) }  \,:=\,  \sum_{x\in B} \eta (x)\, ,
\quad
\cb{\mss m (B) }  \,:=\, \frac{1}{ |B| }\, \mss M(B) \, ,
\quad B \subset  \bb T^d_n\,.
\end{equation}

The next result shows that, in the nonequilibrium density field
fluctuations problem, one may always replace
\begin{equation*}
\sum_{x \in \bb T^d_n} J(x/n) \,
[\eta(x) - \rho(x/n)]
\quad \text{by}\quad
\sum_{x\in \bb T^d_n} J (x+ \Sigma^+_{r_1})
\big\{\, \mss m (x+ \Sigma^+_{r_1})  
-  E_{\nu^n_{\rho(\cdot)}}[ \mss m (x+ \Sigma^+_{r_1})  ]  \,\big\}\,,
\end{equation*}
provided $J$ is of class $C^2(\bb T^d)$ and $r_1 \ll
n^{1-(d/8)}$. This result simplifies the analysis in the subsequent
sections. Here, recall, sums are performed modulo $n$ so that $x+
\Sigma^+_{r_1}$ is a subset of $\bb T^d_n$.

The goal of Sections~\ref{sec4} and onwards is to replace
\begin{gather*}
\sum_{x \in \bb T^d_n} J(x/n) \, \big\{ \tau_x h -  E_{\nu^n_{\rho
(\cdot)}} [\tau_x  h  ]\, \big\}
\quad \text{by}\quad
\sum_{x\in \bb T^d_n} J (x+ \Sigma^+_{r_1})
\Big\{\, \mtt p_{h} \big( \mss m (x+ \Sigma^+_{r_1})  \big)
-  E_{\nu^n_{\rho(\cdot)}} \big[
\mtt p_{h} \big(  \mss m (x+ \Sigma^+_{r_1}) \big) \, \big]
\,\Big\}\,, 
\end{gather*}
where $h$ is a cylinder function and $\mtt p_{h} (\cdot)$ is as
introduced in \eqref{05}. Lemma~\ref{s51b} handles the degree-$1$ part
of the cylinder function $h$, in the sense of
Appendix~\ref{sec8}. Consequently, in Section \ref{sec4}, attention may be
restricted to cylinder functions of degree two or higher.

For a function $J\colon \bb T^d \to \bb R$, of class $C^p(\bb T^d)$,
denote by $\Vert J \Vert _{C^p}$ the norm defined by
\begin{equation*}
\cb{\Vert  J \Vert _{C^p}  } \,:=\,
\sum_{\mtt k \le p} \Vert  \, \partial^{k_1}_{x_1}  \cdots
\partial^{k_d}_{x_d} J \, \Vert _{\infty} \,,
\end{equation*}
where the sum is performed over all vectors
$\mtt k = (k_1, \dots, k_d)$ such that $k_i \ge 0$,
$\sum_i k_i \le p$.

\begin{lemma}
\label{s51b}
Fix a function $J\colon \bb T^d \to \bb R$ of class
$C^2(\bb T^d)$, and let
\begin{equation*}
W_J \,:=\,
\sum_{x\in \bb T^d_n}  J(x/n)\, \big[\, \eta(x) - \rho (x/n) \,\big] \,
- \,
\sum_{x\in \bb T^d_n} J (x+ \Sigma^+_{r_1})
\big\{\, \mss m (x+ \Sigma^+_{r_1})  
-  E_{\nu^n_{\rho(\cdot)}}[ \mss m (x+ \Sigma^+_{r_1})  ]  \,\big\}\,.
\end{equation*}
Then, 
\begin{align*}
& \int \big| \, W_J  \,\big | \, 
f\, d\nu^n_{\rho(\cdot)}
\,\le \,
H_n(f \,|\, \nu^n_{\rho (\cdot)}  )  \,+\, \ln(2) \,+\,
\Vert J \Vert^2_{C^2(\bb T^d)} \,
e^{\Vert J \Vert_{C^2(\bb T^d)} \, (r_1/n)^2} \,
\Big( \frac{r_1}{n}\Big)^4  \, n^d 
\end{align*}
for all density $f$ with respect to $\nu^n_{\rho(\cdot)}$. 
\end{lemma}

\begin{proof}
Fix a function $J\colon \bb T^d \to \bb R$ of class $C^2(\bb T^d)$,
and a density $f$ with respect to $\nu^n_{\rho(\cdot)}$. Recall that
$r_1 = 2b_1 +1$, and let
\begin{equation*}
\cb{ J^{(2)} (x+ \Sigma_{b_1}) } :=
\frac{1}{|\Sigma_{b_1}|}\,
\sum_{y\in \Sigma_{b_1}} J(x - y + \Sigma_{b_1})\,,
\end{equation*}
where $J(x+ \Sigma_{b_1})$ has been introduced in \eqref{156}.  As $J$
is of class $C^2(\bb T^d)$, the absolute value of the difference
$J(x/n)\, -\, J^{(2)}(x+ \Sigma_{b_1}) $ is bounded by
$2\, \Vert J \Vert_{C^2(\bb T^d)} \, (b_1/n)^2$ uniformly in
$x$. Therefore, by Lemma \ref{s47} with $\gamma=1$, $\ell_h=0$,
$\Vert \mf a \Vert_\infty = 2\, \Vert J \Vert_{C^2(\bb T^d)} \,
(b_1/n)^2$,
\begin{align*}
& \int \Big | \, \sum_{x\in \bb T^d_n} 
\big[ \, J(x/n)\, -\,  J^{(2)} (x+ \Sigma_{b_1}) \, \big]
\big[\, \eta(x) - \rho(x/n) \,\big] \,\Big| \, 
f\, d\nu^n_{\rho(\cdot)}
\\
&\quad
\,\le \,  H_n(f \,|\, \nu^n_{\rho (\cdot)}  )  \,+\, \ln(2)\,+\,
C_1(J)\, \Big( \frac{r_1}{n}\Big)^4  \, n^d\, e^{C_2 (J) \, (r_1/n)^2} \,,
\end{align*}
where $C_1(J) = \Vert J \Vert^2_{C^2(\bb T^d)}$,
$C_2(J) = \Vert J \Vert_{C^2(\bb T^d)} $. Here we used the fact
that $b_1 \le r_1/2$.  A summation by parts yields that
\begin{equation*}
\sum_{x\in \bb T^d_n} 
J^{(2)} (x+ \Sigma_{b_1}) \, \big[\, \eta(x) - \rho(x/n) \,\big]
\,=\, \sum_{x\in \bb T^d_n} J (x+ \Sigma_{b_1})
\big\{\, \mss m (x+ \Sigma_{b_1})  
-  E_{\nu^n_{\rho(\cdot)}}[ \mss m (x+ \Sigma_{b_1})  ]  \,\big\} \,, 
\end{equation*}
where $\mss m (x+ B)$ has been introduced in \eqref{156}.  As
$\Sigma_{b_1} = -\, (b_1, \dots, b_1) + \Sigma^+_{r_1}$, this last sum
is equal to
\begin{align*}
\sum_{x\in \bb T^d_n} J (x+ \Sigma^+_{r_1})
\big\{\, \mss m (x+ \Sigma^+_{r_1})  
-  E_{\nu^n_{\rho(\cdot)}}[ \mss m (x+ \Sigma^+_{r_1})  ]  \,\big\} \,.
\end{align*}
This completes the proof of the lemma.
\end{proof}

The next result asserts that if $J$ is the product of two functions,
in the previous lemma we may replace the average of $J$ by the product
of the averages.

\begin{corollary}
\label{s52b}
Fix two functions $F$, $G\colon \bb T^d \to \bb R$ of class
$C^2(\bb T^d)$. Let
\begin{align*}
W_{F,G} \,=\,
& \sum_{x\in \bb T^d_n} 
F(x/n)\, G(x/n)\,  \big[\, \eta(x) - \rho( x/n) \,\big]
\\
& \,-\,
\sum_{x\in \bb T^d_n} F (x+ \Sigma^+_{r_1})\, G (x+ \Sigma^+_{r_1})\, 
\big\{\, \mss m (x+ \Sigma^+_{r_1})  
-  E_{\nu^n_{\rho(\cdot)}}[ \mss m (x+ \Sigma^+_{r_1})  ]  \,\big\}  \,.
\end{align*}
Then,
\begin{equation*}
\int \big |\,  W_{F,G}  \,\big | \,  f\, d\nu^n_{\rho( \cdot)}
\,\le\,
2\, H_n(f \,|\, \nu^n_{\rho (\cdot)}  )  \,+\, 2\, \ln(2) \,+\,
C_1\, \Big( \frac{r_1}{n}\Big)^4  \, n^d
\end{equation*}
for all density $f$ with respect to $\nu^n_{\rho( \cdot)}$.
In this formula,
\begin{align*}
C_1 \,:=\, \Vert F\, G \Vert^2_{C^2(\bb T^d)} \,
e^{\Vert  F\, G \Vert_{C^2(\bb T^d)} \, (r_1/n)^2}
\,+\, 2\,   \Vert F\Vert^2_{C^1(\bb T^d)} \,
\Vert G\Vert^2_{C^1(\bb T^d)} 
\, \exp\{ 2\,  \Vert F\Vert_{C^1(\bb T^d)} \,
\Vert G\Vert_{C^1(\bb T^d)} \, (r_1/n)^2 \}\,.
\end{align*}
\end{corollary}

\begin{proof}
Fix two functions $F$, $G$ in $C^2(\bb T^d)$, and let $J = F\, G$.  We
use the triangle inequality to estimate $W_{F,G} $ summing and
subtracting the term
\begin{equation*}
\widetilde{W}_{F,G} \,:=\,
\sum_{x\in  \bb T^d_n} J (x+ \Sigma^+_{r_1})\,
\big\{\, \mss m (x+ \Sigma^+_{r_1})  
-  E_{\nu^n_{\rho(\cdot)}}[ \mss m (x+ \Sigma^+_{r_1})  ]  \,\big\}  \,.
\end{equation*}
The previous lemma provides an estimate for the difference between the
first term on the right-hand side of the definition of $W_{F,G} $ and
$\widetilde{W}_{F,G}$.

The difference between $\widetilde{W}_{F,G}$ and the second term on
the right-hand side of the definition of $W_{F,G} $ is equal to
\begin{equation}
\label{157}
\sum_{x\in \bb T^d_n} H (x+ \Sigma^+_{r_1})\,
\big\{\, \mss m (x+ \Sigma^+_{r_1})  
-  E_{\nu^n_{\rho(\cdot)}}[ \mss m (x+ \Sigma^+_{r_1})  ]  \,\big\}
\,,
\end{equation}
where
\begin{align*}
H(B)  & =\, \frac{1}{| B|} \,  \sum_{x\in B}
F(x/n) \, G(x/n) \,-\, F(B) \, G(B)
\\
\, &=\, \frac{1}{| B|} \,  \sum_{x\in B}
\big\{\, F(x/n) - F(B)\big\}\,
\big\{\, G(x/n) - G(B)\big\} \,.
\end{align*}
As $F$, $G$ are Lipschitz-continuous,
$|H( B)|\le \Vert F\Vert_{C^1(\bb T^d)} \, \Vert
G\Vert_{C^1(\bb T^d)}\, (r_1/n)^2$.

Write \eqref{157} as
\begin{align*}
& \sum_{x\in \bb T^d_n} H (x+ \Sigma^+_{r_1})\,
\frac{1}{| \Sigma^+_{r_1} |} \,
\sum_{y\in \Sigma^+_{r_1}}  \,[\, \eta_{x+y} - \rho((x+y)/n)\,]
\\
& \quad \,=\,
\sum_{x\in \bb T^d_n} \,[\, \eta_{x} - \rho(x/n)\,] \,
\frac{1}{| \Sigma^+_{r_1} |} \,
\sum_{y\in \Sigma^+_{r_1}}  
H (x - y+ \Sigma^+_{r_1})\,
\,,
\end{align*}
where we performed the change of variables $x'=x+y$ and used the fact
that the torus $\bb T^d_n$ has periodic boundaries.  Let
\begin{align*}
\mtt a (x) \,=\,
\frac{1}{| \Sigma^+_{r_1} |} \,
\sum_{y\in \Sigma^+_{r_1}}  
H (x - y+ \Sigma^+_{r_1})
\quad\text{so that}\;\;
\max_{z\in \bb T^d_n} 
\big|\, \mtt a(z)  \, \big| \,\le\, 
\Vert H\Vert_\infty\,.
\end{align*}
In consequence, by Lemma \ref{s47} with $\gamma =1$, $\ell_h=0$,
$h(\eta)=\eta_0$,  $\mf a (\cdot)= \mtt a(\cdot)$,
\begin{equation*}
\int \Big|\, \sum_{x\in \bb T^d_n} H (x+ \Sigma^+_{r_1})\,
\big\{\, \mss m (x+ \Sigma^+_{r_1})  
-  E_{\nu^n_{\rho(\cdot)}}[ \mss m (x+ \Sigma^+_{r_1})  ]  \,\big\}
\, \Big | \, 
f\, d\nu^n_{\rho( \cdot)} \,\le \, 
H_n(f \,|\, \nu^n_{\rho (\cdot)}  )  \,+\, \ln(2) \,+\,
C_2 \, n^d
\end{equation*}
for all density $f$ with respect to $\nu^n_{\rho (\cdot)} $.  In this
formula, $C_2= 2 \, \Vert H\Vert^2_\infty \, \exp\{ 2\, 
\Vert H\Vert_\infty\}$.  
To complete the proof of the corollary, it remains to recollect the
previous estimates.
\end{proof}

The next result will be used when $\mf a(\cdot)$ is small, that is,
when this function depends on $n$ and $\max_{\mtt i \in \Sigma^+_{r_2}} |\mf
a_n (\mtt i)| \to 0$. Let
\begin{equation}
\label{92}
\cb{\mss M_{\mtt k}} := \sum_{x\in B_{\mtt k}} \eta_x \,,\quad
\cb{\mss m_{\mtt k}} := \frac{\mss M_{\mtt k}} {|B_{\mtt k}| }\,,
\quad
\cb{\bar\rho _{\mtt k} } \,:=\,  \frac{1} { |B_{\mtt k}|} \,
\sum_{x\in B_{\mtt k} } \rho(x/n)  \,, \quad
\mtt k \in \Sigma^+_{r_2} \,.
\end{equation}

\begin{lemma}
\label{s32}
Fix a function $F$ in $C^1([0,1])$, and a function
$\mf a \colon \Sigma^+_{r_2} \to \bb R$.  Then,
\begin{equation*}
\int \Big | \sum_{\mtt i\in \Sigma^+_{r_2}} \mf a (\mtt i) \,
\big\{ \, F(\mss m_{\mtt i})
- E_{\nu^n_{\rho (\cdot)}} [F(\mss m_{\mtt i}) ] \, \big\}\, \Big|
\, f\; d\nu^n_{\rho (\cdot)} 
\,\le \,  \frac{1}{\gamma}\, \Big\{\,
H_n(f \,|\, \nu^n_{\rho (\cdot)}  )  \,+\, \ln (2) \, \Big\} \,+\,
c_0 \gamma\, e^{c_1 \gamma}\, 
\Big( \frac{n}{r_1^{2}} \Big)^d 
\end{equation*}
for all $\gamma>0$, and density $f$ with respect to
$\nu^n_{\rho (\cdot)} $. In this formula,
\begin{equation*}
c_0 \,=\,  \frac{1}{4}   \, 
\Vert \mf a\Vert_\infty^2 \, \Vert F'\Vert^2_\infty\,
\quad\text{and}\quad
c_1 \,=\, 2 \,  \Vert \mf a\Vert_\infty
\, \Vert F\Vert_\infty\,.
\end{equation*}
\end{lemma}

\begin{proof}
The proof is the same as the one of Lemma \ref{s47} up to equation
\eqref{155}. At that point we need to estimate the second term.  As
$\nu^n_{\rho (\cdot)}$ is a product measure, this term is equal to
\begin{equation*}
\frac{1}{\gamma}\, \sum_{\mtt i\in \Sigma^+_{r_2}}
\ln \int 
\exp\Big\{\, b\, \gamma\, \mf a (\mtt i) \,
\big[ \, F(\mss m_{\mtt i}) - E_{\nu^n_{\rho (\cdot)}} [F(\mss m_{\mtt i}) ] \,
\big]  \,\Big\}\, d\nu^n_{\rho (\cdot)}\,.
\end{equation*}
The expression inside braces is bounded by
$2 \, \gamma\, \Vert \mf a\Vert_\infty \, \Vert F\Vert_\infty$. Since
$e^x \le 1 + x + (1/2) x^2 e^{|x|} $ for all $x\in\bb R$, the previous
expression is bounded by
\begin{equation*}
\frac{1}{\gamma}\,  \sum_{\mtt i\in \Sigma^+_{r_2}} \ln \Big\{\, 1
\,+\, \frac{1}{2}   \, \gamma^2\,
\Vert \mf a\Vert_\infty^2 \,
e^{2 \, \gamma\,
\Vert \mf a\Vert_\infty  \, \Vert F\Vert_\infty}\,
\int \big\{ \, F(\mss m_{\mtt i})
- E_{\nu^n_{\rho (\cdot)}} [F(\mss m_{\mtt i}) ] \,\big\}^2
\, d\nu^n_{\rho (\cdot)}\,\Big\}
\end{equation*}
because
$F(\mss m_{\mtt i}) - E_{\nu^n_{\rho (\cdot)}} [F(\mss m_{\mtt i}) ] $
has zero mean. As $\ln (1+a) \le a$, the previous expression is less
than or equal to
\begin{equation*}
\frac{1}{2}   \, \gamma\,
\Vert \mf a\Vert_\infty^2 \,
e^{2 \, \gamma\,
\Vert \mf a\Vert_\infty  \, \Vert F\Vert_\infty}\,
\sum_{\mtt i\in \Sigma^+_{r_2}} 
\int \big\{ \, F(\mss m_{\mtt i})
- E_{\nu^n_{\rho (\cdot)}} [F(\mss m_{\mtt i}) ] \,\big\}^2
\, d\nu^n_{\rho (\cdot)} \,.
\end{equation*}
Introduce $F(\mss m_{\mtt i})$ inside the expectation and apply
Schwarz inequality.  Since $F(\cdot)$ is Lipschitz continuous, the
previous expression is bounded by
\begin{equation*}
\frac{1}{2}   \, \gamma\,
\Vert \mf a\Vert_\infty^2 \, \Vert F'\Vert^2_\infty\,  
e^{2 \, \gamma\,
\Vert \mf a\Vert_\infty  \, \Vert F\Vert_\infty}\,
\sum_{\mtt i\in \Sigma^+_{r_2}} \int \,\int 
\Big\{ \, \frac{1}{|B_{\mtt i}|} \sum_{x\in B_{\mtt i}} \eta_x
\,-\,
\frac{1}{|B_{\mtt i}|} \sum_{x\in B_{\mtt i}} \xi_x
\,\Big\}^2
\, d\nu^n_{\rho (\cdot)} (\eta)  \, d\nu^n_{\rho (\cdot)} (\xi)\,.
\end{equation*}
The double integral is twice the variance of
$|B_{\mtt i}|^{-1}\, \sum_{x\in B_{\mtt i}} \eta_x$. To complete the
proof, recall that $\chi(\rho) = \rho (1-\rho) \le 1/4$.
\end{proof}

\subsection*{Entropy and concentration inequalities}

In this subsection, we estimate expressions of the form
\begin{equation*}
\int \sum_{\mtt i\in \Sigma^+_{r_2}} \mf a (\mtt i) \,
\big\{ \, F(\mss m_{\mtt i})
- E_{\nu^n_{\rho (\cdot)}} [F(\mss m_{\mtt i}) ] \, \big\}
\, f\; d\nu^n_{\rho (\cdot)} 
\end{equation*}
in situations in which, after the entropy inequality, we may not
estimate the expression inside the exponential by its $L^\infty$
norm. We have to use, instead, concentration inequalities to estimate
these terms.

The next result identifies the correct order of magnitude for both the
entropy and the size of the cubes required in the proof of the
nonequilibrium density field fluctuations. To see this, divide both
sides of the inequality in Lemma~\ref{l02} by $n^{d/2}$. For the
right-hand side to vanish as $n \to \infty$, the side length of the
cube must be at least of order $k\sqrt{n}$, where $k \to \infty$ after
$n \to \infty$, while the entropy must be at most $o(n^{d/2})$.

Consequently, in order to derive an entropy bound of order
$o(n^{d/2})$, it is also necessary to consider cubes of side length $k
\sqrt{n}$ with $k \gg 1$.

\begin{lemma}
\label{l02}
Fix $J\colon \Sigma^+_{r_2} \to \bb R$, and let
$\Vert J\Vert_\infty = \max_{\mtt i \in \Sigma^+_{r_2}} |J(\mtt i)|$.
Then,
\begin{equation*}
\int \sum_{\mtt i\in \Sigma^+_{r_2}} J(\mtt i) \,
r_1^d\, ( \mss m_{\mtt i} - \bar \rho_{\mtt i} )^2
\, f\; d\nu^n_{\rho (\cdot)} 
\,\le \, 8\, \|J\|_\infty \, \Big\{   \, H_n(f \,|\, \nu^n_{\rho (\cdot)}  )  \,+\,
\Big( \frac{n}{r_1}\Big)^d\, \Big\}
\end{equation*}
for all density $f$ with respect to $\nu^n_{\rho (\cdot)} $.
\end{lemma}

\begin{proof}
By the entropy inequality, the left-hand side of the displayed
equation appearing in the statement of the lemma is bounded above by 
\begin{equation*}
\frac{1}{\gamma}   \, H_n(f \,|\, \nu^n_{\rho (\cdot)}  )
\,+\,
\frac{1}{\gamma}  \ln \int \exp\Big\{\,
\gamma\, \sum_{\mtt i\in \Sigma^+_{r_2}} J(\mtt i) \,
r_1^d\, ( \mss m_{\mtt i} - \bar \rho_{\mtt i} )^2 \,\Big\}\,
d\nu^n_{\rho (\cdot)}
\end{equation*}
for all $\gamma>0$. As $\nu^n_{\rho (\cdot)}$ is a product measure,
the second term is equal to
\begin{equation*}
\frac{1}{\gamma}  \sum_{\mtt i\in \Sigma^+_{r_2}}
\ln \int 
e^{\gamma\, J(\mtt i) \,
r_1^d\, ( \mss m_{\mtt i} - \bar \rho_{\mtt i} )^2} \,
d\nu^n_{\rho (\cdot)} \,.
\end{equation*}
By Hoeffding's inequality \cite[Lemma 2.2]{BLM13},
$\sum_{x\in B_{\mtt i}}
[\eta_x - \rho(x/n)]$ is an $r_1^d$-sub-Gaussian random variable. Thus,
by \cite[Proposition B.1]{jm1},
\begin{equation*}
\int  e^{\gamma\, J(\mtt i) \,
r_1^d\, ( \mss m_{\mtt i} - \bar \rho_{\mtt i} )^2} \,
d\nu^n_{\rho (\cdot)} \,\le\, e^{8 \gamma \Vert J\Vert_\infty }
\end{equation*}
provided $\gamma \Vert J\Vert_\infty < 1/4$. To complete the proof of
the lemma, it remains to choose $\gamma = 1/ (8\Vert J\Vert_\infty)$.
\end{proof}

Recall that $u\colon [0,T] \times \bb T^d \to (0,1)$ is the solution
of the hydrodynamic equation. The next result is a simple consequence
of the previous lemma and the bound \eqref{255b} on $u(\cdot, \cdot)$.

\begin{corollary}
\label{s19}
Fix $T>0$ and let $(f^{(n)}_t : 0\le t\le T)_{n\ge 1}$ be a sequence
of densities with respect to the measures
$(\nu^{n}_{u(t)} : 0\le t\le T)$ such that
\begin{equation*}
\lim_{n\to \infty}  \int_0^T \frac{1}{n^{d/2}}\,
H_n(f^{(n)}_s | \nu^{n}_{u(s, \cdot)})
\, ds \, = 0\,.
\end{equation*}
Then,
\begin{equation*}
\lim_{n\to \infty}  \int_0^T  ds\,
\int \frac{1}{n^{d/2}} \sum_{\mtt i\in \Sigma^+_{r_2}} J_n(s,\mtt i) \,
r_1^d\, ( \mss m_{\mtt i} - \bar u(s, \mtt i) )^2
\, f^{(n)}_s \; d\nu^n_{u (s, \cdot)} \,=\, 0
\end{equation*}
for every uniformly bounded sequence
$J_n\colon [0,T]\times \Sigma^+_{r_2} \to \bb R$ {\rm
(}$\sup_{n\ge 1} \sup_{0\le s\le T} \max_{\mtt i\in \Sigma^+_{r_2} }
|J_n(s,\mtt i)| < +\infty${\rm )}, and sequence $r_1$ such that
$\sqrt{n}/r_1 \to 0$. In this equation,
$\bar u(s, \mtt i) = |B_{\mtt i}|^{-1} \sum_{x\in B_{\mtt i}}
u(s,x/n)$.
\end{corollary}

When $\mf a (\cdot)$ is large, the alternative bound presented below
might be useful.

\begin{lemma}
\label{s45}
Fix an integer $p\ge 1$, and a function
$\mf a \colon \Sigma^+_{r_2} \to \bb R$.  Suppose that
$\Vert \mf a \Vert_\infty \le \mf c_0 \, r_1^{d/2}$ for some finite
constant $\mf c_0$.  Then, there exists a finite constant $C_p$
depending only on $p$ such that
\begin{equation*}
\int  \Big|\, \sum_{\mtt i\in \Sigma^+_{r_2}} \mf a (\mtt i) \,
\big\{ \, \mss m_{\mtt i}^p
- E_{\nu^n_{\rho (\cdot)}} [\mss m_{\mtt i}^p ] \, \big\} \,\Big|
\, f\; d\nu^n_{\rho (\cdot)} 
\,\le \,
C_p\,  \Big\{\, \mf c_0\, H_n(f \,|\, \nu^n_{\rho (\cdot)}  )
\,+\, \, \Vert \mf a \Vert_\infty\,
\Big(\frac{n}{r_1^{3/2}} \Big)^d \,\Big\} 
\end{equation*}
for all density $f$ with respect to $\nu^n_{\rho (\cdot)} $.
\end{lemma}

\begin{proof}
Suppose that $p=1$. In this case, by Young's inequality, the
expression on the left-hand side of the statement of the lemma is
bounded above by
\begin{equation}
\label{113}
\frac{1}{2\gamma}\,
\sum_{\mtt i\in \Sigma^+_{r_2}} \mf a (\mtt i)^2 \, +\, 
\frac{\gamma}{2}\,
\int
\sum_{\mtt i\in \Sigma^+_{r_2}}
\big\{ \, \mss m_{\mtt i} - E_{\nu^n_{\rho (\cdot)}} [\mss m_{\mtt i} ]  \, \big\}^2
\, f\; d\nu^n_{\rho (\cdot)} 
\end{equation}
for every $\gamma>0$. The first term is bounded by
$\Vert \mf a\Vert_\infty^2 r_2^d / (2\gamma)$. By the entropy
inequality, the second term is bounded by
\begin{equation*}
\frac{1}{\beta}\,  H_n(f \,|\, \nu^n_{\rho (\cdot)}  )
\,+\,
\frac{1}{\beta}\,   \ln \int \exp\Big\{\,
\frac{\beta \,\gamma}{2} \, 
\sum_{\mtt i\in \Sigma^+_{r_2}}
\big\{ \, \mss m_{\mtt i} - E_{\nu^n_{\rho (\cdot)}} [\mss m_{\mtt i} ] \, \big\}^2
\,\Big\}\,
d\nu^n_{\rho (\cdot)}
\end{equation*}
for every $\beta>0$.
Since $\nu^n_{\rho (\cdot)}$ is a product measure and the cubes
$B_{\mtt i}$ are disjoint, this sum is equal to
\begin{equation*}
\frac{1}{\beta}\,  H_n(f \,|\, \nu^n_{\rho (\cdot)}  )
\,+\,
\frac{1}{\beta}\,   \sum_{\mtt i\in \Sigma^+_{r_2}}
\ln \int \exp\Big\{\,
\frac{\beta \,\gamma}{2} \, 
\big\{ \, \mss m_{\mtt i} - E_{\nu^n_{\rho (\cdot)}} [\mss m_{\mtt i} ] \, \big\}^2
\,\Big\}\,
d\nu^n_{\rho (\cdot)}\,.
\end{equation*}
By Hoeffding's inequality \cite[Lemma 2.2]{BLM13} (see \cite[Lemma
4.4]{jl23})
$\mss M_{\mtt i} - E_{\nu^n_{\rho (\cdot)}} [\mss M_{\mtt i} ]$ is
A-subgaussian, where $A= (1/4) |B_{\mtt i}|$. Therefore, by
\cite[Proposition B.1]{jm1}, the previous expression is bounded by
\begin{equation*}
\frac{1}{\beta}\,  H_n(f \,|\, \nu^n_{\rho (\cdot)}  )
\,+\,
\frac{1}{\beta}\,   \sum_{\mtt i\in \Sigma^+_{r_2}}
\frac{\beta \,\gamma}{|B_{\mtt i}|} \,
\,=\,
\frac{1}{\beta}\,  H_n(f \,|\, \nu^n_{\rho (\cdot)}  )
\,+\,
\frac{n^d}{|B_{\mtt i}|} \,
\frac{\gamma}{|B_{\mtt i}|} \, 
\end{equation*}
provided $\beta \,\gamma < 2\, |B_{\mtt i}|$.

Set $\gamma = \Vert \mf a \Vert_\infty\, r_1^{d/2}$ and
$\beta = 1/\mf c_0$, where $\mf c_0$ appeared in the hypotheses
of the lemma.  With these definitions,
$\beta \,\gamma = (1/\mf c_0)\, \Vert \mf a \Vert_\infty\, r_1^{d/2} \le
(1/\mf c_0)\, \mf c_0 \, r_1^{d} = |B_{\mtt i}| < 2\, |B_{\mtt i}|$.
Therefore, $\beta \,\gamma$ is less than $2\, |B_{\mtt i}|$, as
required. With these values the last displayed equation becomes
\begin{equation*}
\mf c_0 \,  H_n(f \,|\, \nu^n_{\rho (\cdot)}  )
\,+\, \Vert \mf a \Vert_\infty\,
\frac{n^d}{|B_{\mtt i}|^{3/2}} \, .
\end{equation*}
Adding the first term in \eqref{113} with
$\gamma =\Vert \mf a \Vert_\infty\, r_1^{d/2}$ yields that
\begin{equation}
\label{114}
\int  \Big|\, \sum_{\mtt i\in \Sigma^+_{r_2}} \mf a (\mtt i) \,
\big\{ \, \mss m_{\mtt i} - E_{\nu^n_{\rho (\cdot)}} [\mss m_{\mtt i} ] \, \big\} 
\, \Big| \, f\; d\nu^n_{\rho (\cdot)} 
\,\le \,
\mf c_0 \,  H_n(f \,|\, \nu^n_{\rho (\cdot)}  )
\,+\, 2\, \Vert \mf a \Vert_\infty\,
\frac{n^d}{|B_{\mtt i}|^{3/2}} \,,
\end{equation}
as claimed.

We turn to the case $p\ge 2$. As $\nu^n_{\rho (\cdot)} $ is a product
measure, $E_{\nu^n_{\rho (\cdot)}} [\mss m_{\mtt i}^p] = E_{\nu^n_{\rho
(\cdot)}} [\mss m_{\mtt i}]^p + \mf R_n$, where $|\mf R_n|\le C_p / |B_{\mtt
i}|$ for some finite constant $C_p$ depending only on $p$. The
integral we wish to bound is thus less than or equal to
\begin{equation*}
\int  \Big|\, \sum_{\mtt i\in \Sigma^+_{r_2}} \mf a (\mtt i) \,
\big\{ \, \mss m_{\mtt i}^p
- E_{\nu^n_{\rho (\cdot)}} [\mss m_{\mtt i} ] ^p \, \big\} \,\Big|
\, f\; d\nu^n_{\rho (\cdot)} 
\, +\, C_p \, \Vert \mf a \Vert_\infty\,
\frac{n^d}{|B_{\mtt i}|^{2}} \,\cdot
\end{equation*}
Let $f_p(x) = x^p$. Add and subtract $f_p'(\,E_{\nu^n_{\rho (\cdot)}}
[\mss m_{\mtt i} ]\,) \, \{\, \mss m_{\mtt i} - E_{\nu^n_{\rho (\cdot)}}
[\mss m_{\mtt i} ]\,\}$ inside the sum.  The linear term has been taken
care of in the first part of the proof. It is enough to remark that
$\sup_{0\le \rho\le 1} | f_p'(\rho)| \le C_p$ to obtain the bound
\eqref{114} with an extra multiplicative factor $C_p$.

On the other hand, by a Taylor expansion, the non-linear term is
bounded by
\begin{equation*}
C_p\,  \int  \sum_{\mtt i\in \Sigma^+_{r_2}}
|\, \mf a (\mtt i) \,| \,
\big\{ \, \mss m_{\mtt i} - E_{\nu^n_{\rho (\cdot)}} [\mss m_{\mtt i} ] \, \big\}^2
\, f\; d\nu^n_{\rho (\cdot)} \,.
\end{equation*}
Repeating exactly the same steps as in the first part of the proof
yields that this expression is less than or equal to
\begin{equation*}
C_p\,  \Big\{\, \frac{1}{\beta}\,  H_n(f \,|\, \nu^n_{\rho (\cdot)}  )
\,+\,
\frac{1}{\beta}\,   \sum_{\mtt i\in \Sigma^+_{r_2}}
\frac{2\, \beta \, |\mf a (\mtt i)|}{|B_{\mtt i}|} \, \Big\} 
\end{equation*}
provided $\beta \, \Vert \mf a \Vert_\infty < |B_{\mtt i}|$. Choose
$\beta = 1/\mf c_0$, and recollect all previous estimates to complete
the proof.
\end{proof}

\begin{remark}
It is the linear term
$\mf a (\mtt i) \, \{ \, \mss m_{\mtt i} - E_{\nu^n_{\rho (\cdot)}}
[\mss m_{\mtt i} ] \, \}$ which gives a bound of order
$\Vert \mf a \Vert_\infty \, (n/r_1^{3/2})^d$. All the other terms
are bounded by $\Vert \mf a \Vert_\infty \, (n/r_1^{2})^d$.
\end{remark}

\subsection*{Entropy and logarithmic Sobolev inequalities}

Let $\cb{ \Omega_{\mtt k} } := \{0,1\}^{B_{\mtt k}}$,
$\mtt k\in \Sigma^+_{r_2}$, and denote by $\Omega_{\mtt k, \mss M}$,
$0\le \mss M \le |B_{\mtt k} |$, the configurations of
$\Omega_{\mtt k }$ with $\mss M $ particles:
$\cb{\Omega_{\mtt k, \mss M} } :=\{\xi\in \Omega_{\mtt k }: \sum_{x\in
B_{\mtt k} } \xi_x = \mss M\}$. Hereafter,
$\nu^{B_{\mtt k}}_{\rho (\cdot)} $ represents the marginal of the
measure $\nu^{n}_{\rho (\cdot)} $ on $\Omega_{\mtt k} $, and 
$\nu^{\rm c}_{\mtt k, \mss M}$ the measure
$\nu^{B_{\mtt k} } _{\rho (\cdot)}(\cdot)$ conditioned to
$\Omega_{\mtt k, \mss M} $:
\begin{equation}
\label{91}
\cb{\nu^{B_{\mtt k}}_{\rho (\cdot)} (\xi) } \,:=\,
\nu^{n}_{\rho (\cdot)} \big( \,
\{\eta \in \Omega_n : \eta_x = \xi_x\;\forall
\,x\in B_{\mtt k}\}\, \big) \,, \quad
\cb{\nu^{\rm c}_{\mtt k, \mss M} (\xi') } \,=\,
\frac{\nu^{B_{\mtt k} }_{\rho(\cdot)} (\xi')}
{\sum_{\zeta \in \Omega_{\mtt k, \mss M}} \nu^{B_{\mtt k}
}_{\rho(\cdot)} (\zeta)} \,, \quad
\end{equation}
for $\xi \,\in\, \Omega_{\mtt k}$,
$\xi' \,\in\, \Omega_{\mtt k, \mss M}$.

The next lemma provides an estimate which will be used throughout the
article. Let $\{ V_{\mtt i} : \mtt i \in \Sigma^+_{r_2}\}$ be a
collection of functions: $V_{\mtt i}  \colon \Omega_n \to \bb R$.

\begin{lemma}
\label{s23}
Assume that $V_{\mtt i}$, $\mtt i\in \Sigma^+_{r_2 }$, is measurable
with respect to $\{\eta_x : x\in B_{\mtt i} \}$, that there exists a
finite constant $C_1$ such that
$ \Vert V_{\mtt i} \Vert_\infty \le C_1(n/r_1)^2$ for all
$\mtt i\in \Sigma^+_{r_2 }$, and that $V_{\mtt i}$ has zero-mean with
respect to all canonical measures:
\begin{equation*}
\int V_{\mtt i} \;  d \nu^{\rm c}_{\mtt i, \mss M} \,=\, 0
\quad
\text{for all $\mtt i\in \Sigma^+_{r_2 }$, $0\le \mss M \le r_1^d$}\,.
\end{equation*}
Then,  
\begin{equation*}
\sum_{\mtt i\in \Sigma^+_{r_2 }} 
\int V_{\mtt i}  \, f \, d  \nu^n_{\rho(\cdot)} \,\le\,
\delta\, n^2 \,  I_n  ( f \,;\, \nu^n_{\rho (\cdot)} ) \,+\,
C_2\, \Big(\frac{r_1}{n}\Big)^2\, \sum_{\mtt i\in \Sigma^+_{r_2 }} 
\max_{0\le \mss M \le r_1^d} \int V^2_{\mtt  i} \;
d \nu^{\rm c}_{\mtt i, \mss M} 
\end{equation*}
for all $\delta>0$, and density $f$ with respect to
$\nu^n_{\rho(\cdot)} $. In this formula,
$C_2 = (\mf c_{\rm LS} / 2 \delta)\, \exp\{C_1 \mf c_{\rm LS}
/\delta\}$, and $\mf c_{\rm LS}$ is the logarithmic Sobolev constant
appearing in Theorem \ref{s83}.
\end{lemma}

\begin{proof}
The proof is divided into several steps. 

\smallskip\noindent{\bf Step 1: Projecting on $\Omega_{\mtt k} $.}
Denote by $\mu^{B_{\mtt k}}_f$, $\mtt k\in \Sigma^+_{r_2}$, the
marginal of the measure $f(\cdot) \, \nu^{n}_{\rho (\cdot)} (\cdot)$
on the set $\Omega_{\mtt k} $:
\begin{gather*}
\cb{ \mu^{B_{\mtt k}}_f (\xi) } \,:=\, \int \mtt 1
\{\eta \in \Omega_n : \eta_x = \xi_x\;\forall
\,x\in B_{\mtt k}\}\, f(\eta) \, d\nu^{n}_{\rho (\cdot)} (
\eta)\,, 
\quad
\text{for all}\; \xi \in \Omega_{\mtt k}  \,.
\end{gather*}
In this formula and below, $\cb{\mtt 1\{A\}}$ represents the indicator
function of the set $A$.

Fix $\mtt k\in \Sigma^+_{r_2}$. Since $V_{\mtt k}$ depends only on the
variables $\{\eta_x: x \in B_{\mtt k}\}$,
\begin{equation}
\label{01c}
\int V_{\mtt k} \, f\, d\nu^n_{\rho (\cdot)} \,=\,
\int V_{\mtt k}  \, f_{\mtt k} \, d\nu^{B_{\mtt k}}_{\rho (\cdot)} \,,
\end{equation}
where $f_{\mtt k} $ represents the density (the Radon-Nikodym
derivative) of the measure $\mu^{B_{\mtt k}}_f$  with respect to
$\nu^{B_{\mtt k} }_{\rho(\cdot)}$:
\begin{equation}
\label{06c}
\cb{ f_{\mtt k}   (\xi)}  \,:=\, \frac{\mu^{B_{\mtt k}}_f  (\xi)}
{\nu^{B_{\mtt k} } _{\rho (\cdot)}(\xi)} \,, 
\quad
\xi \in \Omega_{\mtt k}  \,.
\end{equation}

\smallskip\noindent{\bf Step 2: Projecting on hyperplanes.}  
Decomposing the expectation on the right-hand side of \eqref{01c}
according to the total number of particles yields that
\begin{equation}
\label{04b}
\int V_{\mtt k} \, f_{\mtt k} \, d\nu^{B_{\mtt k}}_{\rho (\cdot)} 
\; =\; 
\sum_{\mss M =0}^{| B_{\mtt k}|}
\nu^{B_{\mtt k}}_{\rho (\cdot)}  (\Omega_{\mtt k, \mss  M})\,
\< f_{\mtt k} \>_{\mtt k, \mss M} \,
\int V_{\mtt k} \, f_{\mtt k, \mss M} \, d\nu^{\rm c}_{\mtt k, \mss M}
\;,
\end{equation}
where $\< f_{\mtt k} \>_{\mtt k, \mss M} $ is the expectation of
$f_{\mtt k}$ with respect to the measure
$\nu^{\rm c}_{\mtt k, \mss M}$, and $f_{\mtt k, \mss M} $ is the
density (with respect to the measure $\nu^{\rm c}_{\mtt k, \mss M}$)
given by
\begin{gather*}
\cb{ \< f_{\mtt k} \>_{\mtt k, \mss M} }  \, :=\,
\int  f_{\mtt k} \, d \nu^{\rm c}_{\mtt k, \mss M}\,,
\quad
\cb{ f_{\mtt k, \mss M} (\xi) } \, :=\, \frac {f_{\mtt k} (\xi)}
{\< f_{\mtt k} \>_{\mtt k, \mss M}} \,, \quad \xi \in
\Omega_{\mtt k, \mss M}\,.
\end{gather*}
With this notation, the expression appearing in the statement of the
lemma is equal to
\begin{equation}
\label{93}
\sum_{\mtt i\in \Sigma^+_{r_2 }} 
\sum_{\mss M =0}^{| B_{\mtt i} |}
\nu^{B_{\mtt i}}_{\rho (\cdot)}  (\Omega_{\mtt i, \mss M})\,
\< f_{\mtt i}  \>_{\mtt i , \mss M} \,
\int V_{\mtt  i}  \,
\, f_{\mtt i, \mss M} \, d\nu^{\rm c}_{\mtt i, \mss M} \,.
\end{equation}

Fix $\mtt i\in \Sigma^+_{r_2 }$, $0\le \mss M \le r_1^d$.  By the
entropy inequality,
\begin{equation}
\label{69}
\int V_{\mtt  i}  \,
\, f_{\mtt i, \mss M} \, d\nu^{\rm c}_{\mtt i, \mss M}
\,\le\,
\frac{1}{\gamma} H_{B_{\mtt i}}
(f_{\mtt i, \mss M} \,|\,  \nu^{\rm c}_{\mtt i, \mss M} )
\,+\, \frac{1}{\gamma}  \ln \int e^{ \gamma\,  V_{\mtt  i} }\;
d\nu^{\rm c}_{\mtt i, \mss M}  
\end{equation}
for all $\gamma>0$. In this formula and below, for a cube $B$
contained in $\bb Z^d$, a probability measure $\nu$ defined on
$\{0,1\}^B$, and a density $g\colon \{0,1\}^B \to \bb R_+$ with
respect to $\nu$, $\cb{H_{B} (g \,|\, \nu)}$ represents the entropy of
the measure $g\, d\nu$ with respect to $\nu$. Similarly,
$\cb{I_{B} (g \,;\, \nu)}$ stands for the carr\'e du champ of $g$
averaged with respect to $\nu$:
\begin{equation*}
{\color{blue} I_{B} (g \,;\, \nu) }
\;:=\; \frac{1}{2}\, \sum_{j=1}^d \sum_{x} 
\int  c_{x,x+e_j} (\eta)\,  \Big[\,  \sqrt{g (T_{x,x+e_j} \eta )} \,-\,
\sqrt{g(\eta)}\, \Big]^2 \;  \nu (d\eta) \;.
\end{equation*}
In this equation, the second sum is performed over all $x\in \bb Z^d$
such that $\{x, x+e_j\}\subset B$.

By the logarithmic Sobolev inequality stated in Theorem \ref{s83}, 
\begin{equation*}
\frac{1}{\gamma} H_{B_{\mtt i}}
(f_{\mtt i, \mss M} \,|\,  \nu^{\rm c}_{\mtt i, \mss M} )
\,\le\, 
\frac{1}{\gamma} \, \mf c_{\rm LS}\, r_1^2\, 
I_{B_{\mtt i}}
( f_{\mtt i, \mss M} \,;\,  \nu^{\rm c}_{\mtt i, \mss M} )\,,
\end{equation*}
where $\mf c_{\rm LS}$ represents the logarithmic Sobolev constant. An
elementary computation yields that
\begin{equation*}
\sum_{\mtt i\in \Sigma^+_{r_2 }} 
\sum_{\mss M =0}^{| B_{\mtt i} |}
\nu^{B_{\mtt i}}_{\rho (\cdot)}  (\Omega_{\mtt i, \mss M})\,
\< f_{\mtt i}  \>_{\mtt i , \mss M} \, I_{B_{\mtt i}}
( f_{\mtt i, \mss M} \,;\,  \nu^{\rm c}_{\mtt i, \mss M}  ) \,\le\,
I_n  ( f \,;\, \nu^n_{\rho (\cdot)} ) \,.
\end{equation*}
Hence, the sum over $\mtt i$ and $\mss M$ of the first term on the
right-hand side of \eqref{69} multiplied by
$\nu^{B_{\mtt i}}_{\rho (\cdot)} (\Omega_{\mtt i, \mss M})\, \<
f_{\mtt i} \>_{\mtt i , \mss M} $ is bounded by
$\mf c_{\rm LS}\, r_1^2\, \gamma^{-1} \, I_n ( f \,;\, \nu^n_{\rho
(\cdot)} )$. Choose $\gamma = \mf c_{\rm LS}\, r_1^2 / (\delta n^2)$ to
get the first term on the right-hand side of the statement of the
lemma.

We turn to the second term on the right-hand side of \eqref{69}.  As
in the previous proofs, since $e^x \le 1 + x + (1/2) x^2 e^{|x|}$ for all
$x\in \bb R$, and since $|V_{\mtt i}| \le C_1 (n/r_1)^2$,
\begin{equation*}
\frac{1}{\gamma}  \ln \int e^{ \gamma\,  V_{\mtt  i} }\;
d\nu^{\rm c}_{\mtt i, \mss M}   \,\le\, 
\frac{1}{\gamma}  \ln \Big\{\, 1
\,+\, \gamma \int V_{\mtt  i} \;
d \nu^{\rm c}_{\mtt i, \mss M}   
\,+\,  \frac{1}{2}\, e^{\gamma C_1 (n/r_1)^2}\,
\gamma^2  \int V^2_{\mtt  i} \;
d \nu^{\rm c}_{\mtt i, \mss M} \, \Big\}\,.
\end{equation*}
As $V_{\mtt i}$ has zero-mean and $\ln (1+a) \le a$, by the choice of
$\gamma$, this expression is less than or equal to
\begin{equation*}
\frac{1}{2}\, \gamma \, 
e^{\gamma C_1 (n/r_1)^2}\,\int V^2_{\mtt  i} \;
d\nu^{\rm c}_{\mtt i, \mss M}
\,=\,
\frac{\mf c_{\rm LS}}{2}\, \frac{1}{\delta}\, 
e^{C_1 \mf c_{\rm LS} /\delta}\,
\Big(\frac{r_1}{n}\Big)^2 \,  \int V^2_{\mtt  i} \;
d\nu^{\rm c}_{\mtt i, \mss M}
\,.
\end{equation*}
To complete the proof, it remains to maximize over $\mss M$, and to sum
over $\mtt i$.
\end{proof}

\subsection*{Large deviations estimates}

We conclude this section with an estimate on atypical densities.

\begin{lemma}
\label{s07b}
Fix $0<\mf a_0<\mf r$. Then, there exists a positive and finite
constant $\mf a_1 = \mf a_1(\mf a_0, \mf r)$, depending only on $\mf a_0$,
$\mf r$, such that
\begin{equation*}
\int \sum_{ \mtt i\in \Sigma^+_{r_2} }  r_1^d\, 
\mtt 1\{ \mss m_{\mtt i} \not\in (\mf a_0, 1-\mf a_0) \}\,
f (\eta) \,  d \nu^n_{\rho(\cdot)} (\eta)
\,\le\,
\mf a_1  \, \Big\{\, H_n (f\,|\, \nu^n_{\rho(\cdot)})
\,+\, \Big( \frac{n} {r_1} \Big)^d\, 
e^{-r_1^d/\mf a_1} \,\Big\}
\end{equation*}
for all density $f$ with respect to $\nu^n_{\rho(\cdot)}$.
\end{lemma}

\begin{proof}
By the entropy inequality the expression on the left-hand side of the
statement of the lemma is bounded above by
\begin{equation*}
\frac{1}{\gamma}\, H_n (f\,|\, \nu^n_{\rho(\cdot)})
\,+\,
\frac{1}{\gamma}\, \ln
E_{\nu^n_{\rho(\cdot)}}\Big[
\, \exp\Big\{ \,\gamma\, 
\sum_{ \mtt i\in \Sigma^+_{r_2} }  r_1^d\, 
\mtt 1 \big\{ \, \mss m_{\mtt i} \not\in (\mf a_0, 1-\mf a_0) \,\big\}\,
\,  \Big\} \,\Big]
\end{equation*}
for every $\gamma>0$. Since $\nu^n_{\rho(\cdot)}$ is a product measure
and the cubes $B_{\mtt i}$ are disjoint, the second term of the
previous expression is equal to
\begin{equation*}
\frac{1}{\gamma}\, \sum_{ \mtt i\in \Sigma^+_{r_2} }  \, \ln
E_{\nu^n_{\rho(\cdot)}}\Big[
\, \exp\Big\{ \,\gamma\, r_1^d\, 
\mtt 1 \big\{ \, \mss m_{\mtt i} \not\in (\mf a_0, 1-\mf a_0) \,\big\}\,
\,  \Big\} \,\Big]\,.
\end{equation*}

Since the density profile $\rho(\cdot)$ takes values in the interval
$[\mf r, 1-\mf r]$ and $0<\mf a_0<\mf r$, there exists a positive
constant $c_1 = c_1(\mf a_0,\mf r)$ such that
$\nu^n_{\rho(\cdot)} \{ \, \mss m_{\mtt i} \not\in (\mf a_0, 1-\mf a_0) \, \}
\le e^{- c_1 r_1^d}$ for all $\mtt i\in \Sigma^+_{r_2}$. The
previous expression is thus less than or equal to
\begin{equation*}
\frac{1}{\gamma}\,  \sum_{ \mtt i\in \Sigma^+_{r_2} } \,
\ln \Big \{ 1 \,+\,  \Big( e^{\gamma r_1^d} - 1\Big) \,
\nu^n_{\rho(\cdot)} \{ \, \mss m_{\mtt i} \not\in (\mf a_0, 1-\mf a_0) \,\}
\Big\}
\,\le\,
\frac{1}{\gamma}\, \Big( \frac{n} {r_1}\Big)^d\,
e^{\gamma r_1^d}  \,  e^{- c_1 r_1^d}
\end{equation*}
because $\ln (1+a) \le a$, $a\ge 0$.  To complete the proof, it remains
to choose $\gamma < c_1$.
\end{proof}

\section{Multiscale analysis: First step.}
\label{sec4}

In this section, we present the first estimate required in the proofs
of Theorems~\ref{t02} and~\ref{t01}. Here and up to the end of the
article, $\cb{\mf a_1}$ stands for a finite and positive constant
which depends only on $\mf r$, $\cb{\mf a_2}$ for one which depends
only on $\mf r$ and a cylinder function $h$ (which sometimes will be
$\eta_A$), $\cb{\mf a_{3,p}}$, $p\ge 1$, for one which depends only on
$\mf r$, a cylinder function $h$, and $\Vert \rho\Vert_{C^p}$. All
these constants are allowed to change from line to line. It may happen
that these constants also depend on an extra variable (say some
$\delta>0$ which prevents the empirical density from being close to
$0$ or $1$). In this case this dependence appears explicitly in the
notation, and we write, for instance, $\mf a_2(\delta)$ instead of
$\mf a_2$.

Suppose that $n=\ell_1\, L_2$ for positive integers $\ell_1\ge 2$,
$L_2\ge 2$. Recall the notation introduced above equation \eqref{156}
and in equation \eqref{92} with $r_1=\ell_1$, $r_2=L_2$.  Fix a
cylinder function $h$ and recall that $\ell_h\ge 0$ represents the
smallest integer $p$ that makes the support of $h$ contained in the
cube $\Sigma_p= \{-p, \dots, p\}^d$. Let $B^o_{\mtt k} $ be the
interior of $B_{\mtt k}$ in the sense that
\begin{equation}
\label{95}
\cb{B^o_{\mtt k} } \,:=\, \big\{ x\in B_{\mtt k} : x + \Sigma_{\ell_h}
\subset B_{\mtt k}\,\big\}\,.
\end{equation}
By definition, the support of $\tau_x h$, $x\in B^o_{\mtt k} $, is
contained in $B_{\mtt k} $. On the other hand, for any continuous
function $J\colon \bb T^d \to \bb R$,
\begin{equation}
\label{115}
\sum_{x\in \bb T^d_n}  J(x/n)\,  \tau_x h
\,=\,
\frac{1}{|B^o_{\mtt 0} |}\, 
\sum_{y \in B_{\mtt 0}}  \sum_{\mtt k\in \Sigma^+_{L_2}}
\sum_{x\in y+ B^o_{\mtt k}}  J(x/n)\,  \tau_x h\,,
\end{equation}
where sums over $x$, as we are on the torus $\bb T^d_n$, are performed
modulo $n$.

\begin{remark}
\label{r01}
By \eqref{115}, any estimate of
\begin{equation}
\label{116}
\sum_{\mtt k\in \Sigma^+_{L_2}}  \int
\sum_{x\in y+ B^o_{\mtt k}} J(x/n)\, \Big(\,
\tau_x h \, -
\, E_{\nu^n_{\rho(\cdot)} } \big[\,
\tau_x h \,\big]\, \Big)
f\, d\nu^n_{\rho(\cdot)}
\end{equation}
uniform over $y\in B_{\mtt 0}$ can be turned into an estimate of the
same quantity with the inner set $B^o_{\mtt k}$ replaced by the full
set $B_{\mtt k}$ at a cost given by the multiplicative constant equal
to $|B_{\mtt k}|/|B^o_{\mtt k}|$.  On the other hand, for a fixed
$y \in \bb Z^d$, the cylinder functions
$h_{y, \mtt k} := \sum_{x\in y+ B^o_{\mtt k}} J(x/n)\, \tau_x h$,
$\mtt k\in \Sigma^+_{L_2}$, have disjoint support and the support of
$h_{y, \mtt k}$ is contained in the cube $y+ B_{\mtt k}$. This
property will be important in the proofs of this section and the next.

In view of this remark, in this section and in the next, we estimate
the sum \eqref{116} for $y=0$ and make sure that the arguments provide
the same bound if the expression is translated by $y\in B_{\mtt 0}$.
\end{remark}

\begin{remark}
\label{r03}
The definition of the cubes $B_{\mtt k}^o$ depends on the cylinder
function $h$. There is a finite number of such functions which appear
in the proofs. One could define an integer $\ell'$ as the maximum of
such integers $\ell_h$ and define $B_{\mtt k}^o$ accordingly.
\end{remark}

For a cylinder function $h$, let $\mtt q_{h}\colon [0,1]\to \bb R$ be
the polynomial, and let $J_{\mtt k}$, $\mtt k\in \Sigma^+_{L_2}$, be
the average given by
\begin{equation}
\label{90}
\mtt q_{h} (\alpha)   \;=\,
-\, \frac{1}{2}\, \mtt p_{h}'' (\alpha)
\, \chi( \alpha) \,, \quad
\cb{J_{\mtt k}}  \, :=\, \frac{1}{ |B^o_{\mtt k}| }\,
\sum_{x\in B^o_{\mtt k}} J(x/n)\,. 
\end{equation}
Recall also the definition of the empirical density $\mss m_{\mtt k}$,
$\mtt k \in \Sigma^+_{L_2}$, introduced in \eqref{92}. The main result
of this section reads as follows. Let
\begin{gather*}
\cb{W^h_{\mtt k, x} } \,:=\,
\tau_x h \, -\,  \mtt p_h(\mss m_{\mtt k})
\,-\,  \frac{1}{\ell^d_1}  \, \mtt q_h(\mss m_{\mtt k}) \,, \quad
\overline{W}^h_{\mtt k, x} \,:=\, W^h_{\mtt k, x} \,-\,
E_{\nu^n_{\rho(\cdot)}} [W^h_{\mtt k, x} ]\,
\end{gather*}
and note that for any continuous function $J\colon \bb T^d \to \bb R$, 
\begin{gather*}
\sum_{x\in B^o_{\mtt k}} J(x/n)\, W^h_{\mtt k, x}
\,=\, 
\sum_{x\in B^o_{\mtt k}} J(x/n)\, \tau_x h \, -\,
|B^o_{\mtt k}|\, J_{\mtt k}  \, \Big\{\, \mtt p_h(\mss m_{\mtt k})
\,+\,  \frac{1}{\ell^d_1}  \, \mtt q_h(\mss m_{\mtt k}) \,\Big\}\,.
\end{gather*}
When $h$ is the cylinder function $\eta_A$, $A$ a finite subset of
$\bb Z^d$, we denote $\ell_h$, $W^h_{\mtt k, x}$, $\mtt p_h$,
$\mtt q_h$, by $\cb{\ell_A}$, $\cb{W^A_{\mtt k, x}}$, $\cb{\mtt p_A}$,
$\cb{\mtt q_A}$, respectively.

Recall that we denote by $\mf c_{\rm SG}$ the constant which appears
in Corollary \ref{s84}, where a spectral gap is stated.

\begin{theorem}
\label{s01}
Let $\rho\colon \bb T^d \to (0,1)$ be a density profile of class
$C^3(\bb T^d)$ satisfying \eqref{83}.  For every cylinder function
$h$, and function $J\colon \bb T^d\to\bb R$ of class $C^1(\bb T^d)$,
there exist constants $\mf a_1$, $\mf a_2 = \mf a_2(\mf c_{\rm SG})$,
$\mf a'_2 = \mf a'_2(\mf c_{\rm SG})$, $\mf a_{3,1}$, $\mf a_{3,3}$
such that
\begin{equation*}
\begin{aligned}
& \sum_{\mtt k\in \Sigma^+_{L_2}}  \int
\Big\{  \sum_{x\in B^o_{\mtt k}} J(x/n)\, \overline{W}^h_{\mtt k, x}  
\Big\}\, f\, d\nu^n_{\rho(\cdot)}
\\
&\quad
\,\le\, \delta \, n^2\, I_n(f \,;\, \nu^n_{\rho(\cdot)})
\,+\, \mf a_{3,1} \, \Vert J\Vert_\infty \,
\Big[  \, H_n (f\,|\, \nu^n_{\rho(\cdot)})
\,+\, \ln(2) \,+\, \Big( \frac{n} {\ell_1}\Big)^d\, 
e^{- \mf a_1 \ell_1^d}\,\Big]
\\
&\quad 
\,+\,  \Big\{\, \mf a_2  \, \frac{1}{\delta} \,  
\Vert J\Vert^2_\infty\, +\,
\mf a_{3,3}\, \Vert J\Vert_{C^1} \,\Big\}
\Big(\frac{\ell_1}{n}\Big)^2\, n^d
\,+\,
\mf a_{3,3}\, \Vert J\Vert_{C^1} \,
\frac{1}{n\, \ell^d_1} \, n^d
\,+\, \mf a_{2} \,  \Vert J\Vert_\infty
\Big( \frac{n}{\ell_1^{3}} \Big)^d
\end{aligned}
\end{equation*}
for all $\delta>0$, integers $\ell_1 \ge \ell_h$ such that
$\Vert h \Vert_\infty\,  \Vert J\Vert_\infty\,  \ell_1^{d+2} \le \mf a'_2 \, \delta \,
n^2$, $L_2\ge 2$, and densities $f$ with respect to
$\nu^n_{\rho(\cdot)}$.
\end{theorem}

\begin{remark}
\label{r02b}
This estimate remains valid if we replace $J(\cdot)$ by
$-J(\cdot)$. In particular, the same bound holds for the absolute
value of the left-hand side.
\end{remark}

\begin{remark}
\label{r02}
In dimension $d\ge 4$, for the term $n^d (\ell_1/n)^2$ to be of order
$n^{d/2}$, $\ell_1$ has to be bounded by a constant. This shows that, in
dimension $d\ge 4$, Theorem \ref{s01} does not allow us to replace
the average of a cylinder function over a large box by its expectation with
respect to the canonical measure with density given by the empirical
density $\mss m_{\mtt k}$.
\end{remark}

\begin{remark}
The condition
$\Vert h \Vert_\infty\, \Vert J\Vert_\infty\, \ell_1^{d+2} \le \mf a'_2 \,
\delta \, n^2$ is imposed by the use in the proof of the standard
perturbation theorem on the largest eigenvalue of a symmetric operator
\cite[Theorem A3.1.1]{kl} which requires a $L^\infty$ bound on the
perturbation.
\end{remark}

The next result follows from Theorem \ref{s01} and Lemma \ref{s32}.
It will be used in dimension $d=1$ and $2$ where $\ell_1$ can be taken
much larger than $n^{1/4}$. Note that the last term on the right-hand
side of the statement of Theorem \ref{s01}, originally equal to
$(n/\ell_1^3)^d$, is changed in Corollary \ref{s46} to
$(n/\ell_1^2)^d$, a larger quantity, due to Lemma \ref{s32}.

\begin{corollary}
\label{s46}
Let $\rho\colon \bb T^d \to (0,1)$ be a density profile of class
$C^3(\bb T^d)$ satisfying \eqref{83}.  For every cylinder function
$h$, and function $J\colon \bb T^d\to\bb R$ of class $C^1(\bb T^d)$,
there exist constants $\mf a_1$, $\mf a_2 = \mf a_2(\mf c_{\rm SG})$,
$\mf a'_2 = \mf a'_2(\mf c_{\rm SG})$, $\mf a_{3,1}$, $\mf a_{3,3}$
such that
\begin{equation*}
\begin{aligned}
& \sum_{\mtt k\in \Sigma^+_{L_2}}  \int
\Big\{  \sum_{x\in B^o_{\mtt k}} J(x/n)\, \Big(\,
\tau_x h  \, -
\, E_{\nu^n_{\rho(\cdot)} } \big[\,
\tau_x h \,\big]\, \Big)
\,-\,   |B^o_{\mtt k}| \,  J_{\mtt k}  \, \Big(\, \mtt p_h(\mss m_{\mtt k})
\,-\,  E_{\nu^n_{\rho(\cdot)} } \big[\, \mtt p_h (\mss m_{\mtt k})
\,\big] \,\Big) \, \Big\}\,
f\, d\nu^n_{\rho(\cdot)}
\\
&\quad
\,\le\, \delta \, n^2\, I_n(f \,;\, \nu^n_{\rho(\cdot)})
\,+\, \mf a_{3,1} \, \Vert J\Vert_\infty \,
\Big[  \, H_n (f\,|\, \nu^n_{\rho(\cdot)})
\,+\, \ln(2) \,+\, \Big( \frac{n} {\ell_1}\Big)^d\, 
e^{- \mf a_1 \ell_1^d}\,\Big]
\\
&\quad 
\,+\,  \Big\{\, \mf a_2  \, \frac{1}{\delta} \,  
\Vert J\Vert^2_\infty\, +\,
\mf a_{3,3}\, \Vert J\Vert_{C^1} \,\Big\}
\Big(\frac{\ell_1}{n}\Big)^2\, n^d
\,+\,
\mf a_{3,3}\, \Vert J\Vert_{C^1} \,
\frac{1}{n\, \ell^d_1} \, n^d
\,+\, \mf a_{2} \,  \Vert J\Vert_\infty
\Big( \frac{n}{\ell_1^2} \Big)^d
\end{aligned}
\end{equation*}
for all $\delta>0$, integer $\ell_1 \ge \ell_h$ such that
$\Vert h\Vert_\infty \,  \Vert J\Vert_\infty \,  \ell_1^{d+2} \le \mf a'_2 \, \delta \,
n^2$, $L_2\ge 2$, and density $f$ with respect to
$\nu^n_{\rho(\cdot)}$.
\end{corollary}

\begin{proof}
By Lemma \ref{s32} with $\gamma^{-1} = \Vert J\Vert_\infty$,
$r_1=\ell_1$, there exists a finite constant $\mf a_2$ such that
\begin{equation*}
\int \sum_{\mtt k\in \Sigma^+_{L_2}}  
J_{\mtt k}  \, 
\big\{  \, \mtt q_h(\mss m_{\mtt k}) 
\,-\,  E_{\nu^n_{\rho(\cdot)} } \big[\, \mtt q_h(\mss m_{\mtt k}) \,\big] \,\big\}\,
f\, d\nu^n_{\rho(\cdot)}
\,\le\, \, \Vert J\Vert_\infty\,
\Big\{\, H_n (f\,|\, \nu^n_{\rho(\cdot)}) \,+\,
\ln (2) \, \Big\} 
\,+\,  \mf a_2\, \Vert J\Vert_\infty \, \Big( \frac{n}{\ell_1^2}\Big)^d
\end{equation*}
for all integers $\ell_1\ge \ell_h$, $L_2\ge 2$, and densities $f$
with respect to $\nu^n_{\rho(\cdot)}$.  This estimate takes care of
the correction term
$|B^o_{\mtt k}|\, J_{\mtt k} \, \ell^{-d}_1\, \mtt q_h(\mss m_{\mtt
k}) $ appearing in the statement of Theorem \ref{s01}.
\end{proof}

\noindent{\bf Proof of Theorem \ref{s01}.} We start with a sketch of
the proof of this result. Recall from \eqref{91} that we denote by
$\nu^{\rm c}_{\mtt k, \mss M} $, $\mtt k\in \Sigma^+_{L_2}$,
$0\le \mss M\le |B_{\mtt k}|$, the canonical measure concentrated on
the configurations of $\Omega_{\mtt k}$ with $\mss M$ particles.  The
proof is divided into two main steps.  In the first one, stated in
Proposition \ref{s02}, we estimate
\begin{equation*}
\int
\sum_{x\in B^o_{\mtt k}} J(x/n)\, \Big(\,
\tau_x \eta_A \, -
\, E_{\nu^{\rm  c}_{\mtt k, \mss M_{\mtt k}} } \big[\,
\tau_x \eta_A \,\big]\, \Big)
f\, d\nu^n_{\rho(\cdot)}\,.
\end{equation*}
In the second one, stated in Proposition \ref{s36}, we replace the
canonical measure $\nu^{\rm c}_{\mtt k, \mss M_{\mtt k}}$ by the grand
canonical $\nu^{\rm gc}_{\mtt k, \mss M_{\mtt k}}$, through the local
central limit theorem (LCLT). As the error is required to be small, we
need to consider the first three terms of the Edgeworth expansion in
the LCLT. This explains the presence of the term $\mtt q_A(\cdot)$
in Theorem \ref{s01}.

\begin{proof}[Proof of Theorem \ref{s01}]
The result is a consequence of Propositions \ref{s02} and \ref{s36}. 
\end{proof}

\subsection*{A: Replacing local averages by canonical expectation}

The main result of this subsection reads as follows. Recall that we
denote by $\mf c_{\rm SG}$ the constant which appears in Corollary
\ref{s84}. 

\begin{proposition}
\label{s02}
Under the hypotheses of Theorem \ref{s01}, for each cylinder function
$h$ there exist constants $\mf a_2 = \mf a_2(\mf c_{\rm SG})$,
$\mf a'_2 = \mf a'_2(\mf c_{\rm SG})$ such that
\begin{equation*}
\begin{aligned}
& \sum_{\mtt k\in \Sigma^+_{L_2}}  \int
\Big\{  \sum_{x\in B^o_{\mtt k}} J(x/n)\,  \big\{ \tau_x h
-  E_{\nu^{\rm c}_{\mtt k, \mss M_{\mtt k}}} [\tau_x h] \,\big\}   \, \Big\}\,
f\, d\nu^n_{\rho(\cdot)}
\\
&\quad 
\,\le\, \delta \, n^2\, I_n(f \,;\, \nu^n_{\rho(\cdot)})
\,+\, \mf a_2 \, \Vert J\Vert_\infty  \, \Big\{\, H_n (f\,|\, \nu^n_{\rho(\cdot)})
\,+\, \Big( \frac{n} {\ell_1} \Big)^d\, 
e^{- \mf a_1 \ell_1^d} \,\Big\}
\,+\,  \mf a_2  \, \frac{1}{\delta} \,  
\Vert J\Vert^2_\infty\, \Big(\frac{\ell_1}{n}\Big)^2\, n^d
\end{aligned}
\end{equation*}
for all $\delta>0$, integers $\ell_1 \ge \ell_h$ such that
$\Vert h\Vert_\infty\, \Vert J\Vert_\infty\, \ell_1^{d+2} \le \mf a'_2
\, \delta \, n^2$, $L_2\ge 2$, and densities $f$ with respect to
$\nu^n_{\rho(\cdot)}$.
\end{proposition}

\begin{proof}
The proof is divided into several steps. We first discard small and
large densities. Let $\delta_1 = \mf r/2$ and rewrite the left-hand
side of the statement of the proposition as
\begin{equation}
\label{109}
\sum_{\mtt k\in \Sigma^+_{L_2}}  \int
\mtt 1 \{\mss m_{\mtt k} \not\in (\delta_1 , 1- \delta_1)\}
V_{\mtt k} \, f\, d\nu^n_{\rho(\cdot)}
\,+\,
\sum_{\mtt k\in \Sigma^+_{L_2}}  \int
\mtt 1 \{\mss m_{\mtt k} \in (\delta_1 , 1- \delta_1)\}
V_{\mtt k} \, f\, d\nu^n_{\rho(\cdot)} \,,
\end{equation}
where
\begin{equation}
\label{161}
\cb{V_{\mtt k} } \,  :=\, 
\sum_{x\in B^o_{\mtt k}} J(x/n)\,  \big\{ \tau_x h
-  E_{\nu^{\rm c}_{\mtt k, \mss M_{\mtt k}}} [\tau_x h] \,\big\} \,.
\end{equation}
Since
$|V_{\mtt k} | \le \, 2\, \Vert J\Vert_\infty\, \Vert h\Vert_\infty \,
\ell_1^d$, by Lemma \ref{s07b} (with $r_1=\ell_1$), the first term is
bounded above by
\begin{equation*}
\mf a_2  \, \Vert J\Vert_\infty  \, \Big\{\, H_n (f\,|\, \nu^n_{\rho(\cdot)})
\,+\, \Big( \frac{n} {\ell_1} \Big)^d\, 
e^{- \mf a_1  \ell_1^d} \,\Big\}
\end{equation*}
where $\mf a_1$ is the constant appearing in Lemma \ref{s07b} with
$\mf a_0 = \mf r/2$.

We turn to the second term in \eqref{109}.  By Lemma \ref{s03}, for
each $\delta>0$ there exist constants $\mf a_2 = \mf a_2(\mf c_{\rm SG})$,
$\mf a'_2 = \mf a'_2(\mf c_{\rm SG})$ such that this term is
bounded above by 
\begin{equation*}
\delta \, n^2\, I_n ( f  \,;\, \nu^n_{\rho (\cdot)}) 
\, +\, \mf a_2\, \frac{1}{\delta} \, \Vert J\Vert^2_\infty \,
\, \Big(\frac{\ell_1}{n}\Big)^2 \, n^d
\end{equation*}
provided $\ell_1$ satisfies the hypothesis of the proposition.
To complete the proof, it remains to recollect the previous estimates.
\end{proof}

\begin{lemma}
\label{s03}
Let $\delta_1=\mf r/2$.  Under the hypotheses of Theorem \ref{s01},
there exist constants $\mf a_2 = \mf a_2(\mf c_{\rm SG})$,
$\mf a'_2 = \mf a'_2(\mf c_{\rm SG})$ such that
\begin{equation*}
\sum_{\mtt k\in \Sigma^+_{L_2}}  \int
\mtt 1\{ \mss m_{\mtt k} \in (\delta_1, 1-\delta_1) \}
\, V_{\mtt k} \, f \, d\nu^n_{\rho(\cdot)} \,\le\, 
\delta \, n^2\, I_n ( f  \,;\, \nu^n_{\rho (\cdot)}) 
\, +\, \mf a_2\, \frac{1}{\delta} \, \Vert J\Vert^2_\infty \,
\, \Big(\frac{\ell_1}{n}\Big)^2 \, n^d
\end{equation*}
for all $\delta>0$, integers $\ell_1\ge \ell_h$ such that
$\Vert h\Vert_\infty \,\Vert J\Vert_\infty \,\ell_1^{d+2} \le \mf a'_2 \, \delta \, 
n^2$, $L_2\ge 2$, and densities $f$ with respect to
$\nu^n_{\rho(\cdot)}$. Here, $V_{\mtt k}$ has been introduced in
\eqref{161}.
\end{lemma}

\begin{proof}
 Repeat steps 1 and 2 in
the proof of Lemma \ref{s23} up to equation \eqref{93} to project the
measure $f(\cdot) \, \nu^n_{\rho(\cdot)}(\cdot)$ on the hyperplanes
$\Omega_{\mtt k, \mss M}$. This yields that the sum on the left-hand
side of the statement of the lemma is equal to
\begin{gather}
\label{94}
\sum_{\mtt k\in \Sigma^+_{L_2 }} 
\sum_{\mss M = \delta_1 | B_{\mtt k} | }^{ (1-\delta_1) | B_{\mtt k} |}
\nu^{B_{\mtt k}}_{\rho (\cdot)}  (\Omega_{\mtt k, \mss M})\,
\< f_{\mtt k}  \>_{\mtt k , \mss M} \,
\int V_{\mtt k} 
\, f_{\mtt k, \mss M} \, d\nu^{\rm c}_{\mtt k, \mss M} \,.
\end{gather}
The sum over $\mss M$ runs from $\delta_1 | B_{\mtt k} |$ to
$(1-\delta_1) | B_{\mtt k} |$ due to the indicator function appearing
inside the integral in the statement of the lemma.

\smallskip\noindent{\bf Step 3: Eigenvalue estimate.} Fix
$\delta_1 | B_{\mtt k} | \le \mss M \le (1-\delta_1) | B_{\mtt k} |$.
Since $V_{\mtt k}$ has zero-mean with respect to
$\nu^{\rm c}_{\mtt k, \mss M}$, by the standard perturbation theorem
on the largest eigenvalue of a symmetric operator \cite[Theorem
A3.1.1]{kl},
\begin{equation}
\label{03}
\int V_{\mtt k} \,  f_{\mtt k, \mss M} 
\, d\nu^{\rm c}_{\mtt k, \mss M} 
\,\le\, \epsilon^{-1} \, I_{B_{\mtt k}} ( f_{\mtt k, \mss M}  \,;\,
\nu^{\rm c}_{\mtt k, \mss M}) \, +\,\epsilon \,
\frac{ \< V_{\mtt k} ,\, (-L_{\mtt k})^{-1} V_{\mtt k} \>_{\mtt k, \mss M}}
{1 - 2 \, \epsilon \, \mf s_{\mtt k, \mss M} \, \Vert V_{\mtt k} \Vert_{\infty}}
\end{equation}
for all
$0<\epsilon < [2 \, \mf s_{\mtt k, \mss M} \, \Vert V_{\mtt k}
\Vert_{\infty}]^{-1}$. In this formula,
$\< \,\cdot\, ,\, \cdot\,\>_{\mtt k, \mss M}$ represents the scalar
product in $L^2(\nu^{\rm c}_{\mtt k, \mss M})$, $L_{\mtt k}$ stands
for the generator of the exclusion dynamics on $B_{\mtt k}$ reversible
with respect to $\nu^{\rm c}_{\mtt k, \mss M}$,
$\mf s_{\mtt k, \mss M}$ for the inverse of the spectral gap of the
operator $L_{\mtt k}$ in $L^2(\nu^{\rm c}_{\mtt k, \mss M})$, and
$I_{B_{\mtt k}} (\,\cdot\,;\, \nu^{\rm c}_{\mtt k, \mss M} )$ for the
Dirichlet form induced by the generator $L_{\mtt k}$:
\begin{gather*}
(L_{\mtt k} g) (\xi) \,=\, \sum_{j=1}^d \sum_{x: \{x,x+e_j\} \subset
B_{\mtt k}} 
c_{x,x+e_j}(\xi ; \rho (\cdot))\, 
[\, g(T_{x,x+e_j} \xi) - g(\xi)\,] \,,
\\
c_{x,x+e_j}(\eta; \rho(\cdot)) \,=\, \exp \big\{ - (1/2) [\, H_{x+e_j} - H_x \,]\,
[\eta_{x+e_j} - \eta_x]\,\big\}\,, \quad H_z \,=\,
\ln \frac{ \rho(z/n)}{1-\rho(z/n)} \,, 
\\
I_{B_{\mtt k}} (g  \,;\, \nu^{\rm c}_{\mtt k, \mss M} )
\,=\, \frac{1}{2}\, \sum_{j=1}^d \sum_{x: \{x,x+e_j\} \subset
B_{\mtt k}}  \int c_{x,x+e_j}(\xi; \rho(\cdot))\,
\Big[\, \sqrt{ g(T_{x,x+e_j} \xi) } - \sqrt{ \vphantom{T_{x,y} } g(\xi) }
\, \Big]^2
\, d \nu^{\rm c}_{\mtt k, \mss M}  \,.
\end{gather*}
Note that our choice of $L_{\mtt k}$ is different from $L_n$ in
\eqref{10}.

By definition, $\Vert V_{\mtt k} \Vert_{\infty}$ is bounded by
$2\, \Vert h\Vert_\infty \, \Vert J\Vert_\infty \, \ell_1^d$. By the
spectral gap, stated in Corollary \ref{s84}, there exists a finite
constant $\mf c_{\rm SG}$ such that
$\mf s_{\mtt k, \mss M} \le \mf c_{\rm SG} \, \ell_1^2$ for all
$\ell_1\ge 2$, $0\le \mss M\le |B_{\mtt k}|$,
$\mtt k\in \Sigma^+_{L_2}$. Therefore, if we choose
$0<\epsilon < [8 \, \mf c_{\rm SG} \, \Vert h\Vert_\infty \, \Vert
J\Vert_\infty \, \ell_1^{d+2}]^{-1}$, since
$\< V_{\mtt k} ,\, (-L_{\mtt k})^{-1} V_{\mtt k} \>_{\mtt k, \mss M}
\le \mf s_{\mtt k, \mss M} \< V_{\mtt k} ,\, V_{\mtt k} \>_{\mtt k,
\mss M} $,
\begin{gather*}
\int V_{\mtt k} \, f_{\mtt k, \mss M}  
\, d \nu^{\rm c}_{\mtt k, \mss M}  
\,\le\, \epsilon^{-1} \, I_{B_{\mtt k}} ( f_{\mtt k, \mss M}   \,;\,
\nu^{\rm c}_{\mtt k, \mss M} ) \, +\, 2\, \epsilon \,
\mf c_{\rm SG} \, \ell_1^2 \,
\int V_{\mtt k}^2  \, d \nu^{\rm c}_{\mtt k, \mss M} \,.
\end{gather*}
By Corollary \ref{s42}, since
$\delta_1 | B_{\mtt k} | \le \mss M \le (1-\delta_1) | B_{\mtt k} |$,
the variance term is less than or equal to
$\mf a_2 \, \Vert J\Vert^2_\infty\, \ell_1^d$. Thus,
\begin{gather*}
\int V_{\mtt k} \, f_{\mtt k, \mss M}  
\, d \nu^{\rm c}_{\mtt k, \mss M}  
\,\le\, \epsilon^{-1} \, I_{B_{\mtt k}} ( f_{\mtt k, \mss M}  \,;\,
\nu^{\rm c}_{\mtt k, \mss M} ) \, +\, \mf a_2\, \epsilon \,
\mf c_{\rm SG} \, \Vert J\Vert^2_\infty \, \ell_1^{d+2}  \,.
\end{gather*}

\smallskip\noindent{\bf Step 4: Conclusion.}  Inserting this bound into
equation \eqref{94} and summing over $\mss M$ yields that
\begin{gather*}
\int \mtt 1 \{\mss m_{\mtt k} \in (\delta_1 , 1- \delta_1)\}\,
V_{\mtt k} \, f_{\mtt k} \, d\nu^{B_{\mtt k}}_{\rho (\cdot)}
\,\le\, \epsilon^{-1} \, I_{B_{\mtt k}}  ( f_{\mtt k}  \,;\,
\nu^{B_{\mtt k}}_{\rho (\cdot)}) \, +\, \mf a_2\, \epsilon \,
\mf c_{\rm SG} \, \Vert J\Vert^2_\infty \, \ell_1^{d+2} \,.
\end{gather*}
As $\nu^n_{\rho(\cdot)}$ is a product measure, by the definition of
$f_{\mtt k}$ introduced in display \eqref{06c}, 
\begin{equation*}
f_{\mtt k} (\xi) \,=\,
\frac{\sum_{\eta} f(\xi, \eta) \, \nu^n_{\rho(\cdot)} (\xi, \eta)}
{\sum_{\eta} \nu^n_{\rho(\cdot)} (\xi, \eta) }
\,=\,
\frac{\sum_{\eta} f(\xi, \eta) \, \nu^{A_{\mtt k}}_{\rho(\cdot)} (\eta)}
{\sum_{\eta} \nu^{A_{\mtt k}}_{\rho(\cdot)} (\eta) }
\,=\,
\sum_{\eta} f(\xi, \eta) \, \nu^{A_{\mtt k}}_{\rho(\cdot)} (\eta)
\,.
\end{equation*}
In this formula, $\xi$ belongs to $\{0,1\}^{B_{\mtt k}}$ and
$\eta$ belongs to $\{0,1\}^{A_{\mtt k}}$, where
$A_{\mtt k} = \bb T_n^d \setminus B_{\mtt k}$, and
$\nu^{A_{\mtt k}}_{\rho(\cdot)} (\cdot)$ is the marginal of the
measure $\nu^n_{\rho(\cdot)} (\cdot)$ on $\{0,1\}^{A_{\mtt k}}$.

It follows from the previous displayed equation and Schwarz inequality
that
\begin{equation*}
I_{B_{\mtt k}} ( f_{\mtt k}  \,;\,
\nu^{B_{\mtt k}}_{\rho (\cdot)})  
\,\le\,
I_{B_{\mtt k}}  (f\,;\, \nu^n_{\rho(\cdot)} )\,,
\end{equation*}
where
\begin{equation*}
I_{B_{\mtt k}}  (f\,;\, \nu^n_{\rho(\cdot)} ) \,=\,
\frac{1}{2}\, 
\sum_{j=1}^d \sum_{x: \{x,x+e_j\} \subset
B_{\mtt k}}  \int c_{x,x+e_j}(\xi)\,
\Big[\, \sqrt{ f(T_{x,x+e_j} \xi) } - \sqrt{ \vphantom{T_{x,y} } f(\xi) }
\, \Big]^2
 \, d \nu^n_{\rho (\cdot)} (\xi) \,.
\end{equation*}
Therefore,
\begin{equation*}
\int \mtt 1 \{\mss m_{\mtt k} \in (\delta_1 , 1- \delta_1)\}\,
V_{\mtt k} \, f_{\mtt k} \, d\nu^{B_{\mtt k}}_{\rho (\cdot)}
\,\le\, \epsilon^{-1} \, I_{B_{\mtt k}} ( f  \,;\,
\nu^n_{\rho (\cdot)}) \, +\, 
\mf a_2\, \epsilon \,
\mf c_{\rm SG} \, \Vert J\Vert^2_\infty \,  \, \ell_1^{d+2}  \,.
\end{equation*}

Summing over $\mtt k\in \Sigma^+_{L_2}$, as the cubes have length
$\ell_1$, we conclude that the sum appearing on the left-hand side of
the statement of the lemma is bounded by
\begin{equation*}
\frac{1}{\epsilon} \,  \sum_{\mtt k\in \Sigma^+_{L_2}}
I_{B_{\mtt k}} ( f  \,;\, \nu^n_{\rho (\cdot)})
\, +\, \mf a_2\, \epsilon \,
\mf c_{\rm SG} \, \Vert J\Vert^2_\infty \,
\Big(\frac{n}{\ell_1}\Big)^d \, \ell_1^{d+2}
\end{equation*}
for all $\epsilon$ satisfying the restrictions imposed in Step 3.  As
the sets $B_{\mtt k}$ are disjoint, the sum is bounded by
$I_n ( f \,;\, \nu^n_{\rho (\cdot)})$. To complete the proof of the
lemma, it remains to set $\epsilon^{-1} = \delta\, n^2$.
\end{proof}

\begin{remark}
\label{rm7}
Here is an alternative proof of Proposition \ref{s02}.
Instead of using the standard perturbation theorem on the largest
eigenvalue of a symmetric operator, one could apply the entropy
inequality and then the logarithmic Sobolev inequality, as in the
proof of Lemma \ref{s23}. 
\end{remark}

\begin{remark}
\label{rm3}
For the last term on the right-hand side of Proposition \ref{s02} to
be of order $o(n^{d/2})$, one needs $\ell_1 \ll n^{1-(d/4)}$. Thus, to
fulfill also the second condition of Lemma \ref{s03} which
requires $\ell_1^{d+2} \le c \, n^2$, in dimension $1$, set
$\ell_1 = \epsilon_n n^{2/3}$, in dimension $2$,
$\ell_1 = \epsilon_n n^{1/2}$, and in dimension $3$,
$\ell_1 = \epsilon_n n^{1/4}$ for some sequence $\epsilon_n\to 0$.
These are the side lengths of the cubes reached after this first step.
However, we will need lengths of order at least $k\, n^{1/2}$ for some
$1\ll k \ll n$. Thus, in dimensions $2$ and $3$ we need to increase the
length further. This is the goal of the next section. 
\end{remark}

\begin{remark}
\label{rm1}
This estimate is analogous to the one presented in \cite[Proposition
5.2]{jl23} with $a_n=1$, whose proof follows the method introduced by
Jara and Menezes \cite{jm1, jm2}. In \cite[Proposition 5.2]{jl23}, the
constant term is of order $g_d(\ell) \, n^{d-2}$, while it is here
$\ell^2 \, n^{d-2}$, which is much larger since $g_d(\ell) \ll \ell^2$
in all dimensions. However, unlike the bound established in
Proposition \ref{s02}, the estimate in \cite[Proposition 5.2]{jl23}
includes on the right-hand side an entropy term multiplied by a
constant which diverges as $n\to\infty$. Also, the
structure of the function $V$ in \cite[Proposition 5.2]{jl23} is
different from $V_{\mtt k}$ here.
\end{remark}

\subsection*{ B: From canonical to grand canonical measures.}

In this subsection, we complete the proof of Theorem \ref{s01}.  Fix a
finite subset $A$ of $\bb Z^d$.  Recall the definition of $\ell_A$,
introduced just before the statement of Theorem \ref{s01}, and of the
polynomials $\cb{\mtt p_A(\cdot)}$, $\cb{\mtt q_A(\cdot)}$: 
\begin{equation}
\label{168}
\cb{\mtt p_A(\alpha)}  = \alpha^{|A|}\,, \quad
\cb{\mtt q_A(\alpha)}  = -\, \frac{1}{2}\, |A|\, (|A|-1)\,
\alpha^{|A|-2} \, \chi(\alpha)\,.
\end{equation}
Let 
\begin{gather*}
\cb{W_{\mtt k, A} } \,:=\, 
\mss m_{\mtt k}^{|A|} \,+\, \frac{1}{ \ell^d_1} \,
\mtt q_A(\mss m_{\mtt k}) \,,
\quad
\cb{\overline{W}_{\mtt k, A} } \,:=\, 
W_{\mtt k, A} 
\,-\,  E_{\nu^n_{\rho(\cdot)} } \big[\, W_{\mtt k, A}  \,\big]
\,,\quad \mtt k\in \Sigma^+_{L_2}\,,
\\
\cb{W_{x, A} } \,:=\,
E_{\nu^{\rm c}_{\mtt k, \mss M_{\mtt k}}} [\tau_x \eta_A] \,-\,
W_{\mtt k, A} \,, \quad
\cb{\overline {W}_{x, A} } \,:=\, W_{x, A}
\,-\,  E_{\nu^n_{\rho(\cdot)} } \big[\, W_{x, A}  \,\big]
\,,\quad x \in B^o_{\mtt k} \,.
\end{gather*}
Note that
\begin{gather}
\label{211}
\sum_{x\in B^o_{\mtt k}} J(x/n)\, W_{x, A}\,=\,
\sum_{x\in B^o_{\mtt k}} J(x/n)\,
E_{\nu^{\rm c}_{\mtt k, \mss M_{\mtt k}} }
\big[\, \tau_x \eta_A \,\big]
\,-\,
|B^o_{\mtt k}|\, J_{\mtt k} \,
\Big\{\, \mss m_{\mtt k}^{|A|} \,+\, \frac{1}{ \ell^d_1 } \,
\mtt q_A(\mss m_{\mtt k}) \,\Big\}\,.
\end{gather}
The next result allows a replacement of the first term on the
right-hand side by the second one.

\begin{proposition}
\label{s36}
Let $\rho\colon \bb T^d \to (0,1)$ be a density profile of class
$C^3(\bb T^d)$ satisfying \eqref{83}.  Assume that
$\sup_{n\ge 1} (\ell^d_1/n) <\infty$.  Fix a function $J$ in
$C^1(\bb T^d)$, and a finite set $A$.  Then, there exist constants
$\mf a_1$, $\mf a_2$, $\mf a_{3,1}$, $\mf a_{3,3}$ such that
\begin{equation*}
\begin{aligned}
&  \int \Big|\, \sum_{\mtt k\in \Sigma^+_{L_2}}  
\sum_{x\in B^o_{\mtt k}} J(x/n)\, \overline {W}_{x, A} \,
\Big |\, f\, d\nu^n_{\rho(\cdot)}
\\
&\quad
\,\le\, \mf a_{3,1} \, \Vert J\Vert_\infty \,
\Big[  \, H_n (f\,|\, \nu^n_{\rho(\cdot)})
\,+\, \ln(2) \,+\, \Big( \frac{n} {\ell_1}\Big)^d\, 
e^{- \mf a_1 \ell_1^d}\,\Big] 
\\
&\quad
\,+\, \mf a_{3,3}\, \Vert J\Vert_{C^1} \,
\Big\{\, \frac{1}{n\, \ell^d_1}
\,+\, \Big( \frac{\ell_1} {n}\Big)^2\,\Big\} \, n^d
\,+\, \mf a_{2} \, \Vert J\Vert_\infty
\Big( \frac{n}{\ell_1^{3}} \Big)^d 
\end{aligned}
\end{equation*}
for all densities $f$ with respect to $\nu^n_{\rho(\cdot)}$,
$\ell_1\ge \ell_A$, $L_2\ge 1$.
\end{proposition}

We first remove small and large densities. Let
$\cb{\delta_1} : = \mf r/2$. Since
$\max_{\mtt k\in \Sigma^+_{L_2}} \max_{x\in B^o_{\mtt k}} |W_{x, A} |
\le \mf a_2$, by Lemma \ref{s07b}, 
\begin{equation}
\label{170}
\begin{aligned}
&  \int  \sum_{\mtt k\in \Sigma^+_{L_2}}
\mtt 1 \{\mss m_{\mtt k} \not\in (\delta_1 , 1-\delta_1)\}
\sum_{x\in B^o_{\mtt k}}  \big|\, J(x/n)\, \overline {W}_{x, A}
\,\big| \, f\, d\nu^n_{\rho(\cdot)}
\\
&\quad
\,\le\, \mf a_2\, \Vert J\Vert_\infty \,
\Big[  \, H_n (f\,|\, \nu^n_{\rho(\cdot)})
\,+\, \ln(2) \,+\, \Big( \frac{n} {\ell_1}\Big)^d\, 
e^{- \mf a_1 \ell_1^d}\,\Big] \,.
\end{aligned}
\end{equation}
In this formula, $\mf a_1$ is the constant $\mf a_1$ appearing in
Lemma \ref{s07b} with $\mf a_0 = \mf r/2$.

Recall from the statement of Proposition \ref{s12b} the definitions of
$\Gamma_N$, $\mtt T_{i, A}$, $i=1$, $2$.  These variables depend on
the cube $\Sigma^+_r$, the density profile $\varrho(\cdot)$, the
empirical density $\mss m$, and the set $A$. Below, $\Gamma_{\mtt k}$,
$\mtt T^{(i)}_{\mtt k, A}$ stand for $\Gamma_N$, $\mtt T_{i, A}$,
taking $\Sigma^+_r = B_{\mtt k}$, $\varrho(\cdot)$ as the restriction
to the cube $B_{\mtt k}$ of the density profile $\rho(\cdot)$, and
$\mss m = \mss m_{\mtt k}$.

Let
\begin{gather*}
W^{\rm exp}_{x,A} (\mss m_{\mtt k}) \,=\,
E_{\nu^{\rm gc} _{\mtt k, \mss M_{\mtt k}}} [\, \tau_x \eta_A\,] \,
\, \Big( 1 \,+\, \frac{1}{\Gamma_{\mtt k}} \, \mtt T^{(1)}_{\mtt k, x+A}
\,+\, \frac{1}{\Gamma_{\mtt k}^2} \, \mtt T^{(2)}_{\mtt k, x+A} \,\Big) \,.
\end{gather*}
Proposition \ref{s12b} asserts that there exists a finite constant
$\mf a_2$ such that
\begin{equation}
\label{88}
\max_{ \mtt k\in \Sigma^+_{L_2} }  
\max_{x\in B^o_{\mtt k}} 
\sup_{\delta_1 \le \mss m_{\mtt k}  \le 1-\delta_1}\,
\, \big|\,  E_{\nu^{\rm c} _{\mtt k, \mss M_{\mtt k}}} [\, \tau_x \eta_A\,]
\,-\, W^{\rm exp}_{x,A}  (\mss m_{\mtt k}) \, \big| 
\,\le\,
\mf a_2  \, \frac{1}{\ell_1^{3d}}
\end{equation}
for all $\ell_1\ge \ell_A$. In this equation,
$\nu^{\rm gc} _{\mtt k, \mss M_{\mtt k}}$ is the product measure
$\nu^n_{\rho(\cdot)}$ tilted for the expectation of
$\sum_{x\in B_{\mtt k}} \eta_x$ to be $\mss M_{\mtt k}$.  As mentioned
before, the set $B^o_{\mtt k}$, introduced in \eqref{95}, guarantees
that the support of $ \tau_x \eta_A$, $x\in B^o_{\mtt k}$, is
contained in $B_{\mtt k}$.

\begin{remark}
\label{rm4}
We need the bound on the right-hand side of Theorem \ref{s01} and
Proposition \ref{s36} to be of order $o(n^{d/2})$. In view of the
definition of $\ell_1$ presented in Remark \ref{rm3}, in dimension $1$
we may end the expansion in equation \eqref{88} at the zero-order
term which is $1$. In dimension $2$, at the first-order term
($ 1 + \Gamma_{\mtt k}^{-1} \, \mtt T^{(1)}_{\mtt k, x+A} )$.  It is
only in dimension $3$ that we need to deal with the second-order term.
\end{remark}

The next result is a straightforward consequence of  \eqref{88}.  Let
$\cb{\mtt 1_{\delta_1} (\mss m_{\mtt k}) } : = \mtt 1 \{\mss m_{\mtt
k} \in (\delta_1 , 1-\delta_1)\} $.

\begin{lemma}
\label{s33}
Let $\rho\colon \bb T^d \to (0,1)$ be a continuous density profile
satisfying \eqref{83}.  Fix a continuous function
$J\colon \bb T^d \to \bb R$, and a finite set $A$. Then, there exist
constants $\mf a_1$, $\mf a_2$ such that
\begin{equation*}
\begin{aligned}
\sum_{\mtt k\in \Sigma^+_{L_2}}  \int
\mtt 1_{\delta_1} (\mss m_{\mtt k}) \,
\sum_{x\in B^o_{\mtt k}}   \Big | \, J(x/n)\,  
\Big\{   \, E_{\nu^{\rm c}_{\mtt k, \mss M_{\mtt k}}} [\tau_x \eta_A] \,
-\, W^{\rm exp}_{x,A}  (\mss m_{\mtt k}) \, \Big\}
\, \Big|  \, 
f\, d\nu^n_{\rho(\cdot)}
\,\le\, \mf a_2  \, \Vert J\Vert_\infty\, \Big(\frac{n}{\ell_1^{3}} \Big)^d 
\end{aligned}
\end{equation*}
for all densities $f$ with respect to
$\nu^n_{\rho(\cdot)}$, $\ell_1\ge \ell_A$, $L_2\ge 1$.
\end{lemma}

At this point we have not yet modified the term
$E_{\nu^n_{\rho(\cdot)}} \big[\, E_{\nu^{\rm c}_{\mtt k, \mss M_{\mtt
k}}} [\tau_x \eta_A] \,-\, W_{\mtt k, A} \big]$ which appears in the
left-hand side of the inequality asserted in Proposition \ref{s36}. We
may proceed as above, inserting first the indicator
$\mtt 1_{\delta_1} (\mss m_{\mtt k})$ inside the expectation
$E_{\nu^n_{\rho(\cdot)} } \big[\, \cdot \,\big]$ (this is needed
because the estimate \eqref{88} holds only for densities
$\mss m_{\mtt k}$ away from $0$ and $1$ and the expressions
$\mtt T^{(i)}_{\mtt k, x+A} $ may explode for densities close to these
values). Then, we may replace
$\mtt 1_{\delta_1} (\mss m_{\mtt k}) \, E_{\nu^{\rm c}_{\mtt k, \mss
M_{\mtt k}}} [\tau_x \eta_A]$ by
$\mtt 1_{\delta_1} (\mss m_{\mtt k}) \, W^{\rm exp}_{x,A} (\mss
m_{\mtt k}) $.  The cost of these two steps is the same as above and
simpler because the expectation is carried out with respect to the
product measure $\nu^n_{\rho(\cdot)} (\eta)$ (that is with
$f(\eta)=1$) and not $f(\eta)\, \nu^n_{\rho(\cdot)} (\eta)$.

After these steps, it remains to estimate
\begin{equation}
\label{171}
\int \Big|\, \sum_{\mtt k\in \Sigma^+_{L_2}}
\mtt 1_{\delta_1} (\mss m_{\mtt k}) \,
\sum_{x\in B^o_{\mtt k}} J(x/n)\, \overline {W}^{(1)}_{x, A} \,
\Big |\, f\, d\nu^n_{\rho(\cdot)}\,,
\end{equation}
where
\begin{equation}
\label{172}
\cb{W^{(1)}_{x, A} } \,:=\, W^{\rm exp}_{x,A}  (\mss m_{\mtt k})
\,-\, W_{\mtt k, A}\,, \quad
\cb{  \overline {W}^{(1)}_{x, A} } \, :=\,
W^{(1)}_{x, A} \,-\, E_{\nu^n_{\rho(\cdot)}}
\big[\, \mtt 1_{\delta_1} (\mss m_{\mtt k}) \, W^{(1)}_{x, A} \,\big]\,.
\end{equation}
Due to the presence of the indicator, we are not
subtracting the mean of $W^{(1)}_{x, A} $, but an expectation very
close to it.

Divide the term $W^{\rm exp}_{x,A} (\mss m_{\mtt k})$ into two pieces: let
\begin{equation}
\label{173}
\cb{W^{{\rm e}, 1}_{x,A} (\mss m_{\mtt k}) } \, :=\,
E_{\nu^{\rm gc} _{\mtt k, \mss M_{\mtt k}}} [\, \tau_x \eta_A\,] \,
\, \Big( 1 \,+\, \frac{1}{\Gamma_{\mtt k}} \,
\mtt T^{(1)}_{\mtt k, x+A} \,\Big) \,, \quad
\cb{W^{{\rm e}, 2}_{x,A} (\mss m_{\mtt k})} \,:=\,
E_{\nu^{\rm gc} _{\mtt k, \mss M_{\mtt k}}} [\, \tau_x \eta_A\,] \,
\, \frac{1}{\Gamma_{\mtt k}^2} \, \mtt T^{(2)}_{\mtt k, x+A}\,, 
\end{equation}
so that
\begin{equation}
\label{175}
W^{\rm exp}_{x,A} (\mss m_{\mtt k}) \,=\,
W^{{\rm e}, 1}_{x,A} (\mss m_{\mtt k})
\,+\, W^{{\rm e}, 2}_{x,A} (\mss m_{\mtt k}) \,.
\end{equation}
We estimate the second term on the right-hand side in the next
subsection, and the first one in the following one.

\subsection*{Negligible terms}

In Lemma \ref{s33}, we replaced the canonical expectation
$E_{\nu^{\rm c}_{\mtt k, \mss M_{\mtt k}}} [\tau_x \eta_A]$ by the
grand canonical one
$E_{\nu^{\rm gc} _{\mtt k, \mss M_{\mtt k}}} [\, \tau_x \eta_A\,]$
multiplied by the sum of three expressions. In this subsection, we
show that the terms
$E_{\nu^{\rm gc} _{\mtt k, \mss M_{\mtt k}}} [\, \tau_x \eta_A\,]\,
\Gamma_{\mtt k}^{-2}\, \mtt T^{(2)}_{\mtt k, x+A}$ are small and can be
disregarded.

Recall the statement of Proposition \ref{s12b}.  Denote by
$\cb{\varrho_{\mtt k} }\colon \bb T^d \to (0,1)$,
$\mtt k \in \Sigma^+_{L_2} $, the density introduced in \eqref{50c}
and satisfying the identity stated in this equation with
$\Sigma^+_{r} = B_{\mtt k}$, $\mss m = \mss m_{\mtt k}$,
$\varrho (\cdot) = \rho(\cdot)$ restricted to $B_{\mtt k}$. Recall the
expression for $\mtt T^{(2)}_{\mtt k, x+A}$. Since all terms forming
$\mtt T^{(2)}_{\mtt k, A}$ are handled in the same way, we focus
on one, say $\mtt s_A^2$.

Let
\begin{equation*}
\cb{W^{\rm neg}_{x, A} } \, :=\, 
E_{\nu^{\rm gc} _{\mtt k, \mss M_{\mtt k}}} [\, \tau_x \eta_A\,] \,
\frac{\mtt s_{x+A}^2} { \Gamma_{\mtt k}^2}\,, \quad
x\in B^o_{\mtt k} \,, \;\; \mtt k\in \Sigma^+_{L_2}\,.
\end{equation*}
Denote by $F\colon [0,1]\to \bb R$ the function that
$W^{\rm neg}_{x, A}$ would be if the density profile
$\rho(\cdot)$ were constant (without the normalizing constant
$|B_{\mtt k}|^2$ reflected in $\Gamma_{\mtt k}^2$):
\begin{equation}
\label{167}
F(\mss m) \,=\, 
E_{\nu^n_{\mss m}} [\, \eta_A\,] \,
\frac{1}{\chi(\mss m)^2} \, |A|^2 \, [1-\mss m]^2
\,=\, |A|^2\, \mss m^{|A|-2}\, .
\end{equation}
In this example $F(\cdot)$ is a nice polynomial. However, due to the
denominator $\chi(\mss m)^2$, in other cases $F(\cdot)$ might explode,
or the derivative might explode at the boundaries of the interval
$[0,1]$. This is irrelevant, however, due to the indicator function
$\mtt 1_{\delta_1} (\mss m_{\mtt k}) $ which prevents densities close
to the boundary. Denote by $F_{\mf r}\colon [0,1]\to \bb R$ a smooth
function which coincides with $F$ on the interval $[\delta_1/2, 1-
(\delta_1/2)]$, and let
\begin{equation*}
\cb{\overline{W}^{\rm neg}_{x, A} } \, :=\, 
W^{\rm neg}_{x, A}
\,-\, \frac{1}{\ell^{2d}_1}\, 
E_{\nu^n_{\rho(\cdot) }  } \big [\, F_{\mf r} (\mss m_{\mtt k}) \,\big] \,,
\quad x\in B^o_{\mtt k} \,, \;\; \mtt k\in \Sigma^+_{L_2}\,.
\end{equation*}

\begin{lemma}
\label{s68}
Let $\rho\colon \bb T^d \to (0,1)$ be a density profile of class
$C^3(\bb T^d)$ satisfying \eqref{83}.  Fix a continuous function
$J\colon \bb T^d \to \bb R$, and a finite set $A$. Then, there exist
finite constants $\mf a_1$, $\mf a_2$, $\mf a_{3,3}$ such that
\begin{equation*}
\begin{aligned}
\int \Big| \, \sum_{\mtt k\in \Sigma^+_{L_2}}  
\mtt 1_{\delta_1} (\mss m_{\mtt k}) \,
\sum_{x\in B^o_{\mtt k}} J(x/n)\,
\overline{W}^{\rm neg}_{x, A} \, \Big | \, f\,
d\nu^n_{\rho(\cdot)}
\, &\le\, \Vert J\Vert_\infty
\Big\{\, \mf a_{2} \, \Big( \frac{n}{\ell_1^{4}} \Big)^d
\,+\, \mf a_{3,3} \, 
\frac{\ell_1} {n} \,  \Big(\frac{n}{\ell_1^{2}} \Big)^d\, 
\,\Big\} 
\\
& + \, \mf a_2 \, \Vert J\Vert_\infty\,
\Big\{   \, H_n (f\,|\, \nu^n_{\rho(\cdot)})
\,+\, \ln(2) \,+\, \Big( \frac{n} {\ell_1}\Big)^d\, 
e^{- \mf a_1 \ell_1^d} \Big\}
\end{aligned}
\end{equation*}
for all densities $f$ with respect to $\nu^n_{\rho (\cdot)}$,
$\ell_1\ge \ell_A$, $L_2\ge 1$.
\end{lemma}

\begin{proof}
The proof is based on Lemmata \ref{s32}, \ref{s07b}, and \ref{s15}. By
definition,
\begin{align*}
W^{\rm neg}_{x, A}\,=\,
E_{\nu^{\rm gc} _{\mtt k, \mss M_{\mtt k}}} [\, \tau_x \eta_A\,] \,
\Big( \, \frac{1} {\sum_{z\in B_{\mtt k} } \chi (\varrho_{\mtt k}  (z/n) ) } \,\Big)^2
\, \Big( \sum_{y\in A} \big[\, 1 - \varrho_{\mtt k}  ([x+y]/n)\, \big] \, \Big)^2 \,,
\end{align*}
where $\varrho_{\mtt k} (\cdot)$ has been introduced at the beginning
of this subsection. As $\delta_1 \le \mss m_{\mtt k} \le 1-\delta_1$,
$\mf r \le \rho (\mtt x) \le 1-\mf r$, and $\delta_1$ has been set to
be $\mf r/2$, by \eqref{56b}, there exists $0<\mf a_1 <1/2$ such that
$\mf a_1 \le \varrho_{\mtt k} (x/n) \le 1-\mf a_1$ for all
$x\in B_{\mtt k}$.  Thus, by Lemma \ref{s15} with $s_1=1$,
$s_2=\ell_1$, $s_3=L_2$, $\delta = \delta_1$, on the set
$\{\mss m_{\mtt k} \in (\delta_1 , 1-\delta_1)\}$, 
\begin{align}
\label{176}
W^{\rm neg}_{x, A}\,=\,
\frac{1} {\ell_1^{2d}} \, F (\mss m_{\mtt k})
\,+\, \mf R_{n} (x) \,, \quad
\text{where}\quad
\max_{\mtt k \in \Sigma^+_{L_2}} \max_{x\in B^o_{\mtt k}} |\mf R_{n}
(x)| \le \mf a_{3, 3}\, \frac{1}{\ell_1^{2d}}  \, \frac{\ell_1}{n}\,, 
\end{align} 
and $F$ has been introduced in \eqref{167}. On the set $\mtt
1_{\delta_1} (\mss m_{\mtt k})$, we may replace $F(\cdot)$ by $F_{\mf
r}(\cdot)$. Thus, by the previous estimate and the definition of
$\overline{W}^{\rm neg}_{x, A} $, 
\begin{align*}
& \int \Big | \, \sum_{\mtt k\in \Sigma^+_{L_2}}
\mtt 1_{\delta_1} (\mss m_{\mtt k}) \,
\sum_{x\in B^o_{\mtt k}} J(x/n)\, \overline{W}^{\rm neg}_{x, A}  \, \Big |
\, f\, d\nu^n_{\rho(\cdot)}
\\
& \quad \leq
\int \Big| \, \sum_{\mtt k\in \Sigma^+_{L_2}}  
\mtt 1_{\delta_1} (\mss m_{\mtt k}) \, 
\frac{|B^o_{\mtt k}|} {\ell_1^{2d}} \,  J_{\mtt k} \, 
\Big\{\, F_{\mf r} (\mss m_{\mtt k}) \,-\,
E_{\nu^n_{\rho(\cdot) }  } \big [\, F_{\mf r} (\mss m_{\mtt k}) \,\big]\,
\Big\}\, \Big | \, 
f\, d\nu^n_{\rho(\cdot)}
\,+\, \mf R_n \,,
\end{align*}
where
$|\mf R_{n}| \le \mf a_{3, 3}\, \Vert J\Vert_\infty\, \ell_1^{-2d} \,
(\ell_1/n) \, n^d$, and $J_{\mtt k}$ is defined in \eqref{90}.

At this point, and since $\Vert F_{\mf r}\Vert_\infty \le \mf a_2$,
applying Lemma \ref{s07b}, we may remove the indicator function
$\mtt 1_{\delta_1} (\mss m_{\mtt k}) \, $ at a cost given by the
right-hand side of \eqref{170}.  The remaining term is estimated
applying Lemma \ref{s32}.  By Lemma \ref{s32} with $r_1=\ell_1$,
$\mf a (\cdot) = (|B^o_{\mtt k}|/\ell_1^{2d} ) \, J_{\cdot}$,
$\gamma = 1/\Vert J\Vert_\infty$, the sum on the left-hand side
(without the indicator) is bounded by
\begin{equation*}
\Vert J\Vert_\infty  \,
\Big\{   \, H_n(f \,|\, \nu^n_{\rho (\cdot)}  )
\,+\, \ln(2) \, \Big\} \,+\, \mf c_2 \, \Vert J\Vert_\infty  \,
\frac{1} {\ell_1^{2d}} \,  \Big( \frac{n}{\ell_1^{2}} \Big)^d
\end{equation*}
for some finite constant $\mf c_2$, depending only on
$\Vert F_{\mf r}\Vert_{C^1}$. However, by definition, this last quantity
depends only on $A$ and $\mf r$. Hence, by convention, $\mf c_2$ can
be represented by $\mf a_2$. This completes the proof of the lemma.
\end{proof}

The proof of Lemma \ref{s68} yields that in the expectation
$ E_{\nu^n_{\rho(\cdot)}} \big[\, \mtt 1_{\delta_1} (\mss m_{\mtt k})
\, W^{\rm neg}_{x, A}\,\big] $ we may replace
$W^{\rm neg}_{x, A}$ by
$\ell^{-2d}_1\, E_{\nu^n_{\rho(\cdot) } } \big [\, F_{\mf r} (\mss
m_{\mtt k}) \,\big]$ at the same cost. Hence,
\begin{equation*}
\begin{aligned}
& \int \Big| \, \sum_{\mtt k\in \Sigma^+_{L_2}}  
\mtt 1_{\delta_1} (\mss m_{\mtt k}) \,
\sum_{x\in B^o_{\mtt k}} J(x/n)\,
\Big\{ W^{\rm neg}_{x, A} \, -\,
E_{\nu^n_{\rho(\cdot)}} \big[   \mtt 1_{\delta_1} (\mss m_{\mtt k})
\, W^{\rm neg}_{x, A} \big] \, \Big\} \, \Big | \, f\,
d\nu^n_{\rho(\cdot)}
\\
& \quad \le\,
\int \Big| \, \sum_{\mtt k\in \Sigma^+_{L_2}}  
\mtt 1_{\delta_1} (\mss m_{\mtt k}) \,
\frac{|B^o_{\mtt k}|} {\ell_1^{2d}} \,  J_{\mtt k} \,
E_{\nu^n_{\rho(\cdot) }  } \big [\, F_{\mf r} (\mss m_{\mtt k}) \,\big] 
\, \Big\{\,
1 \,-\,
E_{\nu^n_{\rho(\cdot) }  } \big [\,
\mtt 1_{\delta_1} (\mss m_{\mtt k}) \,\big] \, \Big\}\, \Big | \, f\,
d\nu^n_{\rho(\cdot)}
\, + \, 2\, \mf R_n\,,
\end{aligned}
\end{equation*}
where $\mf R_n$ stands for the right-hand side in Lemma
\ref{s68}. Since
$\nu^n_{\rho(\cdot) } \big [\, \mss m_{\mtt k} \not \in (\delta_1,
1-\delta_1) \,\big]$ decays as $e^{-\mf a_1 \ell^d_1}$ by the
classical large deviations estimate for independent Bernoulli random
variables, the right-hand side of the previous equation is bounded by
$3\, \mf R_n$.  As this estimate holds for all terms composing
$\mtt T^{(2)}_{\mtt k, x+A}$,
\begin{equation}
\label{174}
\int \Big| \, \sum_{\mtt k\in \Sigma^+_{L_2}}  
\mtt 1_{\delta_1} (\mss m_{\mtt k}) \,
\sum_{x\in B^o_{\mtt k}} J(x/n)\,
\Big\{ W^{{\rm e},2}_{x,A} \, -\,
E_{\nu^n_{\rho(\cdot)}} \big[   \mtt 1_{\delta_1} (\mss m_{\mtt k})
\, W^{{\rm e},2}_{x,A} \big] \, \Big\} \, \Big | \, f\,
d\nu^n_{\rho(\cdot)} \, \le\, \mf a_2 \, \mf R_n\,,
\end{equation}
where $\mf R_n$ stands for the right-hand side in Lemma \ref{s68}, and
$W^{{\rm e},2}_{x,A}$ is defined in \eqref{173}.

\subsection*{Non-negligible terms}

In this subsection, we examine the remaining term
$W^{{\rm e},1}_{x,A}$, defined in \eqref{173}. We show that both
expressions forming $W^{{\rm e},1}_{x,A}$ can be replaced by functions
of the empirical density $\mss m_{\mtt k}$.  We start with
$E_{\nu^{\rm gc} _{\mtt k, \mss M_{\mtt k}}} [\, \tau_x \eta_A\,]$. Let
\begin{equation*}
\cb{W^{{\rm e}, a}_{x, A}} \,:=\,  
E_{\nu^{\rm gc} _{\mtt k, \mss M_{\mtt k}}} [\, \tau_x \eta_A\,] 
\,-\, \mss m_{\mtt k}^{|A|}\,, \quad
\cb{\overline{W}^{{\rm e}, a}_{x,A}} \,:=\, W^{{\rm e}, a}_{x,A}
- E_{\nu^n_{\rho(\cdot)}}  \big[ \,
\mtt 1_{\delta_1} (\mss m_{\mtt k})
\, W^{{\rm e}, a}_{x,A} \,\big] \,.
\end{equation*}
Note that we are not subtracting the mean, but an expression close to
it. Moreover,
\begin{gather*}
\sum_{x\in B^o_{\mtt k}}  J(x/n)\, W^{{\rm e}, a}_{x,A} \,=\,
\sum_{x\in B^o_{\mtt k}}  J(x/n)\,
E_{\nu^{\rm gc} _{\mtt k, \mss M_{\mtt k}}} [\, \tau_x \eta_A\,] \,-\,
|B^o_{\mtt k} | \,  J_{\mtt k}\, \mss m_{\mtt k}^{|A|}\,.
\end{gather*}

\begin{lemma}
\label{s69}
Let $\rho\colon \bb T^d \to (0,1)$ be a density profile of class
$C^3(\bb T^d)$ satisfying \eqref{83}.  Assume that
$\sup_{n\ge 1} (\ell^d_1 /n)<\infty$.  Fix a function $J$ in
$C^1(\bb T^d)$, and a finite set $A$. Then, there exist constants
$\mf a_1$, $\mf a_{3,1}$, $\mf a_{3,3}$ such that
\begin{align*}
& \int \Big|\, \sum_{\mtt k\in \Sigma^+_{L_2}}
\mtt 1_{\delta_1} (\mss m_{\mtt k}) \,
\sum_{x\in B^o_{\mtt k}}  J(x/n)\, \overline{W}^{{\rm e}, a}_{x,A} \, 
\Big |\, f\, d\nu^n_{\rho(\cdot)}
\nonumber
\\
&\quad 
\,\le\, \mf a_{3,1} \, \Vert J\Vert_\infty \,
\Big[  \, H_n (f\,|\, \nu^n_{\rho(\cdot)})
\,+\, \ln(2) \,+\, \Big( \frac{n} {\ell_1}\Big)^d\, 
e^{- \mf a_1 \ell_1^d}\,\Big] 
\,+\, \mf a_{3,3}\, \Vert J\Vert_{C^1} \,
\Big( \frac{\ell_1} {n}\Big)^2\, n^d
\end{align*}
for all densities $f$ with respect to $\nu^n_{\rho (\cdot)} $,
$\ell_1 \ge \ell_A$, $L_2\ge 1$.
\end{lemma}

\begin{proof}
As $\mss m_{\mtt k}$ belongs to the interval $(\delta_1, 1-\delta_1)$,
by Corollary \ref{s35}, with $s_1=1$, $s_2=\ell_1$, $\delta=\delta_1$,
\begin{align*}
\sum_{x\in B^o_{\mtt k}}  J(x/n)\,
W^{{\rm e}, a}_{x,A} 
\,=\, W_J( \mtt k) \, \mss m_{\mtt k}^{|A|-1}\,
\chi(\mss m_{\mtt k}) \,+\, \mf R_n (\mtt k)\,,
\end{align*}
where
\begin{equation*}
W_J( \mtt k) \,=\,
|B_{\mtt k}^o| \, J_{\mtt k}\, 
\frac{1} {\chi(\bar\rho_{\mtt k})}\,
\sum_{y\in A} \frac{y}{n}\cdot (\nabla \rho)(x_{\mtt k}/n)\,,
\end{equation*}
and
\begin{equation*}
\max_{\mtt k \in \Sigma^+_{L_2}}
\max_{\delta_1 \le \mss m_{\mtt k}  \le 1- \delta_1 } 
|\, \mf R_n (\mtt k)\, | \,\le\,  \mf a_{3, 3}\, \Vert J\Vert_{C^1} \,
\, (\ell_1/n)^2\, \ell^d_1 \,.
\end{equation*}

We may also apply Corollary \ref{s35} to the expectation of the
previous sum multiplied by the indicator to get that
\begin{equation*}
E_{\nu^n_{\rho(\cdot)}} \Big[ \,
\mtt 1_{\delta_1} (\mss m_{\mtt k}) \, \sum_{x\in B^o_{\mtt k}}
J(x/n)\, W^{{\rm e}, a}_{x,A}\,\Big]
\,=\, 
W_J( \mtt k) \,
E_{\nu^n_{\rho(\cdot)}} \big[ \,
\mtt 1_{\delta_1} (\mss m_{\mtt k}) \,
\mss m_{\mtt k}^{|A|-1}\, \chi(\mss m_{\mtt k})
\,\big] \,+\, \mf R_n (\mtt k)\,,
\end{equation*}
where the remainder $\mf R_n (\mtt k)$ satisfies the same bound. 

Clearly,
$\max_{\mtt k} |\, W_J (\mtt k) \,| \le \Vert J\Vert_\infty\, \mf
a_{3,1}\, (\ell^d_1/n)$ which is bounded by hypothesis.  Let $G$ be
the polynomial $G(\alpha) = \alpha^{|A|-1} \, \chi(\alpha)$. After
replacing the sum
$\sum_{x\in B^o_{\mtt k}} J(x/n)\, \overline{W}^{{\rm e}, a}_{x, A}$
by
$W_J(\mtt k) \, \big \{ \, G(\mss m_{\mtt k}) \, - \,
E_{\nu^n_{\rho(\cdot)}} [\, \mtt 1_{\delta_1} (\mss m_{\mtt k}) \,
G(\mss m_{\mtt k}) \,] \, \big\}$, a step which required
$\mss m_{\mtt k}$ to be bounded away from $0$ and $1$, as $W_J$ is
uniformly bounded (cf. \eqref{56b}), we may apply Lemma \ref{s07b}, as
in \eqref{170}, to remove all the indicators of the set
$\{\mss m_{\mtt k} \in (\delta_1, 1-\delta_1)\}$ at the cost which
appears on the right-hand side of \eqref{170} with $\mf a_{3,1}$
instead of $\mf a_2$ due to the presence of the gradient
$\nabla \rho$.

Without the indicators, we may apply Lemma \ref{s32} with
$1/\gamma = \Vert J\Vert_\infty$ to conclude that
\begin{align*}
& \int \Big|\, \sum_{\mtt k\in \Sigma^+_{L_2}}  
W_J (\mtt k)  \,
\big \{ \, G(\mss m_{\mtt k}) \, - \,
E_{\nu^n_{\rho(\cdot)}} [\, G(\mss m_{\mtt k}) \,] \, \big\}
\, \Big |\, f\, d\nu^n_{\rho(\cdot)}
\,\le\, \Vert J\Vert_\infty \, \Big\{\, 
\Big[  \, H_n (f\,|\, \nu^n_{\rho(\cdot)})
\,+\, \ln(2) \, \Big] \,+\, n^{d-2}  \,\Big\} \,.
\end{align*}
To complete the proof of the lemma, it remains to recollect all
previous estimates, and to observe that $n^{d-2} \le (\ell_1/n)^2\,
n^d$. 
\end{proof}

The same proof applies to the term
$\Gamma_{\mtt k}^{-1}\, E_{\nu^{\rm gc} _{\mtt k, \mss M_{\mtt k}}}
[\, \tau_x \eta_A\,] \, \mtt T^{(1)}_{\mtt k, x+A}$. It is actually
simpler as this term is of smaller order. As in \eqref{167},
denote by $F\colon [0,1]\to \bb R$ the function that
$\Gamma_{\mtt k}^{-1}\, E_{\nu^{\rm gc} _{\mtt k, \mss M_{\mtt k}}}
[\, \tau_x \eta_A\,] \, \mtt T^{(1)}_{\mtt k, x+A}$ would be if the
density profile $\rho(\cdot)$ were constant (without the normalizing
constant $|B_{\mtt k}|$ which appears in $\Gamma_{\mtt k}$):
\begin{equation*}
F(\mss m) \,=\, 
E_{\nu^n_{\mss m}} [\, \eta_A\,] \,
\frac{1}{\chi(\mss m)} \,
\frac{1}{2} \, \big\{ \, |A| \chi (\mss m) - |A|^2 (1-\mss m)^2
\,+\,  |A| (1-\mss m) \, (1-2\mss m) \, \,\big\} \, .
\end{equation*}
This function is exactly the polynomial $\mtt q_A(\mss m)$ introduced
in \eqref{90}, \eqref{168}.

Let
\begin{gather*}
\cb{V^{{\rm e},b}_{x, A} } \,:=\, \frac{1}{\Gamma_{\mtt k} }\,
E_{\nu^{\rm gc} _{\mtt k, \mss M_{\mtt k}}} [\, \tau_x \eta_A\,] \,
\, \mtt T^{(1)}_{\mtt k, x+A} \,, \quad
\cb{W^{{\rm e},b}_{x, A} } \,:=\, V^{{\rm e},b}_{x, A} 
\,-\, \frac{1}{\ell^d_1} \,\mtt q_A (\mss m_{\mtt k})\,,
\\
\cb{\overline{W}^{{\rm e},b}_{x, A} } \,:=\,
W^{{\rm e},b}_{x, A}  \,-\,
E_{\nu^n_{\rho(\cdot)}} \big[ \, \mtt 1_{\delta_1} (\mss m_{\mtt k})
\, W^{{\rm e},b}_{x, A}  \big]\,.
\end{gather*}
Note that we are not subtracting the mean of $W^{{\rm e},b}_{ x, A}$,
but an expectation close to it. On the other hand,
\begin{gather*}
\sum_{x\in B^o_{\mtt k}} J(x/n)\, 
W^{{\rm e},b}_{x, A} \,=\,
\sum_{x\in B^o_{\mtt k}} J(x/n)\, 
\frac{1}{\Gamma_{\mtt k} }\,
E_{\nu^{\rm gc} _{\mtt k, \mss M_{\mtt k}}} [\, \tau_x \eta_A\,] \,
\, \mtt T^{(1)}_{\mtt k, x+A}
\,-\, \frac{|B^o_{\mtt k}|}{\ell^d_1}
\, J_{\mtt k}\, \mtt q_A (\mss m_{\mtt k})\,.
\end{gather*}

Fix a function $J$ in $C^1(\bb T^d)$, and a finite subset $A$ of
$\bb Z^d$. The proof of Lemma \ref{s69}, using Corollary \ref{s31}
instead of Corollary \ref{s35}, yields the following bound. Note that
while in Corollary \ref{s35} the right-hand side has one term beyond
the remainder, in Corollary \ref{s31}, due to the extra factor
$\ell_1^{-d}$ coming from $\Gamma_{\mtt k}^{-1}$, there is only the
error term. Thus, the expression $W_J(\mtt k)$ which appears in the
previous proof is absent here, and the remainder $\mf R_n(\mtt k)$,
according to Corollary \ref{s31}, is now bounded by
\begin{gather*}
\max_{\mtt k \in \Sigma^+_{L_2}}
\max_{\delta_1 \le \mss m_{\mtt k}  \le 1- \delta_1 } 
|\, \mf R_n (\mtt k)\, | \,\le\,  \mf a_{3, 3}\, \Vert J\Vert_{C^1} \,
\Big\{\, \frac{1}{n} \,+\, \Big(\frac{\ell_1}{n}\Big)^2\, \Big\}\,.
\end{gather*}
Moreover, we do not need to apply Lemma \ref{s07b} as $W_J(\mtt k)$ is
not present here. We get, therefore, that there exist constants
$\mf a'_1$, $\mf a_1$, $\mf a_{3,3}$ such that
\begin{equation}
\label{169}
\int \Big|\, \sum_{\mtt k\in \Sigma^+_{L_2}}
\mtt 1_{\delta_1} (\mss m_{\mtt k}) \,
\sum_{x\in B^o_{\mtt k}} J(x/n)\, \overline{W}^{{\rm e},b}_{x,A} 
\Big |\, f\, d\nu^n_{\rho(\cdot)}
\,\le\, \mf a_{3,3}\, \Vert J\Vert_{C^1} \,
\Big\{\, \frac{1}{n} \,+\, \Big(\frac{\ell_1}{n}\Big)^2\, \Big\}\,
\Big(\frac{n}{\ell_1}\Big)^d
\end{equation}
for all densities $f$ with respect to $\nu^n_{\rho (\cdot)} $,
$\ell_1\ge \ell_A$, $L_2\ge 1$.

\begin{remark}
This remainder should be smaller than that obtained in Lemma \ref{s69}
due to the extra factor $\ell^{-d}_1$ coming from
$\Gamma_{\mtt k}^{-1}$.  Indeed, in \eqref{169} we have
$(\ell_1/n)^2 \ell^{-d}_1$ instead of the $(\ell_1/n)^2$ that appears
in Lemma \ref{s69}. However, there is also the term
$n^{-1} \ell^{-d}_1$, which is not comparable to $(\ell_1/n)^2$. This
happens because in Corollary \ref{s31} we roughly estimated the
first-order term of Lemma \ref{s15}, which is sufficient for our
purposes. We could have taken this term into account and proceeded as
in Lemma \ref{s69} to obtain a better estimate.
\end{remark}

\begin{proof}[Proof of Proposition \ref{s36}]
By \eqref{170}, Lemma \ref{s33} and the paragraph right after the
statement of this lemma, it is enough to estimate \eqref{171}.

In \eqref{175} we decomposed $W^{\rm exp}_{x,A}$ as the sum of
$W^{{\rm e}, 1}_{x,A}$ and $W^{{\rm e}, 2}_{x,A}$. Equation
\eqref{174} provides an estimate for
$W^{{\rm e}, 2}_{x,A} \, -\, E_{\nu^n_{\rho(\cdot)}} \big[ \mtt
1_{\delta_1} (\mss m_{\mtt k}) \, W^{{\rm e},2}_{x,A} \big] $.  Lemma
\ref{s69} and equation \eqref{169}, noting the form of $W_{\mtt k, A}$
in \eqref{168} and $F(\mss m) = \mtt q_A(\mss m)$ in the paragraph
after the proof of Lemma \ref{s69}, present an estimate for
$[\, W^{{\rm e}, 1}_{x,A} - W_{\mtt k, A}\,] \, -\,
E_{\nu^n_{\rho(\cdot)}} \big[ \mtt 1_{\delta_1} (\mss m_{\mtt k}) \,
[\, W^{{\rm e},1}_{x,A} - W_{\mtt k, A}\,] \,\big] $.

To complete the proof of the proposition, it remains to recollect the
remainders obtained in these steps, and to observe that
$(\ell_1/n) \, \ell^{-2d}_1 \le (1/n) \ell^{-d}_1$.
\end{proof}

\section{Multiscale analysis, a recursive bound}
\label{sec9}

In this section we present a recursive argument to augment the length
of the cube obtained in Section \ref{sec4}. Throughout this section
$\rho\colon \bb T^d \to [0,1]$ is a density profile of class $C^3$
satisfying \eqref{83}. We adopt the same convention as in the previous
section regarding the constants $\mf a_1$, $\mf a_2$, $\mf a_{3, p}$.

Fix positive integers $\mtt t_i$, $1\le i\le 3$, and let
$n=\mtt t_ 1\mtt t_2 \mtt t_3$. Divide the torus $\bb T^d_n$ into
$(\mtt t_2\mtt t_3)^d$ disjoint cubes of side length $\mtt t_1$,
represented by $\cb{B_{\mtt k}}$,
$\mtt k =(k_1, \dots, k_d) \in \Sigma^+_{\mtt t_2\mtt t_3}$, and in
$\mtt t_3^d$ disjoint cubes of side length $\mtt t_1\mtt t_2$,
represented by $\cb{B_{2, \mtt j}}$,
$\mtt j =(j_1, \dots, j_d) \in \Sigma^+_{\mtt t_3}$.  Thus
$\bb T^d_n = \cup_{\mtt k} B_{\mtt k} = \cup_{\mtt j} B_{2, \mtt j}
$.  For a cube $B_{2, \mtt j}$, denote by $\Lambda_{\mtt j}$ the set
of indices $\mtt k$ of the sets $B_{\mtt k}$ which are contained in
$B_{2, \mtt j}$:
\begin{equation*}
\cb{\Lambda_{\mtt j}} \,:=\,
\big\{ \mtt k \in \Sigma^+_{\mtt t_2\mtt t_3} :  B_{\mtt k} \subset
B_{2, \mtt j} \, \big\}\,, \;\;\text{so that}\;\;
B_{2, \mtt j}  \,=\, \bigcup_{\mtt k \in  \Lambda_{\mtt j}} B_{\mtt k} \,.
\end{equation*}
Let
\begin{equation*}
\cb{\mss M_{\mtt k}} := \sum_{x\in B_{\mtt k}} \eta_x \,,\quad
\cb{\mss m_{\mtt k}} := \frac{\mss M_{\mtt k}} {|B_{\mtt k}| }\,, \quad 
\cb{\mss M_{2, \mtt j}} := \sum_{x\in B_{2, \mtt j}} \eta_x \,, \quad 
\cb{\mss m_{2, \mtt j}} := \frac{\mss M_{2, \mtt j}} {|B_{2, \mtt j}| }
\,,\quad
\quad
\mtt j \in \Sigma^+_{\mtt t_3}\,,\;\;
\mtt k \in \Lambda_{\mtt j} \,.
\end{equation*}

Let $\mf p\colon [0,1] \to \bb R$ be a polynomial.  The goal of the
section is to replace $\mf p (\mss m_{\mtt k})$ by
$\mf p (\mss m_{2, \mtt j})$ in Theorem \ref{s01}, where
$\mf p(\cdot)$ represents the polynomials $\mtt p_{h}(\cdot)$ and
$\mtt q_h(\cdot)$.  For a polynomial $\mf p$, let $W_{\mf p}$ be
the polynomial given by
\begin{equation}
\label{111}
\begin{aligned}
\cb{ W_{\mf p} } (\mss m) \, :=\, 
\mf p(\mss m)  \,+\, \frac{1}{2} \, \frac{1}{\mtt t_1^d} \,
\mf p'' (\mss m)\, \chi (\mss m) \,.
\end{aligned}
\end{equation}
When $\mf p(\mss m)=\mss m^p$ we represent $W_{\mf p} $ by
$\cb{\mtt V_{p}}$.

\begin{remark}
\label{r09}
Suppose that $\mf p$ is the polynomial $\mf p_A$ introduced in the
statement of Theorem \ref{s01}:
$\cb{\mf p_A (\rho)} := \mtt p_A (\rho) + \ell_1^{-d} \mtt q_A
(\rho)$.  Substituting $\mf p$ by $\mf p_A$ in equation \eqref{111}
with $\mtt t_1 = \ell_1$, a cancellation of the $\ell_1^{-d}$- order
terms occurs and yields that
\begin{equation*}
W_{\mf p_A}  (\mss m) \,=\, \mtt p_A(\mss m) \,+\, \frac{1}{2\,\ell_1^{2d}}\,
\mtt q_{A}''(\mss m)\, \chi(\mss m)\,.
\end{equation*}

This cancellation simplifies considerably the proof in dimension $d=3$
as it eliminates the problematic terms of order $\ell_1^{-d}$ which
would require careful attention. On the one hand, it is a fortuitous
cancellation as the correction terms in Theorems \ref{s01} and
\ref{s26} have different origins.  The corrector term
$\ell_1^{-d} \mtt q_A (\rho)$ in Theorem \ref{s01} is produced by the
first order term in the equivalence of ensembles, while the corrector
term \eqref{111} in Theorem \ref{s26} is produced by the expansion of
$E_{\nu_\rho} [ \mss m_{\mtt k}^{|A|}]$ as a function of $\rho$. The
$0$-order term is $\rho^{|A|}$ and the first order one is
$(1/2) |A| (|A|-1) \rho^{|A|-2} \chi(\rho)\, |B_{\mtt k}|^{-1}$.

On the other hand, in our multi-scale scheme, we replace $\mtt p_A$ by
its canonical expectation on a block of width $\ell_1$.  If now,
instead of approximating by the equivalence of ensembles and then a
further reduction of the associated grand canonical quantity, we would
state the approximation by the conditional expectation on a block of
width $\ell_1\ell_2$ directly, the cancellation would not be seen at
all.
\end{remark}

The main result of this section reads as
follows.

\begin{theorem}
\label{s26}
Let $\rho\colon \bb T^d \to [0,1]$ be a density profile of class $C^3$
satisfying \eqref{83}.  Fix a polynomial
$\mf p\colon [0,1] \to \bb R$, a function $J\colon \bb T^d \to \bb R$
of class $C^1$, and constants $\delta>0$, $\epsilon>0$, $K_1\ge 1 $,
$K_2\ge 1$.  Then, there exist finite constants
$\mf a_{3,3} = \mf a_{3,3} (\mf p)$,
$\mf b_1 = \mf b_1(\mf r, \Vert \rho\Vert_{C^3}, \mf p)$, $C_1$, and
$C_2 = C_2(\mf p)$ such that
\begin{equation*}
\begin{aligned}
& \int
\sum_{\mtt j \in \Sigma^+_{\mtt t_3}}
\sum_{\mtt k \in \Lambda_{\mtt j}}  \mtt t_1^d\, J_{\mtt k}\,
\Big\{\,   \mf p(\mss m_{\mtt k})
- E_{\nu^n_{\rho(\cdot)}} [ \, \mf p(\mss m_{\mtt k})\, ] 
\,-\, 
\big(\, W_{\mf p} (\mss m_{2, \mtt j})
- E_{\nu^n_{\rho(\cdot)}} [ W_{\mf p} (\mss m_{2, \mtt j})] \,\big)
\, \Big\}\,
f\, d\nu^n_{\rho(\cdot)}
\\
&\quad  \le\, \delta\, n^2 \,  I_n  ( f \,;\, \nu^n_{\rho (\cdot)} )
\,+\, 
\mf a_{3,3} \, \Vert J\Vert_{C^1} 
\,\Big\{   \, H_n (f\,|\, \nu^n_{\rho(\cdot)})\, + \, \ln(2)\,\Big\}
\,+\, C_1 \,  n^d\, 
\Big(\frac{\mtt t_1\mtt t_2}{n}\Big)^4
\\
&\quad 
\,+\, \mf a_{3,3}\, \Vert J\Vert^2_{C^1} \,  n^d\, \Big\{ \,
\frac{1}{\delta}\, 
\frac{\mtt t_2^{d+2}}{\mtt t_1^{d-2}}   \, \frac{1}{n^2} 
\,+\,  e^{C_2 K_2^d}\, \frac{1 }{(\mtt t_1\mtt t_2)^{d}} 
\, \frac{1 }{\mtt t_1^{d}}   \, \Big\}
\\
&\quad 
\,+\, \mf a_{3,3}\, \Vert J\Vert_{C^1} \,  n^d\,
\Big\{ 
\Big(\frac{\mtt t_1\mtt t_2} {n}\Big)^2 \, \frac{1}{(\mtt t_1\mtt
t_2)^{d/2}}
\,+\, 
\frac{\mtt t_1\mtt t_2 }{n} \, \frac{1} {(\mtt t_1\mtt t_2)^d } \,
\,+\, \Big(\frac{\mtt t_1 }{n} \Big)^2 \, \frac{1 }{\mtt t_1^{d}}  \Big\}
\end{aligned}
\end{equation*}
for all integers $\mtt t_1\ge 2$, $\mtt t_2\ge 2$ such that
$(\mtt t_1\mtt t_2)^{d+3} \le K_1\, n^3$,
$(\mtt t_1\mtt t_2)^{1+(d/4)} \le n$, $\mtt t_2 \le K_2\, \mtt t_1$,
$\mtt t_2^{2+d} \le (\delta / \mf c_{\rm LS}) (n/\mtt t_1)^2 \, (\mf b_1
\wedge \mtt t_1^{-\epsilon})$, and densities $f$ with respect to
$\nu^n_{\rho(\cdot)}$. In this formula,
$C_1 = \exp\{\mf c_{\rm LS} K_1 [1+\mf a_{3,3}]\,[1+\Vert
J\Vert_{C^1}]/\delta\} \, \exp\{\mf a_{3,3} K_2^d\}$.
\end{theorem}

\begin{remark}
In this article, Theorem \ref{s26} will only be used with $K_2=1$. We
stated it in this generality as it might be useful in other models. On
the other hand, in $d=3$, in the second step $K_1$ is chosen equal to
$1$, and in the third $K_1=\kappa^6$, where $\kappa$ is a constant
which converges to $+\infty$ after $n\to +\infty$.
\end{remark}

\begin{remark}
\label{r07}
We observed in Remark \ref{rm3} that Theorem \ref{s01} and Corollary
\ref{s46} reach sufficiently large cubes in dimension $1$, since they
allow intervals of length of order $\epsilon_n n^{2/3}$, while the
proof of the non-equilibrium fluctuations only requires cubes of
length $K \, n^{1/2}$, for some $1\ll K\ll n$. Thus, Theorem
\ref{s26}, though valid in dimension $1$, will only be used in
dimensions $2$ and $3$. Corollary \ref{s48} states that in dimension
$2$, the correction term in $W_{\mf p} (\mss m_{2, \mtt j})$ appearing
in the statement of Theorem \ref{s26} is not needed.
\end{remark}

\begin{remark}
As observed above, the correction term in
$W_{\mf p} (\mss m_{2, \mtt j}) $ is the same as the one obtained in
Theorem \ref{s01} in the first step of the multiscale analysis, but
with a negative sign. It comes from the replacement of
$E_{\nu^n_{\rho(\cdot)}} [\mf p(\mss m_{\mtt k})]$ by
$\mf p(E_{\nu^n_{\rho(\cdot)}} [\mss m_{\mtt k}])$. We need
$\mf p(\cdot)$ to be a polynomial because we are only able to estimate
powers of $\mss m_{\mtt k}$, that is, $\mss m_{\mtt k}^p$, $p\ge 1$.
\end{remark}

\begin{remark}
One of the main challenges in the proof of the hydrodynamic limit or
the convergence of the density field to the solution of some SPDE
consists in replacing the average of local functions by a function of
the density over a large box. Theorem \ref{s26} provides a method to
increase the width of the box. Its proof does not rely on specific
properties of the model. It might therefore be of independent
interest.
\end{remark}

\begin{remark}
In dimension $d=3$, the condition
$(\mtt t_1\mtt t_2)^{d+3} \le K_1\, n^3$, that is,
$\mtt t_1\mtt t_2 \le K'_1\, n^{1/2}$ prevents us from reaching cubes
of length larger than $K'_1 \, n^{1/2}$.
\end{remark}

As mentioned in Remark \ref{r07}, in dimension $2$, we can estimate
the correction term  $W_{\mf p} (\mss m_{2, \mtt j}) $ applying
Lemma \ref{s45}. This is the content of the next result. In this
result we set $K_2=1$ to avoid the extra constant $K_2$ multiplying
the entropy.

\begin{corollary}
\label{s48}
Fix a polynomial $\mf p\colon [0,1] \to \bb R$, a function
$J\colon \bb T^d \to \bb R$ of class $C^1$, and constants $\delta>0$,
$\epsilon>0$, $K_1\ge 1$.  Then, there exist finite
constants $\mf a_{3,3} = \mf a_{3,3} (\mf p)$,
$\mf b_1 = \mf b_1(\mf r, \Vert \rho\Vert_{C^3}, \mf p)$ such that
\begin{equation*}
\begin{aligned}
& \int
\sum_{\mtt j \in \Sigma^+_{\mtt t_3}}
\sum_{\mtt k \in \Lambda_{\mtt j}}  \mtt t_1^d\, J_{\mtt k}\, \Big\{\, 
\mf p(\mss m_{\mtt k}) 
- E_{\nu^n_{\rho(\cdot)}} [\, \mf p(\mss m_{\mtt k})\, ] 
\,-\, 
\big(\, \mf p(\mss m_{2, \mtt j})
- E_{\nu^n_{\rho(\cdot)}} [\, \mf p(\mss m_{2, \mtt j}) \, ] \,\big)
\, \Big\}\,
f\, d\nu^n_{\rho(\cdot)}
\\
&\quad
\le\, \delta\, n^2 \,  I_n  ( f \,;\, \nu^n_{\rho (\cdot)} )
\,+\, 
\mf a_{3,3} \, \Vert J\Vert_{C^1} 
\,\Big\{   \, H_n (f\,|\, \nu^n_{\rho(\cdot)})\, + \, \ln(2)\,\Big\}
\,+\, \mtt C (J)\,  n^d\, 
\Big(\frac{\mtt t_1\mtt t_2}{n}\Big)^4
\\
&\quad 
\,+\, \mf a_{3,3}\, (1+\Vert J\Vert_{C^1})^2 \,  n^d\, \Big\{ \,
\frac{1}{\delta}\, 
\frac{\mtt t_2^{d+2}}{\mtt t_1^{d-2}}   \, \frac{1}{n^2} 
\,+\,  \frac{1 }{(\mtt t_1\mtt t_2)^{d/2}} 
\, \frac{1 }{\mtt t_1^{d}}   \, \Big\}
\\
&\quad 
\,+\, \mf a_{3,3}\, \Vert J\Vert_{C^1} \,  n^d\,
\Big\{ 
\Big(\frac{\mtt t_1\mtt t_2} {n}\Big)^2 \, \frac{1}{(\mtt t_1\mtt
t_2)^{d/2}}
\,+\, 
\frac{\mtt t_1\mtt t_2 }{n} \, \frac{1} {(\mtt t_1\mtt t_2)^d } \,
\,+\, \Big(\frac{\mtt t_1 }{n} \Big)^2 \, \frac{1 }{\mtt t_1^{d}}  \Big\}
\end{aligned}
\end{equation*}
for all integers $\mtt t_1\ge 2$, $\mtt t_2\ge 2$ such that
$(\mtt t_1\mtt t_2)^{d+3} \le K_1\, n^3$,
$(\mtt t_1\mtt t_2)^{1+(d/4)} \le n$, $\mtt t_2 \le \mtt t_1$,
$\mtt t_2^{2+d} \le (\delta / \mf c_{\rm LS}) (n/\mtt t_1)^2 \, (\mf b_1
\wedge \mtt t_1^{-\epsilon})$, and densities $f$ with respect to
$\nu^n_{\rho(\cdot)}$. In this formula,
$\mtt C(J)  = \exp\{\mf c_{\rm LS} K_1 [1+\mf a_{3,3}]\,[1+\Vert
J\Vert_{C^1}]/\delta\}$.
\end{corollary}

\begin{proof}
In view of Theorem \ref{s26}, we need to estimate the correction term
in $W_{\mf p} (\mss m_{2, \mtt j})$. Recall that 
$\mf q (x) :=  -\, (1/2)\, \mf p'' (x)\, \chi (x)$,
$J_{2, \mtt j} := \sum_{\mtt k \in \Lambda_{\mtt j}} J_{\mtt
k}$. Clearly,
$\Vert J_{2, \mtt j}\Vert_\infty \le \mtt t_2^{d} \, \Vert
J\Vert_\infty$. As $\mtt t_2\le \mtt t_1$,
$\Vert J_{2, \mtt j}\Vert_\infty \le \Vert J\Vert_\infty \, (\mtt
t_1\mtt t_2)^{d/2}$.  Therefore, by Lemma \ref{s45} with
$r_1=\mtt t_1\mtt t_2$ and $\mf a(\mtt j) = J_{2, \mtt j}$,
using  $\mtt t_2\le \mtt t_1$ again,
\begin{align*}
& \int \sum_{\mtt j \in \Sigma^+_{\mtt t_3}}
J_{2, \mtt j}
\, \Big\{\, \mf q( \mss m_{2,\mtt j})
- E_{\nu^n_{\rho(\cdot)}} [\, \mf q( \mss m_{2,\mtt j})\, ]
\, \Big\}\, f \, d\nu^n_{\rho(\cdot)}
\\
&\quad \,\le \,
C(\mf p)\,  \Vert J\Vert_\infty \, \Big\{\, 
\, H_n(f \,|\, \nu^n_{\rho (\cdot)}  )
\,+\, 
\Big(\frac{1}{\mtt t_1\mtt t_2} \Big)^{d/2}\, 
\Big(\frac{n}{\mtt t_1} \Big)^d \,\Big\} \,.
\end{align*}
This additional error is to be added to the right-hand side of Theorem
\ref{s26} statement.  This error is larger than the existing error
$(n/\mtt t_1)^d\, (\mtt t_1\mtt t_2)^{-d}$, and replaces it in the
statement of Corollary \ref{s48}. This completes the proof.
\end{proof}

The proof of Theorem \ref{s26} is divided into three steps.  We start by
excluding large deviations of the density. Then, we compare the sum
$\sum_{\mtt k \in \Lambda_{\mtt j}} \mtt t_1^d\, J_{\mtt k}\, \mf p(\mss
m_{\mtt k})$ to its expectation with respect to the canonical measure
on $B_{2, \mtt j} = \cup_{\mtt k \in \Lambda_{\mtt j}} B_{\mtt
k}$. In the last step, we replace the canonical measure by
the grand-canonical one and show that the expectation with respect to the
grand-canonical measure in $B_{2, \mtt j}$ of $\mf p(\mss m_{\mtt k})$
is close to $\mf p(\mss m_{2, \mtt j})$. \smallskip

We first restrict the density to values bounded away from $0$ and $1$.
By Lemma \ref{s07b} with $r_1 = \mtt t_1 \mtt t_2$ and
$\mf a_0=\mf r/2$, there exists a positive constant
$\mf a_1 = \mf a_1(\mf p)$ such that
\begin{equation}
\label{67}
\begin{aligned}
& \sum_{\mtt j \in \Sigma^+_{\mtt t_3}}
\int \mtt 1\{ \mss m_{2, \mtt j} \not\in (\mf r_-, \mf r_+) \}
\, \sum_{\mtt k \in \Lambda_{\mtt j}} J_{\mtt k}\,
\mtt t^d_1\,   \big[\, \mf p(\mss m_{\mtt k}) \,-\,
\mf p(\mss m_{2, \mtt j})\,\big]  \, f\;
d\nu^n_{\rho(\cdot)}
\\
&\quad \,\le\,
\mf a_1 \, \Vert J\Vert_\infty  \, \Big\{ \, H_n (f\,|\, \nu^n_{\rho(\cdot)})
\,+\, \Big( \frac{n} {\mtt t_1 \mtt t_2} \Big)^d\, 
e^{-   (\mtt t_1 \mtt t_2)^d/\mf a_1} \, \Big\}
\end{aligned}
\end{equation}
for all $\mtt t_i\ge 1$, $i=1$, $2$, $3$. In this formula and below,
$\cb{\mf r_-:=\mf r/2}$, $\cb{\mf r_+:=1- (\mf r/2)}$. This estimate
will apply to the correction term of $W_{\mf p}$ and to the
expectation term $E_{\nu^n_{\rho(\cdot)}}[\mf p(\mss m_{\mtt k})]$
appearing in the statement of Theorem \ref{s26}.

\subsection*{Canonical measures}

In this subsection, we replace $\mf p(\mss m_{\mtt k})$ by its
expectation with respect to the canonical measure on the larger cube
$B_{2,\mtt j}$. This is the content of Proposition \ref{s37}, whose
statement requires some notation.

For $\mtt j \in \Sigma^+_{\mtt t_3}$, let
$R^{(2)}_{\mtt j} \colon \bb R \to (0,1)$ be the strictly increasing
function given by
\begin{equation}
\label{42}
\cb{R^{(2)}_{\mtt j}  (\varphi)} \,:=\, \frac{1}{|B_{2, \mtt j} |}\, 
\sum_{x\in B_{2, \mtt j} }  \frac{e^\varphi \rho(x/n)  }
{e^\varphi \rho(x/n)  + [1- \rho(x/n)] } \,\cdot
\end{equation}
The parameter $\varphi$ is called the chemical potential.  Denote by
$\cb{\Phi^{(2)}_{\mtt j} (\cdot)}$ the inverse of
$R^{(2)}_{\mtt j} (\cdot)$, and let
\begin{equation}
\label{97}
\cb{\varrho^{(2)}_{\mtt j, \vartheta} (x/n)} \,:=\, 
\frac{e^{\Phi^{(2)}_{\mtt j} (\vartheta) }  \rho(x/n)  }
{e^{\Phi^{(2)}_{\mtt j} (\vartheta) }  \rho(x/n)  + [1- \rho(x/n)] }\,,
\quad x\in B_{2, \mtt j}\,, \;\;\text{where}\;\;  \vartheta\in (0,1) \,.
\end{equation}
By definition,
\begin{equation}
\label{50}
\frac{1}{|B_{2, \mtt j} |}\, 
\sum_{x\in B_{2, \mtt j} }  \varrho^{(2)}_{\mtt j, \vartheta } (x/n)
\, =\, \vartheta
\end{equation}
for all $0<\vartheta<1$. We adjusted the chemical potential for the
particle density on the cube $B_{2, \mtt j}$ to be $\vartheta$.

Denote by $\nu^{\rm gc} _{2,\mtt j, \mss M}$ the Bernoulli product
measure on $\cb{\Omega_{2, \mtt j} }:= \{0,1\}^{B_{2, \mtt j}}$
with marginals given by
\begin{equation*}
\cb{ \nu^{\rm gc} _{2,\mtt j, \mss M} } \{ \eta_x = 1\} \,=\,
\varrho^{(2)}_{\mtt j, \mss m} (x/n)\,, \quad x\in B_{2, \mtt j} \,,
\quad \text{where}\;\;
\mss m \,=\,  \mss M/|B_{2, \mtt j}| \,.
\end{equation*}
Let $\cb{\nu^{B_{2,\mtt j} } _{\rho(\cdot)}}$ be the marginal on
$\Omega_{2, \mtt j}$ of the product measure $\nu^{n} _{\rho(\cdot)}$,
and $\nu^{\rm c} _{2,\mtt j, \mss M} $ be the measure
$\nu^{B_{2, \mtt j}} _{\rho(\cdot)}$ conditioned on configurations
with $\mss M$ particles:
\begin{equation}
\label{84}
\cb{ \nu^{\rm c} _{2,\mtt j, \mss M}  (\eta) }\,:=\,
\nu^{B_{2, \mtt j} } _{\rho(\cdot)}
\Big(\, \eta \; \big|\, \sum_{x\in B_{2, \mtt j} } \eta_x = \mss M\,\Big)\,.
\end{equation}
We are now in a position to state the main result of this subsection.
Recall that $\mf a_0 = \mf r/2$.

\begin{proposition}
\label{s37}
Fix $\delta>0$, $\epsilon>0$, a polynomial $\mf p(\cdot)$, a function
$J\colon \bb T^d \to \bb R$ of class $C^1$, and $K_1\ge 1$.  Then,
there exist finite constants $\mf a_{3,1}$, $\mf a_{3,3}$ and
$\mf b_1 = \mf b_1(\mf r, \Vert \rho\Vert_{C^3}, \mf p)$ such that
\begin{equation*}
\begin{aligned}
\sum_{\mtt j \in \Sigma^+_{\mtt t_3}}
\int  \mtt 1\{ \mss m_{2, \mtt j} \in (\mf r_- ,  \mf r_+) \}
\sum_{\mtt k \in \Lambda_{\mtt j}} \mtt t^d_1\,   J_{\mtt k}\,
\big\{\, \mf p(\mss m_{\mtt k}) \,-\,
E_{\nu^{\rm c} _{2,\mtt j, \mss M} } [\, \mf p(\mss m_{\mtt k}) 
\,]  \,\big\} \, f\;
d\nu^n_{\rho(\cdot)}
\\
\,\le\,
\delta\, n^2 \,  I_n  ( f \,;\, \nu^n_{\rho (\cdot)} ) \,+\,
\frac{\mf a_{3,3}}{\delta}\,
\Big\{\, \mf c_{\rm LS} \, \mf C_2
\Big(\frac{\mtt t_1\mtt t_2}{n}\Big)^4\, n^d
\,+\,  \Vert J\Vert^2_\infty\,
\mtt t_2^{d+2}  \, \Big( \frac{n}{\mtt t_1}\Big)^{d-2} \,\Big\}
\end{aligned}
\end{equation*}
for all densities $f$ with respect to $\nu^n_{\rho(\cdot)}$, and all
$\mtt t_1\ge 2$, $\mtt t_2\ge 2$ such that
$(\mtt t_1\mtt t_2)^{d+3} \le K_1\, n^3$,
$\mtt t_2^{2+d} \le (\delta / \mf c_{\rm LS}) \Vert J\Vert_\infty^{-1}\,
(n/\mtt t_1)^2 \, (\mf b_1 \wedge \mtt t_1^{-\epsilon})$.  In this
equation,
$\mf C_2 = \exp\{ \mf c_{\rm LS} C_0(\mf p)\, K_1 \, (\, \mf a_{3,1} +
\Vert J\Vert_{C^1})/\delta\}$.
\end{proposition}

The proof of this proposition is divided into several lemmata. We first
introduce some notation. For $\mtt j \in \Sigma^+_{\mtt t_3}$,
$\mtt k \in \Lambda_{\mtt j} $, write
\begin{align*}
\mf p(\mss m_{\mtt k}) \,-\,
E_{\nu^{\rm c} _{2,\mtt j, \mss M} } [\mf p(\mss m_{\mtt k}) 
\,\big] \,=\,
\widehat{\mf p} (\mss m_{\mtt k}) \,+\,  V_{\mtt j, \mtt k} \,, 
\end{align*}
where
\begin{equation}
\label{82}
\begin{gathered}
\cb {\widehat{\mf p} (\mss m_{\mtt k}) } \,:=\,
\mf p(\mss m_{\mtt k}) -
\mf p_{\mtt j, \mtt k}  (\mss m_{2, \mtt j} )
- \mf p'( \widetilde{\mss m}_{\mtt j, \mtt k} ) ( \mss m_{\mtt k} -
\mss m_{\mtt j, \mtt k}  )\,,
\quad
\cb{V_{\mtt j, \mtt k} } \,:=\,
\mf p'( \widetilde{\mss m}_{\mtt j, \mtt k} ) ( \mss m_{\mtt k} -
\mss m_{\mtt j, \mtt k}  )\,,
\\
\cb{  \mf p_{\mtt j, \mtt k} (\mss m_{2, \mtt j}) }\, :=\,
E_{\nu^{\rm c} _{2,\mtt j, \mss M} } [\mf p(\mss m_{\mtt k}) ]
\,,\quad
\cb{ \mss m_{\mtt j, \mtt k}  }\, :=\,
E_{\nu^{\rm c} _{2,\mtt j, \mss M} } [\mss m_{\mtt k} ]\,,
\quad
\cb{ \widetilde{\mss m}_{\mtt j, \mtt k}  }\, :=\,
E_{\nu^{\rm gc} _{2,\mtt j, \mss M} } [\mss m_{\mtt k} ]\,,
\end{gathered}
\end{equation}
and $\mss M = \mss M_{2, \mtt j}$.  The next lemma provides an
estimate for $\widehat{\mf p} (\mss m_{\mtt k})$.

\begin{lemma}
\label{s16}
Fix $\epsilon>0$, $0<\delta < 1/2$, and a polynomial
$\mf p\colon [0,1]\to \bb R$.  Then, there exist finite constants
$\mf a_{3,3} = \mf a_{3,3}(\mf p(\cdot), \delta)$,
$\mf a^*_{3,3} = \mf a^*_{3,3}(\mf p(\cdot), \delta)$,
$C_0 = C_0(\epsilon)$ such that
\begin{equation*}
\max_{\mtt j \in \Sigma^+_{\mtt t_3}} \max_{\mtt k \in \Lambda_{\mtt j}}
\sup_{ \delta \le  \mss m \le 1-\delta}
\ln  E_{\nu^{\rm c} _{2,\mtt j, \mss M} } 
\big[ \, e^{ \beta \, \mtt t_1^d \, \widehat{\mf p} (\mss m_{\mtt k})
}\,  \big]
\; \le  \;  \mf a_{3,3}\, C_0(\epsilon) \, \beta^2
\end{equation*}
for all $\mtt t_1\ge 2$, and
$ \beta < \min\{ 1, 1/\mf a^*_{3,3} , \mtt t_1^{-\epsilon}\}$.
\end{lemma}

\begin{proof}
Fix $0<\delta<1/2$. In this proof, all constants $\mf a_i$,
$1\le i\le 3$, are allowed to depend on $\delta$ and the polynomial
$\mf p$, but not on $\epsilon$. The proof relies on two bounds of
$\widehat{\mf p} (\cdot)$. 

We first affirm that there exists a finite constant $\mf a_{3,3}$ such
that
\begin{equation}
\label{72}
\max_{\mtt j \in \Sigma^+_{\mtt t_3}} \max_{\mtt k \in \Lambda_{\mtt j}}
\sup_{ \delta \le  \mss m \le 1-\delta}
E_{\nu^{\rm c} _{2,\mtt j, \mss M} }
\big[ \, \big\{\, \mtt t_1^d \; \widehat{\mf p}  (\mss m_{\mtt k})
\big\}^2 \big] \; \le\; \mf a_{3,3} \; .
\end{equation}
The proof of claim \eqref{72} relies on results presented in Appendix
\ref{sec12}. By Lemma \ref{s17}, in the definitions of
$\mf p_{\mtt j, \mtt k} (\cdot)$ and $\mss m_{\mtt j, \mtt k}$, we may
replace the canonical measure by the grand canonical one at a cost
bounded by $\mf a_1\, (\mtt t_1\mtt t_2)^{-d}$. It remains to estimate
the expectation in \eqref{72} with $\widehat{\mf p} (\mss m_{\mtt k})$
replaced by
\begin{align*}
\mf p(\mss m_{\mtt k})
\,-\,  E_{\nu^{\rm gc} _{2,\mtt j, \mss M} } [\mf p(\mss m_{\mtt k}) ]
- \mf p'( \widetilde{\mss m}_{\mtt j, \mtt k} )
\,\big[\, \mss m_{\mtt k} - \widetilde{\mss m}_{\mtt j, \mtt k} 
\,\big]\,.
\end{align*}
By Lemma \ref{s22} with $s_1=\mtt t_1$, $s_2=\mtt t_2$, the absolute
value of the difference between
$E_{\nu^{\rm gc} _{2,\mtt j, \mss M} } [\mf p(\mss m_{\mtt k}) ]$ and
$\mf p(\widetilde{\mss m}_{\mtt j, \mtt k})$ is bounded by
$\mf a_{3,3} \, \mtt t_1^{-d}$.  By Lemma \ref{s10}, we may estimate
the new expression (the previous displayed equation with
$E_{\nu^{\rm gc} _{2,\mtt j, \mss M} } [\mf p(\mss m_{\mtt k}) ]$
replaced by $\mf p(\widetilde{\mss m}_{\mtt j, \mtt k})$) with respect
to the grand canonical measure instead of the canonical one at a
multiplicative cost equal to $\mf a_1$. At this point, and after a
second order Taylor expansion, the bound \eqref{72} becomes a
straightforward estimate on independent inhomogeneous Bernoulli random
variables.

We next assert that
\begin{align}
\label{71}
\big|\, \widehat{\mf p}  (\mss m_{\mtt k}) \,\big| \,\le \,
\mf a_{3,3} \,\big\{ ( \mss m_{\mtt k} - \mss m_{\mtt j, \mtt k}
)^2 \, +\, \mtt t_1^{-d}\,\big\}\,.
\end{align}
To prove this bound, write $\widehat{\mf p} (\mss m_{\mtt k})$
replacing the expectations with respect to the canonical measures by
expectations with respect to the grand canonical ones:
\begin{align*}
& \widehat{\mf p}  (\mss m_{\mtt k}) \,=\,
\mf p(\mss m_{\mtt k})  \,-\,
\mf p(\mss m_{\mtt j, \mtt k}) 
\,-\, \mf p'( \mss m_{\mtt j,
\mtt k} ) ( \mss m_{\mtt k} - \mss m_{\mtt j, \mtt k}   )
\,- \, \mf R_n \,
\\
&\quad \text{where}\quad
\mf R_n \,=\,
E_{\nu^{\rm c} _{2,\mtt j, \mss M} } [\mf p(\mss m_{\mtt k}) ]
\,-\,  \mf p( \mss m_{\mtt j, \mtt k})
\,+\, \big[\,
\mf p'( \widetilde{\mss m}_{\mtt j, \mtt k} )
\,-\, \mf p'( \mss m_{\mtt j, \mtt k} ) \,\big] \, 
(  \mss m_{\mtt k} - \mss m_{\mtt j, \mtt k}   )\,.
\end{align*}
Rewrite the first difference in the definition of $\mf R_n$ as the sum of
\begin{align*}
E_{\nu^{\rm c} _{2,\mtt j, \mss M} } [\mf p(\mss m_{\mtt k}) ]
\,-\,  E_{\nu^{\rm gc} _{2,\mtt j, \mss M} } [\mf p(\mss m_{\mtt k}) ]
\,+\,  E_{\nu^{\rm gc} _{2,\mtt j, \mss M} } [\mf p(\mss m_{\mtt k}) ]
\,-\, \mf p(\widetilde{\mss m}_{\mtt j, \mtt k}) 
\,+\, \mf p(\widetilde{\mss m}_{\mtt j, \mtt k}) 
\,-\, \mf p(\mss m_{\mtt j, \mtt k}) \,.
\end{align*}
By Lemma \ref{s17}, the absolute value of the first difference is
bounded by $\mf a_1 / (\mtt t_1\mtt t_2)^d$. Since $\mf p(\cdot)$ and
$\mf p'(\cdot)$ are Lipschitz continuous on $[0,1]$, the same result
with $p=1$ yields that
$|\, \mf p(\widetilde{\mss m}_{\mtt j, \mtt k}) \,-\, \mf p(\mss
m_{\mtt j, \mtt k}) \,| \le \mf a_1 / (\mtt t_1\mtt t_2)^d$, and
$|\, \mf p'(\widetilde{\mss m}_{\mtt j, \mtt k}) \,-\, \mf p'(\mss
m_{\mtt j, \mtt k}) \,| \le \mf a_1 / (\mtt t_1\mtt t_2)^d$.  On the
other hand, by Lemma \ref{s22},
$|\, E_{\nu^{\rm gc} _{2,\mtt j, \mss M} } [\mf p(\mss m_{\mtt k}) ] -
\mf p(\widetilde{\mss m}_{\mtt j, \mtt k})| \le \mf a_{3,3} / \mtt
t_1^d$. To complete the proof of \eqref{71}, apply  a second order
Taylor expansion and Young's inequality $2ab\le a^2 + b^2$.

We turn to the proof of the lemma. Suppose first that
$\beta \, \mtt t_1^{3d}\le 1$. Expand the exponential up to the third
order to obtain that the expectation appearing in the statement of the
lemma is bounded above by
\begin{equation*}
1 \; +\;  \beta \,  E_{\nu^{\rm c} _{2,\mtt j, \mss M} }
\big[ \mtt t_1^d \, \widehat{\mf p}  (\mss m_{\mtt k}) \big] \;
+\; \beta^2  \, E_{\nu^{\rm c} _{2,\mtt j, \mss M} }
\big[ \big( \mtt t_1^d \, \widehat{\mf p}  (\mss m_{\mtt k})
 \big)^2 \big]
\; + \; \mf a_1\,  \beta^2
\end{equation*}
because $\beta \, \mtt t_1^{3d}\le 1$.  By definition of
$\widehat{\mf p} (\cdot)$, the linear term vanishes. By \eqref{72},
the second-order term is bounded by $\mf a_{3,3}\, \beta^2 $. To
complete the proof in the case $\beta \, \mtt t_1^{3d}\le 1$, it
remains to recall that $\ln (1+a)\le a$ for all $a>0$.

We turn to the case $\beta \, \mtt t_1^{3d} > 1$.  Fix $\epsilon>0$,
and let $\cb{\kappa } := \epsilon/8$.  We estimate the expectation of
$\exp\{ \beta \, \mtt t_1^d\, \widehat{\mf p} (\mss m_{\mtt k})\,\}$
according to the value of $\mss m_{\mtt k} - \mss m_{\mtt j, \mtt k}$.
We first consider the case in which
$|\mss m_{\mtt k} - \mss m_{\mtt j, \mtt k} |\ge \mtt
t_1^{-(d/2)+\kappa}$. In this case, the expectation
\begin{equation*}
E_{\nu^{\rm c} _{2,\mtt j, \mss M} } 
\Big[ \, e^{ \beta \, \mtt t_1^d\, \widehat{\mf p} (\mss m_{\mtt k}) }\,
\mtt 1 \big\{ |\mss m_{\mtt k} - \mss m_{\mtt j, \mtt k} |\ge 
\mtt t_1^{-(d/2)+\kappa}\,\big\}\, \Big]
\end{equation*}
can be written as the sum of two terms. The first one contains the
values of $\mss m_{\mtt k} $ such that
$\mss m_{\mtt k} - \mss m_{\mtt j, \mtt k} \ge
\mtt t_1^{-(d/2)+\kappa}$, and the second one the values such
that
$\mss m_{\mtt k} - \mss m_{\mtt j, \mtt k} \le -\,
\mtt t_1^{-(d/2)+\kappa}$. As both terms are estimated similarly, we
consider only the first one which is equal to
\begin{equation*}
\sum_{\mss n \ge \mss m_{\mtt j, \mtt k} +
\mtt t_1^{-(d/2)+\kappa}}
e^{ \beta \, \mtt t_1^d\, \widehat{\mf p} (\mss n) }\,
\nu^{\rm c} _{2,\mtt j, \mss M} 
\big\{ \mss m_{\mtt k} =\mss n \, \big\}
\,=\,
\sum_{\mss N \ge \mss M_{\mtt j, \mtt k} +
\mtt t_1^{(d/2)+\kappa}}
e^{ \beta \, \mtt t_1^d\, \widehat{\mf p} (\mss N/\mtt t^d_1) }\,
\nu^{\rm c} _{2,\mtt j, \mss M} 
\big\{ \mss M_{\mtt k} =\mss N \, \big\} \,.
\end{equation*}
In this formula,
$\mss M_{\mtt j, \mtt k} = \mtt t^d_1 \, \mss m_{\mtt j, \mtt k}
$. The first sum is performed over elements $\mss n$ in the set
$\{k/\mtt t^d_1: k\in \bb Z\}$, while the second one over integers
$\mss N$.  Estimate
$\nu^{\rm c} _{2,\mtt j, \mss M} \big\{ \mss M_{\mtt k} =\mss N \,
\big\} $ by
$\nu^{\rm c} _{2,\mtt j, \mss M} \big\{ \mss M_{\mtt k} \ge \mss N \,
\big\} $.  By Corollary \ref{s25}, the previous sum is less than or
equal to
\begin{equation*}
\mf a_1\, \sum_{\mss N \ge \mss M_{\mtt j, \mtt k} +
\mtt t_1^{(d/2)+\kappa}}
e^{ \beta \, \mtt t_1^d\, \widehat{\mf p} (\mss N/\mtt t^d_1) }\,
e^{ - (1/2) ( \mss N - \mss M_{\mtt j, \mtt k})^2 \mtt t_1^{-d}}\,.
\end{equation*}
Since $\beta \le 1/(4\, \mf a_{3,3})$, by \eqref{71} and elementary
Gaussian bounds, this expression is less than or equal to
$\mf a_{3,3}\, \exp\{-\, (1/4) \, \mtt t_1^ {2\kappa} \}$ for
$\mtt t_1 \ge 2$. Therefore,
\begin{align}
\label{177}
E_{\nu^{\rm c} _{2,\mtt j, \mss M} } 
\Big[ \, e^{ \beta \, \mtt t_1^d\, \widehat{\mf p} (\mss m_{\mtt k}) }\,
\mtt 1 \big\{ |\mss m_{\mtt k} - \mss m_{\mtt j, \mtt k} |\ge 
\mtt t_1^{-(d/2)+\kappa}\,\big\}\, \Big]
\,\le\, \mf a_{3,3}\, C_0(\kappa)\, e^{-\, \mtt t_1^ {2\kappa} /4}\,.
\end{align}

We turn to the case
$|\mss m_{\mtt k} - \mss m_{\mtt j, \mtt k} |\le \mtt
t_1^{-(d/2)+\kappa}$.  By \eqref{71} and since, by hypothesis
$\beta\le \mtt t_1^{-8\kappa}$, in this region there exists a constant
$\mf a_{3,3}$ such that
\begin{equation*}
\beta \, \mtt t_1^d \, \widehat{ \mf p} (\mss m_{\mtt k}) \,
\le \; \mf a_{3,3}\,  \beta \, \mtt t_1^{2\kappa} 
\; <\;  \mf a_{3,3}\, \mtt t_1^{-6\kappa} \quad\text{and}\quad
\big\{\, \beta \,\mtt t_1^d \, \widehat{ \mf p}
(\mss m_{\mtt k})  \, \big\}^3 \; 
\le \; \mf a_{3,3}\,  \beta^3 \, \mtt t_1^{6\kappa} \; <\;
\mf a_{3,3} \,  \beta^2 \,  \mtt t_1^{-2\kappa}\,.
\end{equation*}
Expanding the exponential up to third order yields that
\begin{align*}
& E_{\nu^{\rm c} _{2,\mtt j, \mss M} } 
\Big[ \, e^{ \beta \, \mtt t_1^d\, \widehat{\mf p} (\mss m_{\mtt k}) }\,
\mtt 1 \big\{ |\mss m_{\mtt k} - \mss m_{\mtt j, \mtt k} |\le 
\mtt t_1^{-(d/2)+\kappa}\,\big\}\, \Big]
\\
&\quad \,\le\,
1 \; +\;  E_{\nu^{\rm c} _{2,\mtt j, \mss M}}
\big[ \, \beta \, \mtt t_1^d \, \widehat{ \mf p} (\mss m_{\mtt k}) \,
\mtt 1 \big \{ | \mss m_{\mtt k} - \mss m_{\mtt j, \mtt k} |
\, \le\,  \mtt t_1^{-(d/2)+\kappa}\big\}\big] 
\, +\; \beta^2  \,
E_{\nu^{\rm c} _{2,\mtt j, \mss M}} \big[ \, \big( 
\mtt t_1^d \, \widehat{ \mf p} (\mss m_{\mtt k})  \big)^2 \big]
\; + \; \mf a_{3,3}  \, \beta^2  \, \mtt t_1^{- 2 \kappa}\; .
\end{align*}

We may rewrite the first order term as
\begin{align*}
E_{\nu^{\rm c} _{2,\mtt j, \mss M}}
\big[ \, \beta \, \mtt t_1^d \, \widehat{ \mf p} (\mss m_{\mtt k}) \, \big] 
\,-\,
E_{\nu^{\rm c} _{2,\mtt j, \mss M}}
\big[ \, \beta \, \mtt t_1^d \, \widehat{ \mf p} (\mss m_{\mtt k}) \,
\mtt 1 \big \{ | \mss m_{\mtt k} - \mss m_{\mtt j, \mtt k} |
\, > \,  \mtt t_1^{-(d/2)+\kappa}\big\}\big] 
\; .
\end{align*}
The first expectation vanishes by definition of
$\widehat{ \mf p} (\cdot)$. By Corollary \ref{s25}, as $\beta \le 1$
and $\widehat{ \mf p} (\cdot)$ is bounded, the second expectation is
bounded by
$\mf a_1\, C_0(\kappa)\, \exp\{ - (1/4) \mtt t_1^{2\kappa}\}$.  On
the other hand, by \eqref{72}, the second order term in the
penultimate displayed equation is bounded by $\mf a_{3,3}\,
\beta^2$. Therefore, 
\begin{align*}
E_{\nu^{\rm c} _{2,\mtt j, \mss M} } 
\Big[ \, e^{ \beta \, \mtt t_1^d\, \widehat{\mf p} (\mss m_{\mtt k}) }\,
\mtt 1 \big\{ |\mss m_{\mtt k} - \mss m_{\mtt j, \mtt k} |\le 
\mtt t_1^{-(d/2)+\kappa}\,\big\}\, \Big]
\,\le\, 1\,+\, 
\mf a_1 \, C_0(\kappa)\, e^{ - \mtt t_1^{2\kappa}/4}
\; +\; \mf a_{3,3}\, \beta^2\,.
\end{align*}
This expression is less than or equal to
$1 + \mf a_{3,3}\, C_0(\kappa) \, \beta^2$ for all $\mtt t_1\ge 2$
because $\beta^{-1} <\mtt t_1^{3d}$.  To complete the proof of the
lemma in the case $\beta \, \mtt t_1^{3d} > 1$, it remains to recall
that $\ln (1+u)\le u$, and to sum the previous estimate with the one
obtained in \eqref{177}.
\end{proof}

Denote by $\mf b_1 = \mf b_1(\mf r, \Vert \rho\Vert_{C^3}, \mf p)$ the
constant obtained from $\mf a^*_{3,3}$ of Lemma \ref{s16} when $\delta = \mf r/2$:
\begin{align*}
\cb{\mf b_1} := \min\{ 1, 1/\mf a^*_{3,3} \} \,.
\end{align*}

\begin{corollary}
\label{s27}
Fix a polynomial $\mf p(\cdot)$, a continuous function
$J\colon \bb T^d \to \bb R$, and $\epsilon>0$.  Then, there exists a
finite constant $\mf a_{3,3} = \mf a_{3,3}(\epsilon, \mf p)$, such
that
\begin{equation*}
\sum_{\mtt j \in \Sigma^+_{\mtt t_3}}
\int \mtt 1\{ \mss m_{2, \mtt j} \in (\mf r_- , \mf r_+) \}
\,  \sum_{\mtt k \in \Lambda_{\mtt j}} \mtt t^d_1\, J_{\mtt k}\,
 \widehat{ \mf p} (\mss m_{\mtt k})  \, f\;
d\nu^n_{\rho(\cdot)}
\,\le\,
\delta\, n^2\, I_n (f\,;\, \nu^n_{\rho(\cdot)})
\,+\, \frac{\mf a_{3,3}}{\delta} \, \Vert J\Vert^2_{\infty}\,
\mtt t_2^{2d}  \, \Big( \frac{n}{\mtt t_1 \mtt t_2}\Big)^{d-2} 
\end{equation*}
for all $\delta>0$, integers $\mtt t_1\ge 2$, $\mtt t_2\ge 2$ such
that
$\mtt t_2^{2+d} \le (\delta / \mf c_{\rm LS}) \Vert J\Vert_\infty^{-1}\,
(n/\mtt t_1)^2 \, (\mf b_1 \wedge \mtt t_1^{-\epsilon})$, and
densities $f$ with respect to $\nu^n_{\rho(\cdot)}$.
\end{corollary}

\begin{proof}
Repeat Steps 1 and 2 in the proof of Lemma \ref{s23} to project on
the hyperplanes
$\cb{\Omega_{2,\mtt j, \mss M}} := \{\eta \in
\{0,1\}^{B_{2,\mtt j}} : \sum_{x\in B_{2,\mtt j}} \eta_x =
\mss M\}$. These projections yield that the left-hand side of the
expression in the statement of the corollary is equal to
\begin{equation}
\label{04c}
\begin{gathered}
\sum_{\mtt j \in \Sigma^+_{\mtt t_3}} \sum_{\mss M }
\nu^{2, \mtt j}_{\rho (\cdot)}  (\Omega_{2,\mtt j, \mss M})\,
\< f_{2,\mtt j}  \>_{2,\mtt j , \mss M} \,
\int W_{\mtt j}  \, f^{(\mss M)}_{2,\mtt j} \, d\nu^{\rm c}_{2, \mtt j, \mss  M}
\;,
\\
\text{where}\quad
\cb{W_{\mtt j} } :=\sum_{\mtt k \in \Lambda_{\mtt j}} J_{\mtt k}\,
\mtt t^d_1\; \widehat{ \mf p} (\mss m_{\mtt k}) \,, \quad 
\mtt j \in \Sigma^+_{\mtt t_3}\,,
\end{gathered}
\end{equation}
and the sum over $\mss M$ is carried over integers $\mss M$ such that
$\mf r_- \le \mss M/ | B_{2, \mtt j} | \le \mf r_+$. By the
entropy inequality, this expression is less than or equal to
\begin{equation*}
\sum_{\mtt j \in \Sigma^+_{\mtt t_3}} \sum_{\mss M }
\nu^{2, \mtt j}_{\rho (\cdot)}  (\Omega_{2, \mtt j, \mss M})\,
\< f_{2,\mtt j}  \>_{2,\mtt j , \mss M} \,
\Big\{ \, \frac{1}{\gamma}\,
H  ( f^{(\mss M)}_{2,\mtt j}  \,|\,  \nu^{\rm c}_{2, \mtt j, \mss  M})
\,+\, \frac{1}{\gamma}\, \ln
\int e^{\gamma W_{\mtt j} }  \, d\nu^{\rm c}_{2, \mtt j, \mss  M}\, \Big\}
\end{equation*}
for all $\gamma>0$. Let $\cb{L} = \mtt t_1\mtt t_2$. By the logarithmic
Sobolev inequality, the first term is bounded by
\begin{equation*}
\sum_{\mtt j \in \Sigma^+_{\mtt t_3}} \sum_{\mss M }
\nu^{2, \mtt j}_{\rho (\cdot)}  (\Omega_{2, \mtt j, \mss M})\,
\< f_{2,\mtt j}  \>_{2,\mtt j , \mss M} \,
\frac{1}{\gamma}\, \mf c_{\rm LS} \, L^2\, 
I_{B_{2, \mtt j} }   ( f^{(\mss M)}_{2,\mtt j}
\,;\,  \nu^{\rm c}_{2, \mtt j, \mss  M})
\;\le\;
\frac{1}{\gamma}\, \mf c_{\rm LS} \, L^2\, 
I_{n}   ( f \,;\,  \nu^n_{\rho(\cdot)})\,.
\end{equation*}
Set $\gamma = \mf c_{\rm LS} \, \delta^{-1}\, (L/n)^2$ to obtain the
first term on the right-hand side of the statement of the corollary.

By H\"older's inequality the second term in the next to last equation
is bounded by
\begin{equation*}
\sum_{\mtt j \in \Sigma^+_{\mtt t_3}} \sum_{\mss M }
\nu^{2, \mtt j}_{\rho (\cdot)}  (\Omega_{2, \mtt j, \mss M})\,
\< f_{2,\mtt j}  \>_{2,\mtt j , \mss M} \,
\Big\{ \, \frac{1}{\gamma}\, \frac{1}{\mtt t_2^d }\,
\sum_{\mtt k \in \Lambda_{\mtt j}}   \ln
\int
e^{\gamma J_{\mtt k}\, \mtt t_2^d\, \mtt t^d_1\; \widehat{ \mf p} (\mss m_{\mtt k})   }
\, d\nu^{\rm c}_{2, \mtt j, \mss  M}\, \Big\}\,.
\end{equation*}
By hypothesis,
$\gamma \, |J_{\mtt k}|\, \mtt t_2^d \le \mf b_1 \wedge \mtt
t_1^{-\epsilon}$. We may therefore apply Lemma \ref{s16}, which
asserts that the previous expression is less than or equal to
\begin{equation*}
\mf a_{3,3} (\epsilon, \mf p) \,
\sum_{\mtt j \in \Sigma^+_{\mtt t_3}} \sum_{\mss M }
\nu^{2, \mtt j}_{\rho (\cdot)}  (\Omega_{2, \mtt j, \mss M})\,
\< f_{2,\mtt j}  \>_{2,\mtt j , \mss M} \,
\gamma\,  J_{\mtt k}^2\, \mtt t_2^{2d}
\,\le\,
\mf a_{3,3} (\epsilon, \mf p) \,
\Vert J\Vert_\infty^2\, \gamma\, n^d \,
\Big(\frac{\mtt t_2}{\mtt t_1}\Big)^d 
\,,
\end{equation*}
which completes the proof of the corollary.
\end{proof}

It remains to estimate $V_{\mtt j, \mtt k}$  introduced in
\eqref{82}.

\begin{lemma}
\label{s28}
Fix a polynomial $\mf p(\cdot)$, a function
$J\colon \bb T^d \to \bb R$ of class $C^1$, $\delta>0$, and $K_1\ge 1$.
Then, there exist finite constants $\mf a_{3,1}$, $C_0 = C_0(\mf p)$
such that
\begin{equation*} 
\sum_{\mtt j\in \Sigma^+_{\mtt t_3 }} 
\int  \mtt 1\{ \mss m_{2, \mtt j} \in (\mf r_- ,\, \mf r_+) \}
\sum_{\mtt k\in \Lambda_{\mtt j}} \mtt t^d_1 \,  
J_{\mtt k}  \, V_{\mtt j, \mtt k} \, f \, d  \nu^n_{\rho(\cdot)} \,\le\,
\delta\, n^2 \,  I_n  ( f \,;\, \nu^n_{\rho (\cdot)} ) \,+\,
\mf a_{3,1} \, \frac{\mf c_{\rm LS}}{\delta} \,
\mf C_2 \, \Big(\frac{\mtt t_1\mtt t_2}{n}\Big)^4\, n^d  
\end{equation*}
for all densities $f$ with respect to $\nu^n_{\rho(\cdot)}$, and
$\mtt t_1\ge 2$, $\mtt t_2\ge 2$ such that
$(\mtt t_1\mtt t_2)^{d+3} \le K_1\, n^3$. In this equation, $\mf C_2 =
\exp\{ \mf c_{\rm LS} C_0(\mf p)\, K_1 \, (\, \mf a_{3,1} + \Vert
J\Vert_{C^1})/\delta\}$. 
\end{lemma}

\begin{proof}
Rewrite the sum over $\mtt k$ on the left-hand side as
\begin{equation}
\label{70}
\begin{aligned}
\sum_{\mtt k\in \Lambda_{\mtt j}} \mtt t^d_1 \,  
J_{\mtt k}  V_{\mtt j, \mtt k}
\, & =\, 
\sum_{\mtt k\in \Lambda_{\mtt j}} \mtt t^d_1 \,  
[\, J_{\mtt k}  -  \, J_{2, \mtt j}  \,]
\,  \mf p'( \widetilde{\mss m}_{\mtt j, \mtt k} ) \,
( \mss m_{\mtt k} - \mss m_{\mtt j, \mtt k} )
\\
\, & +\, 
\mtt t^d_1 \,  J_{2, \mtt j} 
\sum_{\mtt k\in \Lambda_{\mtt j}}
[\, \mf p'( \widetilde{\mss m}_{\mtt j, \mtt k} )
- \mf p'( \mss m_{2, \mtt j}   )\,] \, 
( \mss m_{\mtt k} -
\mss m_{\mtt j, \mtt k} ) \,,
\end{aligned}
\end{equation}
where
$\cb{ J_{2, \mtt j} } = \mtt t_2^{-d} \sum_{\mtt k\in \Lambda_{\mtt
j}} J_{\mtt k}$. In the last identity we used the fact that
\begin{equation*}
\sum_{\mtt k\in \Lambda_{\mtt j}} \mf p'( \mss m_{2, \mtt j}   )\, 
( \mss m_{\mtt k} - \mss m_{\mtt j, \mtt k} ) \,=\, 0 
\end{equation*}
for all $\mtt j\in \Sigma^+_{\mtt t_3}$ because, by definition of $\mss
m_{\mtt k} $ and $\mss m_{\mtt j, \mtt k}$,
\begin{align*}
\sum_{\mtt k\in \Lambda_{\mtt j}}  \mss m_{\mtt k} \,=\, 
\sum_{\mtt k\in \Lambda_{\mtt j}}  \frac{1}{|B_{\mtt k}|}
\sum_{x\in B_{\mtt k}} \eta_x \,=\, \frac{1}{|B_{\mtt k}|} \,
\mss M_{2, \mtt j}
\,=\, \sum_{\mtt k\in \Lambda_{\mtt j}}  \mss m_{\mtt j, \mtt k}\,.
\end{align*}

Now, Lemma \ref{s23} provides bounds for the terms \eqref{70}.  Let
\begin{equation*}
V_{\mtt j} \,=\, \mtt 1\{ \mss m_{2, \mtt j} \in (\mf r_- , \mf r_+) \}
\sum_{\mtt k\in \Lambda_{\mtt j}} \mtt t^d_1 \,  
[\, J_{\mtt k}  -  \, J_{2, \mtt j}  \,]
\,  \mf p'( \widetilde{\mss m}_{\mtt j, \mtt k} ) \,
( \mss m_{\mtt k} - \mss m_{\mtt j, \mtt k} )\,, \qquad
\mtt j\in \Sigma^+_{\mtt t_3} \,.
\end{equation*}
Clearly,
$\Vert V_{\mtt j} \Vert_\infty \le C_0(\mf p) \, \Vert J \Vert_{C^1}\,
(\mtt t_1\mtt t_2)^d \, (\mtt t_1\mtt t_2/n)$. Since, by hypothesis,
$(\mtt t_1\mtt t_2)^{d+3} \le K_1 \, n^3$,
$\Vert V_{\mtt j} \Vert_\infty \le C_1 \, (n/\mtt t_1\mtt t_2)^2$, where
$\cb{C_1} := C_0(\mf p) \, K_1\, \Vert J \Vert_{C^1}\,$. Thus, each
function $V_{\mtt j} $ satisfies the assumptions of Lemma \ref{s23}
with $r_1 = \mtt t_1 \mtt t_2$, $r_2=\mtt t_3$.  Therefore, by this
result,
\begin{equation*}
\sum_{\mtt j\in \Sigma^+_{\mtt t_3 }} 
\int V_{\mtt j}  \, f \, d  \nu^n_{\rho(\cdot)} \,\le\,
\delta\, n^2 \,  I_n  ( f \,;\, \nu^n_{\rho (\cdot)} ) \,+\,
C_2\, \Big(\frac{\mtt t_1\mtt t_2}{n}\Big)^2\,
\sum_{\mtt j\in \Sigma^+_{\mtt t_3}} 
\max_{\mf r_- \le \mss m \le \mf r_+} \int V^2_{\mtt  j} \;
d\nu^{\rm c}_{ 2, \mtt j, \mss M} 
\end{equation*}
for all densities $f$ with respect to $\nu^n_{\rho(\cdot)}$.  The
maximum over densities $\mss m$ is carried on the interval
$(\mf r_- ,\, \mf r_+)$ because $V_{\mtt j} $ vanishes outside it, and
$\cb{ C_2} := (\mf c_{\rm LS} / 2 \delta)\, \exp\{C_1 \mf c_{\rm LS}
/\delta\}$, where $C_1$ has been introduced above in the proof and
$\mf c_{\rm LS}$ is the logarithmic Sobolev constant.

We estimate the remaining integral. As $J$ is of class $C^1$,
$|\, J_{\mtt k} - \, J_{2, \mtt j} \,|\le \Vert J\Vert_{C^1} \, (\mtt
t_1\mtt t_2/n)$. Hence, by Corollary \ref{s71},
\begin{equation*}
\max_{\mf r_-  \le \mss m \le \mf r_+} \int V^2_{\mtt  j}\;
d\nu^{\rm c}_{ 2, \mtt j, \mss M}
\,\le\,
\mf a_1  \, \Vert J\Vert^2_{C^1}
\, (\mtt t_1 \,\mtt t_2)^d\, (\mtt t_1\mtt t_2/n)^2
\quad \text{for all}\;\; \mtt j\in \Sigma^+_{\mtt t_3} \,.
\end{equation*}

To complete the proof of the estimate for the first term in
\eqref{70}, it remains to recollect all previous bounds, noting the
overestimate $\|J\|^2_{C^1} \le 2 \, \exp\{\|J\|_{C^1}\}$.

The same argument applies to the second term in \eqref{70} because, by
\eqref{57},
$|\widetilde{\mss m}_{\mtt j, \mtt k} - \mss m_{2, \mtt j}|\le \mf
a_{3,1} \, (\mtt t_1 \mtt t_2/n)$. Thus, in the previous argument the
estimate
$|\, J_{\mtt k} - \, J_{2, \mtt j} \,|\le \Vert J\Vert_{C^1} \, (\mtt
t_1\mtt t_2/n)$ is replaced by the previous one, and we get the same
bounds with $\mf a_{3,1}$ in place of $\Vert J\Vert_{C^1}$.
\end{proof}

\begin{proof}[Proof of Proposition \ref{s37}]
It follows from Corollary \ref{s27} and Lemma \ref{s28}.
\end{proof}

\subsection*{From canonical to grand canonical measures}

In Proposition \ref{s37}, we replaced $\mf p(\mss m_{\mtt k})$ by its
expectation with respect to the canonical measure
$\nu^{\rm c} _{2,\mtt j, \mss M}$. The next result replaces this
expectation by a sum of polynomials in $\mss m_{2,\mtt j}$.  The main
result of this subsection reads as follows. Recall from \eqref{111}
the definition of $W_{\mf p}$, and that we denote by $\mtt V_{p}$ the
polynomial $W_{\mf p}$ when $\mf p(\mss m) = \mss m^p$, $p\ge 1$. Let
\begin{equation*}
\cb{\mf W_{p, J} (\mss m)}
\,:= \, \sum_{\mtt k \in \Lambda_{\mtt j}}  \mtt t_1^d\, J_{\mtt k}\,
\Big\{\, \big(\, E_{\nu^{\rm c}_{\mss 2, \mtt j, \mss M}}[\mss m_{\mtt k}^p]
- E_{\nu^n_{\rho(\cdot)}} [\mss m_{\mtt k}^p] \,\big)
\,-\, 
\big(\, \mtt V_{p} (\mss m_{2, \mtt j})
- E_{\nu^n_{\rho(\cdot)}} [ \mtt V_{p} (\mss m_{2, \mtt j})] \,\big) \, \Big\}\,,
\quad \mss m = \frac{\mss M}{|B_{2, \mtt j}|}\, \cdot
\end{equation*}

\begin{proposition}
\label{s44}
Fix $p\ge 1$, a function $J\colon \bb T^d \to \bb R$ of class $C^1$,
and constants $K_1\ge 1$, $K_2\ge 1$.  Then, there exist finite
constants $\mf a_1$, $\mf a_{3,3} = \mf a_{3,3}(p)$, $C_2= C_2(p)$
such that
\begin{equation*}
\begin{aligned}
&\int \Big |\,
\sum_{\mtt j \in \Sigma^+_{\mtt t_3}}
\mtt 1\{ \mss m_{2, \mtt j} \in (\mf r_- , \mf r_+) \}
\, \mf W_{p, J} (\mss m_{2, \mtt j}) \,\Big |\, 
f\, d\nu^n_{\rho(\cdot)}
\\
&\quad 
\,\le\, \mf a_{3,3}\,  \, \Vert J\Vert_{\infty} n^d\, \Big\{ \, \Big(
\frac{\mtt t_1\mtt t_2 }{n} \Big)^4
\,+\,  \frac{\mtt t_1\mtt t_2 }{n} \, \frac{1} {(\mtt t_1\mtt t_2)^d } \, 
\,+\,  \frac{1 }{\mtt t_1^{d}}   \, \Big[ \, 
\frac{1 }{\mtt t_1^{2d}} \,+\, \frac{1 }{(\mtt t_1\mtt t_2)^{d}} \,+\,
\Big(\frac{\mtt t_1 }{n} \Big)^2\,\Big] \Big\}
\\
&\quad
\,+\, 
\mf a_{3,3} \, \Vert J\Vert_{C^1}
\,\Big\{   \, H_n (f\,|\, \nu^n_{\rho(\cdot)})\, + \, \ln(2)\,\Big\}
\,+\, \mf a_1 \, \Vert J\Vert_{\infty}\,
\Big( \frac{n}{\mtt t_1\mtt t_2}\Big)^d\,
e^{- (\mtt t_1\mtt t_2)^d/\mf a_1}
\\
&\quad
\,+\, \mf a_{3,3}\, \, \Vert J\Vert_{C^1}\, \Big\{
\Big( \frac{\mtt t_1\mtt t_2}{n}\Big)^2\, 
\Big( \frac{n}{ (\mtt t_1\mtt t_2)^{1/2}} \Big)^d
\,+\,   \Big[\, e^{C_2 K_2^d}\,
+ e^{\mf a_{3,3} K_1}\, (\mtt t_1\mtt t_2)^{2d}\,
\Big( \frac{\mtt t_1\mtt t_2}{n}\Big)^6
\,\Big] \, 
\Big( \frac{n}{ (\mtt t_1\mtt t_2)^{2}} \Big)^d \, \Big\}
\end{aligned}
\end{equation*}
for all integers $\mtt t_1\ge 2$, $\mtt t_2\ge 2$ such that
$(\mtt t_1\mtt t_2)^{d+3} \le K_1\, n^3$,
$ (\mtt t_1\mtt t_2)^{1+(d/4)} \le n$, $\mtt t_2 \le K_2\, \mtt t_1$,
and densities $f$ with respect to $\nu^n_{\rho(\cdot)}$.
\end{proposition}

\begin{remark}
\label{rm6}
One would wish to replace
$\sum_{\mtt k \in \Lambda_{\mtt j}} \mtt t_1^d\, J_{\mtt k}\, \mss m_{\mtt
k}^p$ by $(\mtt t_1\mtt t_2)^d\, J_{2,\mtt j}\, \mss m_{2,\mtt j}^p $. The
previous result states that such a replacement carries a correction
which arises from replacing expectations with respect to
canonical measures by expectations with respect to grand canonical
measures. 
\end{remark}

\begin{proof}[Proof of Proposition \ref{s44}]
We first remove the indicator
$\mtt 1\{\mss m_{2, \mtt j} \in (\mf r_- , \mf r_+) \}$ applying Lemma
\ref{s07b} with $r_1 = \mtt t_1\mtt t_2$. This explains the presence
of the term
$\mf a_1 \, ( n/\mtt t_1\mtt t_2)^d\, \exp\{- (\mtt t_1\mtt t_2)^d/\mf
a_1\}$ on the right-hand side.

The proof relies on Lemma \ref{s32} and Corollary \ref{s30b}. Lemma
\ref{s32} requires subtracting the expectation of the integrand $F$ with
respect to $\nu^n_{\rho(\cdot)}$. To guarantee the presence of this
expectation, write
\begin{equation*}
E_{\nu^n_{\rho(\cdot)}} [\mss m_{\mtt k}^p] \,=\,
E_{\nu^n_{\rho(\cdot)}} \Big[ E_{\nu^n_{\rho(\cdot)}}
\big[\, \mss m_{\mtt k}^p \,|\, \mss
M_{2,\mtt j}\,\big]\, \Big] \,=\,
E_{\nu^n_{\rho(\cdot)}} \Big[ E_{\nu^{\rm c}_{2, \mtt j, \mss M}}
\big[\, \mss m_{\mtt k}^p \,\big]\, \Big] \quad
\text{for}\;\; \mtt k\in \Lambda_{\mtt j} \,.
\end{equation*}
Every operation carried out on
$E_{\nu^{\rm c}_{\mss 2, \mtt j, M}} [\mss m_{\mtt k}^p]$ is also applied to
the last term of the previous displayed equation. For example, to
estimate the sum appearing in the statement of the proposition we expand
the sum
$\sum_{\mtt k \in \Lambda_{\mtt j}} J_{\mtt k}\, E_{\nu^{\rm c}_{\mss 2,
\mtt j, M}} [\mss m_{\mtt k}^p]$ applying Corollary \ref{s30b}. The
same expansion is performed to the corresponding term
$\sum_{\mtt k \in \Lambda_{\mtt j}} J_{\mtt k}\,
E_{\nu^n_{\rho(\cdot)}} [ \, E_{\nu^{\rm c}_{2, \mtt j, \mss M}} \big[\, \mss
m_{\mtt k}^p \,\big] \,]$ without explicitly mentioning it.

Corollary \ref{s30b} expresses
$\sum_{\mtt k \in \Lambda_{\mtt j}} J_{\mtt k}\, E_{\nu^{\rm c}_{\mss 2,
\mtt j, M}} [\mss m_{\mtt k}^p]$ as the sum of five terms and a
remainder. Note that we have an extra multiplicative factor
$\mtt t^d_1$ in comparison with Corollary \ref{s30b}.  The terms
inside the brace in the first line on the right-hand side of the
inequality stated in the proposition correspond to the bound on the
remainder obtained in Corollary \ref{s30b}.

To complete the proof, we have to show that the expectation with
respect to the measure $ f\; d\nu^n_{\rho (\cdot)} $ of the terms
which appear in Corollary \ref{s30b} and do not appear in the
statement of Proposition \ref{s44} can be estimated by the entropy
$H_n(f \,|\, \nu^n_{\rho (\cdot)})$ and small constants.

We start with the last term in \eqref{110}. Let
$F(a) = \chi(a) \, a^{p-2}$.  By Lemma \ref{s32} with
$r_1 = \mtt t_1\mtt t_2$, $\mf a (\mtt i) =J_{2, \mtt i}$ and
$\gamma^{-1} = \Vert J\Vert_\infty$, there exists a finite constant
$C$, depending only on $p$, such that
\begin{equation}
\label{100}
\begin{aligned}
& \int\Big |  \sum_{\mtt j\in \Sigma^+_{\mtt t_3}}
J_{2, \mtt j}\,
\big[ \, F(\mss m_{2,\mtt j}) - E_{\nu^n_{\rho (\cdot)}} [F(\mss m_{2, \mtt j}) ] \, \big] 
\,  \Big| \, f\; d\nu^n_{\rho (\cdot)}
\\
&\quad
\,\le \, \Vert J\Vert_\infty  \, \Big\{   \, H_n(f \,|\, \nu^n_{\rho (\cdot)}  )
\,+\,\ln(2) \,+\,  C(p)\, 
\Big( \frac{n}{(\mtt t_1\mtt t_2)^{2}} \Big)^d\, \Big\}\,.
\end{aligned}
\end{equation}
This estimate provides a bound for the expectation with respect to
the measure $ f\; d\nu^n_{\rho (\cdot)} $ of the last term in
\eqref{110} multiplied by $\mtt t_1^d$. Below, we will simply say that
it provides an estimate for the $i$-th term in \eqref{110}, omitting that
it is in reality an estimate for the expectation of the expression
multiplied by $\mtt t_1^d$. Mind also that, without mentioning it, we
are subtracting the mean of the term with respect to
$\nu^n_{\rho (\cdot)}$.

The same argument applies to all terms of order
$(\mtt t_1\mtt t_2/n)^3$.  By definition, see Lemma \ref{s20b},
$[H^{(1)}_{\mtt j, \mtt k}]^3$,
$H^{(1)}_{\mtt j, \mtt k} \, \mtt Y^{(2)}_{\mtt j, \mtt k}$ and
$ [\, J_{\mtt k} - J_{2, \mtt j}\,] \, \mtt Y^{(2)}_{\mtt j, \mtt k} $
are of order $(\mtt t_1\mtt t_2/n)^3$. More precisely, since by
hypothesis $(\mtt t_1\mtt t_2)^{d+3} \le K_1\, n^3$, Lemma \ref{s32}
with $r_1 = \mtt t_1\mtt t_2$,
$\mf a (\mtt j) = \mtt t_1^d \sum_{\mtt k \in \Lambda_{\mtt j}} J_{\mtt k} \,
[H^{(1)}_{\mtt j, \mtt k} ]^3$ and
$F(x) = (1/6) \, Q_p'''(x) \, \chi(x)^3$ provides a bound for the
fourth term on the right-hand side of \eqref{110} similar to
\eqref{100}. On the right-hand side of the inequality, instead of the
expression which appears in \eqref{100}, we have
\begin{equation}
\label{101}
\Vert J \Vert_\infty  \, \Big\{   \, H_n(f \,|\, \nu^n_{\rho (\cdot)}
)  \,+\, \ln (2) \,+\, \mf a_{3,1} (p)\, e^{\mf a_{3,1} (p) \, K_1}\,
(\mtt t_1\mtt t_2)^{2d}\, \Big( \frac{\mtt t_1\mtt t_2}{n}\Big)^6\, 
\Big( \frac{n}{(\mtt t_1\mtt t_2)^{2}} \Big)^d\, \Big\}\,.
\end{equation}
As $K_1\ge 1$,
$\mf a_{3,1} (p)\, e^{\mf a_{3,1} (p) \, K_1} \le e^{2 \mf a_{3,1}
(p) \, K_1}$, and we can get rid of the multiplicative constant.

An identical bound is obtained for the sums involving
$H^{(1)}_{\mtt j, \mtt k} \, \mtt Y^{(2)}_{\mtt j, \mtt k}$ and
$ [\, J_{\mtt k} - J_{2, \mtt j}\,] \, \mtt Y^{(2)}_{\mtt j, \mtt k}$.
Some of these terms are multiplied by polynomials of $\mss m_{2,\mtt
j}$. This does not interfere in the argument since we are always
subtracting the expectations with respect to $\nu^n_{\rho(\cdot)}$.
In the estimate of the term $[\, J_{\mtt k} - J_{2, \mtt j}\,] \,
\mtt Y^{(2)}_{\mtt j, \mtt k}$, we get the $C^1$ norm of $J$ instead
of the $L^\infty$ one, and a constant of type $\mf a_{3,2}$ instead of
$\mf a_{3,1}$ because of the Hessian of $\rho$ in the term $H^{(4)}$.

Since
$\mtt Y_{\mtt j, \mtt k} = H^{(1)}_{\mtt j, \mtt k} + \mtt
Y^{(2)}_{\mtt j, \mtt k}$, it remains to estimate
\begin{equation}
\label{102b}
\begin{aligned}
& \mtt t_2^d\, J_{2,\mtt j}\, \mss Q_p(\mss m_{2,\mtt j})
\,+\, \mss Q'_p(\mss m_{2,\mtt j}) \,
\chi (\mss m_{2,\mtt j})\,
\sum_{\mtt k\in\Lambda_{\mtt j}}  [\, J_{\mtt k} - J_{2, \mtt j}\,] \,
H^{(1)}_{\mtt j, \mtt k}  \\
& +\, \frac{1}{2}\, \mss Q''_p(\mss m_{2,\mtt j}) \,
\chi(\mss m_{2,\mtt j})^2\,  \sum_{\mtt k\in\Lambda_{\mtt j}}  J_{\mtt k} \,
[H^{(1)}_{\mtt j, \mtt k} ]^2 \,.
\end{aligned}
\end{equation}
Recall the definition of $\mss Q_p(\cdot)$ introduced in \eqref{99}.
We claim that the second and third terms in \eqref{102b} with
$\mss Q_p(\cdot)$ replaced by $\mtt t_1^{-id} \mss q_{p,i} (\cdot)$,
$0\le i\le 2$, are negligible, in the sense that they can be estimated
by the entropy and a small constant. It is enough to consider the term
$\mss q_{p,0}$.  Lemma \ref{s45} with $r_1 = \mtt t_1\mtt t_2$,
$\mf a(\mtt j) = \mtt t_1^d \, \sum_{\mtt k\in\Lambda_{\mtt j}} [\,
J_{\mtt k} - J_{2, \mtt j}\,] \, H^{(1)}_{\mtt j, \mtt k}$ provides a
bound for the second term in \eqref{102b} (with $\mss Q_p(\cdot)$
replaced by $\mss q_{p,0} (\cdot)$) similar to \eqref{100}. Indeed, by
definition of $H^{(1)}_{\mtt j, \mtt k}$,
$|\mf a (\mtt j)| \le \mf a_{3,1}\, \Vert J\Vert_{C^1}\, (\mtt t_1\mtt
t_2)^{d+2}/n^2$. By hypothesis,
$(\mtt t_1\mtt t_2)^{d+2}/n^2 \le  (\mtt t_1\mtt t_2)^{d/2}$. Let
$F(\mss m) = \mss q'_{p,0}(\mss m) \, \chi (\mss m)$.  By Lemma
\ref{s45}, there exists a finite constant $\mf a_{3,1}$ such that
\begin{equation}
\label{103b}
\begin{aligned}
& \int \Big|\, \sum_{\mtt j\in \Sigma^+_{\mtt t_3}}
\big[ \, F(\mss m_{2,\mtt j}) - E_{\nu^n_{\rho (\cdot)}} [F(\mss m_{2, \mtt j})
] \, \big]
\, \mtt t_1^d\,
\sum_{\mtt k\in\Lambda_{\mtt j}}  [\, J_{\mtt k} - J_{2, \mtt j}\,] \,
H^{(1)}_{\mtt j, \mtt k}\,\Big|
\, f\; d\nu^n_{\rho (\cdot)}
\\
&\quad \,\le \,
\mf a_{3,1}\, \Vert J\Vert_{C^1} \,
\Big\{ \, H_n(f \,|\, \nu^n_{\rho (\cdot)}  )  \,+\,
\Big( \frac{\mtt t_1\mtt t_2}{n}\Big)^2\, 
\Big( \frac{n}{ (\mtt t_1\mtt t_2)^{1/2}} \Big)^d\, \Big\}\,.
\end{aligned}
\end{equation}
The third term in \eqref{102b} (with $\mss Q_p(\cdot)$ replaced by
$\mss q_{p,0} (\cdot)$) satisfies the same bound under the same
constraints.

We turn to the first term in \eqref{102b}.  Lemma \ref{s32} with
$r_1 = \mtt t_1\mtt t_2$,
$\mf a (\mtt j) = (\mtt t_2/\mtt t_1)^d\, J_{2, \mtt j}$ and
$F(x) = \mss q_{p,2} (x)$ provides a bound for the first term in
\eqref{102b} (with $\mss Q_p(\cdot)$ replaced by
$\mtt t_1^{-2d} \mss q_{p,2} (\cdot)$) similar to \eqref{100}.  Since, by
hypothesis, $\mtt t_2 \le K_2\, \mtt t_1$, on the right-hand side of the
inequality, instead of the expression which appears in \eqref{100}, we
have
\begin{equation}
\label{104}
\Vert J\Vert_\infty \,
\Big\{   \, H_n(f \,|\, \nu^n_{\rho (\cdot)}  )  \,+\, \ln(2) \,+\,
C(p)\, e^{C(p) K_2^d}\, \Big( \frac{\mtt t_2}{\mtt t_1}\Big)^{2d}\, 
\Big( \frac{n}{ (\mtt t_1\mtt t_2)^{2}} \Big)^d\, \Big\}\,.
\end{equation}
As $\mtt t_2\le K_2 \, \mtt t_1$ and $K_2\ge 1$,
$C(p)\, \exp \{ C(p) K_2^d\}\, (\mtt t_2/\mtt t_1)^{2d} \le C(p)\,
K^{2d}_2 \exp \{ C(p) K_2^d\}$. This expression is less than or equal
to $\exp \{ C(p) K_2^d\}$ for some new constant $C(p)$.

In conclusion, all terms in \eqref{102b} can be estimated with the
exception of
\begin{equation*}
\mtt t_2^d\, J_{2,\mtt j}\, \Big\{\, \mss q_{p,0}(\mss m_{2,\mtt j}) +
\frac{1}{\mtt t_1^{d} } \, \mss q_{p,1}(\mss m_{2,\mtt j})\, \Big\}
\,=\,
\mtt t_2^d\, J_{2,\mtt j}\, \Big\{\, \mss m_{2,\mtt j}^{p} \,+\,
\frac{1}{\mtt t_1^{d} } \, \frac{p(p-1)}{2}\,
\mss m_{2,\mtt j}^{p-2}\, \chi(\mss m_{2,\mtt j}) \Big\}
\,=\, \mtt t_2^d\, J_{2,\mtt j}\,   \mtt V_{p} (\mss m_{2,\mtt j}) \,,
\end{equation*}
where the last identity follows from the definition \eqref{111} of
$W_{\mf p}$ and $\mtt V_p$.  To complete the proof of the proposition,
it remains to recollect all previous estimates.
\end{proof}

\begin{proof}[Proof of Theorem \ref{s26}]
The assertion follows from the bound \eqref{67}, and Propositions
\ref{s37}, \ref{s44}.  We may remove from the right-hand side of the
inequality obtained many terms which are less than or equal to other
terms. Here is a list.  The term
$(n/\mtt t_1 \mtt t_2)^d\, \exp \{- c_1 (\mtt t_1 \mtt t_2)^d\}$ is
dominated by the term $(n/[\mtt t_1 \mtt t_2]^2)^d$.  As
$\mtt t_2 \le K_2 \mtt t_1$, $\mtt t_1^{-2d}$ is dominated by
$(\mtt t_1\mtt t_2)^{-d}$. As $\mtt t_1\mtt t_2 \le n$,
$(\mtt t_1\mtt t_2/n)^6$ is dominated by $(\mtt t_1\mtt
t_2/n)^4$. Finally, the term $(\mtt t_1\mtt t_2)^{-2d}$ is clearly
dominated by $\mtt t_1^{-2d} \mtt t_2^{-d}$.  This completes the proof
of the theorem.
\end{proof}

\section{Proof of Theorem \ref{t01}}
\label{sec1}

In this section, we prove Theorem \ref{t01}.  Hereafter, $C_1$ is a
universal constant which may depend only on the dimension $d$ and the
cylinder functions $h_j$, $g_j$, $1\le j\le d$, changing from line to line.
Recall the convention on the constants $\mf c_{1, p}$
adopted below Lemma \ref{l03}.

Fix $T>0$.  By Lemmata~\ref{l03}, \ref{s41},
\begin{equation}
\label{118}
\begin{aligned}
H_n(f_t | \nu^{n}_{u(t, \cdot) } ) \; & \leq\; H_n(f_0 | \nu^{n}_{u(0, \cdot)})
\, -\,2\,  \int_0^t  n^2\, I_n(f_s ; \nu^{n}_{u(s, \cdot)} )\, ds
\;+\;  \int_0^t ds \; \int \mf W_s \; f_s \, d \nu^{n}_{u(s, \cdot)}  
\\
\, & +\, \int_0^t  ds \int \, \big\{\,  I_n (s) \,+\,
\mf R^{(4)}_{s,n}\, n^{d-2}\, \big\}  \, f_s \, d \nu^{n}_{u(s, \cdot)}  
\end{aligned}
\end{equation}
for all $0\le t\le T$.  Each sum forming $\mf W_s$ in equation
\eqref{118} can be written as
\begin{equation}
\label{117}
\int_0^t  ds\, \int \sum_{x\in \bb T^d_n} 
G(s,x/n)\, \Pi_{u(s, \cdot) , h} (x, \eta)\, 
f_s\, d\nu^n_{u( s, \cdot)}\,,
\end{equation}
where $G(s, \cdot)$ is a smooth function which stands for a function
of $u(s,\cdot)$ and its derivatives, $h(\cdot)$ a cylinder function,
and $\Pi$ is defined in \eqref{179}. In view of the previous
calculations, the main step in the proof of the theorem is stated in
the next result. Recall from \eqref{256} that
$\Vert F\Vert_{T,p} = \sup_{0\le t\le T} \Vert F(t) \Vert_{C^p}$,
$p\ge 0$, $F\colon [0,T]\times \bb T^d \to \bb R$.

\begin{proposition} 
\label{s88}
Assume that the initial condition $u_0$ satisfies \eqref{255a}, and
belongs to $H^{3+\beta}$ for some $0<\beta<1$.  Let $u(t,x)$ be the
solution of the hydrodynamic equation \eqref{127}.  Fix $\delta_0>0$,
a cylinder function $h\colon \{0,1\}^{\bb Z^d} \to \bb R$, and a
function $G$ in $C^{0,2}(\bb R_+ \times \bb T^d)$. In dimension $1$ we
also fix $0<b<1/6$, and in dimension $3$ a constant $\kappa>1$ and
$0< \mtt e < 1/30$.  Then, there exist a finite positive constant
$\mf c_{1,3} = \mf c_{1,3} (h, \delta_0)$ (in dimension $1$,
$\mf c_{1,3}$ also depends on $b$), and an integer
$n_0= n_0(\delta_0, \mf r, h, \Vert u\Vert_{T,3}, \Vert G\Vert_{T,2})$
(in dimension $3$, $n_0$ also depends on $\kappa$) such that
\begin{align*}
& \Big|\, \bb E_{\mu_n} \Big[\, \int_0^t
\sum_{x\in \bb T^d_n} G (s, x/n)\,
\Pi_{u(s, \cdot), h}(x, \eta^n(s)) \, ds \,  \Big] \,\Big|
\\
& \;\le \;
\int_0^t   \delta_0\, n^2\, I_n  (f_s;  \nu^n_{u(s, \cdot) } ) \; ds
\,+\, \mf c_{1,3} \int_0^t  
\,\big\{   \, \mtt K_{d,1} (s) \, H_n (f_s\,|\, \nu^n_{u(s, \cdot)})
\,+\, \mtt K_{d,2} (s)  \,\big\}\, ds
\end{align*}
for all $ 0\le t\le T$, and $n\ge n_0$. In this formula,
\begin{gather*}
\cb{\mtt K_{1,1}(s) } \, :=\, \Vert G(s) \Vert_{\infty} (1 + \Vert G(s)
\Vert_{\infty}) + \Vert G(s) \Vert_{C^2}\,, \quad
\cb{\mtt K_{1,2}(s) }  \,:= \, \mtt K_{1,1}(s) \, n^{(1/3)+b}\,,
\\
\cb{\mtt K_{2,1} (s) } \,:=\, \exp\{ \mf c_{1,3}\, \Vert G(s)\Vert_{C^1}\}
+ \Vert G(s) \Vert_{C^2}\,, \quad
\cb{\mtt K_{2,2} (s)} \,:=\, \mtt K_{2,1} (s) \, n^{4/5}\,,
\\
\cb{\mtt K_{3,1} (s) } \,:=\, \Vert G(s) \Vert_{C^2}\,, \quad
\cb{\mtt K_{3,2} (s) } \,:=\, \mtt K_{3,3} (s)\, n^{3/2 - 2\mtt e}
\,+\, \Vert G(s) \, \Vert_\infty\,  
\frac{n^{3/2}}{\kappa^3}  \,,
\\
\cb{\mtt K_{3,3} (s)} \, := 
\Vert G(s) \Vert_{C^2} \, \kappa^4 \,+\,
\kappa^{1/2}\, e^{\mf c_{1,3} \,\kappa^6\,[1 + \Vert G(s)\Vert_{C^1}]}
\end{gather*}
\end{proposition}

\begin{remark}
It is enough to prove this proposition without the absolute value
and observe that the same argument applies to $-\, G$ in place of $G$.
\end{remark}

\begin{proof}[Proof of Theorem \ref{t01}]
We start with $d=1$, $2$.  Fix $T>0$, $0<b<1/6$. By \eqref{118},
the bound \eqref{108} for $I_n(s)$, the one for
$\mf R^{(4)}_{t,n}$ stated in Lemma \ref{s41}, and Proposition
\ref{s88} with $\delta_0 = 1/(3d)$,
\begin{equation*}
H_n(f_t | \nu^n_{u( t,\cdot)} ) \; \leq\;
H_n(f_0 | \nu^n_{u_0(\cdot)} )
\, - \,  \int_0^t  n^2\, I_n(f_s ; \nu^n_{u( s,\cdot)})\, ds
\, +\, \mf c_{1,4} \, \int_0^t   \big\{\, H_n(f_s | \nu^n_{u( s,\cdot)} )
\,+\, n^{a_d}  \big\} \; ds  
\end{equation*}
for all $n$ sufficiently large and $0\le t\le T$. In this formula,
$a_1= (1/3) + b$ in dimension $1$ and $a_2=4/5$ in dimension $2$.  As
$b<1/6$, apply Gronwall's inequality to complete the proof for $n$
large. For $n$ small the bound follows from equation \eqref{118},
elementary bounds for the terms on the right-hand side of this
inequality, and adjusting the constant $\mf c_{1,4}$.

The proof in dimension $3$ is similar.  Fix $T>0$, $\kappa>1$. Let
$0<\mtt e<1/30$.  By \eqref{118}, the bound \eqref{108} for $I_n(s)$,
the one for $\mf R^{(4)}_{t,n}$ stated in Lemma \ref{s41}, and
Proposition \ref{s88} with $\delta_0 = 1/(3d)$,
\begin{align*}
H_n(f_t | \nu^n_{u( t,\cdot)} ) \; &\leq\;
H_n(f_0 | \nu^n_{u_0(\cdot)} )
\, - \,  \int_0^t  n^2\, I_n(f_s ; \nu^n_{u( s,\cdot)})\, ds
\\
& +\, \mf c_{1,4}   \, \int_0^t   \Big\{\, \, H_n(f_s | \nu^n_{u( s,\cdot)} )
\,+\, e^{\mf c_{1,3} \kappa^6} \, n^{3/2-2\mtt e} \,+\,
\frac{n^{3/2}}{\kappa^3} \Big\} \; ds  
\end{align*}
for all $n\ge n_1 = n_1(\kappa, \mf r, \Vert u \Vert_{T,3}, h)$.  For
$n \ge n_2(\kappa, \mf r, \Vert u \Vert_{T,3}, h)$,
$e^{\mf c_{1,3} \kappa^6} \, n^{3/2-2\mtt e} \le n^{3/2} / \kappa^3$.
Therefore, by Gronwall's inequality, 
\begin{equation*}
H_n (f_t \,|\, \nu^n_{u(t, \cdot)}  ) \,\le\,
\Big\{\, H_n (f_0 \,|\, \nu^n_{u_0 (\cdot)}  ) \,+\,
\mf c_{1,4}  
\frac{n^{3/2}}{\kappa^3} \,\Big\}\, e^{\mf c_{1,4} \,  t}
\end{equation*}
for all $n\ge n_0 = n_1 \vee n_2$, $0\le t\le T$. 
\end{proof}

The remainder of this section is divided into four parts. In the first
one we examine the linear piece of
$\Pi_{u(s, \cdot), h}(x, \eta^n(s)) $. In the following three, we prove
Proposition \ref{s88} in dimensions $1$, $2$ and $3$, respectively.

\subsection*{The linear terms}

In this subsection, we replace the linear terms of the expression
$L_{n,t}^{*} \, \bs 1 \,-\, \partial_t \ln \psi^n_t$ by their space
average. The main result is summarized in the next lemma.

We adopt below the notation introduced in Section \ref{sec3}, see
equation \eqref{92}.  Hence $r_1\ge 2$, $r_2\ge 2$, $n=r_1 r_2$, and
$B_{\mtt k}$, $\mtt k\in \Sigma^+_{r_2}$, is a cube of length $r_1$.
The density profile here is $u( t, \cdot)$ introduced in Section
\ref{sec2}. 

Assume that $r_1$ is an odd number, $r_1 = 2b_1+1$, and recall
the notation introduced in \eqref{156}:
\begin{align*}
\cb{ J(y+  B_{\mtt k})} := \frac{1}{|\Sigma^+_{r_1}|}\,
\sum_{z\in y + B_{\mtt k}} J(z/n)\,, \quad
\cb{\mss m_{y, \mtt k}} := \mss m (y+  B_{\mtt k})
\,=\, \frac{1}{|\Sigma^+_{r_1}|} \,
\sum_{z\in y + B_{\mtt k}} \eta_z\,, \quad
y\in B_{\mtt 0}\,, \; \mtt k \in \Sigma^+_{r_2}\,.
\end{align*}
Note that
$E_{\nu^n_{u( t, \cdot)}} [\mss m_{y, \mtt k} ] = u(t, y+ B_{\mtt
k})$. Let
\begin{equation*}
\cb{W^{(\ref{s51})} _n} \,:=\,
\sum_{x\in \bb T^d_n} 
J(x/n)\, \big[\, \eta(x) - u( t, x/n) \,\big] \,
- \, \sum_{y\in B_{\mtt 0}} \sum_{\mtt k \in \Sigma^+_{r_2}}
J(y+ B_{\mtt k} ) \,
\big[\, \mss m_{y, \mtt k}  - u(t, y+ B_{\mtt k})  \,\big] \,.
\end{equation*}
If needed, to underline the dependence of $W^{(\ref{s51})} _n$ on some
of the parameters $r_1$, $J$ and $u(t)$, we write
$W^{(\ref{s51})} _n(r_1, J,u(t))$, or
$W^{(\ref{s51})} _n(r_1)$. We adopt the same convention for the error
$\mtt R^{(\ref{s51})} _n$ introduced below.

\begin{lemma}
\label{s51}
Fix a function $J\colon \bb T^d \to \bb R$ of class
$C^2(\bb T^d)$. Then, there exists a finite constant $C_1$ such that
\begin{gather*}
\int \big| \, W^{(\ref{s51})} _n \, \big| \, 
f\, d\nu^n_{u( t, \cdot)} \,\le\,
\mtt R^{(\ref{s51})} _n
\\
\text{where}\;\;
\cb{\mtt R^{(\ref{s51})} _n}\,:=\,
 \Vert J\Vert_{C^2} \,
\Big\{\, H_n(f \,|\, \nu^n_{u(t, \cdot)}  )  \,+\, \ln(2) \,+\,
C_1\, \Big( \frac{r_1}{n}\Big)^4  \, n^d\,\Big\}
\end{gather*}
for all densities $f$ with respect to $\nu^n_{u( t, \cdot)}$, and
$0\le t\le T$. 
\end{lemma}

\begin{proof}
Fix $0\le t\le T$, a function $J\colon \bb T^d \to \bb R$ of class
$C^2(\bb T^d)$, and a density $f$ with respect to
$\nu^n_{u( t, \cdot)}$. Recall the definition of the cube $\Sigma_b$
introduced in \eqref{190}, and let 
\begin{align*}
J^{(1)} (x/n)\,=\, \frac{1}{|\Sigma_{b_1}|}\,
\sum_{y\in \Sigma_{b_1}} J \Big( \frac{y+x}{n}\Big) \,, \quad
J^{(2)} (x/n)\,=\, \frac{1}{|\Sigma_{b_1}|}\,
\sum_{y\in \Sigma_{b_1}} J^{(1)} \Big( \frac{y+x}{n}\Big) \,.
\end{align*}
Note that we use here the centered cube $\Sigma_{b_1}$ instead of
$\Sigma^+_{b_1}$. The reason is that we wish the first order term in
the Taylor expansion
$\sum_{y\in \Sigma_{b_1}} [\, J((x+y)/n) - J(x/n)\,]$ to vanish.  As
$J$ is of class $C^2(\bb T^d)$, the absolute value of the difference
$J(x/n)\, -\, J^{(1)}(x/n) $ is bounded by
$C_1\, \Vert J\Vert_{C^2}\, (r_1/n)^2$ uniformly in $x$.  Similarly,
adding and subtracting $J^{(1)}(x/n)$, $J(x/n)\, -\, J^{(2)}(x/n) $ is
also bounded by $C_1 \, \Vert J\Vert_{C^2} \, (r_1/n)^2$. Therefore, by
Lemma \ref{s47} with $h(\eta) = \eta_0$, $\gamma = 1/\Vert J\Vert_{C^2}$
\begin{align*}
& \int \Big| \, \sum_{x\in \bb T^d_n} 
\big\{\, J(x/n)\, -\,  J^{(2)}(x/n) \, \big\}
\big[\, \eta(x) - u( t, x/n) \,\big] \,\Big | \, 
f\, d\nu^n_{u( t, \cdot)}
\\
&\quad
\,\le \,  \Vert J\Vert_{C^2} \,
\Big\{\, H_n(f \,|\, \nu^n_{u (t, \cdot)}  )  \,+\, \ln(2)  \,+\,
C_1\, \Big( \frac{r_1}{n}\Big)^4  \, n^d\,\Big\} \,.
\end{align*}
Changing variables, and by symmetry, 
\begin{align*}
& \sum_{x\in \bb T^d_n}  J^{(2)}(x/n) \, 
\big[\, \eta(x) - u( t, x/n) \,\big]
\, =\, \sum_{x\in \bb T^d_n}  J^{(1)}(x/n) \,
\frac{1}{|\Sigma_{b_1}|}\,
\sum_{y\in \Sigma_{b_1}} 
\big[\, \eta(x+y) - u( t, [x+y]/n) \,\big]
\\
&\quad =\, 
\sum_{x\in \bb T^d_n} \,
\frac{1}{|\Sigma_{b_1}|}\,
\sum_{z\in \Sigma_{b_1}} J \Big( \frac{z+x}{n}\Big)
\frac{1}{|\Sigma_{b_1}|}\,
\sum_{y\in \Sigma_{b_1}} 
\big[\, \eta(x+y) - u( t, [x+y]/n) \,\big]\,.
\end{align*}
As $r_1 = 2b_1+1$, $|\Sigma_{b_1}| = |\Sigma^+_{r_1}|$ and
$x+\Sigma_{b_1} = x - \mb b_1 + \Sigma^+_{r_1} $, where
$\mb b_1 = (b_1, \dots, b_1)$. Thus the previous expression can be
rewritten as
\begin{align*}
\sum_{y\in B_{\mtt 0}} \sum_{\mtt k \in \Sigma^+_{r_2}}
J(y+ B_{\mtt k} ) \,
\big[\, \mss m_{y, \mtt k}  - u(t, y+ B_{\mtt k}) \,\big]\,,
\end{align*}
which completes the proof of the lemma.
\end{proof}

Consider the formula for
$L_{n,t}^{*} \, \bs 1 \,-\, \partial_t \ln \psi^n_t$ provided by Lemma
\ref{s41}. The linear term $\eta_x - u(t,x/n)$ of the first term on
the right-hand side comes multiplied by the product
$F(x/n) \, G(x/n)$, where $F = (\partial_{x_1}^2 u)/\chi(u)$,
$G= \mtt p'_{h_1} (u)$. The previous lemma replaces the product $FG$
by an average of the product of these functions. In the proof of
Theorem \ref{t01}, one needs instead the product of the averages. This
is the content of the next result.

Let
\begin{align*}
\cb{W^{(\ref{s52})}_n} \,: &=\,
\sum_{x\in \bb T^d_n} 
F(x/n)\, G(x/n)\,  \big[\, \eta(x) - u( t, x/n) \,\big]
\\
\, &-\,
\sum_{y\in B_{\mtt 0}} \sum_{\mtt k \in \Sigma^+_{r_2}}
F(y+ B_{\mtt k} ) \, G(y+ B_{\mtt k} ) \,
\big[\, \mss m_{y+ B_{\mtt k}}  -
u(t, y+ B_{\mtt k}) \,\big] \,.
\end{align*}

\begin{corollary}
\label{s52}
Fix two functions $F$, $G\colon \bb T^d \to \bb R$ of class
$C^2(\bb T^d)$. Then, there exists a finite constant $C_1$ such that
\begin{gather*}
\int \big | \, W^{(\ref{s52})}_n \, \big|
\, f\, d\nu^n_{u( t, \cdot)} \,\le \, \mtt R^{(\ref{s52})}_n\,,
\\
\text{where}\;\;
\cb{\mtt R^{(\ref{s52})}_n}\,:=\,
C_1\, \Vert F\Vert_{C^2}\, \Vert G\Vert_{C^2} \,\Big\{
H_n(f \,|\, \nu^n_{u (t, \cdot)}  )  \,+\, \ln(2) \,+\,
\Big( \frac{r_1}{n}\Big)^4  \, n^d\,\Big\}\,,
\end{gather*}
for all densities $f$ with respect to $\nu^n_{u( t, \cdot)}$, and
$0\le t\le T$.
\end{corollary}

\begin{proof}
Fix two functions $F$, $G$ in $C^2(\bb T^d)$.  Let
\begin{align*}
H(y+B_{\mtt k}) \, &=\, \frac{1}{| B_{\mtt k}|} \,  \sum_{x\in y+B_{\mtt k}}
\big\{\, F(x/n) - F(y+ B_{\mtt k})\big\}\,
\big\{\, G(x/n) - G(y+B_{\mtt k})\big\}
\\
&=\, \frac{1}{| B_{\mtt k}|} \,  \sum_{x\in y+B_{\mtt k}}
F(x/n) \, G(x/n) \,-\, F(y+B_{\mtt k}) \, G(y+B_{\mtt k}) \,.
\end{align*}
As $F$, $G$ are Lipschitz-continuous,
$|H( y + B_{\mtt k})|\le \Vert F\Vert_{C^1}\, \Vert G\Vert_{C^1}\,
(r_1/n)^2$. Thus, summing by parts to recover $\eta_x - u(t,x/n)$ from
$\mss m_{y+ B_{\mtt k}} - u(t, y+ B_{\mtt k})$, by Lemma \ref{s47}
with $h=\eta_0$, $\gamma =1/(\Vert F\Vert_{C^1}\, \Vert G\Vert_{C^1})$,
\begin{align*}
& \int \Big| \, \sum_{y\in B_{\mtt 0}} \sum_{\mtt k \in \Sigma^+_{r_2}}
H(y+ B_{\mtt k} ) \,
\big[\, \mss m_{y+ B_{\mtt k}}
- u(t, y+ B_{\mtt k}) \,\big] \, \Big | \, 
f\, d\nu^n_{u( t, \cdot)}
\\
& \quad \,\le \,
\Vert F\Vert_{C^1}\, \Vert G\Vert_{C^1} \,\Big\{
H_n(f \,|\, \nu^n_{u(t, \cdot)}  )  \,+\, \ln(2) \,+\,
C_1\, \Big( \frac{r_1}{n}\Big)^4  \, n^d\,\Big\}
\end{align*}
for all densities $f$ with respect to $\nu^n_{u(t, \cdot)} $.  The
assertion of the corollary follows from Lemma \ref{s51} with $J = FG$
and this bound.
\end{proof}

In the previous corollary, suppose that $G(x) = H(u(t,x))$ for
some function $H$ of class $C^2(\bb R)$. Let
\begin{align*}
J(y+ B_{\mtt k})  
\,:=\, \frac{1}{| B_{\mtt k}|} \,  \sum_{x\in y+B_{\mtt k}} H( u(t,x/n))
\,-\, H \big(   u(t, y+ B_{\mtt k}) \big) \,.
\end{align*}
A Taylor expansion yields that the absolute value of
$ J(y+ B_{\mtt k} )$ is bounded by
$C_1 \, \Vert H\Vert_{C^2} \, \Vert u_t \Vert_{C^2} \, (1+ \Vert u_t
\Vert_{C^2}) \, (r_1/n)^2$ uniformly in $y$, where
$u_t = u(t, \cdot)$. Therefore, by a summation by parts, to recover
the function $\eta_x - u(t,x/n)$, and Lemma \ref{s47} with
$h = \eta_0$ and
$\gamma=1/[\, \Vert F\Vert_{\infty} \, \Vert H\Vert_{C^2} \,\Vert
u_t\Vert_{C^2}\, ( 1+ \Vert u_t\Vert_{C^2})\,]$,
\begin{align*}
& \int \Big | \, \sum_{y\in B_{\mtt 0}}  \sum_{\mtt k \in \Sigma^+_{r_2}}
F(y+ B_{\mtt k} ) \, J(y+ B_{\mtt k} ) \,
\big[\, \mss m_{y+ B_{\mtt k}}
- u(t, y+ B_{\mtt k})  \,\big] \, \Big| \, 
f\, d\nu^n_{u( t, \cdot)}
\\
& \quad \,\le \,
\Big\{ \, \Vert F\Vert_{\infty} \, \Vert H\Vert_{C^2} \, \Vert u_t\Vert_{C^2}
\, ( 1 + \Vert u_t\Vert_{C^2})\, \Big\} \,\Big\{\, 
H_n(f \,|\, \nu^n_{u (t, \cdot)}  )  \,+\, \ln(2)\,+\,
C_1\, \Big( \frac{r_1}{n}\Big)^4  \, n^d\, \Big\}\,.
\end{align*}
Hence, if we define $V^{\eqref{135}}_n$ by
\begin{align*}
\cb{ V^{\eqref{135}}_n} \, &:=\,
\sum_{x\in \bb T^d_n} 
F(x/n)\, H(u(t, x/n)) \,  \big[\, \eta(x) - u( t, x/n) \,\big]
\\
\,&-\,
\sum_{y\in B_{\mtt 0}} \sum_{\mtt k \in \Sigma^+_{r_2}}
F(y+ B_{\mtt k} ) \,
H( u(t, y+ B_{\mtt k})  )  \,
\big[\, \mss m_{y+ B_{\mtt k}}  - u(t, y+ B_{\mtt k}) 
\,\big] \,,
\end{align*}
by Corollary \ref{s52},
\begin{align}
\label{135}
\int \big | \, V^{(\ref{135})}_n\, \big|\,  f\,
d\nu^n_{u( t, \cdot)} \,\le\,  \mtt S^{(\ref{135})}_{n,r_1} (F, H, u_t)
\end{align}
where
\begin{equation*}
\cb{\mtt S^{\eqref{135}}_{n,r_1} (F, H, u_t)} \,:=\,
C_1\, \Big\{ \, \Vert F\Vert_{C^2} \, \Vert H\Vert_{C^2}
\, \Vert u_t\Vert_{C^2}
\, ( 1 + \Vert u_t\Vert_{C^2})\, \Big\} \,
\Big\{\, H_n(f \,|\, \nu^n_{u(t,\cdot)}  )  \,+\,  \ln(2)
\,+\, \Big( \frac{r_1}{n}\Big)^4  \, n^d\,\Big\} \,.
\end{equation*}

Before we proceed to the proof of Theorem \ref{t01}, we may modify the
statements of Theorems \ref{s01}, \ref{s26} and Corollaries \ref{s46},
\ref{s48} replacing
$E_{\nu^n_{u(t, \cdot)} } [\, \mtt p_h (\mss m_{\mtt k})]$ by
$\mtt p_h ( u(t, B_{\mtt k})) $ in order to obtain a second-order
Taylor expansion. Indeed, let $\rho \colon \bb T^d \to (0,1)$ be a
continuous density profile. As $\mtt p_h(\cdot)$ is a polynomial, the
computation performed in the proof of Lemma \ref{s22} yields that
there exists a constant $C_0=C_0(h)$ such that
\begin{equation}
\label{120}
\big|\, E_{\nu^n_{\rho(\cdot)} } [\, \mtt p_h (\mss m_{\mtt k})]
\,-\, \mtt p_h (\rho (B_{\mtt k}) )\,\big|
\,\le\, C_0 \frac{1}{r_1^d}
\end{equation}
for all $n\ge 1$.

\subsection*{Dimension $d=1$}

In this subsection, we denote the lengths $r_1$, $r_2$ introduced in
Section \ref{sec4}, see equation \eqref{95}, by $\ell_1$, $\ell_2$,
respectively.  Hence $\ell_1\ge \ell_h$, $\ell_2\ge 2$,
$n=\ell_1 \ell_2$, and $B_{\mtt k}$, $\mtt k\in \Sigma^+_{\ell_2}$, is
a cube of length $\ell_1$. In the proof below $\ell_1$ is chosen to be
$n^{(2/3)-b}$, where $0<b<1/6$ is an arbitrarily small and fixed
parameter.

\begin{proof}[Proof of Proposition \ref{s88}]
Fix $0<b<1/6$, $\delta_0>0$, and $0\le s\le t \le T$. Recall equation
\eqref{115} to rewrite the sum appearing on the left-hand side of the
statement of the proposition as
\begin{align}
\label{180}
\frac{1}{|B_{\mtt 0} |}\, 
\sum_{y \in B_{\mtt 0}}  \sum_{\mtt k\in \Sigma^+_{\ell_2}}
\sum_{x\in y+ B^o_{\mtt k}} 
\frac{|B_{\mtt 0} |}{|B^o_{\mtt 0} |}\, 
G(s,x/n)\,  \Pi_{u(s, \cdot) , h} (x, \eta)\,.
\end{align}
We consider below the case $y=0$ and make sure at each step that the
estimate is uniform over $y \in B_{\mtt 0}$.  The proof in dimension
$1$ does not require a multi-scale analysis and relies solely on
Corollary \ref{s46} and the estimates presented in the current
section. Recall the notation introduced in that corollary, in
particular the constant $\mf a'_2 = \mf a'_2(\mf c_{\rm SG})$ which will
appear below.

Set $\ell_1 = n^a$, where $a=(2/3) - b$.  There exists a constant
$C_0 = C_0(h)$ such that the ratio $|B_{\mtt 0} |/|B^o_{\mtt 0} |$ is
bounded by $C_0$. As $b>0$, there exists $\cb{n_0} : = n_0(h,G, \delta)$
such that
$\Vert h\Vert_\infty\, \Vert G(s)\Vert_\infty\, \ell^{d+2}_1 \le \mf
a'_2 \,\delta\, n^2$ for all $n\ge n_0$, $0\le s\le T$.  Thus, by
Corollary \ref{s46} with $\rho(\cdot) = u(s, \cdot)$,
$\delta = \delta_0/C_0$, and $L_2 = \ell_2$, and by \eqref{120} with
$r_1=\ell_1$,
\begin{equation*}
\begin{aligned}
& \frac{|B_{\mtt 0} |}{|B^o_{\mtt 0} |} \sum_{\mtt k\in
\Sigma^+_{\ell_2}} \int \sum_{x\in B^o_{\mtt k}} G(s,x/n)\,
\Big(\, \tau_x h \, - \, E_{\nu^n_{u( s,\cdot)} } \big[\, \tau_x h
\,\big] \,-\,  \big[\, \mtt
p_h (\mss m_{\mtt k}) \,-\, \mtt p_h ( u(s, \mtt k)) \,\big] \,
\Big)\, f_s \, d\nu^n_{u(s,\cdot)}
\\
& \quad \,\le\, \delta_0\, n^2\, I_n(f_s ; \nu^n_{u( s,\cdot)} )
\,+\, \mf c_{1,3} \, \mtt K_{1}(G_s)\,
\Big\{ \,  H_n (f_s \,|\, \nu^n_{u( s,\cdot)} ) \,+\,
n^{2a-1} \,+\, n^{1-2a} \,+\, n^{1-a} \,\Big\}
\end{aligned}
\end{equation*}
for all $n\ge n_0$. In this formula,
$\cb{\mtt K_{1}(G_s) } = \Vert G(s) \Vert_{\infty} ( 1+ \Vert G(s)
\Vert_{\infty}) + \Vert G(s) \Vert_{C^1}$.  As $a>0$,
$n^{1-2a} \le n^{1-a}$ and $1\le n^{2a-1} + n^{1-a} $. The constant
$\delta = \delta_0/C_0$ multiplying $I_n(f_s ; \nu^n_{u( s,\cdot)} )$
becomes $\delta_0$ because we overestimate the ratio
$|B_{\mtt 0} |/|B^o_{\mtt 0} |$ by $C_0$.

On the other hand, the inequality \eqref{135}, with $r_1=\ell_1$,
$F(\cdot) = G(s, \cdot)$, $H(\cdot) = \mtt p_h'(\cdot)$, takes care of
the linear term in $\Pi_{u(s, \cdot) , h} (x, \eta)$ at an extra cost
$n^{4a-3}$. This estimate and the previous one yield that
\begin{equation*}
\begin{aligned}
& \int
\frac{|B_{\mtt 0} |}{|B^o_{\mtt 0} |}\,
\sum_{\mtt k\in \Sigma^+_{\ell_2}}
\sum_{x\in B^o_{\mtt k}} 
G(s,x/n)\,  \Pi_{u(s, \cdot) , h} (x, \eta)\, 
f_s \, d\nu^n_{u( s, \cdot)}
\\
& \quad \le
\int \sum_{\mtt k\in \Sigma^+_{\ell_2}}  
\ell_1^d\,  G (s, \mtt k) \, \Big\{ \, \mtt p_h(\mss m_{\mtt k})
\,-\, \mtt p_h( u(s, \mtt k)) \, -\,
\mtt p_h' ( u(s,\mtt k) ) 
\, [\,\mss m_{\mtt k}  - u(s,\mtt k)\,]\,
\Big\}\,
f_s \, d\nu^n_{u(s,\cdot)}
\\
& \quad
\,+\, \delta_0 \, n^2\, I_n(f_s ; \nu^n_{u( s,\cdot)}  )
\,+\, \mf c_{1,3} \, \mtt K_{1,1}(G_s)   \,
\Big\{   \, H_n (f_s\,|\, \nu^n_{u( s,\cdot)}  )
\,+\, n^{2a-1} \,+\, n^{1-a} \,\Big\} \,.
\end{aligned}
\end{equation*}
We used here that $4a -3 < 1-a$ because $a\le 2/3$.  By Taylor's
development, the expression inside braces in the first line is bounded
by $C_0 (h) [\,\mss m_{\mtt k} - u(s,\mtt k)\,]^2$. By Lemma \ref{l02}
with $r_1=\ell_1$, the first line is therefore bounded by
$C_1 \, \Vert G(s) \Vert_\infty \, \{ \, H_n (f_s\,|\, \nu^n_{u(
s,\cdot)} ) + n^{1-a}\}$. To complete the proof of the proposition for
$n$ large, it remains to recollect all previous estimates. 
\end{proof}

\subsection*{Dimension $d=2$}

In dimension $2$ a multiscale analysis is required. We adopt the
notation of Section \ref{sec9} denoting the lengths by $\ell_i$
instead of $\mtt t_i$. In the proof below,  $\ell_1 = n^{2/5}$,
$\ell_2 = n^{1/5}$, so that $\ell_1\ell_2 = n^{3/5}$.

\begin{proof}[Proof of Proposition \ref{s88}]
Fix $\delta_0>0$, $0\le s\le t \le T$. As in dimension $1$, we
estimate the expectation of \eqref{180} for $y=0$ with respect to
$f_s \, d\nu^n_{u(s, \cdot)}$.  The argument is based on Corollaries
\ref{s46}, \ref{s48} and the estimates presented in this section.

Set $\ell_1 = n^{2/5}$.  In dimension $2$, Corollary \ref{s46}
requires $\ell_1$ to be bounded by $c_1 n^{1/2}$ for some constant
$c_1$. This condition is thus satisfied. As in the proof in one
dimension, by Corollary \ref{s46} with $\delta=\delta_0/[2 C_0(h)]$, and
the previous choice of $\ell_1$, setting $n = \ell_1\, L_2$,
\begin{align*}
& \frac{|B_{\mtt 0} |}{|B^o_{\mtt 0} |}  \sum_{\mtt k\in \Sigma^+_{L_2}}  \int
\sum_{x\in B^o_{\mtt k}} G(s,x/n)\,
\Big\{   \, 
\tau_x h \, -
\, E_{\nu^n_{u( s,\cdot)} } \big[\, \tau_x h \,\big]
\,-\,   \big\{ \, \mtt p_h(\mss m_{\mtt k})
\,-\, E_{\nu^n_{u( s,\cdot)} } \big[\,
\mtt p_h(\mss m_{\mtt k}) \,\big] \,\big\} \, \Big\} \,
f_s\, d\nu^n_{u(s,\cdot)}
\nonumber
\\
& \quad
\,\le\, \frac{\delta_0}{2} \, n^2\, I_n(f_s ; \nu^n_{u( s,\cdot)}  )
\,+\, \mf c_{1,3} \, \mtt K_1 (G_s) \,
\Big\{   \, H_n (f_s \,|\, \nu^n_{u( s,\cdot)}  )
\,+\, n^{4/5} \,\Big\} 
\end{align*}
for all $n\ge n_0(h, G, \delta)$. The constant $\mtt K_1 (G_s) $ is
introduced in the previous subsection.

Since in dimension $d=2$ in Corollary \ref{s46} we are forced to take
cubes of side length $\ell_1 \ll \sqrt{n}$, the last step in the proof of
Proposition \ref{s88}, which consists in applying Lemma \ref{l02} (see
the proof in dimension $1$), does not provide a good bound for cubes
of this size. We are thus led to appeal to the multi-scale analysis to
increase the size of the cube. This is the role of Corollary
\ref{s48}.

Set $\ell_2 = n^{1/5}$, $\ell_3 = n^{2/5}$ so that
$\ell_1\ell_2\ell_3 = n$.  With these choices, the hypotheses of
Corollary \ref{s48} are satisfied for $K_1=1$ and $\epsilon <1$. Thus
by this result with $\delta =\delta_0/2$,
\begin{align*}
& \int
\sum_{\mtt j \in \Sigma^+_{\ell_3}}
\sum_{\mtt k \in \Lambda_{\mtt j}}  \ell_1^d\, G(s, B_\mtt k) \, \Big\{\,
\mtt p_h(\mss m_{\mtt k})
\,-\, E_{\nu^n_{u( s,\cdot)} } \big[\,
\mtt p_h(\mss m_{\mtt k}) \,\big]
\,-\, 
\big(\, \mtt p_h (\mss m_{2, \mtt j})
- E_{\nu^n_{u( s,\cdot)} } [\, \mtt p_h (\mss m_{2, \mtt j}) \, ] \,\big)
\, \Big\}\,
f_s \, d\nu^n_{u( s,\cdot)}
\nonumber
\\
&  \le\, \frac{\delta_0}{2} \, n^2 \,  I_n  ( f_s \,;\, \nu^n_{u (s, \cdot)} )
\,+\, \mf c_{1,3} \, \mtt K_2(G(s))
\,\Big\{   \, H_n (f_s \,|\, \nu^n_{u(s, \cdot)}) \,+\,
n^{4/5} \,\Big\}
\end{align*}
for all $n\ge n_1$, where $\cb{n_1} = n_1(\mf c_{1,3}, h,
\delta_0)$. In this formula,
$\cb{\mtt K_2(G(s))} := [1+\Vert G(s)\Vert_{C^1}]^2 \, \exp\{\mf
c_{1,3} [ 1+\Vert G(s)\Vert_{C^1}]\,\}$.

The two previous bounds together with \eqref{120} with
$r_1=\ell_1\ell_2 = n^{3/5}$ [replacing
$E_{\nu^n_{u( s,\cdot)} } [\, \mtt p_h (\mss m_{2, \mtt j}) \, ]$ by
$\mtt p_h(u(s, B_{2,\mtt j}))$] yield that
\begin{align*}
& \frac{|B_{\mtt 0} |}{|B^o_{\mtt 0} |}
\sum_{\mtt j \in \Sigma^+_{\ell_3}}
\sum_{\mtt k \in \Lambda_{\mtt j}}   \int
\sum_{x\in B^o_{\mtt k}} G(s,x/n)\,
\Big\{   \, 
\tau_x h \, -
\, E_{\nu^n_{u( s,\cdot)} } \big[\, \tau_x h \,\big]
\,-\,   \big\{ \, \mtt p_h(\mss m_{2,\mtt j})
\,-\, 
\mtt p_h(u(s,  B_{2,\mtt j}))  \,\big\} \, \Big\} \,
f_s\, d\nu^n_{u(s,\cdot)}
\nonumber
\\
& \quad
\,\le\, \delta_0 \, n^2\, I_n(f_s ; \nu^n_{u( s,\cdot)}  )
\,+\, \mf c_{1,3} \, \mtt K_2(G(s))\, \Big\{   \, H_n (f_s\,|\, \nu^n_{u( s,\cdot)}  )
\,+\, n^{4/5} \,\Big\} 
\end{align*}
for all $n$ sufficiently large.

The inequality \eqref{135}, with $r_1=\ell_1\ell_2$,
$F(\cdot) = G(s, \cdot)$, $H(\cdot) = \mtt p_h'(\cdot)$, takes care of
the linear term in $\Pi_{u(s, \cdot) , h} (x, \eta)$ at an extra cost
$n^{2/5}$. This estimate and the previous one yield that
\begin{equation*}
\begin{aligned}
& \int
\frac{|B_{\mtt 0} |}{|B^o_{\mtt 0} |}\,
\sum_{\mtt k\in \Sigma^+_{L_2}}
\sum_{x\in B^o_{\mtt k}} 
G(s,x/n)\,  \Pi_{u(s, \cdot) , h} (x, \eta)\, 
f_s \, d\nu^n_{u( s, \cdot)}
\\
& \quad \le
\int \sum_{\mtt j\in \Sigma^+_{\ell_3}}  
(\ell_1\ell_2)^d\,  G (s, B_{2,\mtt j}) \, \Big\{ \, \mtt p_h(\mss m_{2,\mtt j})
\,-\, \mtt p_h( u(s, B_{2,\mtt j})) \, -\,
\mtt p_h' ( u(s, B_{2,\mtt j}) ) 
\, [\,\mss m_{2, \mtt j}  - u(s,B_{2,\mtt j} )\,]\,
\Big\}\,
f_s \, d\nu^n_{u(s,\cdot)}
\\
& \quad
\,+\, \delta_0 \, n^2\, I_n(f_s ; \nu^n_{u( s,\cdot)}  )
\,+\, \mf c_{1,3} \, \mtt K_3(G(s)) \, \Big\{   \, H_n (f_s\,|\, \nu^n_{u( s,\cdot)}  )
\,+\, n^{4/5} \,\Big\} \,,
\end{aligned}
\end{equation*}
where $\cb{ \mtt K_3(G(s)) } := \, \mtt K_2(G(s)) + \Vert G(s)\Vert_{C^2}$.
By Lemma \ref{l02} with $r_1=\ell_1\ell_2$, this expression is less
than or equal to
\begin{align*}
\delta_0  \, n^2 \,  I_n  ( f_s \,;\, \nu^n_{u (s, \cdot)} )
\,+\, \mf c_{1,3} \,\mtt K_3(G(s)) \,
\Big\{   \, H_n (f_s\,|\, \nu^n_{u(s, \cdot)}) \,+\, n^{4/5} \,\Big\}\,.
\end{align*}
This completes the proof in dimension $2$.
\end{proof}

\subsection*{Dimension $d=3$}

The proof in dimension $3$ differs substantially from the ones in
dimensions $1$ and $2$. In dimension $3$, we need to increase twice the
size of the boxes, we cannot neglect the correction which appears
when replacing canonical expectations by grand canonical ones, and we
need to add a linear term to compensate correction terms added on
the way.

Assume that $n = \ell_1 \ell_2 \ell_3 L_4$ for some sequences of
integers $\ell_j = \ell_j(n)$, $L_4(n)$. Below we will choose
$\ell_1= n^{1/4-\mtt e}$, $\ell_2= n^{3/20- \mtt e}$,
$\ell_3= \kappa \, n^{1/10+2\mtt e}$ for some $0<\mtt e< 1/30$ and a
constant $\kappa$ which increases to $\infty$ after $n\to\infty$.
With this choice, $\ell_1 \ell_2 \ell_3 = \kappa \, n^{1/2}$.

We adopt the notation introduced at the beginning of Section
\ref{sec9}. Divide the cube $\Sigma^+_n$ into $L_4^d$ disjoint cubes of
length $\ell_1\ell_2\ell_3$, represented by $\cb{B_{3,\mtt i}}$,
$\mtt i =(i_1, \dots, i_d) \in \Sigma^+_{L_4}$.  Thus
$\Sigma^+_n = \cup_{\mtt i\in \Sigma^+_{L_4}} B_{3,\mtt i} =
\cup_{\mtt j\in \Sigma^+_{L_3}} B_{2, \mtt j} $, where
$\cb{L_3 } = \ell_3L_4$.  For a cube $B_{3, \mtt i}$, denote by
$\Lambda_{3, \mtt i}$ the set of indices $\mtt j$ of the sets
$B_{2,\mtt j}$ which are contained in $B_{3, \mtt i}$:
\begin{equation}
\label{119}
\begin{gathered}
\cb{\Lambda_{3, \mtt i}} \,:=\,
\big\{ \mtt j \in \Sigma^+_{L_3} :  B_{2,\mtt j} \subset
B_{3, \mtt i} \, \big\}\,, \;\;\text{so that}\;\;
B_{3, \mtt i}  \,=\, \bigcup_{\mtt j \in  \Lambda_{3,\mtt i}}
B_{2,\mtt j} \,.
\end{gathered}
\end{equation}
Recall the definitions \eqref{156} and set
$\cb{\mss M_{3, \mtt i}} = \mss M(B_{3, \mtt i})$,
$\cb{\mss m_{3, \mtt i}} = \mss m(B_{3, \mtt i})$,
$\cb{\mss m_{3, \mtt i, y}} = \mss m(y+B_{3, \mtt i})$.  Recall,
furthermore, the definition of the polynomial $\mtt q_h(\cdot)$
introduced in \eqref{90}, and the convention on the constants
$\mf a_{3,p}$ established at the beginning of Section \ref{sec4}.

\begin{lemma}
\label{s49}
Fix a density profile $\rho\colon \bb T^d\to (0,1)$ satisfying the
hypotheses of Theorem \ref{s01}, a cylinder function $h$, a function
$J\colon \bb T^d \to \bb R$ of class $C^1$, and constants $\delta'>0$,
$0<\mtt e<1/30$, and $\kappa \ge 1$. Let $\ell_1= n^{1/4-\mtt e}$,
$\ell_2= n^{3/20- \mtt e}$, $\ell_3= \kappa \, n^{1/10+2\mtt e}$.
Then, there exist finite constants $\mf a_{3,3} = \mf a_{3,3} (h)$,
and an integer
$n_0=n_0(\kappa, \delta', \mf r, \Vert \rho\Vert_{C^3}, h, \Vert
J\Vert_\infty)$ such that
\begin{align*}
& \sum_{x  \in \bb T^d_{n}}  \int
J(x/n)\, \big\{   \,  \tau_x h \, -
\, E_{\nu^n_{\rho(\cdot)} } \big[\, \tau_x h \,\big] \, \big\}\,
f\, d\nu^n_{\rho(\cdot)}
\\
& \quad 
\, \le\,
\frac{1}{|B_{\mtt 0} |}\, 
\sum_{y \in B_{\mtt 0}}  \sum_{\mtt i\in \Sigma^+_{L_4}}  \int
(\ell_1\ell_2\ell_3)^d J(B_{3, \mtt i, y} )\, \big\{   \, 
\mtt p_h(\mss m_{3,\mtt i, y})
\,-\, E_{\nu^n_{\rho(\cdot)} } \big[\,
\mtt p_h(\mss m_{3,\mtt i, y}) \,\big] \,\big\} \, 
f\, d\nu^n_{\rho(\cdot)}
\nonumber
\\
& \quad
\,+\, \delta' \, n^2\, I_n(f ; \nu^n_{\rho(\cdot)}  )
\,+\, \mf a_{3,3}\, \Vert J\Vert_{C^1}\, H_n (f\,|\, \nu^n_{\rho(\cdot)}  )
\,+\,  \,\big[\, \mf a_{3,3} (\delta', h) \,  \Vert J\Vert_{C^1}
\,+\, \mtt K_4(\kappa, J) \,\big] \, n^{3/2 - 2\mtt e}  
\end{align*}
for all densities $f$ with respect to $\nu^n_{\rho(\cdot)}$, and
$n\ge n_0$. In this formula, $B_{3, \mtt i, y} = y + B_{3, \mtt i}$,
and
$\cb{\mtt K_4(\kappa, J) } = \kappa^{1/2}\, e^{(\mf c_{\rm
LS}/\delta') \,\kappa^6\, [1+\mf a_{3,3}]\,[1 + \Vert J\Vert_{C^1}]}
$.
\end{lemma}

\begin{proof}
Let $\cb{L_2} := \ell_2\ell_3L_4$, and recall that $L_3 = \ell_3 L_4$.
Recall equation \eqref{180}. We restrict our attention to the term
$y=0$, and estimate the expression
\begin{align}
\label{181}
\frac{|B_{\mtt 0} |}{|B^o_{\mtt 0} |}\,
\int 
\sum_{\mtt k\in \Sigma^+_{L_2}}
\sum_{x\in B^o_{\mtt k}} 
J(x/n)\, \big\{   \,  \tau_x h \, -
\, E_{\nu^n_{\rho(\cdot)} } \big[\, \tau_x h \,\big] \, \big\}
\, f\, d\nu^n_{\rho(\cdot)} \, .
\end{align}

Fix $\delta'>0$ and $0<\mtt e<1/30$. By Theorem \ref{s01} with
$\delta = \delta'/[8C_0(h)]$, $\ell_1 = n^{1/4 - \mtt e}$, as $\mtt e<1/30$, there
exists a finite constant $\mf a_{3,3}$, and an integer
$n_0 = n_0(\delta', h, \Vert J\Vert_\infty)$ such that
\eqref{181} is bounded above by
\begin{align}
\label{121b}
& \sum_{\mtt k\in \Sigma^+_{L_2}}  \int
\ell_1^d \, J(B_{\mtt k})\, \big\{ \, \mtt p_h(\mss m_{\mtt k})
\,-\, E_{\nu^n_{\rho(\cdot)} } \big[\,
\mtt p_h (\mss m_{\mtt k}) \,\big] \,\big\} \, 
f\, d\nu^n_{\rho(\cdot)}
\nonumber
\\
& \quad
\,+ \,  \sum_{\mtt k\in \Sigma^+_{L_2}}  \int
J (B_{\mtt k})  \, \Big\{ \, \mtt q_h(\mss m_{\mtt k}) 
\,-\,  E_{\nu^n_{\rho(\cdot)}  } \big[\,
\mtt q_h(\mss m_{\mtt k}) \,\big] 
\, \Big\}\,
f\, d\nu^n_{\rho(\cdot)}
\\
& \quad
\,+\, \frac{\delta'}{8} \, n^2\, I_n(f ; \nu^n_{\rho(\cdot)}  )
\,+\, \mf a_{3,1}\, \Vert J\Vert_{\infty}
\, H_n (f\,|\, \nu^n_{\rho(\cdot)}  )
\,+\, \mf a_{3,3}\, \Vert J\Vert_{C^1}\,
\frac{1}{\delta'}\, n^{3/2 - 2\mtt e}  
\nonumber
\end{align}
for all $n\ge n_0$.  Note that the term $\mtt q_h$, via $W^h_{\mtt
k,x}$ came with a multiplicative constant $\ell_1^{-d}$ which canceled
the factor $\ell_1^d$ as the sum in $\mtt k$ is performed over
$(n/\ell_1)^d$ terms only. 
As we are far from cubes of side length $n^{1/2}$, we apply Theorem
\ref{s26} to increase their size. We examine in detail the expression
$\mtt p_h(\mss m_{\mtt k}) \,-\, E_{\nu^n_{\rho(\cdot)} } [\, \mtt
p_h(\mss m_{\mtt k}) \,] $. The same arguments apply to
$\mtt q_h(\mss m_{\mtt k}) \,-\, E_{\nu^n_{\rho(\cdot)} } [\, \mtt
q_h(\mss m_{\mtt k}) \,] $, with a simplification due to the extra
implicit factor $\ell_1^{-d}$. Recall the definition of $W_{\mf p, \mtt j}$
introduced in \eqref{111}.

Fix $0<\epsilon<1$, and recall that $\ell_2 = n^{3/20 - \mtt e}$. A direct
calculation shows that there exists an integer
$n_1 = n_1(\epsilon, \delta', \mf r, \Vert \rho\Vert_{C^3}, h)$ such
that for all $n\ge n_1$, the conditions of Theorem \ref{s26} with
$\mtt t_1 = \ell_1$, $\mtt t_2 =\ell_2$ are satisfied with
$K_1=K_2=1$.  Thus, applying Theorem \ref{s26} with these choices and
$\delta = \delta'/8$, we obtain a finite constant $\mf a_{3,3}$ such
that
\begin{align}
\label{122b}
& \int
\sum_{\mtt j \in \Sigma^+_{L_3}}
\sum_{\mtt k \in \Lambda_{\mtt j}}  \ell_1^d\, J(B_{\mtt k}) \, \Big\{\,
\mtt p_h(\mss m_{\mtt k})
\,-\, E_{\nu^n_{\rho(\cdot)} } \big[\,
\mtt p_h(\mss m_{\mtt k}) \,\big]
\,-\, 
\big(\, W_{\mtt p_h}   (\mss m_{2, \mtt j})
- E_{\nu^n_{\rho(\cdot)} } [\, W_{\mtt p_h}   (\mss m_{2, \mtt j})  \, ] \,\big)
\, \Big\}\,
f\, d\nu^n_{\rho(\cdot)}
\nonumber
\\
&  \le\, \frac{\delta'}{8} \, n^2 \,  I_n  ( f ; \nu^n_{\rho (\cdot)} )
\,+\, \mf a_{3,3}\, \Vert J\Vert_{C^1}  \, H_n (f\,|\, \nu^n_{\rho(\cdot)}) \,+\,
e^{(\mf c_{\rm LS}/\delta')
\, [1+\mf a_{3,3}]\,[1 + \Vert J\Vert_{C^1}]} \,n^{3/2 - 4\mtt e} 
\end{align}
for all $n\ge n_1$. The right-hand side of Theorem \ref{s26} has many
terms with different powers of $n$. To bound all of them by
$n^{3/2 - 4\mtt e}$ we used the assumption that we selected a constant
$\mtt e<1/30$.

As $\ell_1 \ell_2 = n^{2/5 - 2\mtt e} \ll n^{1/2}$ we need to apply
Theorem \ref{s26} once again to increase the size of the cube to reach
lengths of order $n^{1/2}$. By \eqref{111},
$W_{\mtt p_h} (\mss m_{2, \mtt j})$ is the sum of two terms. We
examine in detail the expression
$\mtt p_h(\mss m_{2, \mtt j}) \,-\, E_{\nu^n_{\rho(\cdot)} } [\, \mtt
p_h(\mss m_{2, \mtt j}) \,] $. The same arguments apply to
$\mtt q_h(\mss m_{2,\mtt j}) \,-\, E_{\nu^n_{\rho(\cdot)} } [\, \mtt
q_h(\mss m_{2, \mtt j}) \,] $, with a simplification due to the extra
factor $\ell_1^{-d}$.

For a polynomial $\mf p\colon [0,1]\to \bb R$, 
let
\begin{equation*}
\begin{gathered}
\cb{ W^{(3)}_{\mf p, \mtt i} } \, :=\,
\mf p(\mss m_{3,\mtt i})  \,+\, \frac{1}{2} \, \frac{1}{(\ell_1\ell_2)^d} \,
\mf p'' (\mss m_{3,\mtt i})\, \chi (\mss m_{3,\mtt i})
\, :=\,
\mf p (\mss m_{3,\mtt i})  \,-\, \frac{1}{(\ell_1\ell_2)^d} \,
\cb{\mtt q_{\mf p}} (\mss m_{3,\mtt i}) \,.
\end{gathered}
\end{equation*}

Recall that
$\sum_{\mtt k\in\Lambda_{\mtt j}}J(B_{\mtt k}) = \ell_2^d \, J(B_{2, \mtt
j})$, and that $\ell_3 = \kappa\, n^{1/10 + 2\mtt e}$ for some
$\kappa\ge 1$
so that $\ell_1\ell_2\ell_3 = \kappa \,n^{1/2}$. Assume that
$\epsilon <1$. A direct calculation shows that there exists an integer
$n_2 = n_2(\kappa, \epsilon, \delta', \mf r, \Vert \rho\Vert_{C^3},
\mf p)$ such that for all $n\ge n_2$, the conditions of Theorem
\ref{s26} with $\mtt t_1 = \ell_1\ell_2$, $\mtt t_2 =\ell_3$ are
satisfied with $K_1=\kappa^ 6$, $K_2=1$.  Thus, applying Theorem
\ref{s26} with these choices and $\delta = \delta'/8$, we obtain a
finite constant $\mf a_{3,3}$ such that
\begin{align*}
& \int
\sum_{\mtt i \in \Sigma^+_{L_4}}
\sum_{\mtt j \in \Lambda_{3,\mtt i}}  (\ell_1\ell_2)^d\, J(B_{2,\mtt j}) \, \Big\{\,
\mtt p_h(\mss m_{2,\mtt j})
\,-\, E_{\nu^n_{\rho(\cdot)} } \big[\,
\mtt p_h(\mss m_{2,\mtt j}) \,\big]
\,-\, 
\big(\, W^{(3)}_{\mtt p_h,\mtt i}  
- E_{\nu^n_{\rho(\cdot)} } [\,  W^{(3)}_{\mtt p_h,\mtt i}  \, ] \,\big)
\, \Big\}\,
f\, d\nu^n_{\rho(\cdot)}
\nonumber
\\
&  \le\, \frac{\delta'}{8} \, n^2 \,  I_n  ( f ; \nu^n_{\rho (\cdot)} )
\,+\, \mf a_{3,3}\, \Vert J\Vert_{C^1}  \, H_n (f\,|\, \nu^n_{\rho(\cdot)}) \,+\,
C(\mf p)\, e^{(\mf c_{\rm LS}/\delta')\kappa^6\,
\, [1+\mf a_{3,3}]\,[1 + \Vert J\Vert_{C^1}]}  \, \kappa^{1/2}\,
n^{5/4} 
\end{align*}
for all $n\ge n_2$.

Each time Theorem \ref{s01} or \ref{s26} is applied, a correction
term is produced. As we apply Theorem \ref{s01} once and Theorem
\ref{s26} twice there are at the end $8$ terms arising from
$h(\eta)$. We denote these terms by $T_{abc}$ where $a$, $b$,
$c\in \{0,1\}$. The index $a$ refers to Theorem \ref{s01}, and the
indices $b$, $c$ to the first and second application of Theorem
\ref{s26}, respectively.  For the main terms, we set the index $a$,
$b$ or $c$ to be $0$, and $1$ for the correction terms. Hence, the
term $T_{000}$ is the one obtained at the end of the three steps
described above, choosing at each step the main contribution, so that
\begin{equation*}
T_{000} \,=\, \int
\sum_{\mtt i \in \Sigma^+_{L_4}}
(\ell_1\ell_2\ell_3)^d\, J(B_{3,\mtt i}) \, \Big\{\,
\mtt p_h(\mss m_{3,\mtt i})
\,-\, E_{\nu^n_{\rho(\cdot)} } \big[\,
\mtt p_h(\mss m_{3,\mtt i}) \,\big]  
\, \Big\}\,
f\, d\nu^n_{\rho(\cdot)}
\end{equation*}
On the other hand, the term $T_{100}$ corresponds to the correction
produced in Theorem \ref{s01} and contains expressions of the form
$\ell_1^{-d} \mtt q_h(\mss m_{3,\mtt i})$.

With the exception of $T_{100}$ and $T_{010}$, examined at the end of
this proof, the correction terms can be estimated with Lemma
\ref{s45}. We present the bound for $T_{001}$, the other ones being of
smaller order and easier to estimate. $T_{001}$ is given by
\begin{equation*}
-\, \int \sum_{\mtt i \in \Sigma^+_{L_4}}
\ell_3^d\, J(B_{3,\mtt i}) \, \big\{\,
\mtt q_h(\mss m_{3,\mtt i})
\,-\, E_{\nu^n_{\rho(\cdot)} } \big[\,
\mtt q_h(\mss m_{3,\mtt i}) \,\big] \,\big\}
\, f\, d\nu^n_{\rho(\cdot)}\,.
\end{equation*}
By Lemma \ref{s45} with $r_1 = \ell_1 \ell_2 \ell_3$,
$\mf a(\mtt i) = \ell_3^d\, J(B_{3,\mtt i})$, and
$\mf c_0=\Vert J \Vert_\infty$ [note that
$\Vert \mf a\Vert_\infty \le \Vert J \Vert_\infty \, \ell_3^d \le
\Vert J \Vert_\infty\, (\ell_1\ell_2\ell_3)^{d/2}$ for all
$n\ge n_0(\kappa)$], this expression is less than or equal to
\begin{equation*}
C(h)\, \Vert J\Vert_\infty \,\Big\{   \, H_n (f\,|\, \nu^n_{\rho(\cdot)}) +
\Big(\frac{\ell_3}{\ell_1\ell_2} \Big)^{d/2}\,
\frac{n^{d/2}}{\kappa^d} \, \Big\} 
\end{equation*}
for some finite constant $C(h)$ depending only on $h$. By the previous
choice of $\ell_1$, $\ell_2$, $\ell_3$, the second term inside braces
is bounded by $n^{(21/20)+6\mtt e} \le n^{(3/2) - 2\mtt e}$ because
$\mtt e<1/30$, $\kappa\ge 1$.

We turn to $T_{100}$ and $T_{010}$. The first one appears as we apply
Theorem \ref{s01} and is given by
\begin{equation*}
\frac{1}{|B_{\mtt 0} |}\, 
\sum_{y \in B_{\mtt 0}}  \sum_{\mtt i\in \Sigma^+_{L_4}}  \int
(\ell_2\ell_3)^d \,
J(B_{3, \mtt i, y} ) \, \big\{ \, \mtt q_h(\mss m_{3,\mtt i, y} ) 
\,-\,  E_{\nu^n_{\rho(\cdot)}  } \big[\,
\mtt q_h(\mss m_{3,\mtt i, y} ) \,\big] 
\, \big\}\,
f\, d\nu^n_{\rho(\cdot)}\,.
\end{equation*}
The second one, $T_{010}$, comes from Theorem \ref{s26}. By
\eqref{111} it is equal to the sum appearing in the previous displayed
equation but with a minus sign.  Hence $ T_{100} +T_{010} =0$ and we
do not need to estimate these expressions.  Recollecting all previous
estimates completes the proof of the lemma.
\end{proof}

\begin{remark}
\label{r06}
Along the proof of the replacement of a cylinder function $h$ by its
canonical expectation on a cube of width $n^{1/2}$, presented in
Sections \ref{sec4} and \ref{sec9}, two identical correction terms
appeared. The first one in Proposition \ref{s36} is due to the
equivalence of ensembles. As the width of the cube is small, we need
to take into consideration the first order term in the expansion of
$E_{\nu^{\rm c}_{\mtt k, \mss M}} [h]$ as a function of the density
$\mss M/|B_{\mtt k}|$. More precisely, the local central limit theorem
yields that
\begin{equation*}
E_{\nu^{\rm c}_{\mtt k, \mss M}} [h] \,=\,
\mtt p_h(\mss M/|B_{\mtt k}|)
\,+\, \frac{1}{|B_{\mtt k}|} \, \mtt q_h(\mss M/|B_{\mtt k}|)
\,+\, O\Big( \frac{1}{|B_{\mtt k}|^2} \Big) \,\cdot
\end{equation*}
As $|B_{\mtt k}| = \ell_1^d$ is not large enough in $d=3$, we need to
carry the polynomial $\mtt q_h(\mss M/|B_{\mtt k}|)$.

The second correction term, though identical to the first one up to
the sign, has a completely different origin. In Section \ref{sec9}, at
Proposition \ref{s44}, when we replace
$E_{\nu^{\rm c}_{2, \mtt j, \mss M}} [\mtt p_h (\mss m_{\mtt k})] $ by
$\mtt p_h (\mss m_{2, \mtt j})$, we have to exchange
$E_{\nu^n_\rho} [\mtt p_h (\mss m_{\mtt k})]$ by
$\mtt p_h ( E_{\nu^n_\rho} [\mss m_{\mtt k}])$. This swap is based on
Corollary \ref{s30b}, which relies on Lemma \ref{s22}. In that lemma
it is shown that
\begin{equation*}
E_{\nu^n_\rho(\cdot)} [\mtt p_h (\mss m_{\mtt k})] \,=\,
\mtt p_h ( E_{\nu^n_\rho(\cdot)} [\mss m_{\mtt k}]) \,-\,
\frac{1}{|B_{\mtt k}|} \,  \mtt q_h ( E_{\nu^n_\rho(\cdot)} [\mss m_{\mtt k}])
\,+\, O\Big( \frac{1}{|B_{\mtt k}|^2} \Big) \,\cdot 
\end{equation*}

As noted in Remark \ref{r09}, it is fortuitous that both terms are the
same with opposite signs, and therefore cancel. Without this
cancellation, we would have on the right-hand side of Lemma \ref{s49}
an extra term of the form
\begin{equation*}
\frac{1}{|B_{\mtt 0} |}\, 
\sum_{y \in B_{\mtt 0}}  \sum_{\mtt i\in \Sigma^+_{L_4}}  \int
(\ell_2\ell_3)^d \,
J(B_{3, \mtt i, y} ) \, \big\{ \, \mtt q_h(\mss m_{3,\mtt i, y} ) 
\,-\,  E_{\nu^n_{\rho(\cdot)}  } \big[\,
\mtt q_h(\mss m_{3,\mtt i, y} ) \,\big] 
\, \big\}\,
f\, d\nu^n_{\rho(\cdot)}\,.
\end{equation*}

This cancellation is however not crucial in the proof.  It is possible
to complete the proof of Proposition \ref{s88} with this extra term.
\end{remark}

Proposition \ref{s88} is obtained combining Lemma \ref{s49} and
\eqref{135}.

\begin{proof}[Proof of Proposition \ref{s88}]
Fix $0\le s\le t\le T$, $\kappa \ge 1$, $0<\mtt e<1/30$, and $\delta_0>0$.
By Lemma \ref{s49} with $J(\cdot) = G(s, \cdot)$,
$\delta' = \delta_0$, and \eqref{135}, the expression appearing in the
statement of the proposition is bounded by
\begin{align*}
& \frac{1}{|B_{\mtt 0} |}\, 
\sum_{y \in B_{\mtt 0}}  \int\,  \ell^d\, \sum_{\mtt i\in \Sigma^+_{L_4}}  
G(s, B_{3, \mtt i, y} )\,
\mtt E( \mtt p_{h} ,\mss m_{3, \mtt i,y} ) \,  f\, d\nu^n_{u(s, \cdot)}
\\
& \quad
\,+\, \delta_0 \, n^2\, I_n(f ; \nu^n_{u(s, \cdot)}  )
\,+\, \mf c_{1,3} \, \Vert G(s)\Vert_{C^2}\, H_n (f\,|\, \nu^n_{u(s, \cdot)}  )
\,+\,  \mtt K_5 \, n^{3/2 - 2\mtt e}
\end{align*}
for all integers $n \ge n_0$, and densities $f$ with respect to
$\nu^n_{u(s, \cdot)}$. In this formula,
$\mtt K_5 = \mf c_{1,3} (\delta_0, h) \, \Vert G(s) \Vert_{C^1} \,+\,
\mtt K_4(\kappa, G(s)) + \mf c_{1,2} (h) \, \Vert G(s) \Vert_{C^2}
\kappa^4$, where $\mtt K_4(\kappa, G(s))$ is defined in the statement
of the previous lemma, $\ell = \ell_1 \ell_2 \ell_3 = \kappa n^{1/2}$,
and, for a polynomial $\mf p$,
\begin{equation*}
\mtt E( \mf p ,\mss m )  \, =\,
\mf p(\mss m) \,-\, E_{\nu^n_{u(s, \cdot)}} \big[\, \mf p(\mss m)  \,\big]
\,-\, \mf p' \big(\, E_{\nu^n_{u(s, \cdot)}} [\, \mss m \,]
\,\big) \,
\big \{\, \mss m - E_{\nu^n_{u(s, \cdot)}} [\, \mss m \,] \, \big\} \,.
\end{equation*}

By \eqref{120} with $r_1=\ell$, the expression in the penultimate
equation is less than or equal to
\begin{align*}
& \frac{1}{|B_{\mtt 0} |}\, 
\sum_{y \in B_{\mtt 0}}
\int\,  \ell^d\, \sum_{\mtt i\in \Sigma^+_{L_4}}  
G(s, B_{3, \mtt i,y} )\, 
\mtt E_c( \mtt p_{h} ,\mss m_{3, \mtt i,y} ) \, 
f\, d\nu^n_{u(s, \cdot)}
\\
& \quad
\,+\, \delta_0 \, n^2\, I_n(f ; \nu^n_{u(s, \cdot)}  )
\,+\, \mf c_{1,3} \, \Vert G(s)\Vert_{C^2}\, H_n (f\,|\, \nu^n_{u(s, \cdot)}  )
\,+\,  \mtt K_5 \, n^{3/2 - 2\mtt e}
\,+\, C_0\, \Vert G(s) \Vert_\infty\,  \frac{n^d}{\ell^d} \,,
\end{align*}
where $\mtt E_c( \mf p ,\mss m)$ is defined as
$\mtt E( \mf p ,\mss m)$ replacing
$E_{\nu^n_{u(s, \cdot)}} [\, \mf p(\mss m) \,]$ by
$\mf p( E_{\nu^n_{u(s, \cdot)}} [\, \mss m \,])$. Since
$\ell = \kappa n^{1/2}$, the last term can be written as
$n^{3/2}/\kappa^3$.  Note that $\mtt E_c( \mf p ,\mss m_{3, \mtt i} )$
is a second-order Taylor expansion. Hence, by Lemma \ref{l02} with
$r_1=\ell$, which yields an extra term
$C_0\, \Vert G(s) \Vert_\infty\, [H_n (f\,|\, \nu^n_{u(s, \cdot)} ) +
(n/\ell)^d]$, the previous expression is bounded by
\begin{align*}
\delta_0 \, n^2\, I_n(f ; \nu^n_{u(s, \cdot)}  )
\,+\, \mf c_{1,3} \, \Vert G(s)\Vert_{C^2}\, H_n (f\,|\, \nu^n_{u(s, \cdot)}  )
\,+\,  \mtt K_5 \, n^{3/2 - 2\mtt e}
\,+\, C_0\, \Vert G(s) \Vert_\infty\,  
\frac{n^{3/2}}{\kappa^3}  \,.
\end{align*}
This completes the proof of the proposition.
\end{proof}

\section{The Boltzmann--Gibbs principle}
\label{sec7}

In this section, we present a quantitative version of the
Boltzmann--Gibbs principle stated in Theorem \ref{t02}. This is the
last ingredient in the proof of Theorem \ref{t1}. We adopt the
convention established at the beginning of Section \ref{sec2} on the
constants $\mf r$, $\mf c_{1,p}$, $p\ge 0$. The constants $C_1$
depend only on the dimension $d$, on the time-horizon $T$, and on the
cylinder functions $h_j$, $g_j$, $1\le j\le d$, and may change from
line to line. Recall the definition of $\Pi_{\rho, g}$ in
\eqref{179}. Recall the definition of the constants $\mtt
K_{i,j}(\cdot)$ introduced in the statement of Proposition \ref{s88},
and the one of $\Vert G\Vert_{T,p} $ defined in \eqref{256}.

\begin{theorem}[Boltzmann--Gibbs principle]
\label{t02b}
Assume the hypotheses of Theorem \ref{t02}. Fix $T>0$, a function $G$
in $C^{0,2}(\bb R_+ \times \bb T^d)$, and a cylinder function
$h\colon \{0,1\}^{\bb Z^d} \to \bb R$. In dimension $1$ fix also
$0<b<1/6$, and in dimension $3$ a constant $\kappa>1$.  Let $\mu_n$ be
a sequence of probability measures on $\Omega_n$ such that
$H_n(\mu_n \,|\, \nu^n_{u_0(\cdot)} ) =o(n^{d/2}) $. Then, there exist
a finite constant $\mf c_{1,4} = \mf c_{1,4} (h)$, and an integer
$n_0 = n_0 (\mf r, h, \Vert u\Vert_{T,4})$ such that
\begin{equation*}
\bb E_{\mu_n} \Big[\,\Big| \int_s^t
\frac{1}{n^{d/2}} \sum_{x\in \bb T^d_n} G (r, x/n)\,
\Pi_{u(r, \cdot), h}(x, \eta^n(r) ) \, dr \, \Big| \, \Big] \;\le \;
\mtt A^{(d)}_n (s) \,+\, \int_s^t \mtt B^{(d)}_n(r)\, dr
\end{equation*}
for all $0\le s\le t\le T$, and $n\ge n_0$.  In this formula,
\begin{gather*}
\mtt A^{(d)}_n (s) \, = \, \frac{d}{n^{d/2}}\,  \Vert G\Vert_{T,2}\, 
\big\{ \, H_n(f_s \,|\, \nu^n_{u(s, \cdot)}) \,+\, \ln (2)\,\big\}\,,
\end{gather*}
\begin{align*}
\mtt B^{(d)}_n(r) \,= \, \mf c_{1,4}\, \Vert G\Vert_{T,2}\,
\frac{1}{n^{d/2}}
\Big\{ n^2 \, I_n(f_r ;\, \nu^n_{u(r,\cdot)}) \, +\,
H_n(f_r\,|\,\nu^n_{u(r, \cdot)}) \,+\, \varkappa^{(d)}_n \, \Big\}\,,
\end{align*}
where $\varkappa^{(1)}_n = n^{(1/3) + b}$, $\varkappa^{(2)}_n =
n^{4/5}$, $\varkappa^{(3)}_n = n^{3/2}/\kappa^3$.
In dimension $1$, the constants $\mf c_{1,4}$, $n_0$ may depend on $b$ as
well. In dimension $3$ the constant $\mf c_{1,4}$ does not depend on
$\kappa$, and $n_0$ may depend on $\kappa$.
\end{theorem}

In contrast with Proposition \ref{s88}, $n_0$ depends only on the
solution of the hydrodynamic equation and on the cylinder function
$h$, but not on the test function $G$.  We obtain this improvement by
replacing the test function $G$ by
$\widehat{G} = G/\Vert G\Vert_{T,2}$. This latter function has
$\Vert \cdot \Vert_{T,2}$ norm bounded by $1$. In consequence, the
dependence of the constants on the test function disappears.  In
particular, in the case where $h=h_j$ or $g_j$, $1\le j\le d$, $n_0$
depends only on the parameters of the model, and the bound can be
extended to all $n\ge 1$, as explained next.

Denote by $n_\kappa$ the integer $n_0$ starting from which the
inequality holds for all functions $h_j$, $g_j$. The sequence is
clearly increasing. Let $\cb{\mtt i_n}:= 1/\kappa^3$ for
$n_\kappa \le n < n_{\kappa+1}$. With this definition, in dimension
$d=3$ we may replace on the right-hand side $n^{3/2}/\kappa^3$ by
$n^{3/2}\, \mtt i_n$, where $\mtt i_n\to 0$ as $n\to\infty$.

More precisely, let $\cb{\kappa_n} := \varkappa^{(d)}_n$ in dimension
$1$ and $2$, and $\cb{\kappa_n }:= n^{3/2}\, \mtt i_n$ in dimension
$3$. Then, the inequality holds for $h=h_j$, $g_j$, with $\kappa_n$
replacing $\varkappa^{(d)}_n$, and all
$n\ge n_0 = n_0(\mf r, h, \Vert u\Vert_{T,4})$. As the left-hand side
is bounded by $\mf c_{1,0}\, T\, \Vert G\Vert_{T,2}\, n^{d/2}$, by
modifying the constant $\mf c_{1,4}$, the inequality holds for all
$n\ge 1$. We summarize in the next result these conclusions, where we
set $b=1/10$.

\begin{corollary}
\label{s94}
Assume the hypotheses of Theorem \ref{t02}.  Fix $T>0$.  Then, there
exists a finite constant $\mf c_{1,4} $ such that
\begin{equation*}
\bb E_{\mu_n} \Big[\,\Big| \int_s^t
\frac{1}{n^{d/2}} \sum_{x\in \bb T^d_n} G (r, x/n)\,
\Pi_{u(r, \cdot), h}(x, \eta^n(r) ) \, dr \, \Big| \, \Big] \;\le \;
\mf c_{1,4} (h_j) \, (1+T)\, \Vert G\Vert_{T, 2}\,
\big\{\, 1 + (t-s) \,\big\} \, \mtt h_n
\end{equation*}
for all functions $G$ in $C^{0,2}(\bb R_+ \times \bb T^d)$,
$0\le s\le t\le T$, cylinder function $h=h_j$, $1\le j\le d$, and
$n\ge 1$.  In this formula,
\begin{equation}
\label{225}
\cb{\mtt h_n} \;:=\;
\frac{1}{n^{d/2}}\, \Big\{\,
\sup_{0\le r\le T} H_n\big(f_r \,\big|\,
\nu^n_{u(r, \cdot)}\big) \;+\; \int_0^T n^2\,
I_n\big(f_r ;  \nu^n_{u(r, \cdot)}\big) \, dr \, +\kappa_n \,\Big\} \,.
\end{equation}
\end{corollary}

Theorem \ref{t02} follows from the entropy estimate stated in Theorem
\ref{t01} and from this result letting, in dimension $3$,
$\kappa\to\infty$ after $n\to\infty$.

The result asserts that the local field associated with a cylinder
function $h$ is, on the fluctuation scale and after integration in
time, entirely determined by its projection on the density field. It
is the identity which closes the equation satisfied by $X^n$.
\smallskip

The expression to be estimated carries an absolute value inside the
expectation. To remove it, we appeal to the entropy inequality and to
the inhomogeneous-in-time Feynman--Kac formula, \cite[Lemma
A.2]{jm2}, which yield that, for every $V\colon [0,T] \times \Omega_n
\to \bb R$ and every $\gamma>0$,
\begin{equation}
\label{eq:FK}
\begin{aligned}
&\bb E_{\mu_n}\Big[ \, \Big|\, \int_0^t
V (s,\eta^n(s))  \, ds\, \Big| \, \Big]
\, \le \, \frac{1}{\gamma}\, \big\{ H_n(\mu_n \,|\,
\nu^n_{u_0(\cdot)} ) + \ln (2)\, \big\}
\\
&\quad 
\, +\, \max_{\mtt b=\pm 1} \int_0^t 
\sup_f \Big\{\, \int  \Big[ \, 
\mtt b\, V(s) + \frac{1}{2\, \gamma} \, \big(L^*_{n,s} {\bf  1}
-\partial_s \ln \psi^n_s\big)\, \Big]
\, f\, d\nu^n_{u(s,\cdot)} - \frac{n^2}{\gamma}
\, I_n(f ; \nu^n_{u(s,\cdot)}) \,
\Big\} \; ds \,,
\end{aligned}
\end{equation}
the supremum being carried out over all densities $f$ with respect to
$\nu^n_{u(s,\cdot)}$.

This creates two difficulties. The first one is that the supremum is
taken over \emph{all} densities $f$, and not over the density $f_s$ of
the state of the process. The entropy bounds provided by Theorem
\ref{t01}, which were the main tool of Sections
\ref{sec4}--\ref{sec1}, are therefore no longer available: the only
quantity left to control the expression between square brackets is the
Dirichlet form $I_n(f ; \nu^n_{u(s,\cdot)})$. Returning to the
previous sections, only Lemma \ref{s03} and Proposition \ref{s37}
provide bounds of this type. They allow us to replace a local function
$h$ by its expectation with respect to the canonical measure on a
cube, and a monomial $\mss m_{\mtt k}^p$ by its expectation with
respect to the canonical measure on a larger cube.

The second difficulty is that the Feynman--Kac formula introduces in
the supremum the term $L^*_{n,s} {\bf 1} - \partial_s \ln \psi^n_s$,
which has to be estimated as well, and which is not of the restricted
form just described.

Both difficulties are resolved similarly. The expression
$V(s,\eta) = n^{-d/2} \sum_x G(s,x/n) \Pi_{u(s,\cdot), h}(x,\eta)$
splits into two pieces. The first one can be estimated in absolute
value, by the results of Section \ref{sec3} and by the multi-scale
analysis of Sections \ref{sec4} and \ref{sec9}; the second one is
precisely of the form handled by Lemma \ref{s03} and Proposition
\ref{s37}, and is therefore amenable to \eqref{eq:FK}.  Accordingly,
we estimate the absolute value of the first piece directly, and apply
the Feynman--Kac formula only to the second.

The same dichotomy takes care of the term
$L^*_{n,s} {\bf 1} - \partial_s \ln \psi^n_s$, but a further device is
needed because this term already sits inside the supremum. By Lemma
\ref{s41}, up to negligible errors, it is a sum of expressions
$\sum_x F(s,x/n)\, \Pi_{u(s,\cdot), g}(x,\eta)$, and it therefore
splits in exactly the same way. Denote by $\Psi_n(s)$ the part of
$L^*_{n,s} {\bf 1} - \partial_s \ln \psi^n_s$ which cannot be
estimated by the Dirichlet form alone, and write
\begin{equation*}
\bb E_{\mu_n}\Big[ \, \Big|\, \int_0^t
V (s,\eta^n(s))  \, ds\, \Big| \, \Big] \;=\;
\bb E_{\mu_n}\Big[ \, \Big|\, \int_0^t
V (s,\eta^n(s))  \, ds\, \Big| \, - \int_0^t \Psi_n (s) \, ds\,\, \Big]
\,+\,
\bb E_{\mu_n}\Big[ \, \int_0^t \Psi_n (s) \, ds\, \, \Big]\,.
\end{equation*}
Applying the entropy inequality and \eqref{eq:FK} to the first
expectation carries $-\Psi_n(s)$ into the supremum, where it cancels
exactly the ``bad'' part of $L^*_{n,s} {\bf 1} - \partial_s \ln
\psi^n_s$, leaving inside the supremum only terms of the admissible
form. The second expectation comes without an absolute value and is
handled directly by the results of the previous sections.

\subsection*{Organisation of the section}

Subsection \ref{ssec9.1} contains the estimates which are valid in all
dimensions $d\le 3$. It decomposes $\Pi_{u(s,\cdot),\eta_A}$ into four
pieces along the lines just described, estimates three of them, and
reduces the proof of Theorem \ref{t02b} to a single expression,
Conclusion \ref{s87}. From that point on the argument is dimension
dependent, because it requires cubes of width $n^{1/2}$ and the width
reached at this stage is $\ell_1$. In dimension $1$ the cubes are
already large enough and the proof is completed in Subsection
\ref{ssec9.2}. In dimensions $2$ and $3$ the multi-scale analysis of
Section \ref{sec9} has to be invoked: Subsection \ref{ssec9.3}
performs one further step, of width $\ell_1\ell_2$, and Subsection
\ref{ssec9.4} concludes the proof in dimension $2$. In dimension $3$ a
third scale, of width $\ell_1\ell_2\ell_3 = \kappa\, n^{1/2}$, is
needed; this is carried out in Subsection \ref{ssec9.5}. Throughout,
we focus on cylinder functions of the form $\eta_A$, $A$ a finite
subset of $\bb Z^d$, since every cylinder function is a finite linear
combination of such terms, see \eqref{13}.

\subsection{First estimates}
\label{ssec9.1}

We start the proof with some estimates which hold in all dimensions.
Fix $T>0$, a test function $G$ in $C^{0,2}(\bb R_+ \times \bb T^d)$,
and a cylinder function $h$.  We focus on cylinder functions of the
form $\eta_A$, $A$ a finite subset of $\bb Z^d$, given that $h$ is a
finite linear combination of such terms.

Recall the multi-scale setup in Section \ref{sec4}, for $d\leq 3$,
where the space is divided into $L^d_2$ cubes of width $\ell_1$ so
that $n=\ell_1L_2$.  In dimensions $1$ and $2$, the choice of $\ell_1$
is the same as taken in the proof of Theorem \ref{t01}:
$\ell_1=n^{(2/3) - b}$, for some $0<b<1/6$, in $d=1$, and
$\ell_1=n^{2/5}$ in $d=2$.  In $d=3$, we specify differently.  Let
$\upsilon_n$ be a sequence which decreases slowly to $0$, say
$\upsilon^{-1}_n = n^{3/100}$, and set $\ell_1=\upsilon_n \, n^{1/4}$.
This condition is often unnecessary, and weaker ones would
suffice. But in Lemma \ref{s81c} we apply Proposition \ref{s37} and
produce the error $\mtt R^{( \ref{s37} ), 3}_n$. For this error to be
small, we need $\ell^{d+2}_3\, (n/\ell_1\ell_2)^{d-2} =
o(n^{d/2})$. By our choice of widths, this term is of order
$\kappa^5\, \upsilon^{-12}_n\, n^{11/10}$. For this expression to be
$o(n^{3/2})$, we need $\upsilon^{-1}_n = o (n^{1/30})$.

Let $\cb{\mf a_0 = \mf r/2}$. Recall from the beginning of Section
\ref{sec1} the definition of $\mss m_{y,\mtt k}$, $y\in B_{\mtt 0}$,
$\mtt k\in \Sigma^+_{L_2}$, and write
\begin{align}
\label{203}
\frac{1}{n^{d/2}} \sum_{x\in \bb T^d_n} G (s, x/n)\,
\Pi_{u(s, \cdot), \eta_A} (x, \eta)
\,=\,
\frac{1}{n^{d/2}} \frac{1}{|B_{\mtt 0}^o|} 
\sum_{y\in B_{\mtt 0} } \sum_{\mtt k\in \Sigma^+_{L_2}} \,
\sum_{i=1}^4 \Pi^{(i)}_{y,\mtt k, s}(\eta) \,,
\end{align}
where $\Pi^{(1)}_{y,\mtt k, s}$ replaces the linear part of $\Pi_{u(s,
\cdot), \eta_A} (x, \eta)$ by an average over $y+B_{\mtt k} $:
\begin{gather*}
\cb{\Pi^{(1)}_{y,\mtt k, s}(\eta)} \, :=\, \Lambda_{y,\mtt k, s}(\eta)
-\, \sum_{x\in y+B^o_{\mtt k} }
G (s, x/n)\, \mtt p'_A (u(s, x/n)) \, [\eta_x - u(s, x/n)]\,,
\\
\,  \cb{\Lambda_{y,\mtt k, s}(\eta)} \,:= \,  |B_{\mtt k}^o| \,
G(s, y+ B_{\mtt k} ) \, \mtt p_A'( u(s, y+ B_{\mtt k}) )  \,
[\mss m_{y,\mtt k}  - u(s, y+B_{\mtt k})] \,.
\end{gather*}

To introduce the other terms of the decomposition, denote by
$\cb{\nu^{\rm c}_{y, \mtt k, \mss M}}$ the canonical measure on the
cube $y + B_{\mtt k}$ with $\mss M$ particles.  Recall from Lemma
\ref{s33} that
$\mtt 1_{\delta} (\mss m_{\mtt k}) = \mtt 1 \{\mss m_{\mtt k} \in
(\delta , 1-\delta)\} $, and let
$\cb{\mtt 1^{\rm c}_{\delta} (\mss m_{\mtt k}) }:= \mtt 1 \{\mss
m_{\mtt k} \not \in (\delta , 1-\delta)\} $. Let
\begin{gather*}
\cb{\Pi^{(2)}_{y,\mtt k, s}(\eta)} \,:=\,
\mtt 1_{\mf a_0} (\mss m_{y, \mtt k}) \,
\sum_{x\in y+B^o_{\mtt k} } G (s, x/n)\,
\big\{ \tau_x \eta_A
\,-\, E_{\nu^{\rm c}_{y, \mtt k, \mss M_{y, \mtt k} }}
[\tau_x \eta_A]\,\big\}  \,,
\\
\cb{\Pi^{(3)}_{y,\mtt k, s}(\eta) } \,: =\,
\mtt 1_{\mf a_0} (\mss m_{y, \mtt k}) \,
\Big\{\, \sum_{x\in y+B^o_{\mtt k} } G (s, x/n)\,
\big\{\, E_{\nu^{\rm c}_{y, \mtt k, \mss M_{y, \mtt k} }} [  \tau_x \eta_A ] 
\,-\, E_{\nu^n_{u(s,\cdot)}} [  \tau_x \eta_A ] \ \big\}
\, - \, \Lambda_{y,\mtt k, s}(\eta)  \, \Big\} \,,
\\
\cb{\Pi^{(4)}_{y,\mtt k, s}(\eta)} \,:=\,
\widehat{\Pi}^{(2)}_{y,\mtt k, s}(\eta) \,+\,
\widehat{\Pi}^{(3)}_{y,\mtt k, s}(\eta) \,,
\end{gather*}
where $\widehat{\Pi}^{(q)}_{y, \mtt k, s}(\eta)$, $q=2,3$, is obtained
from $\Pi^{(q)}_{y,\mtt k, s}(\eta)$ by replacing the indicator
$\mtt 1_{\mf a_0} (\mss m_{y, \mtt k})$ by
$\mtt 1^{\rm c}_{\mf a_0} (\mss m_{y, \mtt k})$.  As $G$, $u$ and $A$
are fixed, we omitted them from the notation. When we consider $y=0$,
we will drop the subscript $y$ from the notation.

In this subsection, we estimate all terms of the previous
decomposition. We leave $\Pi^{(2)}_{y,\mtt k, s}(\eta)$ to the end
as it requires a finer analysis.

We start with $\Pi^{(1)}_{y,\mtt k, s}(\eta)$. By \eqref{135},
\begin{align}
\label{192b}
\bb E_{\mu_n} \Big[\,\Big| \int_0^t
\frac{1}{n^{d/2}} \frac{1}{|B_{\mtt 0}^o|} 
\sum_{y\in B_{\mtt 0} } \sum_{\mtt k\in \Sigma^+_{L_2}}
\Pi^{(1)}_{y,\mtt k, s}(\eta^n(s)) \, ds \, \Big| \, \Big]
\,\le\, \frac{1}{n^{d/2}} \, \int_0^t \mtt S^{(\ref{135})}_{n,\ell_1}
(G_s,\mtt p'_A, u_s)\, ds\,,
\end{align}
where $\mtt S^{(\ref{135})}_{n,\ell_1} (G_s,\mtt p'_A, u_s)$ is defined
in \eqref{135}.

We turn to $\Pi^{(4)}_{y,\mtt k, s}$. Note that
$|\Pi^{(4)}_{y,\mtt k, s} (\eta)| \le C_1 \, \ell_1^d \,
\|G(s)\|_\infty$, where, recall, $C_1$ is a finite constant which
depends only on $A$ and $d$, whose value may change from line to
line. Hence, by Lemma \ref{s07b} with $r_1=\ell_1$,
\begin{equation*}
\bb E_{\mu_n} \Big[\,\Big| \int_0^t
\frac{1}{n^{d/2}} \frac{1}{|B_{\mtt 0}^o|} 
\sum_{y\in B_{\mtt 0} } \sum_{\mtt k\in \Sigma^+_{L_2}}
\Pi^{(4)}_{y,\mtt k, s} (\eta^n(s)) \, ds
\, \Big| \, \Big]
\, \le  \, \mtt R^{(\ref{s07b})}_{n,\ell_1,G} \,,
\end{equation*}
where for a continuous function
$H\colon [0,T]\times \bb T^d_n \to \bb R$,
\begin{align}
\label{202}
\cb{\mtt R^{(\ref{s07b})}_{n,\ell, H}} \,:=\,
\mf c_{1,0} \, \frac{1}{n^{d/2}}\,  
\int_0^t \|H(s)\|_\infty \,
\Big\{ H_n(f_s|\nu^n_{u(s,\cdot)})\,
+ \Big(\frac{n}{\ell}\Big)^d e^{-\ell^d/\mf c_{1,0}}  \,\Big\} \, ds
 \,.
\end{align}
We use the previous bound many times in the argument either
to include the indicator or to remove it, producing each time the
same error.

We turn to $\Pi^{(3)}_{y,\mtt k, s}$. Recall the definition of
$W_{x,A}$, $W_{\mtt k, A}$ introduced just before Proposition
\ref{s36}.  Let
\begin{equation}
\label{195b}
\begin{gathered}
\cb{F_{y,\mtt k, s} (\mss m_{y,\mtt k})}  \,:=\,
\cb{\overline{W}_{y, \mtt k, A, s}}
\,-\, \mtt p_A'(u(s, y+B_{\mtt k}))\, 
[\mss m_{y,\mtt k}  - u(s, y+B_{\mtt k})]\,,
\\
W_{\mtt k, A}\,=\, \mf p_A(\mss m_{\mtt k})\,, \quad\text{where}\;\;
\cb{\mf p_A(\rho) } \, :=\,  \mtt p_A(\rho) + \ell_1^{-d} \mtt q_A
(\rho)\,, 
\\
\cb{\overline{W}_{y, \mtt k, A, s}} \,:=\,
\mf p_A(\mss m_{y, \mtt k}) \,-\, E_{\nu^n_{u(s,\cdot)}} [\mf p_A(\mss
m_{y, \mtt k}) ]\,.
\end{gathered}
\end{equation}
We added the subscripts $s$ and $y$ to $\overline{W}_{\mtt k, A} $
because the expectation is taken here with respect to
$\nu^n_{u(s,\cdot)}$ [instead of $\nu^n_{\rho(\cdot)}$ as in the
displayed formula below \eqref{168}], and the particles' density is
taken over the translated cubes $y+B_{\mtt k}$ instead of $B_{\mtt
k}$. Note that the derivative $\mtt p_A'$ which appears in the first
line concerns the polynomial $\mtt p_A$ and not $\mf p_A$.

\begin{lemma}
\label{s78}
For all $n\ge 1$,
\begin{align*}
& \bb E_{\mu_n} \Big[\,\Big| \int_0^t
\frac{1}{n^{d/2}} \frac{1}{|B_{\mtt 0}^o|} 
\sum_{y\in B_{\mtt 0} } \sum_{\mtt k\in \Sigma^+_{L_2}}
\Pi^{(3)}_{y,\mtt k, s} (\eta^n(s))  \, ds \, \Big| \,
\Big]
\\
& \le\,
\bb E_{\mu_n} \Big[\,\Big| \int_0^t
\frac{1}{n^{d/2}} \frac{1}{|B_{\mtt 0}^o|} 
\sum_{y\in B_{\mtt 0} } \sum_{\mtt k\in \Sigma^+_{L_2}}
|B^o_{\mtt k} |\, G (s, y+B_{\mtt k} )\,
F_{y,\mtt k, s} (\mss m_{y, \mtt k}(s)) \,
\mtt 1_{\mf a_0} \big(\mss m_{y,\mtt k}(s) \big)  \, ds \, \Big| \,
\Big] \, + \,  \mtt R^{(\ref{s78})}_n  \,,
\end{align*}
where
\begin{gather*}
\cb{\mtt R^{(\ref{s78})}_n} \,:=\,
\mf c_{1,3}\, 
\frac{1}{n^{d/2}}\, 
\int_0^t \|G(s) \|_{C^1}\, \Big\{ \,  H_n(f_s|\nu^n_{u(s,\cdot)})\,
\,+ \, \ln (2)\,+\,
\Big[ \frac{1}{\ell^{3d}_1} + \frac{1}{n\, \ell^d_1}
+ \Big(\frac{\ell_1}{n}\Big)^2\Big]\, n^d\, \Big\}\, ds\,.
\end{gather*}
\end{lemma}

\begin{proof}
The proof relies on Proposition \ref{s36} and uses the notation
introduced in that result. The linear terms of
$\Pi^{(3)}_{y,\mtt k, s} (\eta)$ and
$F_{y,\mtt k, s} (\mss m_{y, \mtt k})$ coincide.  We are left to estimate
the non-linear term of $\Pi^{(3)}_{y,\mtt k, s}$, that is,
\begin{align*}
\frac{1}{n^{d/2}} \frac{1}{|B_{\mtt 0}^o|} 
\sum_{y\in B_{\mtt 0} } \sum_{\mtt k\in \Sigma^+_{L_2}}
\mtt 1_{\mf a_0} (\mss m_{y, \mtt k}) \,
& \sum_{x\in y+B^o_{\mtt k} } G (s, x/n)\,
\overline{W}_{y, x, A, s} \,.
\end{align*}
As above, we added the subscripts $s$ and $y$ to $\overline{W}_{x,A}$
because the expectation is taken here with respect to
$\nu^n_{u(s,\cdot)}$, and averages in the term $W_{\mtt k, A}$ are
performed over the cubes $y+B_{\mtt k}$ instead of $B_{\mtt k}$.

This expression is equal to the one appearing in the statement of
Proposition \ref{s36}. We have to estimate the expectation of the
absolute value of the time integral of this quantity. We move the sum
over $y$ out of the expectation and estimate the expectation for each
fixed $y$. As the bounds are uniform in $y$, we set $y=0$.  By
Proposition \ref{s36},
\begin{equation*}
\bb E_{\mu_n} \Big[\,\Big|
\int_0^t \frac{1}{n^{d/2}} \sum_{\mtt k\in \Sigma^+_{L_2}}
\sum_{x\in B^o_{\mtt k}} G (s, x/n)\,
\overline{W}_{x,A,s}  \,
\mtt 1_{\mf a_0} (\mss m_{\mtt k}(s)) \, ds \, \Big| \, \Big]
\,\le \, \mtt R^{(\ref{s78})}_n \,,
\end{equation*}
where $\mtt R^{(\ref{s78})}_n$ is introduced in the statement of the
lemma.  We do not have on the right-hand side the term
$(n/\ell_1)^d\, \exp\{-\ell^d_1/\mf a_1\}$ which appears in
Proposition \ref{s36} due to the presence of the indicator
$\mtt 1_{\mf a_0} (\mss m_{\mtt k})$.
\end{proof}

Recollecting all previous estimates yields the first main result of
this subsection stated below. Recall from \eqref{135} the definition
of the error $\mtt S^{\eqref{135}}_{n,\ell_1} (G_s,\mtt p'_A, u_s)$,
from \eqref{202} the one of $\mtt R^{(\ref{s07b})}_{n,\ell_1, G}$, and from Lemma
\ref{s78} the one of $\mtt R^{(\ref{s78})}_n$.

\begin{proposition}
\label{s79}
Fix a finite subset $A$ of $\bb Z^d$, and a function $G$ in
$C^{0,2}(\bb R_+ \times \bb T^d)$. Let $\mf a_0 = \mf r/2$. Then,
\begin{align*}
& \bb E_{\mu_n} \Big[\,\Big| \int_0^t
\frac{1}{n^{d/2}} \frac{1}{|B_{\mtt 0}^o|}
\sum_{i\in \{1,3,4\}} \sum_{y\in B_{\mtt 0} } 
\sum_{\mtt k\in \Sigma^+_{L_2}}
\Pi^{(i)}_{y,\mtt k, s}  (\eta^n(s)) \, ds \, \Big| \, \Big]
\\
&\quad \le \,
\bb E_{\mu_n} \Big[\,\Big| \int_0^t
\frac{1}{n^{d/2}} \frac{1}{|B_{\mtt 0}^o|}  \sum_{y\in B_{\mtt 0} }
\sum_{\mtt k\in \Sigma^+_{L_2}}
|B^o_{\mtt k} |\, G (s, y+B_{\mtt k} )\,
F_{y, \mtt k, s} (\mss m_{y, \mtt k}(s)) \,
\, ds \, \Big| \,
\Big] \,+\, \mtt R^{(\ref{s79}) }_n \,,
\end{align*}
where $\mtt R^{(\ref{s79}) }_n \,=\,  \mtt R^{(\ref{s79}) }_n
(G,\eta_A)$ is given by
\begin{align*}
\cb{\mtt R^{(\ref{s79}) }_n} \,=\,\frac{1}{n^{d/2}}\,
\int_0^t  \mtt S^{\eqref{135}}_{n,\ell_1}
(G_s,\mtt p'_A, u_s) \, ds \,+\, 2\, \mtt R^{(\ref{s07b})}_{n,\ell_1, G} \,+\,
\mtt R^{(\ref{s78})}_n \,.
\end{align*}
\end{proposition}

Note that we removed the indicator $\mtt 1_{\mf a_0} \big(\mss m_{y,
\mtt k}(s) \big)$ which was multiplying $F_{y, \mtt k, s} (\mss m_{y,
\mtt k}(s))$ at an extra cost $\mtt R^{(\ref{s07b})}_{n,\ell_1, G} $. This
explains the factor $2$ in front of $\mtt R^{(\ref{s07b})}_{n,\ell_1, G} $.

\begin{remark}
\label{r04}
Inspecting each term yields that
\begin{equation*}
\mtt R^{(\ref{s79})}_n \,\le\, \mf c_{1,3}\, \frac{1}{n^{d/2}}
\int_0^t \Vert G(s)\Vert_{C^2}\, \big[\,
H_n(f_s \,|\, \nu^n_{u(s,\cdot)})
+  \mtt e^{(\ref{s79})}_{n,d} \,\big] \; ds  \,,
\end{equation*}
where $\mtt e^{(\ref{s79})}_{n,1} = n^{(1/3)-2b}$ in dimension $1$,
$\mtt e^{(\ref{s79})}_{n,2} = n^{4/5}$ in dimension $2$, and
$\mtt e^{(\ref{s79})}_{n,3} = \upsilon_n^{2} \, n^{3/2}$ in dimension
$3$. Since $\upsilon_n\to 0$, by the a priori entropy bounds, the
remainder $\mtt R^{(\ref{s79})}_n$ vanishes as $n\to\infty$ in
dimension $d\le 3$.
\end{remark}

We now turn to $\Pi^{(2)}_{y, \mtt k, s}$.  Recall from Lemma
\ref{s41} that
\begin{align*}
L^*_{n,s} \mb 1 - \partial_s \ln \psi^n_s
\ = \ \sum_{j=1}^d \sum_{a=1}^3 \sum_{x\in \bb T^d_n} G_{a,j}
(s,x/n)\, \Pi_{u(s,\cdot), f_{a,j}} (x, \eta)
+ I_n(s) + \mf R^{(4)}_{s,n} \, n^{d-2}\,,
\end{align*}
where $f_{a,j}$ are cylinder functions, and $G_{a,j}$ are the
associated coefficient functions. By \eqref{13}, we may rewrite the
previous identity as
\begin{align}
\label{209}
L^*_{n,s} \mb 1 - \partial_s \ln \psi^n_s
\ = \ \sum_{D\in \mc A} \sum_{x\in \bb T^d_n} H_{D}
(s,x/n)\, \Pi_{u(s,\cdot), \eta_D} (x, \eta)
+ I_n(s) + \mf R^{(4)}_{s,n} \, n^{d-2}\,,
\end{align}
where the sum is carried over a finite collection $\mc A$ of finite
subsets $D$ of $\bb Z^d$, corresponding to the supports of
$\{ f_{a,j} \}$.

Denote by $\cb{\Pi^{i,D}_{y, \mtt k, s} (\eta)}$, $1\le i\le 4$, the
associated functions introduced in \eqref{203} [that is, with $\eta_A$
replaced by $\eta_D$ and $G(s, \cdot)$ by $H_{D}(s, \cdot)$]. Let
\begin{equation*}
\cb{\Psi^{D}_n(s)} \, :=\,  \sum_{x\in \bb T^d_n} H_{D}
(s,x/n)\, \Pi_{u(s,\cdot),\eta_D} (x, \eta)
\,  -\,  \frac{1}{|B_{\mtt 0}^o|} 
\sum_{y\in B_{\mtt 0} } \sum_{\mtt k\in \Sigma^+_{L_2}} \,
\Pi^{2,D}_{y, \mtt k, s} (\eta) \,.
\end{equation*}
By \eqref{203},
\begin{equation}
\label{204}
\Psi^{D}_n(s) \, =\,  \sum_{i\in \{1,3,4\}}
\frac{1}{|B_{\mtt 0}^o|} 
\sum_{y\in B_{\mtt 0} } \sum_{\mtt k\in \Sigma^+_{L_2}} \,
\Pi^{i,D}_{y, \mtt k, s} (\eta) \,.
\end{equation}
Define $\Psi_n(s)$ by
\begin{equation*}
\begin{gathered}
\cb{\Psi_n(s)} \,:=\, \sum_{D\in \mc A}
\Psi^{D}_n(s)  + I_n(s) + \mf R^{(4)}_{s,n} \, n^{d-2} \,,
\end{gathered}
\end{equation*}
to write $L^*_{n,s} \mb 1 - \partial_s \ln \psi^n_s$ as
\begin{align}
\label{196}
L^*_{n,s} \mb 1 - \partial_s \ln \psi^n_s \,=\, \Psi_n(s)
\ + \
\frac{1}{|B_{\mtt 0}^o|} \sum_{D\in \mc A}
\sum_{y\in B_{\mtt 0} } \sum_{\mtt k\in \Sigma^+_{L_2}} \,
\Pi^{2,D}_{y, \mtt k, s} (\eta) \,.
\end{align}
The terms $\Psi^{D}_n(\cdot)$ have been treated in Proposition
\ref{s79}. The other ones are small remainders. This observation is
stated precisely in the next two results.  We first estimate the
remaining term $\Pi^{(2)}_{y,\mtt k,s}$ in the decomposition
\eqref{203} of $\Pi_{u(s, \cdot), \eta_A}(x,\eta)$ in terms of
$\Psi_n$, and then we estimate $\Psi_n$.

\begin{lemma}
\label{s85}
Under the hypotheses of Proposition \ref{s79},
\begin{align*}
& \frac{1}{|B_{\mtt 0}^o|} 
\sum_{y\in B_{\mtt 0} } 
\bb E_{\mu_n} \Big[\,\Big| \int_0^t
\frac{1}{n^{d/2}} \sum_{\mtt k\in \Sigma^+_{L_2}}
\Pi^{(2)}_{y,\mtt k, s}(\eta^n(s)) \, ds \, \Big| \, \Big]
\,\le \, \frac{1}{2\, n^{d/2}}  \, \bb E_{\mu_n} \Big[\, \int_0^t
\Psi_n(s)  \, ds \, \Big]
\,+\,   \mtt R^{( \ref{s85} )}_n\,,
\end{align*}
where
\begin{gather*}
\cb{ \mtt R^{( \ref{s85} )}_n } \,:=\,
\frac{1}{n^{d/2}}\,
\big\{ \, H_n(f_0 \,|\, \nu^n_{u_0(\cdot)}) \,+\, \ln (2)\,\big\}\,
+\, \mtt R^{( \ref{s03} )}_n\,,
\\
\cb{\mtt R^{( \ref{s03} )}_n}\, =\, \mf c_{1,0}\,
n^{d/2}\, \left(\frac{\ell_1}{n}\right)^2  \,
\int_0^t  \big\{\, \|G(s)\|_\infty^2
+ \sum_{D\in\mc A} \|H_D(s)\|_\infty^2\,\big\} \, ds\,.
\end{gather*}
\end{lemma}

\begin{proof}
As before, we fix $y=0$ and make sure that all estimates are uniform
over $y$. Let
\begin{align*}
\cb{V_{s} }  \,:= \, \sum_{\mtt k\in \Sigma^+_{L_2}}
\Pi^{(2)}_{\mtt k, s}(\eta^n(s)) \,.
\end{align*}
Recall the definition of $\Psi_n(s)$ introduced just before \eqref{196}.
Write the expectation appearing in the statement of the lemma with
$y=0$ as
\begin{equation*}
\bb E_{\mu_n} \Big[\, \Big| \int_0^t \frac{1}{n^{d/2}} \, V_s \, ds
\, \Big| \, - \, \frac{1}{2}\, \frac{1}{n^{d/2}}
\int_0^t \, \Psi_n(s) \, ds \, \Big]
\, + \, \bb E_{\mu_n} \Big[\, \frac{1}{2}\, \frac{1}{n^{d/2}}
\int_0^t \, \Psi_n(s) \,ds\, \Big]\,.
\end{equation*}
The second term appears on the right-hand side of the statement of the
lemma. We turn to the first one. By the entropy inequality, this
expectation is bounded by
\begin{equation}
\label{197}
\begin{aligned}
\frac{1}{n^{d/2}}H_n(f_0 \,|\, \nu^n_{u_0(\cdot)})
\,+\, \frac{1}{ n^{d/2}}
\ln \bb E_{\nu^n_{u_0(\cdot)}}\Big[\exp\Big\{ \, \Big|\, \int_0^t
V_{s} \, ds\, \Big|
\,-\, \int_0^t \frac{1}{2} \, \Psi_n(s) \, ds \, \Big\}\, \Big] \,.
\end{aligned}
\end{equation}

Since $e^{|x|} \leq  e^x + e^{-x}$, for any random variables $Y$, $Z$,
\begin{equation*}
\ln E[e^{|Z| + Y}] \, \le\, \ln \Big\{  E[e^{Z+Y}]  + E[e^{-Z+Y}] \Big\}
\, \le\, \ln \Big\{ \, 2 \max_{b=\pm 1}  E[e^{b Z+Y}]  \Big\}
\, =\, \ln 2 + \max_{b=\pm  1}  \ln E[e^{bZ+Y}]  \,.
\end{equation*}
Hence, invoking the inhomogeneous-in-time Feynman--Kac formula
\cite[Lemma A.2]{jm2}, the second term in \eqref{197} is bounded
by
\begin{align*}
\frac{\ln(2)}{n^{d/2}} \,+\, \max_{b=\pm 1}
\frac{1}{n^{d/2}} \, \int_0^t 
\sup_f \Big\{\, \int  \Big[ \, 
b V_s - \frac{1}{2} \,  \Psi_n(s) + \frac{1}{2} \, \big(L^*_{n,s} {\bf  1}
-\partial_s \ln \psi^n_s\big)\, \Big]
\, f\, d\nu^n_{u(s,\cdot)} - n^2I_n(f ; \nu^n_{u(s,\cdot)}) \,
\Big\} \; ds \,,
\end{align*}
where the supremum is carried out over all densities $f$ with respect to
$\nu^n_{u(s,\cdot)}$. 

By \eqref{196}, the second term of this sum is equal to
\begin{gather}
\label{210}
\max_{b=\pm 1}
\frac{1}{n^{d/2}} \, \int_0^t 
\sup_f \Big\{\, \int  \mtt X^b_n
\, f\, d\nu^n_{u(s,\cdot)} - n^2I_n(f ;\, \nu^n_{u(s,\cdot)}) \, \Big\} \; ds \,,
\end{gather}
where
\begin{align*}
\mtt X^b_n (\eta) \, & =\, b\,
\sum_{\mtt k\in \Sigma^+_{L_2}} \, \Pi^{(2)}_{\mtt k, s} (\eta)
\,+\, \frac{1}{2}\, 
\frac{1}{|B_{\mtt 0}^o|} \sum_{D\in \mc A}
\sum_{y\in B_{\mtt 0} } \sum_{\mtt k\in \Sigma^+_{L_2}} \,
\Pi^{2,D}_{y, \mtt k, s} (\eta) \,.
\end{align*}
By definition of $\Pi^{(2)}_{\mtt k, s} (\eta)$,
$\Pi^{2,D}_{y, \mtt k, s} (\eta) $, and by Lemma \ref{s03} with
$\delta = 1/(1+|\mc A|)$, the expression \eqref{210} is less than or
equal to $\mtt R^{( \ref{s03} )}_n$. To complete the proof, it remains
to recollect all previous estimates.
\end{proof}

It remains to estimate $\Psi_n(s)$. The advantage is that it comes
without the absolute value, and can therefore be estimated with the
results of the previous sections.

Recall from \eqref{196}  the definition of $\Psi_n(s)$ to
write it as
\begin{align*}
\sum_{j=1}^d \sum_{a=1}^3
\sum_{x\in \bb T^d_n} G_{a,j}
(s,x/n)\, \Pi_{u(s,\cdot), f_{a,j}} (x, \eta)
\,  -\,  \sum_{D\in \mc A} \frac{1}{|B_{\mtt 0}^o|} 
\sum_{y\in B_{\mtt 0} } \sum_{\mtt k\in \Sigma^+_{L_2}} \,
\Pi^{2,D}_{y, \mtt k, s} (\eta) 
+ I_n(s) + \mf R^{(4)}_{s, n} \, n^{d-2} \,.
\end{align*}
The first term has been estimated in Proposition \ref{s88}, the second
one in Lemma \ref{s03}, the third one in \eqref{108}. The last one is
bounded by $\mf c_{1,4} n^{d-2}$. Adding all estimates yields the next
result.

\begin{lemma}
\label{s82}
Under the hypotheses of Proposition \ref{s79},
\begin{align*}
\frac{1}{n^{d/2}}  \bb E_{\mu_n} \Big[\, \int_0^t
\Psi_n(s)  \, ds \, \Big]
\, \le\, \mtt R^{( \ref{s82} )}_{n,d}(\delta)
\end{align*}
for all $\delta>0$. In this formula,
\begin{align*}
\cb{\mtt R^{(\ref{s82})}_{n,d} (\delta)} \,=\, \mtt R^{( \ref{s88} )}_{n,d} \,+\, 
\frac{1}{n^{d/2}}\, \int_0^t \Big\{ \delta\, n^2 \,
I_n(f_s ;\, \nu^n_{u(s,\cdot)}) \,+\, C_1 \,
H_n( f_s \,|\, \nu^n_{u(s, \cdot)} ) \,+\,
\mf c_{1,4}\, \frac{1}{\delta}\, \ell_1^2 \, n^{d-2}\, \Big\}\, ds\,,
\end{align*}
and $\cb{ \mtt R^{( \ref{s88} )}_{n,d}}$ is the sum in $a$ and $j$ of the
expression appearing on the right-hand side of Proposition \ref{s88}
with $G=G_{a, j}$, $h=f_{a,j}$ divided by $n^{d/2}$.
\end{lemma}

\begin{remark}
\label{r08}
A simple inspection provides a bound for the error terms
$\mtt R^{( \ref{s85} )}_n$ and
$\mtt R^{( \ref{s82} )}_{n,d} (1/2)$:
\begin{gather*}
\mtt R^{( \ref{s85} )}_n \, \le\,
\frac{1}{n^{d/2}}\,
\big\{ \, H_n(f_0 \,|\, \nu^n_{u_0(\cdot)}) \,+\, \ln (2)\,\big\}\,
+\, \mf c_{1,4}\,  \frac{1}{n^{d/2}}\, \mtt e^{( \ref{s85} )}_{n,d}\,
\int_0^t \big\{ \,  1+ \Vert G(s) \Vert^2_\infty \big\}\, ds 
\\
\mtt R^{( \ref{s82} )}_{n,d} (1/2) \, \le\, \, \frac{1}{n^{d/2}}\, 
\int_0^t \Big\{ \, n^2I_n(f_s ;\, \nu^n_{u(s,\cdot)}) \, +\,
\mf c_{1,3}\,  \big[\, (1 + \mtt K_{d,1}(s)) \, H_n( f_s \,|\, \nu^n_{u(s, \cdot)} )
\,+\, \mtt K_{d,2} (s) \,\big]\, \Big\} \, ds \,.
\end{gather*}
In this equation, the variables $\mtt K_{d,1} (s)$, $\mtt K_{d,2} (s)$
are defined in Proposition \ref{s88}, and
$\mtt e^{( \ref{s85} )}_{n,1} = n^{(1/3)-2b}$,
$\mtt e^{( \ref{s85} )}_{n,2} = n^{4/5}$,
$\mtt e^{( \ref{s85} )}_{n,3} = \upsilon_n^{2} \, n^{3/2}$.  In
Proposition \ref{s88} we chose $\delta_0 = 1/(6d)$ and in Lemma
\ref{s82} $\delta=1/2$ to end up with
$n^2I_n(f ;\, \nu^n_{u(s,\cdot)})$ multiplied by $1$, a purely
cosmetic decision.
\end{remark}

\begin{conclusion}
\label{s87}
In view of the decomposition \eqref{203}, Proposition \ref{s79} (and
Remark \ref{r04}), Lemmata \ref{s85}, \ref{s82} (and the previous
remark), to prove Theorem \ref{t02b}, it remains to estimate the
expressions involving $F_{y, \mtt k, s}$. The argument is dimension
dependent.
\end{conclusion}

\subsection{ Dimension $d=1$}
\label{ssec9.2}

In this subsection, we complete the proof of the Boltzmann--Gibbs
principle in dimension $1$. In view of the previous subsection, we
need to estimate the term on the right-hand side of Proposition
\ref{s79}.

Recall from Proposition \ref{s36} that $\overline{W}_{\mtt k, A}$ is
the sum of two terms, the second one, denoted below by
$\overline{W}^{(1)}_{y, \mtt k, A, s} $, being of order $\ell_1^{-d}$:
\begin{gather}
\label{207}
\cb{\overline{W}^{(1)}_{y, \mtt k, A, s}} \,:=\,
\frac{1}{ \ell^d_1} \,
\mtt q_A(\mss m_{y, \mtt k})
\,-\, E_{\nu^n_{u(s,\cdot)}} \big[\,
\frac{1}{ \ell^d_1} \,
\mtt q_A(\mss m_{y, \mtt k})\,\big]\,.
\end{gather}

We claim that
\begin{equation}
\label{208}
\begin{aligned}
&\bb E_{\mu_n}  \Big[\,\Big| \int_0^t
\frac{1}{n^{d/2}} \sum_{\mtt k\in \Sigma^+_{L_2}} \,
|B^o_{\mtt k}|\, 
G(s,\mtt k) \, \overline{W}^{(1)}_{\mtt k, A, s} \,
\mtt 1_{\mf a_0} (\mss m_{\mtt k}(s)) \, ds\,\Big|\,\Big ]
\leq\, \frac{1}{\ell^{d}_1}\, \mtt R^{(\ref{s07b})}_{n,\ell_1, G}  \,+\,
\mtt R^{(\ref{s32})}_{n, \ell_1} \,,
\end{aligned}
\end{equation}
where $\mtt R^ {(\ref{s07b})}_{n,\ell_1, G}$ is the remainder defined
in \eqref{202}, and
\begin{align*}
\cb{\mtt R^{(\ref{s32})}_{n,\ell_1} }\, &:=\,
\frac{\mf c_{1,0}}{n^{d/2}} \, 
\, \Big(\frac{n}{\ell_1^2}\Big)^d 
\int_0^t  \|G(s) \|_\infty \; ds
\\
& + \, 
\frac{\mf c_{1,0}}{n^{d/2}}\, 
\int_0^t  \Vert G(s)\Vert_\infty  
\big[\, H_n(f_s \,|\, \nu^n_{u(s,\cdot)})
+ \ln (2)\,\big] \; ds  \,.
\end{align*}

To prove \eqref{208}, apply Lemma \ref{s07b} to get rid of the
indicator.  This lemma produces an error bounded by
$\ell^{-d}_1\, \mtt R^ {(\ref{s07b})}_{n,\ell_1, G}$.  Note that we
have an extra factor $1/\ell_1^d$ due to the presence of this factor
in the definition of $\overline{W}^{(1)}_{\mtt k, A, s}$.  Then, apply
Lemma \ref{s32} with $r_2=L_2$, $\mf a(\mtt k) = G(s, \mtt k)$ and
$\gamma= \gamma(s) = 1/\|G(s) \|_\infty$.  This proves \eqref{208}.

By \eqref{195b} and \eqref{208}, it remains to estimate
\begin{equation}
\label{195c}
\cb{F^o_{y,\mtt k, s} (\mss m_{y,\mtt k})}  \,:=\,  \mss m_{y,\mtt k}^{|A|}
\,-\, E_{\nu^n_{u(s,\cdot)}} \big [ \,\mss m_{y,\mtt k}^{|A|} \,\big]
\,-\, \mtt p_A'(u(s, y+B_{\mtt k}))\, 
[\mss m_{y,\mtt k}  - u(s, y+B_{\mtt k})]\,.
\end{equation}

\begin{lemma}
\label{s80}
Under the hypotheses of Proposition \ref{s79},
\begin{align*}
\frac{1}{|B_{\mtt 0}^o|}  \sum_{y\in B_{\mtt 0} }
\bb E_{\mu_n} \Big[\,\Big| \int_0^t
\frac{1}{n^{d/2}} \sum_{\mtt k\in \Sigma^+_{L_2}}
|B^o_{\mtt k} |\, G (s, y+B_{\mtt k} )\,
F^o_{y, \mtt k, s} (\mss m_{y, \mtt k}(s)) \,
\mtt 1_{\mf a_0} \big(\mss m_{y, \mtt k}(s) \big)  \, ds \, \Big|\,\Big]
\,\le\, \mtt R^{(\ref{s80})}_n\,,
\end{align*}
where
$\cb{\mtt R^{(\ref{s80})}_n} \,=\, \mtt S^{(\ref{120})}_n \,+\, \mtt
R^{ ( \ref{l02} ) }_n$,
\begin{gather*}
\cb{\mtt S^{(\ref{120})}_n} \,:=\, C_1 \, \frac{n^{d/2}}{\ell_1^d}
\, \int_0^t \|G(s)\|_\infty\, ds\,,
\quad
\cb{ \mtt R^{ ( \ref{l02} ) }_n} \,:=\, 
\frac{C_1}{n^{d/2}}
\int_0^t \|G(s) \|_\infty \, \Big\{\, H_n(f_s \,|\, \nu^n_{u(s,\cdot)})
+ \Big(\frac{n}{\ell_1}\Big)^d\, \Big\}\, ds \,.
\end{gather*}
\end{lemma}

\begin{proof}
As before, we set $y=0$ and make sure that the estimates are uniform
over $y$.  By \eqref{120}, in the definition of $F^o_{\mtt k, s}$ we may
replace $E_{\nu^n_{u(s,\cdot)}}\big[\mss m_{\mtt k}^{|A|}\big]$ by
$u(s, \mtt k)^{|A|}$ at a cost $\mtt S^{(\ref{120})}_n$, defined in
the statement of the lemma.

At this point, and since $\mtt p_A(\rho) = \rho^{|A|}$, to complete
the proof of the lemma, it remains to estimate
\begin{align*}
\cb{\mtt p^{(2)}_{\mtt k, s} (\mss m_{\mtt k}) } \,:=\,
\mtt p_A(\mss m_{\mtt k}) - \mtt p_A(u(s,\mtt k))-\mtt
p_A'(u(s,\mtt k))\big[\mss m_{\mtt k} - u(s, \mtt k)\big] \,.
\end{align*}
By the Taylor estimate, the absolute value of this expression is less
than or equal to $C_1 \big[\mss m_{\mtt k} - u(s, \mtt k)\big]^2$.
Thus, by Lemma \ref{l02},
\begin{align*}
\bb E_{\mu_n}\Big[ \,\Big| \int_0^t \frac{1}{n^{d/2}}
\sum_{\mtt k\in \Sigma^+_{L_2}} |B^o_{\mtt k}|\, 
G(s, \mtt k) \,  \mtt p^{(2)}_{\mtt k, s} (\mss m_{\mtt k})
\, ds\,\Big|\,\Big] \,\le\, \mtt R^{ ( \ref{l02} ) }_n \,,
\end{align*}
where $\mtt R^{ ( \ref{l02} ) }_n$ is defined in the statement of the
lemma. This completes the proof.
\end{proof}

It will be helpful, for later tightness estimates in Section
\ref{sec11}, to bound the integral appearing in the Boltzmann--Gibbs
principle from times $s$ to $t$.

\begin{proof}[Proof of Theorem \ref{t02b} in $d=1$]  
Fix $T>0$ and $0\le s<t\le T$.
Recall the definition of the norm
$\Vert G\Vert_{T,2}$ introduced just before the statement of the
theorem and let $\widehat G = G / \Vert G\Vert_{T,2} $, so that
$\sup_{0\le r\le T} \Vert \widehat G (r)\Vert_{C^2} \le 1$. Rewrite
the expectation appearing in the statement of the theorem as
\begin{align}
\label{235}
\Vert G \Vert_{T,2} \;
\bb E_{\mu_n}\Big[\, \Big| \, \int_s^t \frac{1}{n^{d/2}}
\sum_{x\in \bb T^d_n} \widehat G(r, x/n) \, \Pi_{u(r,\cdot),h}(x, \eta^n(r))
\, dr\, \Big|\, \Big] \,.
\end{align}

By Proposition \ref{s79}, Lemmata \ref{s85}, \ref{s82}, the estimate
\eqref{208}, and Lemma \ref{s80},
\begin{align*}
\bb E_{\mu_n}\Big[\, \Big| \, \int_s^t \frac{1}{n^{d/2}}
\sum_{x\in \bb T^d_n} \widehat G(r, x/n) \, \Pi_{u(r,\cdot),h}(x, \eta^n(r))
\, dr\, \Big|\, \Big]
\,\leq\,
\mtt R^{( \ref{s79} )}_n \,+\, \mtt R^{( \ref{s85} )}_n
\,+\, \frac{1}{2}\, \mtt R^{( \ref{s82} )}_{n,d} (1/2)
\,+\, \mtt R^{(\ref{s32})}_{n, \ell_1} \,+\, \mtt R^{( \ref{s80} )}_n \,. 
\end{align*}
In the formulae for the remainders we need to replace the
time-integrals between $0$ and $t$ by integrals between $s$ and $t$,
the entropy $H_n(f_0 \,|\, \nu^n_{u_0(\cdot)})$ by
$H_n(f_s \,|\, \nu^n_{u(s, \cdot)})$, and $G$ by $\widehat G$.  The
term $\ell_1^{-d}\, \mtt R^{(\ref{s07b})}_{n,\ell_1, G}$ does not
appear explicitly on the right-hand side because it is part of the
remainder $\mtt R^{( \ref{s79} )}_n$ (without the factor
$\ell_1^{-d}$).

By Remarks \ref{r04} and \ref{r08} and the explicit formula for
$\mtt R^{( \ref{s80} )}_n$, the right-hand side is bounded by
\begin{align*}
& \frac{1}{n^{1/2}}\,
\big\{ \, H_n(f_s \,|\, \nu^n_{u(s, \cdot)}) \,+\, \ln (2)\,\big\}\,
+\, \mf c_{1,4}\,  \frac{1}{n^{(1/6) - b}}\, (t-s)
\\
\, &+\, \frac{1}{n^{1/2}}\, 
\int_s^t \Big\{ \, n^2I_n(f_r ;\, \nu^n_{u(r,\cdot)}) \, +\,
\mf c_{1,3}\,  \big[\, H_n( f_r \,|\, \nu^n_{u(r, \cdot)} )
\,+\, n^{(1/3)+b} \,\big]\, \Big\} \, dr \,.
\end{align*}
This completes the proof in dimension $1$.
\end{proof}

\subsection{Estimates in $d=2$, $3$}
\label{ssec9.3}

In this subsection we present the estimates used in the proof of
Theorem \ref{t02b} in dimensions $d=2$ and $d=3$.  Recall from
Subsection \ref{ssec9.1} the choices of $\ell_1$, and let
$\ell_2 = n^{1/5}$ in $d=2$, $\ell_2 = \upsilon_n\, n^{3/20}$ in
$d=3$.  Set $n = \ell_1\ell_2 L_3$.

Recall from \eqref{195b} the definition of $F_{y,\mtt k, s}$.  By
Conclusion \ref{s87}, to complete the proof of Theorem \ref{t02b} it
remains to estimate the expression
\begin{equation}
\label{206}
\bb E_{\mu_n} \Big[\, \Big|\, \int_0^t
\frac{1}{n^{d/2}} \frac{1}{|B_{\mtt 0}|} 
\sum_{y\in B_{\mtt 0} } 
\sum_{\mtt k\in \Sigma^+_{L_2}}
|B_{\mtt k} |\, G  (s, y+B_{\mtt k} )\,
F_{y,\mtt k, s} (\mss m_{y, \mtt k}(s))
\, ds  \,\Big|\, \Big]\,.
\end{equation}
Mind that we replaced the factors $|B_{\mtt 0}^o| = |B^o_{\mtt k}|$,
which appeared multiplying $F_{y,\mtt k, s}$ in the statement of
Proposition \ref{s79}, by $|B_{\mtt 0}| = |B_{\mtt k}|$.  

In dimension $d=1$ the expression \eqref{206} was estimated directly,
in Lemma \ref{s80}, because the cubes $B_{\mtt k}$ are already large
enough.  In dimensions $d=2$, $3$ this is no longer the case, and we
have to increase the width of the cubes.  This is achieved with the
results of Section \ref{sec9}, and the argument is a repetition of the
one carried out in Subsection \ref{ssec9.1}, one scale higher.

Recall the notation introduced at the beginning of Section \ref{sec9}
with $\mtt t_1 = \ell_1$, $\mtt t_2 = \ell_2$, $\mtt t_3 = L_3$: the
cubes $B_{\mtt k}$, $\mtt k \in \Sigma^+_{L_2}$, have width $\ell_1$,
the cubes $B_{2, \mtt j}$, $\mtt j \in \Sigma^+_{L_3}$, have width
$\ell_1 \ell_2$, and $\Lambda_{\mtt j}$ is the set of indices $\mtt k$
such that $B_{\mtt k} \subset B_{2, \mtt j}$.  Let
$\mss M_{y, 2, \mtt j}$ stand for the total number of particles in the
cube $y+ B_{2, \mtt j}$, and $\mss m_{y, 2, \mtt j}$ for the density:
\begin{align*}
\cb{\mss M_{y, 2, \mtt j}} := 
\sum_{z\in y + B_{2, \mtt j}} \eta_z\,, \quad
\cb{\mss m_{y, 2, \mtt j}} := \mss m (y+  B_{2, \mtt j})
\,=\, \frac{1}{|B_{2, \mtt j}|} \, \mss M_{y, 2, \mtt j}
\,, \quad
y\in B_{\mtt 0}\,, \; \mtt j \in \Sigma^+_{L_3}\,.
\end{align*}
Denote by $\cb{\nu^{\rm c}_{2, \mtt j, \mss M}}$ the canonical measure
on the cube $y+B_{2, \mtt j}$ with $\mss M$ particles.

The polynomial which appears at this scale is not $\mf p_A$ but the
one produced by Theorem \ref{s26}.  Recall from \eqref{195b} that
$\mf p_A = \mtt p_A + \ell_1^{-d}\, \mtt q_A$, and recall the
cancellation stated in Remark \ref{r09}: with $\mtt t_1 = \ell_1$,
\begin{equation}
\label{212}
\cb{\mf p^{(2)}_A(\rho)} \,:=\, W_{\mf p_A} (\rho) \,=\,
\mtt p_A(\rho) \,+\, \frac{1}{\ell_1^{2d}}\, \mtt q^{(2)}_A(\rho)\,,
\qquad
\cb{\mtt q^{(2)}_A(\rho)} \,:=\, \frac{1}{2}\,
\mtt q_{A}''(\rho)\, \chi(\rho)\,.
\end{equation}
The corrector of order $\ell_1^{-d}$ has disappeared, and the one
which survives is of order $\ell_1^{-2d}$.  As observed in Remark
\ref{r09}, this is the reason why the second scale is inexpensive.

As in Subsection \ref{ssec9.1}, write
\begin{equation}
\label{213}
\frac{1}{n^{d/2}} 
\frac{1}{|B_{\mtt 0}|} 
\sum_{y\in B_{\mtt 0} } \sum_{\mtt k\in \Sigma^+_{L_2}}
|B_{\mtt k} |\, G  (s, y+B_{\mtt k} )\,
F_{y,\mtt k, s} (\mss m_{y,\mtt k}) \,
=\,
\frac{1}{n^{d/2}} \frac{1}{|B_{\mtt 0}|} 
\sum_{y\in B_{\mtt 0} }  \sum_{\mtt j\in \Sigma^+_{L_3}} \,
\sum_{i=1}^4 \Pi^{2,i}_{y, \mtt j, s}(\eta) \,,
\end{equation}
where $\Pi^{2,1}_{y,\mtt j, s}$ replaces the linear part of the
left-hand side by an average over cubes of width $\ell_1\ell_2$:
\begin{gather*}
\cb{\Pi^{2,1}_{y,\mtt j, s}(\eta)} \, :=\,
\Lambda^{(2)}_{y,\mtt j, s}(\eta) 
-\, \sum_{\mtt k \in \Lambda_{\mtt j} }
|B_{\mtt k}| \, G(s, y+B_{\mtt k} ) \, \mtt p_A'( u(s, y+B_{\mtt k}) )  \,
[\,\mss m_{y, \mtt k}  - u(s, y+B_{\mtt k}) \,]  \,,
\\
\,  \cb{\Lambda^{(2)}_{y,\mtt j, s}(\eta)} \,:= \,  |B_{2,\mtt j}| \,
G(s, y+B_{2, \mtt j} ) \, \mtt p_A'( u(s, y+B_{2, \mtt j}) )  \,
[\mss m_{y,2,\mtt j}  - u(s, y+B_{2, \mtt j})] \,.
\end{gather*}
The term $\Pi^{2,2}_{y, \mtt j, s}(\eta)$ compares the function
$\mf p_A (\mss m_{y, \mtt k})$ to its average with respect to the
canonical measure on the box $y+ B_{2, \mtt j}$:
\begin{gather*}
\cb{\Pi^{2,2}_{y, \mtt j, s}(\eta)} \,:=\,
\mtt 1_{\mf a_0} (\mss m_{y, 2, \mtt j}) \,
\sum_{\mtt k\in \Lambda_{\mtt j}}
|B_{\mtt k} |\, G  (s, y + B_{\mtt k} )\,
\big\{ \mf p_A (\mss m_{y, \mtt k}) 
\,-\, E_{\nu^{\rm c}_{2, \mtt j, \mss M_{y, 2, \mtt j} }}
[\mf p_A (\mss m_{y, \mtt k}) ]\,\big\}  \,.
\end{gather*}
The definition of $\Pi^{2,3}_{y, \mtt j, s}(\eta)$ requires a comment.
Let
\begin{gather*}
\cb{\Pi^{2,3}_{y, \mtt j, s}(\eta) } \,: =\,
\mtt 1_{\mf a_0} (\mss m_{y, 2, \mtt j}) \, 
\Big\{\, \sum_{\mtt k\in \Lambda_{\mtt j}}
|B_{\mtt k} |\, G  (s, y+ B_{\mtt k} )\,
\big(\, E_{\nu^{\rm c}_{2, \mtt j, \mss M_{y, 2, \mtt j} }}
[\mf p_A (\mss m_{y, \mtt k} )]
\,-\, E_{\nu^n_{u(s,\cdot)}} [  \mf p_A (\mss m_{y, \mtt k} )  ] \ \big)
\, - \, \Lambda^{(2)}_{y, \mtt j, s}(\eta) \, \Big\} \,.
\end{gather*}
Note that in the canonical expectation we have the polynomial
$\mf p_A = \mtt p_A + \ell_1^{-d}\, \mtt q_A$, while in the linear
term $\Lambda^{(2)}_{y, \mtt j, s}(\eta)$ it is the derivative of
$\mtt p_A$ which appears.  This is explained by the identity
\eqref{212}.  Finally,
\begin{equation*}
\cb{\Pi^{2,4}_{y, \mtt j, s}(\eta)} \,:=\,
\widehat{\Pi}^{2,2}_{y, \mtt j, s}(\eta) \,+\,
\widehat{\Pi}^{2,3}_{y, \mtt j, s}(\eta) \,,
\end{equation*}
where $\widehat{\Pi}^{2,q}_{y, \mtt j, s}(\eta)$, $q=2$, $3$, is
obtained from $\Pi^{2,q}_{y, \mtt j, s}(\eta)$ by replacing the
indicator $\mtt 1_{\mf a_0} (\mss m_{y, 2, \mtt j})$ by
$\mtt 1^{\rm c}_{\mf a_0} (\mss m_{y, 2, \mtt j})$.  As in Subsection
\ref{ssec9.1}, when we consider $y=0$ we drop the subscript $y$ from
the notation.

At this point we may repeat the reasoning followed in Subsection
\ref{ssec9.1}.  We start with $\Pi^{2,1}_{y, \mtt j, s}(\eta)$.  By
\eqref{135}, adding and subtracting the associated sum of terms
proportional to $[\eta_x - u(s,x/n)]$,
\begin{align}
\label{192c}
\bb E_{\mu_n} \Big[\,\Big| \int_0^t
\frac{1}{n^{d/2}} \frac{1}{|B_{\mtt 0}|} 
\sum_{y\in B_{\mtt 0} }
\sum_{\mtt j\in \Sigma^+_{L_3}}
\Pi^{2,1}_{y,\mtt j, s}(\eta^n(s)) \, ds \, \Big| \, \Big]
\,\le\, \frac{2}{n^{d/2}} \, \int_0^t \mtt S^{(\ref{135})}_{n,\ell_1\ell_2}
(G_s,\mtt p'_A, u_s)\, ds\,,
\end{align}
where $\mtt S^{(\ref{135})}_{n,r_1}$ is defined in \eqref{135}.  Here
we used the fact that
$\mtt S^{(\ref{135})}_{n,\ell_1} \le \mtt
S^{(\ref{135})}_{n,\ell_1\ell_2}$, which produces the factor $2$
multiplying $\mtt S^{(\ref{135})}_{n,\ell_1\ell_2}$.

We turn to $\Pi^{2,4}_{y,\mtt j, s}$. Note that
$|\Pi^{2,4}_{y,\mtt j, s} (\eta)| \le C_1 \, (\ell_1\ell_2)^d \,
\|G(s)\|_\infty$.  Hence, by Lemma \ref{s07b} with
$r_1=\ell_1\ell_2$,
\begin{equation}
\label{214}
\bb E_{\mu_n} \Big[\,\Big| \int_0^t
\frac{1}{n^{d/2}} \frac{1}{|B_{\mtt 0}|} 
\sum_{y\in B_{\mtt 0} } \sum_{\mtt j\in \Sigma^+_{L_3}}
\Pi^{2,4}_{y,\mtt j, s} (\eta^n(s)) \, ds
\, \Big| \, \Big]
\, \le  \, \mtt R^{(\ref{s07b})}_{n,\ell_1\ell_2, G} \,,
\end{equation}
where $\mtt R^{(\ref{s07b})}_{n,\ell, H}$ is defined in \eqref{202}.

To examine $\Pi^{2,3}_{y,\mtt j, s}$, let
\begin{equation}
\label{195d}
\begin{gathered}
\cb{F^{(2)}_{y,\mtt j, s} (\mss m_{y,2,\mtt j})}  \,:=\,
\mtt p_A (\mss m_{y, 2,\mtt j})
\,-\, E_{\nu^n_{u(s,\cdot)}} \big [ \,\mtt p_A (\mss m_{y, 2,\mtt j}) \,\big]
\,-\, \mtt p_A'(u(s, y+B_{2,\mtt j}))\, 
[\mss m_{y, 2,\mtt j} - u(s, y+B_{2,\mtt j})]\,,
\\
\cb{\overline{W}^{(2)}_{y, \mtt j, A, s}} \,:=\,
\frac{1}{ \ell^{2d}_1} \,
\Big\{\, \mtt q^{(2)}_A(\mss m_{y, 2, \mtt j})
\,-\, E_{\nu^n_{u(s,\cdot)}} \big[\,
\mtt q^{(2)}_A(\mss m_{y, 2, \mtt j})\,\big]\, \Big\}\,,
\end{gathered}
\end{equation}
so that, by \eqref{212},
$\mf p^{(2)}_A (\mss m_{y,2,\mtt j}) - E_{\nu^n_{u(s,\cdot)}}[\mf
p^{(2)}_A (\mss m_{y,2,\mtt j})]$ is the sum of
$F^{(2)}_{y,\mtt j, s}$, of $\overline{W}^{(2)}_{y, \mtt j, A, s}$ and
of the linear term appearing in $\Lambda^{(2)}_{y, \mtt j, s}$.
Denote by
$\cb{\mtt R^{(\ref{s44})}_{\mtt t_1, \mtt t_2, K_1, K_2} (s)}$,
$0\le s\le t$, the expression appearing on the right-hand side of
Proposition \ref{s44} with the monomial $\mss m^p$ replaced by
$\mf p_A$, $J = G(s)$, $\rho (\cdot) = u(s, \cdot)$, $f=f_s$,
$\mf a_1 = \mf c_{1,0}$, $\mf a_{3,3} = \mf c_{1,3}$.

\begin{lemma}
\label{s78b}
For all $n\ge 1$,
\begin{align*}
& \bb E_{\mu_n} \Big[\,\Big| \int_0^t
\frac{1}{n^{d/2}} \frac{1}{|B_{\mtt 0}|} 
\sum_{y\in B_{\mtt 0} } \sum_{\mtt j\in \Sigma^+_{L_3}}
\Pi^{2,3}_{y,\mtt j, s} (\eta^n(s))  \, ds \, \Big| \,
\Big]
\\
& \le\,
\bb E_{\mu_n} \Big[\,\Big| \int_0^t
\frac{1}{n^{d/2}} \frac{1}{|B_{\mtt 0}|} 
\sum_{y\in B_{\mtt 0} } \sum_{\mtt j \in \Sigma^+_{L_3}}
|B_{2, \mtt j} |\, G (s, y+B_{2, \mtt j} )\,
F^{(2)}_{y,\mtt j, s} (\mss m_{y, 2, \mtt j}(s)) \, ds \, \Big| \,
\Big] \, + \,  \mtt R^{(\ref{s78b})}_n  \,,
\end{align*}
where
\begin{equation*}
\cb{\mtt R^{(\ref{s78b})}_n } \,:=\, \frac{1}{n^{d/2}}\,
\int_0^t \mtt R^{(\ref{s44})}_{\ell_1, \ell_2, 1, 1} (s)\, ds
\, +\, 2\,\mtt R^{(\ref{s07b})}_{n,\ell_1\ell_2, G}
\,+ \, \mtt R^{(\ref{s32})}_{n, \ell_1\ell_2, \ell_1}\,,
\end{equation*}
\begin{align*}
\cb{\mtt R^{(\ref{s32})}_{n,\ell_1\ell_2, \ell_1} }\, &:=\,
\frac{\mf c_{1,0}}{n^{d/2}} \,
\Big(\frac{\ell_2}{\ell_1}\Big)^{d}\,
\Big(\frac{n}{(\ell_1\ell_2)^2}\Big)^d 
\int_0^t  \|G(s) \|_\infty \; ds
\\
& + \, 
\frac{\mf c_{1,0}}{n^{d/2}}\, 
\Big(\frac{\ell_2}{\ell_1}\Big)^{d} \, 
\int_0^t  \Vert G(s)\Vert_\infty  
\big[\, H_n(f_s \,|\, \nu^n_{u(s,\cdot)})
+ \ln (2)\,\big] \; ds  \,.
\end{align*}
\end{lemma}

\begin{proof}
The proof follows closely the one of Lemma \ref{s78}, with two
modifications.  As before, we move the sum over $y$ out of the
expectation, fix $y=0$, and make sure that the estimates are uniform
over $y$.

First, instead of applying Proposition \ref{s36} we apply Proposition
\ref{s44} with $\mf p = \mf p_A$, $\mtt t_1 = \ell_1$,
$\mtt t_2 = \ell_2$, $\mtt t_3 = L_3$, $K_1 = K_2 = 1$.  This is
legitimate: as $\ell_1\ell_2 = n^{3/5}$ in $d=2$ and
$\ell_1\ell_2 = \upsilon_n^{2}\, n^{2/5}$ in $d=3$, the conditions
$(\ell_1\ell_2)^{d+3} \le n^3$, $(\ell_1\ell_2)^{1+(d/4)} \le n$ and
$\ell_2 \le \ell_1$ are satisfied for all $n$ large enough.
Proposition \ref{s44} replaces
$\sum_{\mtt k\in \Lambda_{\mtt j}} |B_{\mtt k}| G(s, B_{\mtt k}) \{\mf
p_A(\mss m_{\mtt k}) - E_{\nu^{\rm c}_{2,\mtt j, \mss M}} [\mf
p_A(\mss m_{\mtt k})]\}$ by
$|B_{2,\mtt j}|\, G(s, B_{2, \mtt j})\, \{ \mf p^{(2)}_A(\mss m_{2,
\mtt j}) - E_{\nu^n_{u(s,\cdot)}} [ \mf p^{(2)}_A(\mss m_{2, \mtt j})
]\}$, at a cost bounded by
$n^{-d/2}\int_0^t \mtt R^{(\ref{s44})}_{\ell_1, \ell_2, 1, 1} (s)\,
ds$.

Second, by \eqref{212}, $\mf p^{(2)}_A(\mss m_{2, \mtt j}) -
E_{\nu^n_{u(s,\cdot)}} [ \mf p^{(2)}_A(\mss m_{2, \mtt j})]$ is the
sum of $F^{(2)}_{\mtt j, s}$, of $\overline{W}^{(2)}_{\mtt j, A, s}$
and of the linear term of $\Lambda^{(2)}_{\mtt j, s}$, which cancels
the one appearing in the definition of $\Pi^{2,3}_{\mtt j, s}$.  It
remains to remove $\overline{W}^{(2)}_{\mtt j, A, s}$.  This is done
as in the proof of \eqref{208}: apply first Lemma \ref{s07b} to get
rid of the indicator, at a cost bounded by
$\mtt R^{(\ref{s07b})}_{n,\ell_1\ell_2, G}$ (the extra factor
$\ell_1^{-2d}$ present in the definition of
$\overline{W}^{(2)}_{\mtt j, A, s}$ only improves this bound), and
then Lemma \ref{s32} with $r_1 = \ell_1\ell_2$, $r_2 = L_3$,
$F = \mtt q^{(2)}_A$,
$\mf a(\mtt j) = (\ell_2/\ell_1)^{d}\, G(s, B_{2,\mtt j})$ and
$\gamma = \gamma(s) = (\ell_1/\ell_2)^{d} / \|G(s)\|_\infty$.  Note
that
$|B_{2,\mtt j}|\, \ell_1^{-2d} = (\ell_2/\ell_1)^{d}$, so that
$\Vert \mf a\Vert_\infty \gamma \le 1$.  This yields
$\mtt R^{(\ref{s32})}_{n, \ell_1\ell_2, \ell_1}$.

Finally, we remove the indicator
$\mtt 1_{\mf a_0} \big(\mss m_{y, 2, \mtt j}(s) \big)$ multiplying
$F^{(2)}_{y,\mtt j, s}$, as in Proposition \ref{s79}, at the cost of a
second term $\mtt R^{(\ref{s07b})}_{n,\ell_1\ell_2, G}$.
\end{proof}

It remains to estimate $\Pi^{2,2}_{y, \mtt j, s}$.  Recall equations
\eqref{209}--\eqref{196} and the definition of $\Psi_n(s)$.  The
expression $\Pi^{2,2}_{y, \mtt j, s}$ is handled exactly as
$\Pi^{(2)}_{y, \mtt k, s}$ in Lemma \ref{s85}: the absolute value is
removed by the entropy inequality and the Feynman--Kac formula, and
the resulting variational problem is estimated with the results which
bound expressions in terms of the functional
$I_n(f ; \nu^n_{u(s,\cdot)})$ alone.  At the second scale the relevant
result is Proposition \ref{s37} rather than Lemma \ref{s03}.  Let
\begin{align*}
\cb{\mtt R^{( \ref{s37} )}_n} \, :=\, \frac{1}{n^{d/2}}\,
\int_0^t \mtt R^{( \ref{s37} )}_n (G(s)) \, ds \,,
\end{align*}
where $\mtt R^{( \ref{s37} )}_n (J)$ is the second term on the
right-hand side of the statement of Proposition \ref{s37} with
$\delta^{-1} = 2$, $\mf a_{3,3} = \mf c_{1,3}$,
$\mf a_{3,1} = \mf c_{1,1}$, $\mtt t_1 = \ell_1$,
$\mtt t_2 = \ell_2$, $K_1 = 1$, that is,
\begin{equation*}
\mtt R^{( \ref{s37} )}_n (J) \,=\,
\mf c_{1,3}\, \Big\{\, \mf c_{\rm LS} \, \mf C_2
\Big(\frac{\ell_1\ell_2}{n}\Big)^4\, n^d
\,+\,  \Vert J\Vert^2_\infty\,
\ell_2^{d+2}  \, \Big( \frac{n}{\ell_1}\Big)^{d-2} \,\Big\}\,.
\end{equation*}

\begin{lemma}
\label{s81b}
Under the hypotheses of Proposition \ref{s79},
\begin{align*}
& \frac{1}{|B_{\mtt 0}|} 
\sum_{y\in B_{\mtt 0} } 
\bb E_{\mu_n} \Big[\,\Big| \int_0^t
\frac{1}{n^{d/2}} \sum_{\mtt j\in \Sigma^+_{L_3}}
\Pi^{2,2}_{y,\mtt j, s}(\eta^n(s)) \, ds \, \Big| \, \Big]
\,\le \, \frac{1}{2\, n^{d/2}}  \, \bb E_{\mu_n} \Big[\, \int_0^t
\Psi_n(s)  \, ds \, \Big]
\,+\, \mtt R^{( \ref{s81b} )}_n \,,
\end{align*}
where
\begin{equation*}
\cb{ \mtt R^{( \ref{s81b} )}_n } \,:=\,
\frac{1}{n^{d/2}}\,
\big\{ \, H_n(f_0 \,|\, \nu^n_{u_0(\cdot)}) \,+\, \ln (2)\,\big\}\,
\,+\, \mtt R^{( \ref{s03} )}_n
\,+\, \mtt R^{( \ref{s37} )}_n \,,
\end{equation*}
and $\mtt R^{( \ref{s03} )}_n$ is defined in Lemma \ref{s85}.
\end{lemma}

\begin{proof}
As before, we fix $y=0$ and make sure that all estimates are uniform
over $y$.  Let
\begin{align*}
\cb{V^{(2)}_{s} }  \,:= \, \sum_{\mtt j\in \Sigma^+_{L_3}}
\Pi^{2,2}_{\mtt j, s}(\eta^n(s)) \,,
\end{align*}
and repeat the proof of Lemma \ref{s85} up to equation \eqref{210},
with $V^{(2)}_s$ in place of $V_s$.  We obtain the same bound with
\begin{align*}
\mtt X^b_n (\eta) \, =\, b\,
\sum_{\mtt j\in \Sigma^+_{L_3}} \, \Pi^{2,2}_{\mtt j, s} (\eta)
\,+\, \frac{1}{2}\, 
\frac{1}{|B_{\mtt 0}^o|} \sum_{D\in \mc A}
\sum_{y\in B_{\mtt 0} } \sum_{\mtt k\in \Sigma^+_{L_2}} \,
\Pi^{2,D}_{y, \mtt k, s} (\eta) \,.
\end{align*}
In the first sum we have cubes of width $\ell_1\ell_2$, and in the
second cubes of width $\ell_1$.  It remains to apply Proposition
\ref{s37} with $\mf p = \mf p_A$, $\mtt t_1 = \ell_1$,
$\mtt t_2 = \ell_2$, $K_1=1$ and $\delta = 1/2$ to estimate the first
sum, and Lemma \ref{s03} with $\delta = 1/(2|\mc A|)$ to estimate the
second one.  The hypotheses of Proposition \ref{s37} are fulfilled for
all $n$ large enough by the choices of $\ell_1$, $\ell_2$ made at the
beginning of this subsection.
\end{proof}

Adding the estimates \eqref{192c}, \eqref{214} and Lemmata
\ref{s78b}, \ref{s81b} to the decomposition \eqref{213} yields the
main result of this subsection.

\begin{conclusion}
\label{s87b}
In view of the decomposition \eqref{213}, of the estimates
\eqref{192c}, \eqref{214}, and of Lemmata \ref{s78b}, \ref{s81b},
\ref{s82}, the expression \eqref{206} is bounded by
\begin{align*}
\bb E_{\mu_n} \Big[\,\Big| \int_0^t
\frac{1}{n^{d/2}} \frac{1}{|B_{\mtt 0}|} 
\sum_{y\in B_{\mtt 0} } \sum_{\mtt j \in \Sigma^+_{L_3}}
|B_{2, \mtt j} |\, G (s, y+B_{2, \mtt j} )\,
F^{(2)}_{y,\mtt j, s} (\mss m_{y, 2, \mtt j}(s)) \, ds \, \Big| \,
\Big]
\,+\, \frac{1}{2}  \, \mtt R^{( \ref{s82} )}_{n,d}(\delta)
\,+\, \mtt R^{(\ref{s87b})}_n\,,
\end{align*}
for all $\delta>0$, where
\begin{equation*}
\cb{\mtt R^{(\ref{s87b})}_n} \,:=\,
\frac{2}{n^{d/2}} \, \int_0^t \mtt S^{(\ref{135})}_{n,\ell_1\ell_2}
(G_s,\mtt p'_A, u_s)\, ds
\,+\, \mtt R^{(\ref{s07b})}_{n,\ell_1\ell_2, G}
\,+\, \mtt R^{(\ref{s78b})}_n \,+\, \mtt R^{( \ref{s81b} )}_n \,.
\end{equation*}
It therefore remains to estimate the first one.  The argument is
dimension dependent.
\end{conclusion}

\begin{remark}
\label{r10}
At this stage, the expression left to estimate has the same form as
the one in Subsection \ref{ssec9.1}, with
the cubes $B_{\mtt k}$ of width $\ell_1$ replaced by the cubes
$B_{2, \mtt j}$ of width $\ell_1\ell_2$, and the polynomial $\mf p_A$
replaced by $\mtt p_A$.  In dimension $d=2$ the cubes are now large
enough for the argument of Lemma \ref{s80} to apply, and the proof is
complete.  In dimension $d=3$ one further scale is needed, and the
construction of this subsection is repeated once more.
\end{remark}

\subsection{Dimension $d=2$}
\label{ssec9.4}

In this subsection we complete the proof of the Boltzmann--Gibbs
principle in dimension $2$.  Recall that $\ell_1 = n^{2/5}$,
$\ell_2 = n^{1/5}$, so that $\ell_1\ell_2 = n^{3/5}$ and
$L_3 = n^{2/5}$.  In view of Conclusion \ref{s87b}, we have to
estimate the term involving $F^{(2)}_{y,\mtt j, s}$.  This is the
content of the next result, which is the analogue of Lemma \ref{s80}
at the second scale.

\begin{lemma}
\label{s80b}
Under the hypotheses of Proposition \ref{s79},
\begin{align*}
\bb E_{\mu_n} \Big[\,\Big| \int_0^t
\frac{1}{n^{d/2}} \frac{1}{|B_{\mtt 0}|} 
\sum_{y\in B_{\mtt 0} } \sum_{\mtt j\in \Sigma^+_{L_3}}
|B_{2, \mtt j} |\, G (s, y+B_{2, \mtt j} )\,
F^{(2)}_{y, \mtt j, s} (\mss m_{y, 2, \mtt j}(s)) \, ds \, \Big|\,\Big]
\,\le\, \mtt R^{(\ref{s80b})}_n\,,
\end{align*}
where
$\cb{\mtt R^{(\ref{s80b})}_n} \,=\, \mtt S^{(\ref{120})}_{n, \ell_1\ell_2}
\,+\, \mtt R^{ ( \ref{l02} ) }_{n, \ell_1\ell_2}$,
\begin{gather*}
\cb{\mtt S^{(\ref{120})}_{n,r_1}} \,:=\, C_1 \, \frac{n^{d/2}}{r_1^d}
\, \int_0^t \|G(s)\|_\infty\, ds\,,
\quad
\cb{ \mtt R^{ ( \ref{l02} ) }_{n,r_1}} \,:=\, 
\frac{C_1}{n^{d/2}}
\int_0^t \|G(s) \|_\infty \, \Big\{\, H_n(f_s \,|\, \nu^n_{u(s,\cdot)})
+ \Big(\frac{n}{r_1}\Big)^d\, \Big\}\, ds \,.
\end{gather*}
\end{lemma}

\begin{proof}
As before, we set $y=0$ and make sure that the estimates are uniform
over $y$.  By \eqref{120} with $r_1 = \ell_1\ell_2$, in the definition
\eqref{195d} of $F^{(2)}_{\mtt j, s}$ we may replace
$E_{\nu^n_{u(s,\cdot)}}\big[\mtt p_A (\mss m_{2, \mtt j})\big]$ by
$\mtt p_A ( u(s, B_{2, \mtt j}))$ at a cost bounded by
\begin{align*}
\frac{1}{n^{d/2}}\int_0^t \|G(s)\|_\infty\,
\sum_{\mtt j\in \Sigma^+_{L_3}} (\ell_1\ell_2)^d\,
\big|E_{\nu^n_{u(s,\cdot)}}\big[\mtt p_A(\mss m_{2, \mtt j}(s))\big]
- \mtt p_A( u(s, B_{2, \mtt j}))\big|\; ds
\ \leq \ \mtt S^{(\ref{120})}_{n, \ell_1\ell_2}\, .
\end{align*}
At this point, it remains to estimate
\begin{align*}
\cb{\mtt p^{(2)}_{\mtt j, s} (\mss m_{2, \mtt j}) } \,:=\,
\mtt p_A(\mss m_{2, \mtt j}) - \mtt p_A(u(s, B_{2, \mtt j}))-\mtt
p_A'(u(s, B_{2, \mtt j}))\big[\mss m_{2, \mtt j} - u(s, B_{2,\mtt j})\big] \,.
\end{align*}
By a Taylor expansion, the absolute value of this expression is less
than or equal to
$C_1 \big[\mss m_{2, \mtt j} - u(s, B_{2, \mtt j})\big]^2$.  Thus, by
Lemma \ref{l02} with $r_1 = \ell_1\ell_2$, $r_2 = L_3$,
\begin{align*}
\bb E_{\mu_n}\Big[ \,\Big| \int_0^t \frac{1}{n^{d/2}}
\sum_{\mtt j\in \Sigma^+_{L_3}} (\ell_1\ell_2)^d\, 
G(s, B_{2, \mtt j}) \,  \mtt p^{(2)}_{\mtt j, s} (\mss m_{2, \mtt j}(s))
\, ds\,\Big|\,\Big] \,\le\, \mtt R^{ ( \ref{l02} ) }_{n, \ell_1\ell_2} \,,
\end{align*}
which completes the proof.
\end{proof}

As in dimension $d=1$, it is convenient, for the tightness estimates of
Section \ref{sec11}, to bound the integral appearing in the
Boltzmann--Gibbs principle between two arbitrary times $s$ and $t$.

\begin{proof}[Proof of Theorem \ref{t02b} in $d=2$]
Fix $T>0$ and $0\le s<t\le T$. Repeat the proof of the theorem in
dimension $1$ up to equation \eqref{235}.  By the decomposition
\eqref{203}, Proposition \ref{s79}, Lemmata \ref{s85}, \ref{s82},
Conclusion \ref{s87b}, and Lemma \ref{s80b}, the expression appearing
in \eqref{235} is bounded by
\begin{align*}
\Vert G \Vert_{T,2} \; \big\{\, 
\mtt R^{( \ref{s79} )}_n \,+\, \mtt R^{( \ref{s85} )}_n
\,+\, \mtt R^{( \ref{s82} )}_{n,d} (1/2)
\,+\, \mtt R^{(\ref{s87b})}_n \,+\, \mtt R^{( \ref{s80b} )}_n \,\big\}
\,.
\end{align*}
In the formulae for the remainders we replace the time-integrals
between $0$ and $t$ by integrals between $s$ and $t$, the entropy
$H_n(f_0 \,|\, \nu^n_{u_0(\cdot)})$ by
$H_n(f_s \,|\, \nu^n_{u(s, \cdot)})$, and $G$ by $\widehat{G}$.

We estimate the previous displayed equation. Since
$\sup_{0\le r\le T} \Vert \widehat G (r)\Vert_{C^2} \le 1$, inspecting
all terms yields that the expression inside curly braces is bounded by
\begin{align*}
&\frac{2}{n^{d/2}}\,
\big\{ \, H_n(f_s \,|\, \nu^n_{u(s, \cdot)}) \,+\, \ln (2)\,\big\}\,
+\, \mf c_{1,3}\,  \frac{1}{n^{d/2}}\, 
\int_s^t  \, H_n( f_r \,|\, \nu^n_{u(r, \cdot)} ) \, dr
\\
\, &+\, \frac{1}{n^{d/2}}\, 
\int_s^t \Big\{ \, 2\, n^2 \, I_n(f_r ;\, \nu^n_{u(r,\cdot)}) \, +\,
\mf c_{1,4}\,  n^{4/5} \,  \Big\} \, dr \,,
\end{align*}
as claimed.
\end{proof}

\subsection{Dimension $d=3$}
\label{ssec9.5}

In dimension $3$ the cubes of width $\ell_1\ell_2$ reached in the
previous subsection are still too small, and we have to perform one
further multi-scale step. Recall from Subsection \ref{ssec9.3} the
choices $\ell_1 = \upsilon_n\, n^{1/4}$,
$\ell_2 = \upsilon_n\, n^{3/20}$, and let
$\ell_3 \,=\, \kappa\, \upsilon_n^{-2}\, n^{1/10}$,
$n \,=\, \ell_1 \ell_2 \ell_3 L_4$, where $\kappa\ge 1$ is a constant
which increases to $\infty$ after $n\to\infty$.  With these choices,
$\ell_1\ell_2 = \upsilon_n^{2}\, n^{2/5}$ and
$\ell_1\ell_2\ell_3 = \kappa\, n^{1/2}$, the width required by the
non-equilibrium fluctuations.

We repeat the construction of Subsection \ref{ssec9.3} with the
cubes $B_{2, \mtt j}$ of width $\ell_1\ell_2$ in the role of the cubes
$B_{\mtt k}$, and cubes of width $\ell_1\ell_2\ell_3$ in the role of
the cubes $B_{2, \mtt j}$.  That is, we apply the results of Section
\ref{sec9} with $\mtt t_1 = \ell_1\ell_2$, $\mtt t_2 = \ell_3$,
$\mtt t_3 = L_4$.  Denote by $\cb{B_{3, \mtt i}}$,
$\mtt i\in \Sigma^+_{L_4}$, the cubes of width $\ell_1\ell_2\ell_3$,
by $\cb{\Lambda_{3, \mtt i}}$ the set of indices
$\mtt j \in \Sigma^+_{L_3}$ such that
$B_{2,\mtt j}\subset B_{3, \mtt i}$, and let
\begin{align*}
\cb{\mss M_{y, 3, \mtt i}} := 
\sum_{z\in y + B_{3, \mtt i}} \eta_z\,, \quad
\cb{\mss m_{y, 3, \mtt i}} := \mss m (y+  B_{3, \mtt i})
\,=\, \frac{1}{|B_{3, \mtt i}|} \, \mss M_{y, 3, \mtt i}
\,, \quad
y\in B_{\mtt 0}\,, \; \mtt i \in \Sigma^+_{L_4}\,.
\end{align*}
Denote by $\cb{\nu^{\rm c}_{3, \mtt i, \mss M}}$ the canonical measure
on the cube $y+B_{3, \mtt i}$ with $\mss M$ particles.

Recall the definition of the polynomial $\mf p^{(2)}_A (\rho)$
introduced in \eqref{212}.  The polynomial produced by Theorem
\ref{s26} at this scale is $W_{\mf p^{(2)}_A}$ computed with
$\mtt t_1 = \ell_1\ell_2$.  By \eqref{111} and \eqref{212},
\begin{equation}
\label{216}
\cb{\mf p^{(3)}_A (\rho)} \,:=\, W_{\mf p^{(2)}_A} (\rho)
\,=\,  \mtt p_A(\rho)
\,-\, \frac{1}{(\ell_1\ell_2)^{d}}\, \mtt q_A(\rho)
\,+\, \frac{1}{\ell_1^{2d}}\, \mtt q^{(2)}_A(\rho)
\,+\, \frac{1}{\ell_1^{2d}\, (\ell_1\ell_2)^{d}}\, \mtt q^{(3)}_A(\rho)\,,
\end{equation}
where
$\cb{\mtt q^{(3)}_A(\rho)} := (1/2)\, [\mtt q^{(2)}_A]''(\rho)\,
\chi(\rho)$.  In contrast with what happens at the second scale, there
is no cancellation.  It is, however, harmless, as shown in Lemma
\ref{s78c} below.

Repeating verbatim the decomposition \eqref{213}, write
\begin{equation}
\label{217}
\frac{1}{n^{d/2}} 
\frac{1}{|B_{\mtt 0}|} 
\sum_{y\in B_{\mtt 0} } \sum_{\mtt j\in \Sigma^+_{L_3}}
|B_{2, \mtt j} |\, G  (s, y+B_{2, \mtt j} )\,
F^{(2)}_{y,\mtt j, s} (\mss m_{y, 2, \mtt j}) \,
=\,
\frac{1}{n^{d/2}} \frac{1}{|B_{\mtt 0}|} 
\sum_{y\in B_{\mtt 0} }  \sum_{\mtt i\in \Sigma^+_{L_4}} \,
\sum_{q=1}^4 \Pi^{3,q}_{y, \mtt i, s}(\eta) \,,
\end{equation}
where the functions $\Pi^{3,q}_{y, \mtt i, s}$, $1\le q\le 4$, are
defined as the functions $\Pi^{2,q}_{y, \mtt j, s}$, replacing
everywhere the cubes $B_{\mtt k}$, $\mtt k\in \Sigma^+_{L_2}$, by the
cubes $B_{2, \mtt j}$, $\mtt j\in \Sigma^+_{L_3}$, the cubes
$B_{2, \mtt j}$, $\mtt j\in \Sigma^+_{L_3}$, by the cubes
$B_{3, \mtt i}$, $\mtt i\in \Sigma^+_{L_4}$, the set
$\Lambda_{\mtt j}$ by $\Lambda_{3, \mtt i}$, the canonical measure
$\nu^{\rm c}_{2, \mtt j, \mss M}$ by $\nu^{\rm c}_{3, \mtt i, \mss M}$,
and the polynomial $\mf p_A$ by $\mf p^{(2)}_A$.

The estimates of $\Pi^{3,1}_{y, \mtt i, s}$ and
$\Pi^{3,4}_{y, \mtt i, s}$ are identical to \eqref{192c} and
\eqref{214}, with $\ell_1\ell_2$ replaced by $\ell_1\ell_2\ell_3$:
\begin{equation}
\label{218}
\begin{gathered}
\bb E_{\mu_n} \Big[\,\Big| \int_0^t
\frac{1}{n^{d/2}} \frac{1}{|B_{\mtt 0}|} 
\sum_{y\in B_{\mtt 0} }
\sum_{\mtt i\in \Sigma^+_{L_4}}
\Pi^{3,1}_{y,\mtt i, s}(\eta^n(s)) \, ds \, \Big| \, \Big]
\,\le\, \frac{2}{n^{d/2}} \, \int_0^t
\mtt S^{(\ref{135})}_{n,\ell_1\ell_2\ell_3}
(G_s,\mtt p'_A, u_s)\, ds\,,
\\
\bb E_{\mu_n} \Big[\,\Big| \int_0^t
\frac{1}{n^{d/2}} \frac{1}{|B_{\mtt 0}|} 
\sum_{y\in B_{\mtt 0} } \sum_{\mtt i\in \Sigma^+_{L_4}}
\Pi^{3,4}_{y,\mtt i, s} (\eta^n(s)) \, ds
\, \Big| \, \Big]
\, \le  \, \mtt R^{(\ref{s07b})}_{n,\ell_1\ell_2\ell_3, G} \,.
\end{gathered}
\end{equation}

To examine $\Pi^{3,3}_{y,\mtt i, s}$, let
\begin{equation}
\label{195e}
\cb{F^{(3)}_{y,\mtt i, s} (\mss m_{y,3,\mtt i})}  \,:=\,
\mtt p_A (\mss m_{y, 3,\mtt i})
\,-\, E_{\nu^n_{u(s,\cdot)}} \big [ \,\mtt p_A (\mss m_{y, 3,\mtt i}) \,\big]
\,-\, \mtt p_A'(u(s, y+B_{3,\mtt i}))\, 
[\mss m_{y, 3,\mtt i} - u(s, y+B_{3,\mtt i})]\,,
\end{equation}
and let $\overline{W}^{(3)}_{y, \mtt i, A, s}$ be the sum of the three
correctors appearing in \eqref{216}, centred with respect to
$\nu^n_{u(s,\cdot)}$:
\begin{equation*}
\cb{\overline{W}^{(3)}_{y, \mtt i, A, s}} \,:=\,
\Big\{ \mf p^{(3)}_A(\mss m_{y, 3, \mtt i}) - \mtt p_A(\mss m_{y, 3, \mtt i})
\Big\}
\,-\, E_{\nu^n_{u(s,\cdot)}} \Big[\,
\mf p^{(3)}_A(\mss m_{y, 3, \mtt i}) - \mtt p_A(\mss m_{y, 3, \mtt i})
\,\Big]\,.
\end{equation*}

\begin{lemma}
\label{s78c}
For all $n\ge 1$,
\begin{align*}
& \bb E_{\mu_n} \Big[\,\Big| \int_0^t
\frac{1}{n^{d/2}} \frac{1}{|B_{\mtt 0}|} 
\sum_{y\in B_{\mtt 0} } \sum_{\mtt i\in \Sigma^+_{L_4}}
\Pi^{3,3}_{y,\mtt i, s} (\eta^n(s))  \, ds \, \Big| \,
\Big]
\\
& \le\,
\bb E_{\mu_n} \Big[\,\Big| \int_0^t
\frac{1}{n^{d/2}} \frac{1}{|B_{\mtt 0}|} 
\sum_{y\in B_{\mtt 0} } \sum_{\mtt i \in \Sigma^+_{L_4}}
|B_{3, \mtt i} |\, G (s, y+B_{3, \mtt i} )\,
F^{(3)}_{y,\mtt i, s} (\mss m_{y, 3, \mtt i}(s)) \, ds \, \Big| \,
\Big] \, + \,  \mtt R^{(\ref{s78c})}_n  \,,
\end{align*}
where
\begin{equation*}
\cb{\mtt R^{(\ref{s78c})}_n } \,:=\, \frac{1}{n^{d/2}}\,
\int_0^t \mtt R^{(\ref{s44})}_{\ell_1\ell_2, \ell_3, \kappa^6, 1} (s)\, ds
\, +\, 2\,\mtt R^{(\ref{s07b})}_{n,\ell_1\ell_2\ell_3, G}
\,+ \, \mtt R^{(\ref{s45})}_{n, \ell_1\ell_2\ell_3}\,,
\end{equation*}
\begin{equation*}
\cb{\mtt R^{(\ref{s45})}_{n,\ell_1\ell_2\ell_3} }\, :=\,
\frac{C_1}{n^{d/2}} \,
\int_0^t  \Vert G(s)\Vert_\infty  \,\Big\{ \,
H_n(f_s \,|\, \nu^n_{u(s,\cdot)})
\,+\, \ell_3^d\, \Big(\frac{n}{(\ell_1\ell_2\ell_3)^{3/2}}\Big)^d
\,\Big\} \; ds  \,.
\end{equation*}
\end{lemma}

\begin{proof}
The proof repeats the one of Lemma \ref{s78b}.  We fix $y=0$ and make
sure that the estimates are uniform over $y$.

We apply Proposition \ref{s44} with $\mf p = \mf p^{(2)}_A$,
$\mtt t_1 = \ell_1\ell_2$, $\mtt t_2 = \ell_3$, $\mtt t_3 = L_4$,
$K_1 = \kappa^6$, $K_2 = 1$.  An inspection shows that the hypotheses
of that result are fulfilled for all $n$ large enough. For instance,
$\ell_3 \le \ell_1\ell_2$ because
$\ell_3/(\ell_1\ell_2) = \kappa\, \upsilon_n^{-4}\, n^{-3/10}$, which
vanishes since $\upsilon_n^{-1}\le n^{3/100}$.  This
proposition allows us to replace
\begin{align*}
\sum_{\mtt j\in \Lambda_{3, \mtt i}} |B_{2, \mtt j}| G(s, B_{2,\mtt
j}) \{\mf p^{(2)}_A(\mss m_{2, \mtt j}) - E_{\nu^{\rm c}_{3,\mtt i,
\mss M}} [\mf p^{(2)}_A(\mss m_{2, \mtt j})]\}
\end{align*}
by
\begin{align*}
|B_{3,\mtt i}|\, G(s, B_{3, \mtt i})\, \{ \mf p^{(3)}_A(\mss m_{3,
\mtt i}) - E_{\nu^n_{u(s,\cdot)}} [ \mf p^{(3)}_A(\mss m_{3, \mtt i})
]\}
\end{align*}
at a cost bounded by
$n^{-d/2}\int_0^t \mtt R^{(\ref{s44})}_{\ell_1\ell_2, \ell_3,
\kappa^6, 1} (s)\, ds$.

By \eqref{216}, the resulting expression is the sum of
$F^{(3)}_{\mtt i, s}$, of $\overline{W}^{(3)}_{\mtt i, A, s}$ and of
the linear term, which cancels the one appearing in the definition of
$\Pi^{3,3}_{\mtt i, s}$.  It remains to remove
$\overline{W}^{(3)}_{\mtt i, A, s}$.  We first apply Lemma \ref{s07b}
to get rid of the indicator, at a cost bounded by
$\mtt R^{(\ref{s07b})}_{n,\ell_1\ell_2\ell_3, G}$.  We then estimate
each of the three correctors of \eqref{216} separately.  The dominant
one is $(\ell_1\ell_2)^{-d}\, \mtt q_A$. Its contribution is given by
\begin{equation*}
\frac{1}{n^{d/2}} \int_0^t \int \Big| \sum_{\mtt i\in \Sigma^+_{L_4}}
\ell_3^d\, G(s, B_{3,\mtt i}) \,
\big\{\, \mtt q_A (\mss m_{3, \mtt i}) -
E_{\nu^n_{u(s,\cdot)}} [ \mtt q_A (\mss m_{3, \mtt i})]\, \big\}\,
\Big|\, f_s\, d\nu^n_{u(s,\cdot)}\, ds \,.
\end{equation*}
As $\mtt q_A$ is a polynomial, Lemma \ref{s45} applies to each of its
monomials with $r_1 = \ell_1\ell_2\ell_3$, $r_2 = L_4$,
$\mf a (\mtt i) = \ell_3^d\, G(s, B_{3,\mtt i})$ and
$\mf c_0 = \Vert G(s)\Vert_\infty$.  The hypothesis
$\Vert \mf a\Vert_\infty \le \mf c_0\, r_1^{d/2}$ of that lemma reads
$\ell_3^d \le (\ell_1\ell_2\ell_3)^{d/2}$, that is,
$\ell_3 \le (\ell_1\ell_2\ell_3)^{1/2} = \kappa^{1/2} n^{1/4}$, which
holds for all $n$ large enough because
$\ell_3 = \kappa\, \upsilon_n^{-2}\, n^{1/10}$,
$\upsilon_n^{-1}\le n^{3/100}$.  This yields
$\mtt R^{(\ref{s45})}_{n, \ell_1\ell_2\ell_3}$.  The two remaining
correctors of \eqref{216} carry the extra factor $\ell_1^{-2d}$ and
are estimated in the same way, at a smaller cost.

Finally, we remove the indicator
$\mtt 1_{\mf a_0} \big(\mss m_{y, 3, \mtt i}(s) \big)$ multiplying
$F^{(3)}_{y,\mtt i, s}$ at the cost of a second factor
$\mtt R^{(\ref{s07b})}_{n,\ell_1\ell_2\ell_3, G}$.
\end{proof}

The term $\Pi^{3,2}_{y, \mtt i, s}$ is estimated as
$\Pi^{2,2}_{y, \mtt j, s}$ in Lemma \ref{s81b}.  Let
\begin{align*}
\cb{\mtt R^{( \ref{s37} ), 3}_n} \, :=\, \frac{1}{n^{d/2}}\,
\int_0^t 
\mtt R^{( \ref{s37} ), 3}_n (G(s)) \, ds \,,
\end{align*}
where $\mtt R^{( \ref{s37} ), 3}_n (J)$ is the second term on the
right-hand side of the statement of Proposition \ref{s37} with
$\delta^{-1} = 2$, $\mtt t_1 = \ell_1\ell_2$,
$\mtt t_2 = \ell_3$, $K_1 = \kappa^6$, that is,
\begin{equation*}
\mtt R^{( \ref{s37} ), 3}_n (J) \,=\,
\mf c_{1,3}\,
\Big\{\, \mf c_{\rm LS} \, \mf C_2
\Big(\frac{\ell_1\ell_2\ell_3}{n}\Big)^4\, n^d
\,+\,  \Vert J\Vert^2_\infty\,
\ell_3^{d+2}  \, \Big( \frac{n}{\ell_1\ell_2}\Big)^{d-2} \,\Big\}\,.
\end{equation*}

\begin{lemma}
\label{s81c}
Under the hypotheses of Proposition \ref{s79},
\begin{align*}
& \frac{1}{|B_{\mtt 0}|} 
\sum_{y\in B_{\mtt 0} } 
\bb E_{\mu_n} \Big[\,\Big| \int_0^t
\frac{1}{n^{d/2}} \sum_{\mtt i\in \Sigma^+_{L_4}}
\Pi^{3,2}_{y,\mtt i, s}(\eta^n(s)) \, ds \, \Big| \, \Big]
\,\le \, \frac{1}{2\, n^{d/2}}  \, \bb E_{\mu_n} \Big[\, \int_0^t
\Psi_n(s)  \, ds \, \Big]
\,+\, \mtt R^{( \ref{s81c} )}_n \,,
\end{align*}
where
\begin{equation*}
\cb{ \mtt R^{( \ref{s81c} )}_n } \,:=\,
\frac{1}{n^{d/2}}\,
\big\{ \, H_n(f_0 \,|\, \nu^n_{u_0(\cdot)}) \,+\, \ln (2)\,\big\}\,
\,+\, \mtt R^{( \ref{s03} )}_n
\,+\, \mtt R^{( \ref{s37} ), 3}_n \,.
\end{equation*}
\end{lemma}

\begin{proof}
The proof is the one of Lemma \ref{s81b}, applying Proposition
\ref{s37} with $\mf p = \mf p^{(2)}_A$, $\mtt t_1 = \ell_1\ell_2$,
$\mtt t_2 = \ell_3$, $K_1 = \kappa^6$, $\delta = 1/2$ to estimate the
first sum, and Lemma \ref{s03} with $\delta = 1/(2|\mc A|)$ to
estimate the second one.
\end{proof}

It remains to estimate the term involving $F^{(3)}_{y, \mtt i, s}$.
This is where the width $\ell_1\ell_2\ell_3 = \kappa\, n^{1/2}$ is
used, and where the parameter $\kappa$ appears in the final bound.

\begin{lemma}
\label{s80c}
Under the hypotheses of Proposition \ref{s79},
\begin{align*}
\bb E_{\mu_n} \Big[\,\Big| \int_0^t
\frac{1}{n^{d/2}} \frac{1}{|B_{\mtt 0}|} 
\sum_{y\in B_{\mtt 0} } \sum_{\mtt i\in \Sigma^+_{L_4}}
|B_{3, \mtt i} |\, G (s, y+B_{3, \mtt i} )\,
F^{(3)}_{y, \mtt i, s} (\mss m_{y, 3, \mtt i}(s)) \, ds \, \Big|\,\Big]
\,\le\, \mtt R^{(\ref{s80c})}_n\,,
\end{align*}
where
$\cb{\mtt R^{(\ref{s80c})}_n} \,=\, \mtt S^{(\ref{120})}_{n,
\ell_1\ell_2\ell_3} \,+\, \mtt R^{ ( \ref{l02} ) }_{n,
\ell_1\ell_2\ell_3}$, these remainders being defined in the statement
of Lemma \ref{s80b}.  In particular, as
$\ell_1\ell_2\ell_3 = \kappa\, n^{1/2}$ and $d=3$,
\begin{equation*}
\mtt R^{(\ref{s80c})}_n \,\le\,
\frac{C_1}{\kappa^{3}} \, \int_0^t \|G(s)\|_\infty\, ds
\,+\, \frac{C_1}{n^{d/2}}
\int_0^t \|G(s) \|_\infty \, H_n(f_s \,|\, \nu^n_{u(s,\cdot)})\, ds \,.
\end{equation*}
\end{lemma}

\begin{proof}
The proof is the one of Lemma \ref{s80b}, with $\ell_1\ell_2\ell_3$ in
place of $\ell_1\ell_2$ and $L_4$ in place of $L_3$. 
\end{proof}

\begin{proof}[Proof of Theorem \ref{t02b} in $d=3$]
Fix $T>0$ and $0\le s<t\le T$. Recall the definition of the norm
$\Vert G\Vert_{T,2}$ introduced just before the statement of the
theorem and let $\widehat G = G / \Vert G\Vert_{T,2} $, so that
$\sup_{0\le r\le T} \Vert \widehat G (r)\Vert_{C^2} \le 1$. Rewrite
the expectation appearing in the statement of the theorem as in
\eqref{235}. 

Collecting Proposition \ref{s79}, Lemmata \ref{s85}, \ref{s82},
Conclusion \ref{s87b}, the decomposition \eqref{217}, the estimates
\eqref{218} and Lemmata \ref{s78c}, \ref{s81c}, \ref{s80c} with $G$
replaced by $\widehat G$, we obtain
\begin{align*}
& \bb E_{\mu_n}\Big[\, \Big| \, \int_s^t \frac{1}{n^{d/2}}
\sum_{x\in \bb T^d_n} \widehat G(r, x/n) \, \Pi_{u(r,\cdot),h}(x, \eta^n(r))
\, dr\, \Big|\, \Big]
\,\leq\,
\mtt R^{( \ref{s79} )}_n \,+\,  \mtt R^{( \ref{s85} )}_n
\,+\, \frac{3}{2}\, \mtt R^{( \ref{s82} )}_{n,d}(1/2)
\,+\,   \mtt R^{(\ref{s87b})}_n
\\
&\qquad
\,+\, \frac{2}{n^{d/2}} \, \int_s^t
\mtt S^{(\ref{135})}_{n,\ell_1\ell_2\ell_3} (\widehat G_r,\mtt p'_A, u_r)\, dr
\,+\, \mtt R^{(\ref{s07b})}_{n,\ell_1\ell_2\ell_3, \widehat G}
\,+\, \mtt R^{(\ref{s78c})}_n \,+\, \mtt R^{( \ref{s81c} )}_n
\,+\, \mtt R^{( \ref{s80c} )}_n \,,
\end{align*}
where, as in dimension $2$, in the formulae for the remainders the
time-integrals between $0$ and $t$ are replaced by integrals between
$s$ and $t$, and the entropy $H_n(f_0 \,|\, \nu^n_{u_0(\cdot)})$ by
$H_n(f_s \,|\, \nu^n_{u(s, \cdot)})$.  Note that the third scale
produced one further term
$(1/2)\, n^{-d/2}\, \bb E_{\mu_n}[\int_s^t \Psi_n(r)\, dr]$, coming
from Lemma \ref{s81c}, so that the total contribution of the term
$\Psi_n$ is at most $(3/2)\, \mtt R^{( \ref{s82} )}_{n,3}(1/2)$, by
Lemma \ref{s82}.  

It remains to examine each term on the right-hand side. Recall that
$\ell_1 = \upsilon_n n^{1/4}$,
$\ell_1\ell_2 = \upsilon_n^{2} n^{2/5}$,
$\ell_1\ell_2\ell_3 = \kappa\, n^{1/2}$, and that
$\upsilon_n^{-1}\le n^{3/100}$. 

The remainders $\mtt R^{( \ref{s79} )}_n$, $\mtt R^{( \ref{s85} )}_n$
and $\mtt R^{( \ref{s82} )}_{n,d}(1/2)$ have been estimated in Remarks
\ref{r04} and \ref{r08}. The remainder
$\mtt S^{(\ref{135})}_{n,\ell_1\ell_2\ell_3} (\widehat G_r,\mtt p'_A,
u_r)$ is easily bounded, as well as
$\mtt R^{(\ref{s07b})}_{n,\ell_1\ell_2\ell_3, \widehat G}$.

A meticulous inspection yields that
\begin{align*}
\mtt R^{(\ref{s87b})}_n \,\le\, \frac{1}{n^{d/2}} 
\Big\{ H_n(f_s\,|\,\nu^n_{u(s,\cdot)}) + \ln(2) \Big\}
\,+\, \frac{\mf c_{1,3}}{n^{d/2}} \int_s^t
\big\{\, H_n(f_r\,|\,\nu^n_{u(r, \cdot)})\,+\,
\upsilon^{2}_n\, n^{3/2}\,\big\}\,
dr\,,
\end{align*}
since $\upsilon^{- 1}_n =  n^{3/100}$. Similarly,
\begin{gather*}
\mtt R^{(\ref{s78c})}_n 
\,\le\, 
\frac{\mf c_{1,3}}{n^{d/2}} \int_s^t
\big\{\, H_n(f_r\,|\,\nu^n_{u(r, \cdot)})\,+\,
\kappa^{1/2}\, n^{5/4}\,\big\}\,
dr\,,
\\
\mtt R^{( \ref{s81c} )}_n
\,\le\, \frac{1}{n^{d/2}}\,
\big\{ \, H_n(f_s \,|\, \nu^n_{u_s, \cdot)}) \,+\, \ln (2)\,\big\}\,
\,+\, \mf c_{1,3} \, \int_s^t \big\{\, \upsilon_n^{2}
\,+\, \kappa^5 \, \upsilon_n^{-12}\, n^{-2/5}\,\big\} \, dr \,,
\\
\mtt R^{( \ref{s80c} )}_n  \,\le\,
C_1\, 
\int_s^t  \big\{\, \frac{1}{n^{d/2}} \, H_n(f_r \,|\, \nu^n_{u(r,\cdot)})\,
\,+\, \frac{1}{\kappa^{3}} \, \big\} \, dr \,.
\end{gather*}
To complete the proof of the theorem, it remains to recollect all
previous estimates.
\end{proof}

\section{Proof of Theorem \ref{t1}}
\label{sec11}

In this section, we prove the tightness of the sequence $\bb Q_{\mu_n}$,
and we characterize its limit points. We also prove the convergence of
the finite dimensional distributions of the process $X^n_t$, and prove
Theorem \ref{t1}. The constants introduced in this section are denoted
by $\mf c_{\rm f}$. They may depend on $T$, $d$ and $h_j$, $g_j$,
$1\le j\le d$, and on $u$ in $[0,T]\times \bb T^d$. If the context
requires us to make explicit the dependence of the constant
$\mf c_{\rm f}$ on some variable $\alpha$, we write
$\mf c_{\rm f} (\alpha)$. As before, the constant may change value
from line to line.

\subsection*{The linear equation}

We present in this subsection some properties of the solutions to a
linear partial differential equation needed in the proof of the
tightness of the sequence $\bb Q_{\mu_n}$.  Exceptionally, in this
subsection constants are represented by $\mf c^{(\rm k)}_j$, $j\ge 1$.

Throughout this subsection, $0<t\le T$ is fixed.  Recall from the
paragraph below \eqref{127} the definition of the H\"older spaces
$H^\beta = H^\beta(\bb T^d)$,
$H^{\beta/2,\beta} = H^{\beta/2,\beta} ([0,t]\times \bb T^d)$.

Fix a non-integer, positive $\beta$, and $F\in H^{\beta+2}$.  Recall
from \eqref{41} the definition of the linear operator $\mc A_r$,
$r\ge 0$. Theorem \ref{s65} below, a particular case of \cite[Theorem
IV.5.1]{lsu68}, establishes the existence and uniqueness of a solution
to the linear parabolic equation
\begin{equation}
\label{159b}
\left\{
\begin{aligned}
& \partial_s \mf v \,=\, \mc A_{t-s}  \mf v \,,
\quad 0\le s\le t\,, 
\\
& \mf v(0) \,=\, F\,.
\end{aligned}
\right.
\end{equation}

\begin{theorem}
\label{s65}
Recall that $0<t\le T$ is fixed. Fix a non-integer, real number
$\beta>0$, $F\in H^{\beta+2}$. Assume that the solution $u$ of the
hydrodynamic equation \eqref{127} belongs to
$H^{\beta/2,\beta} ([0,t]\times \bb T^d)$ and is bounded away from $0$
and $1$: there exists $0<\mf r<1/2$ such that
$\mf r\le u(s,x) \le 1-\mf r$ for all $(s,x) \in [0,T]\times \bb T^d$.
Then, there exists a unique classical solution
$\mf v \colon [0,t] \times \bb T^d \to \bb R$ in
$H^{(\beta+2)/2,\beta+2}$ to the Cauchy problem
\eqref{159b}. Moreover, there exists a finite constant
$\mf c^{(\rm k)}_{1}$, depending only on $T$, $\beta$ and $u$ in
$[0,t]\times \bb T^d$, such that
\begin{equation*}
\| \mf v \|_{(\beta +2)} \,\leq \, 
\mf c^{(\rm k)}_{1}\, \| F \|_{(\beta +2)} \,.
\end{equation*}
\end{theorem}

Mind that \cite[Theorem IV.5.1]{lsu68} requires the strict ellipticity
of the operator $\mc A_s$ which is guaranteed here from the fact that
$D_j(u)\ge c_0>0$ because $\mf r\le u(s,x) \le 1-\mf r$ for all
$(s,x) \in [0,T]\times \bb T^d$.

To recover, from the estimate stated in the theorem, bounds with $\beta$
replaced by integers, let $p$ be an odd integer, $p=2q-1$, $q\ge
1$. Assume that $F$ belongs to $C^{p+2}(\bb T^d)$ and $u$ to
$C^{(p+1)/2,p} ([0,t]\times \bb T^d)$. Fix a non-integer
$\beta \in (p-1, p)$.  It follows from the previous result that
\begin{equation}
\label{220}
C(p,\beta) \, \| \mf v \|_{C^{(p+1)/2,p+1}} \le \| \mf v \|_{(\beta +
2)}\,\leq \,  \mf c^{(\rm k)}_{1} \, \| F \|_{(\beta+2)} \le
\mf c^{(\rm k)}_{1} (\beta, p)
\, \| F \|_{C^{p+2}}\,.
\end{equation}

We will also need to compare the solutions for different time
intervals. Fix $0<t_1<t_2$, a non-integer, real number $\beta>0$, and
$F\in H^{\beta+2}$. Assume that $u$ belongs to
$H^{\beta/2,\beta} ([0,t_2]\times \bb T^d)$. Denote by $\mf v_i$, $i=1$,
$2$, the solutions to equation \eqref{159b} with $t=t_i$.  By
definition, on the time interval $[0,t_1]$,
$\mf w = \mf v_2 - \mf v_1$ solves the Cauchy problem
\begin{equation*}
\left\{
\begin{aligned}
&  \partial_s \mf w \, =\, \mc A_{t_2-s} \mf w \, +\, \mf f \,,
\quad 0\le s\le t_1\,, 
\\
& \mf w(0) \,=\, 0 \,.
\end{aligned}
\right.
\quad \text{where} \;\;
\mf f (s, \cdot) \, =\, [\,  (\mc A_{t_2-s} - \mc A_{t_1-s}) \, \mf
v_1\,] (s, \cdot) \,.
\end{equation*}
By \cite[Theorem IV.5.1]{lsu68},  there exists a finite constant
$\mf c^{(\rm k)}_{2}$, depending only on $T$, $\beta$ and $u$ in
$[0,t_2]\times \bb T^d$, such that
\begin{equation*}
\| \mf w \|_{(\beta +2)} \,\leq \, 
\mf c^{(\rm k)}_{2} \, \| \mf f \|_{(\beta)} \,.
\end{equation*}

Let $p$ be an odd integer, $p=2q-1$, $q\ge 1$, and set $\beta=p-1/2$.
By the previous estimate,
\begin{equation*}
C(p)\, \| \mf w \|_{C^{(p+1)/2,p+1}} \,\leq \, 
\| \mf w \|_{(\beta +2)} \,\leq \, 
\mf c^{(\rm k)}_{2} (p)  \, \| \mf f \|_{(\beta)}
\,\leq \, 
\mf c^{(\rm k)}_{2} (p) \, \| \mf f \|_{C^{(p+1)/2,p}} \,.
\end{equation*}
Here, we lost the ratio $(r/2,r)$ between time and space regularity 
because, as $p$ is odd and $\beta = p - 1/2$, the first integer larger
than $\beta/2$ is $(p+1)/2$ and the first integer larger than $\beta$
is $p$.

As $D_j (\cdot)$, $1\le j\le d$, are polynomials, an
elementary computation yields that for all $k\ge 0$, $\ell\ge 0$,
\begin{equation}
\label{244}
\| \mf f \|_{C^{k, \ell}}
\,\leq \, 
\mf c^{(\rm k)}_{3} 
\big(\, k, \ell, \| u  \|_{C^{k+1,\ell}} \,\big) 
\, \| \mf v_1  \|_{C^{k,\ell+2}} \, (t_2-t_1) \,,
\end{equation}
where, as indicated,
$\mf c^{(\rm k)}_{3} \big(\, k, \ell, \| u \|_{C^{k+1,\ell}} \,\big)$
is a finite constant which depends on $k$, $\ell$ and
$\| u \|_{C^{k+1,\ell}}$.  All these norms are computed on
$[0,t_2] \times \bb T^d$.

By Theorem
\ref{s65}, for $\beta' = p+2 + (1/2)$, if $F$ belongs to $H^{\beta'}$
and $u$ to $H^{(\beta'-2)/2, \beta'-2} ([0,t_1]\times \bb T^d)$,
\begin{equation*}
C(p)\, \| \mf v_1 \|_{C^{(p+1)/2,p+2} } \,\le\,  \| \mf v_1 \|_{(\beta') }
\,\le\,  \mf c ^{(\rm k)}_{1} \, \| F \|_{(\beta') } \,\le\,
\mf c ^{(\rm k)}_{1}  \, \| F \|_{C^{p+3} }\,.
\end{equation*}

Recollecting the previous estimates yields the following result.  Fix
$0<t_1<t_2$, an odd integer $p$, $p=2q-1$ for some $q\ge 1$.  Assume
that $F\in C^{p+3} (\bb T^d)$,
$u\in C^{(p+3)/2,p+1} ([0,t_2]\times \bb T^d)$.  Set in the previous
argument $k=(p+1)/2$, $\ell=p$. Then, there exists a finite constant
$\mf c ^{(\rm k)}_{4}$ depending only on $T$, $p$ and $u$ in
$[0,t_2]\times \bb T^d$ such that
\begin{equation}
\label{162}
\| \mf w \|_{C^{(p+1)/2,p+1} } \,\leq \, 
\mf c ^{(\rm k)}_{4} \, \| F  \|_{C^{p+3}} \, (t_2-t_1)\,.
\end{equation}
The time-regularity imposed on $u$ is required in \eqref{244}, while
the space regularity comes from the hypothesis that $u$ belongs to
$H^{(\beta'-2)/2, \beta'-2} ([0,t_1]\times \bb T^d)$, imposed in the
penultimate paragraph.

Denote by $\cb{(P_{s,t}: 0\le s \le t)}$ the semigroup associated to
the linear operators $(\mc A_r: r\ge 0)$ introduced in \eqref{41}.
With this notation,
\begin{equation}
\label{220c}
\mf v(s) \,=\,  P_{t-s,t} F\,.
\end{equation}
In terms of the semigroup, Theorem \ref{s65} states that for each
non-integer $\beta>0$, there exists a constant
$\mf c^{({\rm k})}_{1}$, depending on $T$, $\beta$ and
$\Vert u\Vert_{(\beta)}$, such that
\begin{equation*}
\Vert \,  P_{t- \,\cdot,\, t}\, F \,  \Vert_{(2+\beta)}
\,=\, \Vert \,  \mf v  \,  \Vert_{(2+\beta)}
\;\le\; \mf c^{({\rm k})}_{1} \,  \Vert F \Vert_{(2+\beta)}\,.
\end{equation*}
Here, $P_{t- \,\cdot,\, t}\, F =\mf v$ is a function of $s$ and $x$,
in contrast with $F$ or $P_{t-s,t}\, F$ which are functions of $x$
only.  Since
$\Vert \, P_{t-s,t}\, F \, \Vert_{(2+\beta)} \le
\Vert \, P_{t- \,\cdot,\, t}\, F \, \Vert_{(2+\beta)}$ for all
$0\le s \le t$, 
\begin{equation}
\label{220d}
\sup_{0\le s \le t} \Vert \, P_{s,t}\, F \, \Vert_{(2+\beta)}
\,=\,
\sup_{0\le s \le t} \Vert \, P_{t-s,t}\, F \, \Vert_{(2+\beta)}
\;\le\;
\mf c^{({\rm k})}_{1} \,  \Vert F \Vert_{(2+\beta)}\,.
\end{equation}

On the other hand, since
$\Vert \mf v(t) \Vert_{C^{p+1}} \le \Vert \mf v \Vert_{C^{(p+1)/2,
p+1}}$, by \eqref{220} applied with the odd integer $p=q-1$, for any
even integer $q\ge 2$, there exists a finite constant $C(q)$
such that
\begin{equation}
\label{220b}
\Vert P_{0,t} F \Vert_{C^{q}}  \,=\, 
\Vert \mf v(t) \Vert_{C^{q}}
\, \le\,  \Vert \mf v \Vert_{C^{q/2, q}}
\,\le\, C(q) \, \Vert F \Vert_{C^{q+1}} \,,
\end{equation}
provided $u$ belongs to $C^{q/2, q-1}([0,t]\times \bb T^d)$.

Fix $0<t_1<t_2$, an odd integer $p$, $p=2q-1$ for some $q\ge 1$.
Assume that $F\in C^{p+3} (\bb T^d)$ and that $u$ belongs to
$C^{(p+3)/2, p+1} ([0,T]\times \bb T^d)$.  With the notation
introduced above,
\begin{equation*}
\Vert P_{r,t_2} F - P_{r,t_1} F\Vert_{C^{p+1}} \,=\,
\Vert \mf v_2 (t_2-r)  -   \mf v_1 (t_1-r)  \Vert_{C^{p+1}} 
\end{equation*}
for all $0\le r\le t_1$.  By the triangle inequality, \eqref{220} and
\eqref{162}, since $p\ge 1$, this expression is less than or equal to
\begin{align*}
& \Vert \mf v_2 (t_2-r)  -   \mf v_2 (t_1-r)  \Vert_{C^{p+1}}  \,+\,
\Vert \mf v_2 (t_1-r)  -   \mf v_1 (t_1-r)  \Vert_{C^{p+1}}
\\
&\quad \le
\Vert \mf v_2 \Vert_{C^{1,p+1}}  \, (t_2-t_1) \,+\,
\Vert \mf w \Vert_{C^{(p+1)/2,p+1}}
\le \mf c^{(\rm k)}_{5} (p)
\, \| F \|_{C^{p+3}}  \, (t_2-t_1)  \,.
\end{align*}
Therefore, for every even integer $p'\ge 2$ ($p'=p+1$), there exists a
constant $\mf c^{(\rm k)}_{5} (p')$ such that
\begin{equation}
\label{221}
\sup_{0\le r\le t_1} \Vert P_{r,t_2} F - P_{r,t_1} F\Vert_{C^{p'}}
\,\le\, \mf c^{(\rm k)}_{5} (p')\, 
\| F \|_{C^{p'+2}}  \, (t_2-t_1) 
\end{equation}
for all $F\in C^{p'+2} (\bb T^d)$.

\subsection*{Tightness in $L^1(0,T;\mc H_{-\mtt r})$}
\label{sec6b}

In this subsection, we prove the tightness of the sequence
$\bb Q_{\mu_n} = \bb P_{\mu_n}\circ (X^n)^{-1}$ in
$L^1(0,T; \mc H_{-\mtt r})$, for $\mtt r$ large enough.

The proof requires the following elementary bounds. We claim that
there exist positive constants $\mf c_{\rm f} = \mf c_{\rm f}(p, \mtt
q)$, $\mf c'_{\rm f} = \mf c'_{\rm
f} (p, \mtt q)$ such that
\begin{equation}
\label{243}
\Vert G \Vert_{\mc H_{\mtt q}} \;\le\; \mf c_{\rm f} \,
\Vert G \Vert_{C^{p}} \;\;\text{ if } p > \mtt q + \tfrac d2\;,
\qquad
\Vert G \Vert_{C^{p}} \;\le\; \mf c'_{\rm f} \,
\Vert G \Vert_{\mc H_{\mtt q}}\;\;\text{ if } \mtt q > p + \tfrac d2\;.
\end{equation}

To prove the first assertion, integrate by parts $p$ times to get that
$|\<G,\phi_m\>| \le C_0(p) (1+\Vert m\Vert)^{-p}\Vert G\Vert_{C^p}$.
Thus,
$\Vert G\Vert^2_{\mc H_{\mtt q}} \le C_0(p)\, \Vert G\Vert^2_{C^p}
\sum_m (1+\Vert m\Vert)^{2\mtt q - 2p}$. This sum is finite if
$2p-2\mtt q>d$.  The second inequality is the Sobolev embedding
theorem.
\smallskip

Recall the decomposition \eqref{145b} of $X^n_t(\cdot)$.  All terms on
the right-hand side are bounded by $C_0 n^{d/2}\Vert F\Vert_{C^4}$,
uniformly in $\eta^n(\cdot)$. Hence, by \eqref{243}, they belong to
$\mc H_{-\mtt r}$ for $\mtt r> 4 + (d/2)$. Since they are continuous,
with the exception of $M^n$ which is c\`adl\`ag, they are Bochner
measurable and belong to $L^1(0,T;\mc H_{-\mtt r})$ almost
surely. The main result of this subsection reads as follows.

\begin{proposition}
\label{s77}
Let $d\leq 3$.  Assume that $\mu_n$ is a sequence of probability
measures on $\Omega_n$ satisfying
$H_n(\mu_n \,|\, \nu^n_{u_0(\cdot)}) = o(n^{d/2})$. Suppose that there
exists $\mtt r_0 > d/2$ such that the sequence of laws of $X^n_0$ is
tight in $\mc H_{-\mtt r_0}$. Suppose, finally, that the solution $u$
of \eqref{127} satisfies $\mf r \le u \le 1 - \mf r$ and belongs to
$H^{\beta/2,\beta}([0,T]\times \bb T^d)$ for some non-integer
$\beta > \mtt r_0 + (d/2) + 4$. Then, the sequence of probability
measures $\bb Q_{\mu_n}$ is tight in $L^1(0,T;\mc H_{-\mtt r})$ for
$\mtt r \;>\; \mtt r^\star $, where $\mtt r^\star$ has been introduced
in \eqref{242b}.
\end{proposition}

Recall the decomposition \eqref{145b} of $X^n_t$. The proof of
Proposition \ref{s77} consists in showing that each term on the
right-hand side is tight in $L^1(0,T; \mc H_{-\mtt r})$, for
$\mtt r> \mtt r^\star$.  Define the $\mc H_{-\mtt r}$-valued processes by
\begin{equation*}
\begin{gathered}
\mf B^n_t \,:=\, \int_0^t \mtt B^{(n)}_s \, ds \;,
\quad
\mf D^n_t \,:=\, \int_0^t \mtt D^{(n)}_s \, ds \;, \quad
\mf J^n_t \,:=\, \int_0^t \mc A^*_s X^n_s \, ds \;.
\end{gathered}
\end{equation*}
The last term acts on a smooth function $F\in \mc H_{\mtt r}$ as
$\mf J^n_t (F) = \int_0^t X^n_s (\mc A_s F)\, ds$.  We also recall, as
$L^1(0,T;\mc H_{-\mtt r})$ is a complete, separable space, that if a
sequence of processes in $L^1(0,T;\mc H_{-\mtt r})$ converges weakly
then it is tight by Prokhorov's theorem.

\smallskip\noindent{\bf The process $X^n_0$ is tight in
$L^1(0,T;\mc H_{-\mtt r})$ for every $\mtt r \ge \mtt r_0$.}  The map
$\Phi \colon \mc H_{-\mtt r_0} \to L^1(0,T;\mc H_{-\mtt r})$ defined
by $\Phi(\mtt x) (t) = \mtt x$ for all $0\le t\le T$, is linear and
bounded:
$\Vert \Phi (\mtt x)\Vert_{L^1(0,T;\mc H_{-\mtt r})} = T \Vert \mtt x
\Vert_{\mc H_{-\mtt r}} \le T \Vert \mtt x\Vert_{\mc H_{-\mtt r_0}}$,
because $\mtt r\ge \mtt r_0$.

By hypothesis, the sequence of laws of $X^n_0$ is tight in
$\mc H_{-\mtt r_0}$. As $\Phi$ is bounded and continuous,
$\Phi(X^n_0)$ is tight in $L^1(0,T;\mc H_{-\mtt r})$, as claimed.

\smallskip\noindent{\bf The martingale $M^n$ is tight in
$L^1(0,T;\mc H_{-\mtt r})$ for $\mtt r > 1 +  (3d+1)/2$.}
We verify the two conditions of Corollary \ref{s67}.  Recall from
\eqref{144}, \eqref{143} that
$M^n_t(F)^2 - \int_0^t \Gamma^n(\eta^n(s),F)\, ds$ is a martingale,
and that $|\Gamma^n(\eta,F)| \le C_0 \Vert \nabla F\Vert^2_\infty$
uniformly in the configuration.  Consequently, for $0\le s<t\le T$,
\begin{equation*}
\bb E_{\mu_n} \big[\, \big(M^n_t(F) - M^n_s(F)\big)^2 \,\big]
\;=\; \bb E_{\mu_n}\Big[ \int_s^t \Gamma^n(\eta^n(v),F)\, dv \Big]
\;\le\; C_0\, (t-s)\, \Vert \nabla F\Vert^2_\infty \;.
\end{equation*}
This estimate applied to the real and the imaginary parts of $\phi_m$,
whose gradients are bounded by $2\pi \Vert m \Vert$, and 
Schwarz inequality yield that
\begin{equation*}
\bb E_{\mu_n}\big[\,\big|  \, M^n_t(\phi_m) - M^n_s(\phi_m)
\,\big| \,\big] \;\le\; \mf c_{\rm f} \, (1+\Vert m \Vert)\, |t-s|^{1/2}
\end{equation*}
for all $n\ge 1$, which is condition (B) of Corollary \ref{s67} with
$b=1/2$ and $\beta_m = \mf c_{\rm f} \, (1+\Vert m\Vert)$. The sum
\eqref{239} converges since $\mtt r > 1 + (3d+1)/2$.

Taking $s=0$ and integrating in $t$,
\begin{equation*}
\bb E_{\mu_n}\Big[ \int_0^T \big| \, M^n_t(\phi_m)  \,\big| \, dt
\Big] \;\le\; \mf c_{\rm f} \,   T^{3/2} \, (1+\Vert m\Vert)\;,
\end{equation*}
so that condition (A) of Corollary \ref{s67} holds for
$1 + (3d+1)/2< \mtt r' < \mtt r$.

\smallskip\noindent{\bf The process $\mf B^n$ converges to $0$ in
$L^1(0,T;\mc H_{-\mtt r})$, in $L^1(\bb P_{\mu_n})$, for
$\mtt r > 4 + (3d+1)/2$.}
Recall from \eqref{241} the definition of $\mtt B^{(n)}_r$, in which the
discrete Laplacian $\Delta^n_j$ appears.  For $H\in C^4(\bb T^d)$, let
\begin{equation*}
\cb{\widehat {\mtt B}^{(n)}_r(H)} \; :=\; \frac{1}{n^{d/2}} \sum_{j=1}^d
\sum_{x\in\bb T_n^d} (\partial^2_{x_j} H)  (x/n)\, 
\Pi_{u(r, \cdot), h_j}  (x, \eta^n(r))  \,.
\end{equation*}
By a fourth order Taylor expansion,
$| (\Delta^n_j H) - (\partial^2_{x_j} H )| \le C_0 \Vert H\Vert_{C^4}
n^{-2}$, while $|\Pi_{u(r, \cdot), h_j}| \le C_0$.  Hence,
\begin{equation}
\label{230}
\big|\, \mtt B^{(n)}_r(H) \,-\, \widehat {\mtt B}^{(n)}_r(H) \,\big| \;\le\;
C_0 \, \Vert H\Vert_{C^4}\ n^{(d/2)-2}
\end{equation}
for all $H\in C^4(\bb T^d)$, $0\le r\le T$ and $n\ge 1$.  As $d\le 3$,
the right-hand side vanishes as $n\to\infty$.

Fix $m\in \bb Z^d$ and $0\le s<t\le T$.  By Corollary \ref{s94} with
$G = \partial^2_{x_j} \phi_m$, and summing over $1\le j\le d$, 
\begin{equation*}
\bb E_{\mu_n}\Big[\, \Big| \int_s^t \widehat{\mtt B}^{(n)}_r (\phi_m)\, dr\,
\Big| \, \Big] \;\le\; \mf c_{\rm f} \, (1+\Vert m\Vert)^4 \,
\mtt h_n\, \big\{\, (t-s) \,+\, 1 \,\big\}
\end{equation*}
for all $n\ge 1$.  Adding this bound to the one obtained by
integrating \eqref{230} over $[s,t]$,
\begin{equation*}
\bb E_{\mu_n}\Big[\, \big| \< \mf B^n_t - \mf B^n_s , \phi_m\>
\big| \, \Big] \;\le\;
\mf c_{\rm f} \,  (1+\Vert m\Vert)^4 \, \varrho_n\,
\big\{\,  (t-s) \,+\, 1  \,\big\}\;,
\quad
\varrho_n \,:=\, \mtt h_n + n^{(d/2)-2}\;.
\end{equation*}

By Corollary \ref{c-t01}, the definition of the sequence $\kappa_n$,
and since $d\le 3$, $\varrho_n \to 0$, so that
$\sup_{n\ge 1} \varrho_n < \infty$.  Hence, condition (B) of Corollary
\ref{s90} holds with
$\beta_m = \mf c_{\rm f} \, \sup_{k\ge 1} \varrho_k \, (1+\Vert
m\Vert)^4$, $b_n = \varrho_n/\sup_{k \ge 1} \varrho_k$, and
$\mtt r > 4 + (3d+1)/2$.  Condition \eqref{158} follows from the
previous displayed equation provided we choose
$4 + (3d+1)/2 < \mtt r' <\mtt r$. This proves that $\mf B^n$ is tight
in $L^1(0,T;\mc H_{-\mtt r})$.

Actually, in view of \eqref{236}, the previous estimate with $s=0$
gives the stronger statement
\begin{equation*}
\bb E_{\mu_n}\Big[\, \Vert \mf B^n \Vert_{L^1(0,T;\mc H_{-\mtt r})}
\,\Big]
\;\le\; \mf c_{\rm f} \, T\, (T+ 1)\, \varrho_n \sum_{m\in \bb Z^d}
(1+\Vert m\Vert)^{4}\, \gamma_m^{-\mtt r/2}\, \omega_m^{-1/2} \;,
\end{equation*}
which vanishes as $n\to\infty$, since  $\mtt r > 4 + (3d+1)/2$.

\smallskip\noindent{\bf The process $\mf D^n$ converges to $0$ in
$L^1(0,T;\mc H_{-\mtt r})$ almost surely, provided
$\mtt r > 4 + (3d+1)/2$.}
By \eqref{240},
$|\< \mf D^n_t, \phi_m\>| \le \mf c_{\rm f} \, T \, \Vert \phi_m
\Vert_{C^4}\, n^{(d/2)-2} \le \mf c_{\rm f} \, T \, (1+\Vert m\Vert)^4
\, n^{(d/2)-2}$, uniformly in $t\le T$ and in the trajectory.  Hence,
by \eqref{236}, 
\begin{equation*}
\Vert \mf D^n \Vert_{L^1(0,T;\mc H_{-\mtt r})} \;\le\;
\mf c_{\rm f} \, T^2\, n^{(d/2)-2} \sum_{m\in \bb Z^d}
(1+\Vert m\Vert)^{4}\, \gamma_m^{-\mtt r/2}\, \omega_m^{-1/2}\;.
\end{equation*}
This sum is finite for $\mtt r > 4 + (3d+1)/2$, and the expression
vanishes as $n\to\infty$ in $d\le 3$.

\smallskip\noindent {\bf The term $\mf J^n$.}  Fix
$F \in C^\infty(\bb T^d)$ and $0\le s \le T$.  Recall from
\eqref{220c} that we denote by $(P_{s,t} : 0\le s\le t)$ the semigroup
associated to $\mc A_s$.  By \eqref{146} with $G = \mc A_s F$,
\begin{equation*}
\mf J^n_t(F) \;=\; \mf J^{1,n}_t(F) \,+\, \mf J^{2,n}_t(F)
\,+\, \mf J^{3,n}_t(F) \,+\, \mf J^{4,n}_t(F)\;,
\end{equation*}
where
\begin{gather*}
\cb{\mf J^{1,n}_t(F) } \,:=\, \int_0^t X^n_0\big(P_{0,s}\, \mc A_s
F\big)\, ds \;, \qquad
\cb{\mf J^{2,n}_t(F) } \,:=\, \int_0^t M^{n,s}_s\big(\mc A_s F\big)\, ds\;,
\\
\cb{\mf J^{3,n}_t(F) } \,:=\, \int_0^t \!\! ds \int_0^s \mtt B^{(n)}_v\big(P_{v,s}\,
\mc A_s F\big)\, dv \;, \qquad
\cb{\mf J^{4,n}_t(F)} \,:=\, \int_0^t \!\! ds \int_0^s \mtt D^{(n)}_v
\big(P_{v,s}\, \mc A_s F\big)\, dv \;.
\end{gather*}

\smallskip\noindent{\bf The process $\mf J^{1,n}$ is tight in
$L^1(0,T;\mc H_{-\mtt r})$ for $\mtt r > \mtt r_0 + d + 4$.}
Since $\partial_s P_{s,t} G = -\, \mc A_s P_{s,t}G$ and
$P_{r,t} = P_{r,s}P_{s,t}$, $r \le s \le t$,
$\partial_t P_{r,t} G = P_{r,t}\, \mc A_t G$, so that
\begin{equation*}
\int_r^t P_{r,s}\, \mc A_s G \; ds \;=\; P_{r,t} G \,-\, G \;,
\qquad 0\le r \le t \le T\;, \;\; G \in C^\infty(\bb T^d)\;,
\end{equation*}
and $\mf J^{1,n}_t = P^*_{0,t}X^n_0 - X^n_0$.

Let $\Phi_1 \colon \mc H_{-\mtt r_0} \to L^1(0,T;\mc H_{-\mtt r})$ be
given by $\Phi_1 (\mtt x) (t) := P^*_{0,t} \mtt x - \mtt x$. We claim
that there exists a finite constant $\mf c_{\rm f} $ such that
\begin{equation}
\label{238}
\big\Vert \Phi_1 (\mtt x) (t) \big\Vert_{\mc H_{-\mtt r}}\,=\,
\big\Vert P^*_{0,t} \mtt x - \mtt x \big\Vert_{\mc H_{-\mtt r}}
\;\le\; \mf c_{\rm f}
\, \Vert \mtt x \Vert_{\mc H_{-\mtt r_0}}
\end{equation}
for all $0\le t\le T$, $ \mtt x \in \mc H_{-\mtt r_0}$.

To prove this claim, let $p$ be the smallest even integer larger than
$\mtt r_0 + (d/2)$; thus
$\mtt r_0 + (d/2) < p \le \mtt r_0 + (d/2) + 2$. By \eqref{242b},
$\mtt r > p + 2 + (d/2)$.  By definition of the $\mc H_{-\mtt r}$
norm,
\begin{equation*}
\big\Vert P^*_{0,t} \mtt x - \mtt x \big\Vert_{\mc H_{-\mtt r}}
\;=\; \sup_{\Vert F \Vert_{\mc H_{\mtt r}} \le 1}
\big| \, \mtt x\big(P_{0,t}F - F\big)\, \big|
\;\le\; \Vert \mtt x \Vert_{\mc H_{-\mtt r_0}}
\sup_{\Vert F \Vert_{\mc H_{\mtt r}} \le 1}
\Vert  P_{0,t}F - F \Vert_{\mc H_{\mtt r_0}}\,.
\end{equation*}
Since $\mtt r_0 + (d/2) < p$, by \eqref{243}, the previous expression
is less than or equal to
\begin{equation*}
\mf c_{\rm f}
\, \Vert \mtt x \Vert_{\mc H_{-\mtt r_0}}
\, \sup_{\Vert F \Vert_{\mc H_{\mtt r}} \le 1}
\Vert P_{0,t}F - F \Vert_{C^{p}}\,.
\end{equation*}
By \eqref{220b} and \eqref{243}, as $p + 1+ (d/2) < \mtt r$,
$\Vert P_{0,t}F \Vert_{C^{p}} \le \mf c_{\rm f} \, \Vert F
\Vert_{C^{p+1}} \le \mf c_{\rm f} \, \Vert F \Vert_{\mc H_{\mtt
r}}$. As the supremum is carried over all functions $F$ such that
$\Vert F \Vert_{\mc H_{\mtt r}} \le 1$, to complete the proof of
\eqref{238}, it remains to recollect all previous estimates.

Since $p + 2 + (d/2) < \mtt r$, the same computation, with \eqref{221}
in place of \eqref{220b}, gives
\begin{equation*}
\big\Vert \Phi_1(\mtt x)(t) - \Phi_1(\mtt x)(s) \big\Vert_{\mc H_{-\mtt r}}
\;=\; \big\Vert P^*_{0,t} \mtt x - P^*_{0,s} \mtt x \big\Vert_{\mc H_{-\mtt r}}
\;\le\; \mf c_{\rm f}\, \Vert \mtt x\Vert_{\mc H_{-\mtt r_0}}\, (t-s)\;.
\end{equation*}
This step, which relies on \eqref{221}, requires $u$ to belong to
$C^{(p+2)/2,p} ([0,T]\times \bb T^d)$. As
$p+2 \le \mtt r_0 + (d/2) + 4$, this explains the requirement in the
hypotheses that $u\in H^{\beta/2,\beta}$ for
$\beta > \mtt r_0 + (d/2) + 4$.

Thus, $t\mapsto \Phi_1(\mtt x)(t)$ is Lipschitz continuous with values
in $\mc H_{-\mtt r}$. In particular, it is Bochner integrable and
bounded.  Therefore, $\Phi_1$ is a bounded linear operator with
$\Vert \Phi_1( \mtt x)\Vert_{L^1(0,T;\mc H_{-\mtt r})} \le \mf c_{\rm
f}\, T\, \Vert \mtt x\Vert_{\mc H_{-\mtt r_0}}$.  Since $X^n_0$ is tight in
$\mc H_{-\mtt r_0}$ and $\Phi_1$ is continuous,
$\mf J^{1,n} = \Phi_1(X^n_0)$ is tight in $L^1(0,T;\mc H_{-\mtt r})$.

\smallskip\noindent{\bf The process $\mf J^{2,n}$ is tight in
$L^1(0,T;\mc H_{-\mtt r})$ provided $\mtt r > 4 +  (3d+1)/2$.}
Fix a function $G\in C^\infty(\bb T^d)$, and recall from \eqref{51b}
the definition of the martingale $M^{n,s}_v(G)$, $0\le v\le s$. Its
quadratic variation is given by
\begin{equation*}
\int_0^v
\frac{1}{n^{d} } \sum_{j=1}^d \sum_{x\in \bb T^d_n}
\tau_x g_j (\eta^n(w))\,
[\, (\nabla^n_j P_{w,s} G) (x/n)\,]^2 \, dw\,,
\end{equation*}
where $g_j$ is the cylinder function introduced just above
\eqref{126}.  Therefore, by Schwarz inequality,
\begin{equation*}
\bb E_{\mu_n}\big[\, \big| M^{n,s}_s(G) \, \big| \,\big]
\;\le\; \Big( \mf c_{\rm f} \, s \, \sup_{0\le w\le s}
\Vert \nabla P_{w,s} G \Vert^2_\infty \Big)^{1/2}
\;\le\; \mf c_{\rm f} \, \sup_{0\le w\le s} \Vert P_{w,s}G
\Vert_{C^{1}} \;.
\end{equation*}

Fix $0<\beta_0<1$.  By \eqref{220c}, $\mf v(w) = P_{s-w,s} G$,
$0\le w\le s$, is the solution of equation \eqref{159b} with initial
condition $\mf v(0)=G$. Hence, by \eqref{220d} and by definition
of $\mc A_s$, for $G= \mc A_s \phi_m$,
\begin{equation*}
\sup_{0\le w \le s} \Vert P_{w,s}\, \mc A_s\phi_m \Vert_{C^{2}}
\;\le\;
\sup_{0\le w \le s} \Vert P_{w,s}\, \mc A_s\phi_m \Vert_{(2+\beta_0)}
\;\le\; \mf c_{\rm f}\,  \Vert \mc A_s \phi_m \Vert_{(2+\beta_0)}
\;\le\; \mf c_{\rm f} \, (1+\Vert m \Vert)^{4+\beta_0}\;.
\end{equation*}

Thus, for $0\le s<t\le T$,
\begin{equation*}
\bb E_{\mu_n}\big[\, \big| \< \mf J^{2,n}_t - \mf J^{2,n}_s,
\phi_m\> \big| \,\big]
\;\le\; \int_s^t \bb E_{\mu_n}\big[\, \big| M^{n,v}_v(\mc A_v\phi_m)
\big| \,\big]\, dv
\;\le\; \mf c_{\rm f} \, (1+\Vert m\Vert)^{4+\beta_0}\, (t-s)\;.
\end{equation*}
This is condition (B) of Corollary \ref{s67} with $b=1$ and
$\beta_m = \mf c_{\rm f} \, (1+\Vert m\Vert)^{4+\beta_0}$.

We turn to Condition (A) of this result. The left-hand side of the
previous equation with $s=0$, integrated over $[0,T]$, is bounded by
$\mf c_{\rm f} \, T^2\, (1+\Vert m\Vert)^{4+\beta_0}$.  Both sums in
Corollary \ref{s67} converge provided
$4 + \beta_0 + (3d+1)/2 < \mtt r' < \mtt r$. Since, by hypothesis,
$\mtt r > 4 + (3d+1)/2$, it is possible to choose $\beta_0>0$ small
enough for such an $\mtt r'$ to exist.

\smallskip\noindent{\bf The process $\mf J^{3,n}$ converges to $0$ in
$L^1(0,T;\mc H_{-\mtt r})$, in $L^1(\bb P_{\mu_n})$, if
$\mtt r > 6 + (3d+1)/2$.}
Fix $0<\beta_0<1$, $m \in \bb Z^d$ and $0\le v\le s\le T$. As above,
by \eqref{220d} and by definition of $\mc A_s$,
\begin{align}
\label{245}
\sup_{0\le v\le s} \big\Vert P_{v,s}\, \mc A_s \phi_m
\big\Vert_{C^{4}}
\;\le\;\sup_{0\le v\le s} \big\Vert P_{v,s}\, \mc A_s \phi_m
\big\Vert_{(4+\beta_0)}
\;\le\; \mf c_{\rm f}\,  \Vert \mc A_s \phi_m \Vert_{(4+\beta_0)}
\;\le\; \mf c_{\rm f} \, (1+\Vert m \Vert)^{6+\beta_0}\;.
\end{align}
Recall from \eqref{230} the definition of
$\widehat{\mtt B}^{(n)}_v(H)$.  In the formula for $\mf J^{3,n}$
replace $\mtt B^{(n)}_v(P_{v,s}\, \mc A_s \phi_m)$ by
$\widehat{\mtt B}^{(n)}_v(P_{v,s}\, \mc A_s \phi_m)$ at a cost bounded
by
$C_0\, \sup_{0\le v\le s} \big\Vert P_{v,s}\, \mc A_s \phi_m
\big\Vert_{C^{4}} n^{-2}$. By the previous equation, this expression
is less than or equal to
$\mf c_{\rm f} \, (1+\Vert m \Vert)^{6+\beta_0} \,n^{-2}$. Summing
over $x$, dividing by $n^{d/2}$ yields an
additional factor $n^{d/2}$. Thus, the cost of this replacement in
$\int_0^s \mtt B^{(n)}_v(P_{v,s}\, \mc A_s\phi_m)\, dv$ is bounded by
$\mf c_{\rm f}\, (1+T)\, (1+\Vert m \Vert)^{6+\beta_0}\, n^{(d/2)-2}$.

Apply Corollary \ref{s94} to
$\widehat{\mtt B}^{(n)}_v(P_{v,s}\, \mc A_s \phi_m)$ and add the
previous bound to get that
\begin{equation*}
\bb E_{\mu_n} \Big[\, \Big| \int_0^s \mtt B^{(n)}_v\big(P_{v,s}\mc A_s
\phi_m\big)\, dv \, \Big| \, \Big]
\;\le\; \mf c_{\rm f} \,  (1+\Vert m\Vert)^{6+\beta_0}\,
(1+T)^2\, \varrho_n
\end{equation*}
for every $m\in\bb Z^d$, $0\le s\le T$. In this formula,
$(\varrho_n:n \ge 1)$ is a sequence such that $\varrho_n \to 0$.
Therefore,
\begin{equation*}
\bb E_{\mu_n}\big[\, \big| \< \mf J^{3,n}_t ,
\phi_m\> \big| \,\big]
\;\le\; \mf c_{\rm f} \,   (1+\Vert m\Vert)^{6+\beta_0}\,
(1+T)^2\, \varrho_n\, t
\end{equation*}
for all $0\le t\le T$.

Since, by hypothesis, $\mtt r > 6 + (3d+1)/2$, there exists
$\beta_0>0$ such that $\mtt r> 6+\beta_0 + (3d+1)/2$. Hence, by
\eqref{236} and the previous estimate,
\begin{equation*}
\bb E_{\mu_n}\big[\, \big\Vert \mf J^{3,n}_t\big\Vert_{\mc H_{- \mtt r}} \,\big]
\;\le\; \mf c_{\rm f} \,  (1+T)^2\,  \varrho_n\, t
\end{equation*}
for all $0\le t\le T$. Thus,
$\bb E_{\mu_n}[\, \Vert \mf J^{3,n}\Vert_{L^1(0,T;\mc H_{-\mtt r})} \,] \le
\mf c_{\rm f} \, (1+T)^4\,  \varrho_n \to 0$.  

\smallskip\noindent{\bf The process $\mf J^{4,n}$ converges to $0$ in
$L^1(0,T;\mc H_{-\mtt r})$ almost surely, if
$\mtt r > 6 + (3d+1)/2$.}
By \eqref{240},
$|\mtt D^{(n)}_v(P_{v,s}\mc A_s\phi_m)| \le \mf c_{\rm f}\, \Vert\,
P_{v,s}\mc A_s\phi_m\Vert_{C^4} \, n^{(d/2)-2}$. By \eqref{245}, this
expression is bounded by
$\mf c_{\rm f}\, (1+\Vert m\Vert)^{6+\beta_0} n^{(d/2)-2}$ for all
$0\le v\le s\le T$, and uniformly in $\eta^n(\cdot)$. Hence,
$|\<\mf J^{4,n}_t ,\phi_m\>| \le \mf c_{\rm f}\, T^2 \, (1+\Vert
m\Vert)^{6+\beta_0}n^{(d/2)-2}$.

As $\mtt r > 6 + (3d+1)/2$, there exists $\beta_0>0$ such that
$\mtt r> 6+\beta_0 + (3d+1)/2$.  By the previous estimate and
\eqref{236},
\begin{equation*}
\Vert \mf J^{4,n} \Vert_{L^1(0,T;\mc H_{-\mtt r})} \;\le\;
\mf c_{\rm f}\,  T^3\, n^{(d/2)-2} \sum_{m\in\bb Z^d} (1+\Vert
m\Vert)^{6+\beta_0} \gamma_m^{-\mtt r/2}\omega_m^{-1/2}\,,
\end{equation*}
uniformly in $\eta^n(\cdot)$. The right-hand side vanishes as
$n\to\infty$, which completes the proof of the claim and the one of
the tightness of the process $\mf J^{n}(\cdot)$.

\begin{proof}[Proof of Proposition \ref{s77}]
The proof follows from the previous estimates, from the decomposition
\eqref{145b} and from the fact that a finite sum of tight sequences of
random elements of a separable Banach space is tight.
\end{proof}
  
\begin{remark}
\label{s91b}
The threshold $\mtt r > 6 + (3d+1)/2$ comes from the terms
$\mf J^{3,n}$, $\mf J^{4,n}$, and the requirement
$\mtt r > \mtt r_0 + d+ 4$ appears in the proof of the tightness of
the process $\mf J^{1,n}$.
\end{remark}

\subsection*{Convergence of the finite-dimensional distributions}

The main result of this subsection reads as follows. 

\begin{proposition}
\label{s55}
Let $\mu_n$ be a sequence of probability measures on $\Omega_n$
satisfying the hypotheses of Theorem \ref{t1}.  Fix $p\ge 1$,
$0\le t_1\le t_2 \le \cdots \le t_p$ and functions
$F_j\in C^\infty(\bb T^d)$, $1\le j\le p$. Under the measure
$\bb P_{\mu_n}$, the random vector
$(X^{n}_{t_1}(F_1), \dots, X^{n}_{t_p}(F_p))$ converges in
distribution to
$(X_{0}(P_{0,t_1} F_1), \dots, X_{0}( P_{0,t_p} F_p)) + (\mss
M_{t_1}(F_1), \dots, \mss M_{t_p}(F_p))$, where the two vectors are
independent, $X_0$ is the limit of the sequence of
$\mc H_{-\mtt r_0}$-valued fields $X^n_0$, and
$(\mss M_{t_1}(F_1), \dots, \mss M_{t_p}(F_p))$ is a zero-mean
Gaussian vector with covariances given by
\begin{equation}
\label{152}
\bb E \big[\, \mss M_{t_i}(F_i) \, \mss M_{t_j}(F_j) \,\big]
\,=\,
\int_0^{t_i\wedge t_j} \int_{\bb T^d} 2\, \chi(u(s,x))\,
(\nabla P_{s,t_i} F_i)(x) \cdot  D(u(s,x)) \, (\nabla P_{s,t_j}
F_j) (x) \; dx \; ds \,.
\end{equation}
\end{proposition}

The proof of this result is divided into several steps. We start with a
simple consequence of Lemma \ref{s47}.

\begin{lemma}
\label{s54}
Fix a cylinder function $f$. Then, there exists a finite constant $C_0
= C_0(f)$ such that
\begin{align*}
& \bb E_{\mu_n} \Big[\,
\int_0^t \,\Big|\, \sum_{x\in \bb T^d_n} J(s,x) \,
\Big\{ \, f(\tau_x \eta^n(s)) -  E_{\nu^n_{u(s, \cdot)} } [ \tau_x f]  \,\Big\}
\,\Big| \,  \, ds \, \Big] \\
&\qquad \le\;
\int_0^t \Big\{ \, \frac{1}{\gamma}\, 
\big[\, H_n(f^n_s\,|\, \nu^n_{u(s, \cdot)}  ) \,+\, \ln (2) \,\big] \;+\;
C_0 \, \gamma\, n^d\, \Vert  J_s \Vert^2_\infty\,
e^{C_0 \, \gamma\,  \Vert  J_s \Vert_\infty} \,\Big\} \; ds
\end{align*}
for every function $J\colon \bb R_+ \times \bb T^d_n \to \bb R$,
$t>0$, $n\ge 1$, $\gamma>0$. In this formula, $J_s$ is the function
which takes the value $J(s,x)$ at $x\in \bb T^d_n$.
\end{lemma}

\begin{lemma}
\label{s62}
For each $F$, $G \in C^\infty(\bb T^d)$, $t>0$, the random variables
$X^n_0(G)$ and $M^{n,t}_t(F)$ are asymptotically independent:
for every $\theta\in\bb R$,
\begin{equation*}
\begin{aligned}
& \lim_{n\to\infty} \bb E_{\mu_n} \Big[\,
\exp i\, \theta\,  \big\{  M^{n,t}_t(F) + X^n_0 (G) \big\}
\,\Big]
\\
&\quad \,=\,
\lim_{n\to\infty} \bb E_{\mu_n} \Big[\,
\exp \big\{   i\, \theta\,  M^{n,t}_t(F) \big\}
\,\Big] \lim_{n\to\infty} \bb E_{\mu_n} \Big[\,
\exp \big\{   i\, \theta\,  X^n_0 (G) \big\}
\,\Big] \,.
\end{aligned}
\end{equation*}
\end{lemma}

\begin{proof}
As $M^{n,t}_t(F)=0$ for $t=0$, we may assume that $t>0$.
For $\theta\in \bb R$, $0\le s\le t$, consider the exponential
\begin{equation*}
e^{ i\, \theta\, M^{n,t}_s (F)} \,=\,
\exp\Big\{ i\,\theta\, X^n_s(H_s) \,-\, i\,\theta\, X^n_0(H_0)
\,-\, i\,\theta\, \int_0^s (\partial_r + L_n) X^n_r(H_r)\,
dr\Big\}\,, 
\end{equation*}
where $\cb{H_r = P_{r,t} F}$. Let
\begin{equation*}
\cb {W^n_s(H_\cdot) } \,:=\, e^{- i\,\theta\, X^n_s (H_s)  }
(\partial_s + L_n)  e^{i\,\theta\, X^n_s (H_s)  }
\,-\,  i\,\theta \, (\partial_s + L_n) X^n_s(H_s)\,,
\end{equation*}
and write the previous exponential as
\begin{equation}
\label{151}
\bb M^n_s (H_\cdot)\, \exp\Big\{ \int_0^s W^n_r (H_\cdot) \, dr \Big\}
\end{equation}
where $\bb M^n_s (H_\cdot)$ is the exponential martingale given by 
\begin{equation*}
\cb{ \bb M^n_s (H_\cdot)} \, :=\, 
\exp\Big\{ i\,\theta\, X^n_s(H_s) \,-\, i\,\theta\, X^n_0(H_0)
\,-\, \int_0^s e^{- i\,\theta\, X^n_r (H_r)  }
(\partial_r + L_n)  e^{i\,\theta\, X^n_r (H_r)  }\, dr\Big\}\,.
\end{equation*}
Note that $\bb M^n_s(H_\cdot)$ is equal to $1$ at time $0$.

In the definition of $W^n_s(H_\cdot)$, the time-derivatives
cancel. Perform a Taylor expansion in the first term. The first-order
term in the expansion cancels with the second term in the definition
of $W^n_s(H_\cdot)$ so that
\begin{equation*}
W^n_s (H_\cdot)\,=\, \frac{-\, \theta^2}{2}\, 
\frac{1}{n^d} \sum_{j=1}^d \sum_{x\in \bb T^d_n}
c_{0,e_j}  (\tau_x \eta^n(s)) \, [\eta^n_{x+e_j} (s) - \eta^n_x(s)]^2\,
[(\nabla^n_j H_s)(x/n)]^2 \,+\, O(n^{-(d/2)-1})\,.
\end{equation*}
By Lemma \ref{s54} with $J(s,x)$ replaced by
$n^{-d} \, [\nabla^n_j P_{s,t}F (x/n)]^2$, the time integral over the
interval $[0,s]$ of $W^n_r(H_\cdot)$ converges in $L^1(\bb P_{\mu_n})$
to  $- (1/2) \,\theta^2\, \mss V_s (H_\cdot)$, where
\begin{equation}
\label{153}
\cb{ \mss V_s (H_\cdot)}  \,:=\, 
\int_0^s \int_{\bb T^d} 2\, \chi(u(r,x))\,
(\nabla H_r)(x) \cdot  D(u(r,x)) \, (\nabla H_r) (x) \; dx \; dr \;.
\end{equation}
In conclusion, $\exp\{ i\, \theta\, M^{n,t}_s (F)\}$ is the product of
an exponential martingale with
$\exp \{ \int_0^s W^n_r (H_\cdot) \, dr \}$, where $W^n_r(H_\cdot)$ is
a uniformly bounded expression which converges in probability to the
deterministic value $- (1/2) \,\theta^2\, \mss V_s (H_\cdot)$.

Fix $\theta\in \bb R$. By the previous considerations,
\begin{align*}
& \lim_{n\to\infty} \bb E_{\mu_n} \Big[\,
\exp i\, \theta\,  \big\{  M^{n,t}_t(F) + X^n_0 (G) \big\}
\,\Big]
\\
&\quad 
\,=\,
\lim_{n\to\infty} \bb E_{\mu_n} \Big[\,
\exp \big\{ i\, \theta\, X^n_0 (G) \big\}\, \bb M^n_t (H_\cdot)
\, \exp\Big\{ \int_0^t W^n_r (H_\cdot)\, dr \Big\} \,\Big]
\end{align*}
Since the time-integral of $W^n_r(H_\cdot)$ converges in probability to $-
(1/2) \,\theta^2\, \mss V_t(H_\cdot)$ and $\bb M^n_t $ is a
martingale equal to $1$ at time $0$, the previous expression is equal
to
\begin{align*}
&\quad \lim_{n\to\infty} \bb E_{\mu_n} \Big[\,
\exp \big\{ i\, \theta\, X^n_0 (G) \big\}\, \bb M^n_t (H_\cdot)
\, \,\Big]\, \exp\big\{  - (1/2) \,\theta^2\, \mss V_t(H_\cdot)
\big\}
\\
&\quad  \,=\,
\lim_{n\to\infty} \bb E_{\mu_n} \Big[\,
\exp \big\{ i\, \theta\, X^n_0 (G) \big\} \,\Big]\,
\exp\big\{  - (1/2) \,\theta^2\, \mss V_t(H_\cdot) \big\}
\end{align*}
The same computation shows that the second term is equal to
$\lim_{n\to\infty} \bb E_{\mu_n} [\, \exp \{i\, \theta\,
M^{n,t}_t(F)\} \,]$. This completes the proof of the lemma.
\end{proof}

\begin{remark}
\label{s60}
The previous proof yields also that $M^{n,t}_t(F)$ converges in
distribution to a Gaussian random variable with zero mean and variance
equal to $\mss V_t(H_\cdot) = \mss V_t(P_{\cdot, t}F)$.
\end{remark}

The previous lemma can be extended to linear combinations of variables
$M^{n,t}_{t}(F)$.

\begin{lemma}
\label{s59}
For each $p\ge 1$, $0\le  t_1\le \cdots \le t_p$, $F_1, \dots, F_p$,
$G \in C^\infty(\bb T^d)$, the random variables $X^n_0(G)$ and
$\sum_{1\le j\le p} M^{n,t_j}_{t_j}(F_j)$ are asymptotically
independent: for every $\theta\in\bb R$,
\begin{equation*}
\begin{aligned}
& \lim_{n\to\infty} \bb E_{\mu_n} \Big[\,
\exp i \theta\,  \Big\{ X^n_0 (G) + \sum_{j=1}^p 
M^{n,t_j}_{t_j}(F_j) \Big\} \,\Big]
\\
&\quad \,=\,
\lim_{n\to\infty} \bb E_{\mu_n} \Big[\,
\exp \big\{   i\, \theta\,  X^n_0 (G) \big\}
\,\Big] 
\lim_{n\to\infty} \bb E_{\mu_n} \Big[\,
\exp i \theta\, \sum_{j=1}^p 
M^{n,t_j}_{t_j}(F_j)  \,\Big]
\,.
\end{aligned}
\end{equation*}
\end{lemma}

\begin{proof}
The proof is similar to the previous one. We may also assume here that
$t_1>0$.  Write the expression inside braces on the left-hand side as
\begin{align*}
X^n_0 ( G) 
\,+\, \sum_{k=0}^{p-1} \Big\{ \, \mss M^{n} _{t_{k+1}}
( \sum_{j=k+1}^p P_{\cdot, t_j} F_j) - \mss M^{n}_{t_{k}}
( \sum_{j=k+1}^p P_{\cdot,t_j} F_j)\, \Big\}\,,
\end{align*}
where $t_0=0$. Recall the definitions of the exponential martingale
$\bb M^n_s(H_\cdot)$, and the process $W^n_s (H_\cdot)$ introduced in
\eqref{151}. In view of the previous displayed equation rewrite the
expectation appearing in the statement of the lemma as
\begin{equation*}
\bb E_{\mu_n} \Big[\, \exp \big\{ i  \theta\,  X^n_0 (G) \big\} \,
\prod_{k=0}^{p-1} \,
\frac{\bb M^n_{t_{k+1}} (G^{(k)}_\cdot) } {\bb M^n_{t_k}(G^{(k)}_\cdot)}
\, \exp\Big\{ \int_{t_k}^{t_{k+1}} W^n_r  (G^{(k)}_\cdot) \, dr \Big\} \,
\Big]\,, 
\end{equation*}
where $G^{(k)} = \sum_{j=k+1}^p P_{\cdot, t_j} F_j$.

As in the previous lemma, $\int_{[t_k, t_{k+1}]} W^n_r
(G^{(k)}_\cdot) \, dr$ converges in $L^1(\bb P_{\mu_n})$ to the
deterministic value
\begin{equation*}
-\, \theta^2\,
\int_{t_k}^{t_{k+1}}   \int_{\bb T^d} \chi(u(r,x))\,
(\nabla G^{(k)}_r)(x) \cdot  D(u(r,x)) \, (\nabla G^{(k)}_r) (x) \; dx \; dr \;.
\end{equation*}
In the penultimate displayed equation, we may, therefore, replace
$\int_{[t_k, t_{k+1}]} W^n_r (G^{(k)}_\cdot) \, dr$ by the previous
quantity paying a cost which vanishes as $n\to\infty$. After this
replacement, taking successive conditional expectations we show that
the remaining expectation is equal to
$\bb E_{\mu_n} [\, \exp \{ i \theta\, X^n_0 (G) \} \,]$. This proves
the asymptotic independence claimed.
\end{proof}

\begin{remark}
\label{s61}
The previous proof yields that the sequence of random variables
$\sum_{1\le j\le p} M^{n,t_j}_{t_j}(F_j)$ converges in distribution to
a zero-mean Gaussian variable with variance equal to
\begin{equation*}
\sum_{k=0}^{p-1} \int_{t_k}^{t_{k+1}}   \int_{\bb T^d} 2\, \chi(u(r,x))\,
(\nabla G^{(k)}_r)(x) \cdot  D(u(r,x)) \, (\nabla G^{(k)}_r) (x) \; dx \; dr \;.
\end{equation*}
\end{remark}

\begin{proof}[Proof of Proposition \ref{s55}]
Let $\mu_n$ be a sequence of probability measures on $\Omega_n$
satisfying the hypotheses of the proposition.  Fix $p\ge 1$,
$0 \le  t_1\le t_2 \le  \cdots \le t_p$, functions $F_j\in C^\infty(\bb T^d)$,
$1\le j\le p$, and $T>t_p$.

Since $M^{n,t}_t (F)=0$ if $t=0$, we may assume that $t_1>0$. 
We first prove the proposition for $p=1$. By Lemma \ref{s56},
$X^n_{t_1}(F_1) - X^n_0(P_{0,t_1}F_1) -M^{n,t_1}_{t_1}(F_1)$ converges
to $0$ in $L^1(\bb P_{\mu_n})$. By hypothesis, $X^n_0(P_{0,t_1}F_1)$
converges in law to $X_0(P_{0,t_1}F_1)$, and by Remark \ref{s60},
$M^{n,t_1}_{t_1}(F_1)$ converges in law to $\mss M_{t_1}(F_1)$. By
Lemma \ref{s62}, $X^n_0(P_{0,t_1}F_1)$, $M^{n,t_1}_{t_1}(F_1)$ are
asymptotically independent. Thus,
$X^n_0(P_{0,t_1}F_1) +M^{n,t_1}_{t_1}(F_1)$ converges in law to the sum of
two independent random variables which are distributed as
$X_0(P_{0,t_1}F_1)$ and $\mss M_{t_1}(F_1)$. This completes the proof of
the proposition for $p=1$.

The proof for the general case is similar. By Lemma \ref{s56}, it is
enough to consider the limit of the vector
$(X^n_0(P_{0,t_1}F_1) +M^{n,t_1}_{t_1}(F_1), \dots, X^n_0(P_{0,t_p}F_p)
+M^{n,t_p}_{t_p}(F_p))$. By Lemma \ref{s59}, the vectors
$(X^n_0(P_{0,t_1}F_1) , \dots, X^n_0(P_{0,t_p}F_p))$,
$(M^{n,t_1}_{t_1}(F_1), \dots, M^{n,t_p}_{t_p}(F_p))$ are
asymptotically independent [By Fourier transform two vectors are
independent, or asymptotically independent, if all linear combinations
of the coordinates are independent]. By hypothesis,
$(X^n_0(P_{0,t_1}F_1) , \dots, X^n_0(P_{0,t_p}F_p))$ converges in
distribution to $(X_0(P_{0,t_1}F_1) , \dots, X_0(P_{0,t_p}F_p))$. By
Remark \ref{s61}, and a straightforward computation, the vector
$(M^{n,t_1}_{t_1}(F_1), \dots, M^{n,t_p}_{t_p}(F_p))$ converges in
distribution to a zero-mean Gaussian vector with covariances given by
\eqref{152}. This completes the proof of the proposition.
\end{proof}

\subsection*{Characterization of the limit points}

In the first part of this section, we proved that the sequence of
probability measures $\bb Q_{\mu_n}$ on $L^1(0,T; \mc H_{-\mtt r})$ is
tight provided $\mtt r$ is large enough. In this subsection we
characterize the limit points.  According to Corollary \ref{s74}, it is
enough to determine the distribution of the $p$-tuples
\begin{equation}
\label{188}
\Big( \int_{s_1}^{t_1} \< X_{r} , H_1\> \, dr,
\dots,
\int_{s_p}^{t_p} \< X_{r} , H_p\> \, dr \, \Big)
\end{equation}
where $H_1, \dots, H_p$ belong to $\mc H_{\mtt r}$, and
$0\le s_1 < t_1 \le s_2 < \cdots < t_{p-1} \le s_p <t_p\le T$. The
proof of the convergence of these $p$-tuples is similar to the one of
the finite-dimensional distribution convergence.

Throughout this subsection, besides the hypotheses of Proposition
\ref{s77}, we assume hypothesis {\rm (b)} of Theorem \ref{t1}, in the
form

\begin{itemize}
\item[(H)] the sequence of probability measures
$\bb P_{\mu_n} \circ (X^n_0)^{-1}$ on $\mc H_{-\mtt r_0}$ converges
weakly to a probability measure $\bb Q_0$.
\end{itemize}

Denote by $X_0$ a random element of $\mc H_{-\mtt r_0}$ distributed
according to $\bb Q_0$.  Clearly, $X^n_0$ also converges weakly to
$X_0$ in $\mc H_{-\mtt r}$ if $\mtt r>\mtt r_0$.

Fix $G \in \mc H_{\mtt r}$, $0\le s<t\le T$, and let
\begin{equation}
\label{246}
\cb{\Gamma^{s,t}_r G} \;:=\; \int_{r\vee s}^{t} P_{r,r'}\, G \;dr'
\;, \qquad 0\le r \le t\;.
\end{equation}

\begin{lemma}
\label{s95}
Fix $\mtt r > \mtt r^\star$, where $\mtt r^\star$ has been introduced
in \eqref{242b}.  Fix $G\in \mc H_{\mtt r}$ and $0\le s<t\le T$.  Then,
$\Gamma^{s,t}_0 G$ belongs to $\mc H_{\mtt r_0}$ and
\begin{equation}
\label{247}
\Vert \Gamma^{s,t}_0 G \Vert_{\mc H_{\mtt r_0}} \;\le\;
\mf c_{\rm f} \, T\, \Vert G \Vert_{\mc H_{\mtt r}}\;,
\qquad
\sup_{0\le r\le t} \Vert \Gamma^{s,t}_r G \Vert_{C^{2}} \;\le\;
\mf c_{\rm f} \, T\, \Vert G \Vert_{C^{3}} \;.
\end{equation}
Moreover, for all $0\le r\le t$ and $0\le a\le b\le t$,
\begin{equation}
\label{248}
\Big\Vert \, \Gamma^{s,t}_r G \,-\, \int_{a}^{b} P_{r,r'}\, G \; dr'\,
\Big\Vert_{C^{2}} \;\le\; \mf c_{\rm f} \, 
\Vert G \Vert_{C^{4}} \, \big\{\, |a - (r\vee s)| \,+\, |b-t| \,\big\}\;.
\end{equation}
\end{lemma}

\begin{proof}
Let $p$ be the smallest even integer larger than $\mtt r_0 + (d/2)$,
so that $p\le \mtt r_0 + (d/2) + 2$ and, by \eqref{242b},
$\mtt r > p + 2 + (d/2)$.  By \eqref{220b} with $q=p$,
$\Vert P_{0,r'} G\Vert_{C^p} \le \mf c_{\rm f} \Vert G
\Vert_{C^{p+1}}$. Thus, since $p> \mtt r_0 + (d/2)$ and
$\mtt r> p+1+ (d/2)$, by \eqref{243},
$\Vert P_{0,r'} G\Vert_{\mc H_{\mtt r_0}} \le \mf c_{\rm f} \Vert G
\Vert_{\mc H_{\mtt r}}$ for all $0\le r'\le T$.

The same argument, with \eqref{221} in place of \eqref{220b}, yields
that the map $r'\mapsto P_{0,r'} G$ is Lipschitz-continuous with
values in $\mc H_{\mtt r_0}$. Hence, the integral in \eqref{246} is a
Bochner integral in $\mc H_{\mtt r_0}$, and the first bound in
\eqref{247} follows.  The second one is proved in the same way.  By
\eqref{220d} with $0<\beta<1$,
$\Vert P_{r,r'} G\Vert_{C^2} \le \Vert P_{r,r'} G\Vert_{(2+\beta)} \le
\mf c_{\rm f} \Vert G\Vert_{(2+\beta)} \le \mf c_{\rm f} \Vert
G\Vert_{C^3}$.

We turn to \eqref{248}.  By \eqref{221} with $p'=2$,
\begin{equation}
\label{254}
\sup_{0\le r\le a\wedge b} \Vert P_{r,b}\, G - P_{r,a}\, G
\Vert_{C^{2}} \;\le\; \mf c_{\rm f}\, \Vert G\Vert_{C^4}\, |b-a| \;,
\end{equation}
so that $r'\mapsto P_{r,r'} G$ is Lipschitz continuous in $C^2$,
uniformly in $r$, with Lipschitz constant $\mf c_{\rm f}\, \Vert
G\Vert_{C^4}$.  As the two integrals in \eqref{248} are carried over
the intervals $[r\vee s, t]$ and $[a,b]$, whose symmetric difference
has Lebesgue measure bounded by $|a-(r\vee s)| + |b-t|$, and as the
integrand is bounded by $\mf c_{\rm f}\Vert G\Vert_{C^3}$ in $C^2$,
the assertion follows.
\end{proof}

\begin{remark}
\label{s96}
The first bound in \eqref{247} is the only place where the
condition $\mtt r > \mtt r_0+d+4$ of \eqref{242b} is needed in this
subsection; it ensures that $P_{0,r}$ maps $\mc H_{\mtt r}$ into
$\mc H_{\mtt r_0}$.  It may be replaced by $\mtt r\ge \mtt r_0$ if one
disposes of the $\mc H_{\mtt r_0}$-boundedness of the semigroup, which
follows from an elementary energy estimate.

This proof requires the same regularity of $u$ as in Proposition
\ref{s77} as it uses \eqref{221} with the same exponents.
\end{remark}

For $G\in \mc H_{\mtt r}$ and $0\le s<t\le T$, let
\begin{equation}
\label{249}
\cb{\sigma^2(G; s,t)} \;:=\;
\int_{0}^{t} \int_{\bb T^d} 2\, \chi(u(r,x))\,
(\nabla \Gamma^{s,t}_r G)(x) \cdot  D(u(r,x)) \,
(\nabla \Gamma^{s,t}_r G)(x) \; dx \; dr \;,
\end{equation}
which is finite by the second bound in \eqref{247}.

\begin{lemma}
\label{s75}
Fix $H\in \mc H_{\mtt r}$, $0\le s < t \le T$. Then, the random
variable $\int_{s}^{t} X^n_r (H) \, dr$ converges in law to
$X_0 ( \Gamma^{s,t}_0 H) + \mf G$, where $X_0$ and $\mf G$ are
independent, and $\mf G$ is a centered Gaussian random variable with
variance $\sigma^2(H;s,t)$.
\end{lemma}

\begin{proof}
Fix $H\in \mc H_{\mtt r}$, $0\le s < t \le T$, and let
\begin{equation*}
\cb{\varepsilon^n_r} \,:= \, X^n_r (H) - X^n_0 (P_{0,r}\, H) - M^{n,r}_r
(H)\,, \quad 0\le r\le T\,.
\end{equation*}
By Fubini's theorem,
\begin{equation}
\label{250}
\int_{s}^{t} X^n_r (H) \, dr \,=\,
X^n_0\big( \Gamma^{s,t}_0 H \big) \;+\;
\int_{s}^{t} M^{n,r}_r (H) \, dr \;+\;
\int_{s}^{t} \varepsilon^n_r \, dr \,.
\end{equation}
By Lemma \ref{s56}, for each $r>0$, $\varepsilon^n_r$ converges
to $0$ in $L^1 (\bb P_{\mu_n})$. The same proof yields that the time
integral also converges to $0$ in $L^1 (\bb P_{\mu_n})$.

It remains to examine the first two terms on the right-hand side.  Fix
an integer $m\ge 1$, and let $t_i = s + (i /m)\, (t-s)$,
$0\le i\le m$. By linearity, we may rewrite these terms as
\begin{align}
\label{187}
X^n_0 \Big(\, \int_{s}^{t} P_{0,r} H\, dr\, \Big)
+ \frac{1}{m}\, (t-s)\, \sum_{i=0}^{m-1} M^{n, t_i}_{t_i} (H) 
\, +\, \sum_{i=0}^{m-1} \int_{t_i}^{t_{i+1}} 
\big[\, M^{n,r}_{r} (H) - M^{n,t_i}_{t_i} (H)\,\big] \, dr \,.
\end{align}
Denote by $\mss S^n_m$ the second term of \eqref{187} and by
$\mss E^n_m$ the third one, and recall that the first one is
$X^n_0(\Gamma^{s,t}_0H)$.

\smallskip\noindent {\it Step 1: The main term.}  Fix $m\ge 1$, and
apply Lemma \ref{s59} and Remark \ref{s61} with $p=m$, with the test
functions $F_j := [(t-s)/m]\, H$, and with $G := \Gamma^{s,t}_0H$.
Mind that $t_0 = s$, while, in Lemma \ref{s59}, $t_0=0$.  Mind also
that the proof of Lemma \ref{s59} uses the function $G$ only through
the random variable $X^n_0(G)$, and never differentiates it: the lemma
therefore holds for every $G$ in $\mc H_{\mtt r_0}$, and
$\Gamma^{s,t}_0G$ belongs to this space by Lemma \ref{s95}.
Similarly, only finitely many derivatives of the functions $F_j$ are
used, so that $F_j\in C^5(\bb T^d)$ suffices, which holds by
\eqref{243} since $\mtt r > 5 + (d/2)$.

By hypothesis {\rm (H)}, since $\Gamma^{s,t}_0 H \in \mc H_{\mtt
r_0}$ and since the map $\mtt x \mapsto \mtt x(\Gamma^{s,t}_0 H)$ is
continuous from $\mc H_{-\mtt r_0}$ to $\bb R$,
$X^n_0(\Gamma^{s,t}_0H)$ converges in law to $X_0(\Gamma^{s,t}_0H)$.
Hence, by Lemma \ref{s59}, Remark \ref{s61} and the Cram\'er--Wold
device, as $n\to\infty$, the sum of the first two terms of \eqref{187}
converges in law to $X_0 ( \Gamma^{s,t}_0H ) + \mf G_m$, where
$X_0 ( \Gamma^{s,t}_0 H)$ and $\mf G_m$ are independent and $\mf G_m$
is a centered Gaussian random variable with variance
\begin{gather*}
\sigma^2_m \;=\; \sum_{k=0}^{m-1} \int_{\tau_k}^{\tau_{k+1}}
\int_{\bb T^d} 2\, \chi(u(r,x))\,
(\nabla H^{(k)}_r)(x) \cdot  D(u(r,x)) \, (\nabla H^{(k)}_r) (x) \; dx \; dr \;.
\end{gather*}
In this formula, $\tau_0=0$, $\tau_j = t_{j-1}$ for $1\le j\le m$, and 
$H^{(k)}_r= (1/m) (t-s) \sum_{i=k}^{m-1} P_{r, t_i} H$,
$\tau_k \le r\le \tau_{k+1}$.

\smallskip\noindent{\it Step 2: $\sigma^2_m\to\sigma^2(H;s,t)$.}  By
\eqref{248} and the Lipschitz continuity of $r'\mapsto P_{r,r'} G$ in
$C^2$, established in \eqref{254}, we obtain that
\begin{equation*}
\sup_{0\le k<m} \; \sup_{\tau_k\le r\le \tau_{k+1}}
\big\Vert \, H^{(k)}_r \,-\, \Gamma^{s,t}_r H \, \big\Vert_{C^2}
\;\le\; \mf c_{\rm f}\,  \Vert H\Vert_{C^4}\,
\frac{t-s}{m}\;\cdot
\end{equation*}
Since $\tau_m = t_{m-1} = t - [(t-s)/m] \to t$, and since the
integrand in \eqref{249} is bounded, replacing $H^{(k)}_r$ by
$\Gamma^{s,t}_rH$ in the previous displayed formula for $\sigma^2_m$,
and $[0,\tau_m]$ by $[0,t]$, shows that $\sigma^2_m$ converges to
$\sigma^2(H;s,t)$ as $m\to\infty$.

Thus, as $m\to\infty$, the sequence of random variables
$\mf G_{m}$ converges in distribution to a centered Gaussian
random variable $\mf G$ with variance $\sigma^2(H;s,t)$,
independent of $X_0$.

\smallskip\noindent{\it Step 3: The remainder $\mss E^n_m$.}
It remains to show that the third term in \eqref{187} converges to $0$
in $L^1(\bb P_{\mu_n})$ as $n\to\infty$ and then $m\to\infty$.  We
first obtain a bound for $r$ fixed. We affirm that there exists a
finite constant $C_0$, which depends only on the cylinder functions
$h_j$, such that  
\begin{equation}
\label{186}
\begin{aligned}
& \bb E_{\mu_n} \big[\, \big|\, M^{n,r}_{r} (H) - M^{n,t_i}_{t_i}
(H)\,\big|\, \big] \,\le\,
C_0\, \Vert H \Vert_{C^1}\, (r-t_i)^{1/2}
\\
&\quad \, +\,
C_0 \, T\, \int_0^{t_i}\,
\Vert \mc A_u \,  P_{u,r}\,H - \mc A_u \,  P_{u,t_i} \,H \Vert_{C^1}
\, du \, +\,
C_0\, T\, \sup_{0\le u\le r}
\Vert \mc A_u \,  P_{u,r}\,H \Vert_{C^1}\, (r-t_i)
\end{aligned}
\end{equation}
for all $t_i \le r\le t_{i+1}$, $0\le i<m$.  To prove this estimate,
recall from \eqref{56} the formula for $M^{n,r}_{r} (H)$. By the
triangle inequality, the expectation on the left-hand side is
bounded by
\begin{equation*}
\begin{aligned}
& \bb E_{\mu_n} \big[\, \big|\, M^{n}_{r} (H) - M^{n}_{t_i}
(H)\,\big|\, \big] \,+\,
\int_0^{t_i} \bb E_{\mu_n} \Big[\, \big|\,  M^{n}_{u} \big( \mc A_u \, P_{u,r}
\,H - \mc A_u \,  P_{u,t_i} \,H \big) \,\big|\, \Big] \, du 
\\
&\quad \,+\, \int_{t_i}^r
\bb E_{\mu_n} \big[\, \big|\, M^{n}_{u} ( \, \mc A_u \,  P_{u,r} \, H
\,)
\,\big|\, \big] \, du \,.
\end{aligned}
\end{equation*}
To complete the proof of the claim \eqref{186}, it remains to bound
the $L^1(\bb P_{\mu_n})$ norm by the $L^2(\bb P_{\mu_n})$ norm, to recall
from \eqref{143} the explicit formula for the quadratic variation of
the martingale $M^n_r(F)$, and the bound \eqref{134}.

The $C^1$-norm appearing in the second term of \eqref{186} is bounded
by $\mf c_{\rm f}\, \Vert P_{u,r} H - P_{u,t_i} H\Vert_{C^{4}}$.  By
\eqref{221} with $p'=4$, this expression is less than or equal to
$\mf c_{\rm f}\, \Vert H\Vert_{C^{6}}\, (r-t_i)$. On the other hand,
by \eqref{220b} with $q=4$, the supremum in the third one is bounded
by
$\mf c_{\rm f} \sup_u \Vert P_{u,r} H \Vert_{C^{4}} \le \mf c_{\rm f}
\Vert H \Vert_{C^{5}}$.  Since $r-t_i\le (t-s)/m$ and, by \eqref{243},
$\Vert H \Vert_{C^6}\le \mf c_{\rm f} \Vert H\Vert_{\mc H_{\mtt r}}$
because $\mtt r > 6 + (d/2)$, we conclude that
\begin{equation}
\label{251}
\sup_{n\ge 1} \bb E_{\mu_n} \big[ \, \big| \mss E^n_m \big|\,
\big] \;\le\; \mf c_{\rm f} \, (1+T)^2\, \Vert H \Vert_{\mc H_{\mtt
r}}\, \Big( \frac{t-s}{m} \Big)^{1/2}\;\cdot
\end{equation}

\smallskip\noindent{\it Step 4: Conclusion.}
It remains to recollect all previous estimates to complete the proof
of the lemma.
\end{proof}

We can extend this result, with the same proof, to vectors.

\begin{lemma}
\label{s76}
Fix $p\ge 2$, elements $H_1, \dots, H_p$ of $\mc H_{\mtt r}$, and
$0\le s_1 < t_1 \le s_2 < \cdots < t_{p-1} \le s_p <t_p\le T$. As
$n\to\infty$, the random vector \eqref{188} converges in law to the
sum of the random vector
$( X_0(\Gamma^{s_1,t_1}_0 H_1), \dots, X_0(\Gamma^{s_p,t_p}_0H_p))$ and of a
centered Gaussian random vector $\mf G$ in $\bb R^p$, independent of
$X_0$, whose covariances are given by
\begin{equation}
\label{252}
\bb E\big[\, \mf G_i \, \mf G_j\,\big] \;=\;
\int_{0}^{t_i\wedge t_j} \int_{\bb T^d} 2\, \chi(u(r,x))\,
\big(\nabla \Gamma^{s_i,t_i}_r H_i\big)(x) \cdot  D(u(r,x)) \,
\big(\nabla \Gamma^{s_j,t_j}_r H_j\big)(x) \; dx \; dr \;.
\end{equation}
\end{lemma}

\begin{proof}
By the Cram\'er--Wold device, it is enough to prove the convergence in
law of $\sum_{i=1}^p \theta_i \int_{s_i}^{t_i} X^n_r(H_i)\, dr$ for
every $\theta\in\bb R^p$.  The proof of Lemma \ref{s75} applies
with minor modifications.
\end{proof}

\begin{remark}
\label{s97}
In contrast with what happens at a fixed time, the components of
$\mf G$ are not independent, even when the intervals $[s_i,t_i]$ are
disjoint: they are all driven by the same space-time white noise on
$[0,t_p]\times \bb T^d$, and the covariance \eqref{252} does not
vanish for $i\neq j$.
\end{remark}

We conclude this section with the identification of the limit
point, which completes the proof of Theorem \ref{t1}.  Denote by
$\mss M$ the centered Gaussian field with covariances given by
\eqref{152}, independent of $X_0$, and by $\bb Q$ the law, on
$L^1(0,T;\mc H_{-\mtt r})$, of the process
\begin{equation}
\label{253}
Y_t \;=\; P^*_{0,t}\, X_0 \;+\; \mss M_t \;, \qquad 0\le t\le
T\;,
\end{equation}
that is, of the mild solution of the equation \eqref{58} with
initial condition $X_0$.

\begin{proposition}
\label{s98}
Under the hypotheses of Proposition \ref{s77} and {\rm (H)}, the
sequence $\bb Q_{\mu_n}$ converges weakly, in
$L^1(0,T;\mc H_{-\mtt r})$, to $\bb Q$.
\end{proposition}

\begin{proof}
By Proposition \ref{s77}, the sequence $\bb Q_{\mu_n}$ is tight.
By Corollary \ref{s74}, a probability measure on
$L^1(0,T;\mc H_{-\mtt r})$ is determined by the laws of the
$p$-tuples \eqref{188}.  By Lemma \ref{s76}, these laws are the same
for all the limit points of the sequence $\bb Q_{\mu_n}$, so that this
sequence converges.

It remains to identify the limit with $\bb Q$.  Let $Y$ be the
process \eqref{253}.  By \eqref{152}, $\mss M_r(H)$ may be represented
as
\begin{equation*}
\sum_{j=1}^d \int_0^r \int_{\bb T^d} \sqrt{2\, \chi(u(v,x))\,
D_j(u(v,x))}\, (\partial_{x_j} P_{v,r} H)(x)\, \xi^j(dv,dx)\,,
\end{equation*}
where
$\xi = (\xi^1,\dots,\xi^d)$ is the space-time white noise which
appears in \eqref{58}.  Hence, by the stochastic Fubini theorem and
since
$\int_{v\vee s}^t \partial_{x_j} P_{v,r}\, H\, dr = \partial_{x_j}
\Gamma^{s,t}_v H$,
\begin{equation*}
\int_s^t \mss M_r(H)\, dr \;=\; \sum_{j=1}^d \int_0^t
\int_{\bb T^d} \sqrt{2\, \chi(u(v,x))\, D_j(u(v,x))}\,
\big(\partial_{x_j} \Gamma^{s,t}_v H \big)(x)\, \xi^j(dv,dx)\;,
\end{equation*}
a centered Gaussian random variable with variance
$\sigma^2(H;s,t)$.  On the other hand,
$\int_s^t \< P^*_{0,r}X_0 , H\>\, dr = X_0(\Gamma^{s,t}_0H)$.  Hence,
under $\bb Q$, the $p$-tuple \eqref{188} has exactly the law described
in Lemma \ref{s76}, which completes the proof.
\end{proof}

\appendix

\section{Equivalence of Ensembles}
\label{sec5}

In this section we obtain an expansion, up to the second order, for
the expectation of a cylinder function with respect to the tilted
canonical measure.

Fix a cube $\cb{\Sigma^+_r}$, $r\ge 2$, and set $\cb{N=r^d}$.  Let
$\cb{\Omega^+_{r}} := \{0,1\}^{\Sigma^+_r}$,
$\cb{ \Omega^+_{r, \mss M}} := \{\eta\in \Omega^+_{r} : \sum_{x\in
\Sigma^+_r} \eta_x = \mss M\}$, $0\le \mss M\le N$, be the state space
of the exclusion process on $\Sigma^+_{r}$.  We first examine the
equivalence of ensembles in the homogeneous case. Denote by
$\cb{\nu_{r,\mss M}}$ the stationary state of the symmetric simple
exclusion process on $\Omega^+_{r, \mss M}$. This is the uniform
measure over all configurations in $\Omega^+_{r, \mss M}$.  The result
reads as follows.

\begin{lemma}
\label{s05}
For each finite subset $A$ of $\Sigma^+_r$,
\begin{equation*}
E_{\nu_{r,\mss M}} [ \eta_A ] \,=\, 
\sum_{k=0}^{|A|-1} \Big(\frac{-1}{N}\Big)^k\,
\rho^{|A|-k} \, (1-\rho)^k\, \sum_{|B|=k}
\prod_{j\in B} j \, \Big(  1 + \frac{j} {N-j} \Big) \,.
\end{equation*}
In this formula, $\rho = \mss M/N$ and the sum over $B$ is carried out
over all subsets $B$ of $\{1, \dots, |A|-1\}$ with $k$ elements. In
particular,
\begin{align*}
E_{\nu_{r,\mss M}} [ \eta_A ] \,=\,  \rho^{|A|}\,
+\, \frac{1}{N} \, \mtt p_1(\rho) \,
+\, \frac{1}{N^2} \, \mtt p_2\, (\rho)
+\mf R_N
\end{align*}
where
\begin{equation*}
\begin{gathered}
\mtt p_1(\rho) \,=\,  -\, 
\frac{1}{2}\, |A|\, (|A|-1) \, \rho^{|A|-1}\, (1-\rho)\,
\\
\mtt p_2(\rho) \,=\,  \frac{1}{2}\, \chi(\rho)\, \rho^{|A|-3}\,
\big\{ \, \mf n^2_{1} \, (1-\rho) - \mf n_{2} \,(1+\rho)\,
\big\} \,,
\end{gathered}
\quad\text{and}\quad 
|\mf R_N| \, \le \, C_A \, \frac{1}{N^3} 
\end{equation*}
for some finite constant $C_A$.
In this formula, $\mf n_{k} = 1^k + \cdots + (|A|-1)^k$.
\end{lemma}

\begin{proof}
Since $\nu_{r,\mss M}$ is the uniform measure,
\begin{equation*}
E_{\nu_{r, \mss M}} [ \eta_A] \,=\, \frac{\mss M}{N}  \cdots  \frac{\mss M-|A|+1}
{N-|A|+1}\,\cdot
\end{equation*}
Rewrite $(\mss M-j)/(N-j)$ as $(\mss M/N) + \epsilon_j$, where
$\epsilon_j = -\, (1-\rho) \, (j/N)  [ 1 + j/(N-j)]$ to conclude that
\begin{align*}
E_{\nu_{r, \mss M}} [ \eta_A] \, & =\, \rho\, \prod_{j=1}^{|A|-1}
\Big\{ \rho -  (1-\rho) \, \frac{j}{N}\,  \Big(  1 + \frac{j} {N-j}
\Big)\, \Big\}\,.
\end{align*}
A straightforward computation yields the first assertion of the lemma.

We may now compute the first three terms in the expansion of
$E_{\nu_{r, \mss M}} [ \eta_A]$. The first two terms are easy to
compute and correspond to the first two terms in the sum appearing in
the statement of the lemma. The third one, which has order $1/N^2$, is
obtained considering the term $k=2$ in the sum and the term $k=1$ of
order $1/N^2$. It is given by
\begin{align*}
\frac{1}{N^2} \, \Big\{
\rho^{|A|-2} \, (1-\rho)^2\, \frac{1}{2}\,
\Big[ \, \Big( \sum_{j=1}^{|A|-1} j \Big)^2
\,-\, \sum_{j=1}^{|A|-1} j^2 \,\Big]
\,-\, \rho^{|A|-1} \, (1-\rho)
\sum_{j=1}^{|A|-1} j^2
\Big\}\,.
\end{align*}
This completes the proof of the lemma.
\end{proof}

As observed in  \cite[Lemma A2.2.1]{kl} the polynomials appearing in
the expansion can be represented in terms of covariances and
cumulants. For example,
\begin{equation}
\label{37b}
\mtt p_1(\rho) \,=\, \frac{1}{2}\, \Big\{
\frac{\gamma_3(\rho)}{\gamma_2(\rho)^2}\,
E_{\nu_\rho} \big[ \eta_A \, ;\, \sum_{x\in A}  (\eta_x - \rho) \, \big]
\,-\,
\frac{1}{\gamma_2(\rho)}\,
E_{\nu_\rho} \Big[ \eta_A \, ;\,
\Big( \sum_{x\in A}  (\eta_x - \rho) \,\Big)^2\, \Big]\,\Big\}\,.
\end{equation}
In this formula, $\gamma_p(\rho)$, $p\ge 1$, represents the $p$-th
cumulant of a Bernoulli distribution of parameter $\rho$, and
$E_{\nu_\rho} [ g \, ;\, h\,]$ the covariance between two functions
$g$ and $h$. Recall that the first three cumulants are given by
$\gamma_1(\rho) = \rho$, $\gamma_2(\rho) = \chi(\rho)$,
$\gamma_3(\rho) = \chi(\rho)(1-2\rho)$.

The proof of the previous lemma is based on the fact that the
canonical measures of the exclusion dynamics are the uniform
measures. In more general contexts one relies on the local central
limit theorem as shown below in the inhomogeneous case. \medskip

\subsection*{Locally uniform local central limit theorem}

We turn to the inhomogeneous case.  Fix a density profile
$\varrho\colon \bb R^d \to [0,1]$ of class $C^2$ with compact support
in $ (-2,2)^d$ and for which there exists $0<\mf r<1/2$ such that
\begin{equation*}
\mf r  \,\le\, \varrho(\mtt x )\,\le 1-\mf r\quad
\text{for all $\mtt x \in [-1,1]^d$.}
\end{equation*}
Although all quantities in this section depend on the density profile
$\varrho(\cdot)$, it is omitted from the notation.

Fix integers $n\ge 8$, $2\le r \le n/4$. In the applications, $r = r_n$ is
a sequence such that $r_n/n\to 0$.  For $\varphi\in \bb R$, let
$\varrho_\varphi \colon (-2,2)^d \to [0,1]$ be given by
\begin{equation}
\label{21b}
\cb{\varrho_\varphi(\mtt x)} 
\,:=\, \frac{e^\varphi \varrho(\mtt x)  }
{e^\varphi \varrho(\mtt x)  + [1- \varrho(\mtt x)] }\,\cdot
\end{equation}
The parameter $\varphi$ is called the chemical potential.
Denote by $\cb{\nu^{n, r}_{\varrho_\varphi(\cdot)}}$ the Bernoulli
product measure on $\Omega^+_r$ with marginals given by
$\nu^{n, r}_{\varrho_\varphi (\cdot)} \{\eta: \eta_x =1 \} =
\varrho_\varphi(x/n)$, $x\in \Sigma^+_r$.

For the statement of the main result of this subsection, we recall the
Hermite polynomials up to order $9$, for the reader's convenience:
\begin{gather*}
H_0(x)=1\,, \quad H_1(x)=x\,, \quad
H_2(x)=x^2-1\,, \quad H_3(x)=x^3-3x \,, \quad
H_4(x)= x^4-6x^2+3\,,
\\
H_5(x) = x^5-10x^3+15x\,, \quad
H_6(x)=x^6-15x^4+45x^2-15\,, \quad
H_7(x) = x^7-21x^5+105x^3-105x\,,
\\
H_8(x)=x^8-28x^6+210x^4-420x^2+105\,, \quad
H_9(x)=x^9-36x^7+378x^5-1260x^3+945x\,.
\end{gather*}
Fix a finite subset $A$ of $\Sigma^+_r$.  Recall the definition of the
functions $\cb{q_{\nu, n}}$, $\nu\ge 1$, defined in equation (1.10),
Chapter VI in \cite{pet}, page 138. In the statement of \cite[Theorem
VII.12]{pet}, $n$ represents the number of variables, which is here
$\cb{N_A} := |\Sigma^+_r \setminus A| = N- |A|$, where, recall,
$N=|\Sigma^+_r|$, so that
\begin{equation}
\label{78}
q_{j,N,A} (z)\; =\; \frac{1}{\sqrt{2\pi}} \, e^{-z^2/2} \, \sum H_{j+2s} (z)
\prod_{m=1}^j \frac{1}{k_m!} \left( \frac{\lambda_{m+2,N, A}(\varphi)}
{(m+2)! }\right)^{k_m}\; ,
\end{equation}
where the sum is carried out over all non negative integer solutions
of $k_1+2k_2+\cdots + j k_j =j$, and $s=k_1+k_2+\cdots + k_j$.  In
this formula, $H_m(\cdot)$ is the Hermite polynomial of degree $m$,
\begin{equation}
\label{79}
\lambda_{j, N, A} (\varphi) \,=\, \frac{N_A^{(j-2)/2}}  {\Gamma_{N, A, \varphi}
^{j/2}} \sum_{x\in \Sigma^+_r \setminus A} \gamma_{j,x} (\varphi)
\,, 
\quad 
\cb {\Gamma_{N,  A, \varphi}  } \,:=\, \sum_{x\in \Sigma^+_r \setminus A} 
\chi (\varrho_\varphi (x/n))\,,
\end{equation}
and $\gamma_{j,x} (\varphi)$ is the cumulant of order $j$ of a
Bernoulli random variable with mean $\varrho_\varphi (x/n)$.  When $A$
is the empty set, we write $N_\varnothing = |\Sigma^+_r| = N$,
and $\cb{q_{j,N} (\cdot)}$, $\cb{\Gamma_{N, \varphi}}$ for
$q_{j,N,\varnothing} (\cdot)$, $\Gamma_{N, \varnothing, \varphi}$, respectively.

\begin{theorem}
\label{s29}
For each $p\ge 2$, finite subset $A$ of $\Sigma^+_r$, and bounded
interval $I \subset \bb R$, there exist a finite constant
$\mf c_{_{\rm CLT}} = \mf c_{_{\rm CLT}}(\varrho(\cdot) , p, I)$
and a finite integer $r_0=r_0(\varrho(\cdot) , p, I, |A|)$
such that
\begin{equation*}
\sup_{\varphi \in I} \, \sup_{0\le M\le |\Sigma^+_r \setminus A|} \,
\Big\vert
\, \sqrt{ \Gamma_{N, A, \varphi} }  \;
\nu^{n, r}_{\varrho_\varphi(\cdot)} \big\{
\sum_{x\in \Sigma^+_{r} \setminus A} \eta_x = M  \, \big\}
\; -\;
\frac{1}{\sqrt{2\pi}}\, e^{-z^2/2} \,-\, 
\sum_{j=1}^{p-2} \frac{1}{N_A ^{j/2}} q_{j, N,  A} (z) \, \Big\vert
\; \le\;
\frac{ \mf c_{_{\rm CLT}} }{ N^{(p-1)/2}}
\end{equation*}
for all $r \ge r_0$. In this formula, 
\begin{equation}
\label{74}
z \,=\, \frac{ M -  \sum_{x\in \Sigma^+_{r} \setminus A}
\varrho_\varphi (x/n)}
{\sqrt{ \Gamma_{N,  A, \varphi} } } \,.
\end{equation}
\end{theorem}

\begin{proof}
For a fixed $\varphi\in \bb R$, this result is the local central limit
theorem presented in \cite[Theorem VII.12]{pet}. Going through the
proof, it is easily seen that the estimates are uniform for $\varphi$
in bounded sets. For $\varrho(\cdot)$ constant, this is \cite[Theorem
A2.1.3]{kl}, see \cite[Remark A2.1.2]{kl}.
\end{proof}

Note that we stated the result with an error of order
$O(1/N^{(p-1)/2})$ instead of an error of order $o(1/N^{(p-2)/2})$
which appears in \cite[Theorem VII.12]{pet}. To get this apparently
better bound, apply \cite[Theorem VII.12]{pet} with $p+1$ and insert
the last term, $ N_A ^{-(p-1)/2} q_{p-1, N, A} (z)$ in the error.

\subsection*{Locally uniform equivalence of ensembles}

Denote by $\cb{\nu^{\rm c}_{r,\mss M} }$,
$0\le \mss M\le |\Sigma_r^+|$, the product measure
$\nu^{n, r}_{\varrho(\cdot)}$ conditioned on $\Omega^+_{r,\mss M}$, the
set of configurations of $\Omega^+_r$ with $\mss M$ particles.  The next
result is an extension of the so-called equivalence of ensembles. It
encompasses the case where the density profile is slowly varying, it
is uniform over the density on compact subsets of $(0,1)$, and it
provides an expansion up to the second order.

Let $R_r\colon \bb R \to (0,1)$ be the function
given by
\begin{equation}
\label{42c}
\cb{R_r (\varphi)} \,:=\, \frac{1}{|\Sigma^+_r|}\, 
\sum_{x\in \Sigma^+_r}   \varrho_\varphi (x/n) \,=\, 
\frac{1}{|\Sigma^+_r|}\, 
\sum_{x\in \Sigma^+_r}  \frac{e^\varphi \varrho(x/n)  }
{e^\varphi \varrho(x/n)  + [1- \varrho(x/n)] } \,\cdot
\end{equation}
It is easily seen that $R_r(\cdot)$ is strictly increasing, and
that $\lim_{\varphi\to+\infty} R_r (\varphi) = 1$,
$\lim_{\varphi\to -\infty} R_r (\varphi) = 0$.  Fix
$\mss m \in (0,1)$ and choose
$\varphi = \varphi (\varrho(\cdot), r, \mss m)$ to get
\begin{equation}
\label{50c}
\frac{1}{|\Sigma^+_r|}\, 
\sum_{x\in \Sigma^+_r}  \varrho_{\varphi} (x/n )
\, =\, \mss m \,.
\end{equation}
The difference between this subsection and the previous one is that,
while in the previous subsection $\varphi$ was a free parameter, it
here depends on $\mss m$ and is fixed by the previous
equation. Moreover, to lighten the notation we omit $\varphi$ from it
and write, for example, $\cb{\Gamma_N}$, $\cb{\Gamma_{N,A}}$ instead
of $\Gamma_{N,\varphi}$, $\Gamma_{N,A,\varphi}$, respectively.

For a finite set $A$ contained in $\Sigma^+_r$, let
\begin{equation}
\label{77}
\quad
\cb{\mtt s_A } \,:=\, 
\sum_{y\in A} [\varrho_\varphi (y/n) - 1] \,,
\quad 
\cb{\mtt t_{k,A} } \,:=\,
\sum_{y\in  A} 
\gamma_k(\varrho_\varphi  (y/n))\,,
\quad k\ge 2\,.
\end{equation}
In this formula, $\cb{\gamma_k(\alpha)}$, $\alpha\in (0,1)$, represents
the $k$-th cumulant of a Bernoulli distribution of parameter $\alpha$.
Let, furthermore,
\begin{equation}
\label{43b}
\begin{gathered}
\cb{\mtt T_{k} } \,:=\, 
\sum_{y\in \Sigma^+_r}  \gamma_k(\varrho_\varphi (y/n)) 
\,, \quad
\cb{\mtt r_{k}} \,:=\,
\frac{1}{\Gamma_{N}}
\sum_{y\in \Sigma^+_r}  \gamma_k(\varrho_\varphi (y/n))\,, \quad k\ge
3\,,
\\
\cb{\mtt T_{1,A}}\,:=\,  \frac{1}{2} \, \big\{ \, \mtt t_{2,A} - \mtt s_A^2
\,-\,  \mtt r_{3} \, \mtt s_A\,\big\} \,,
\quad
\cb{\mtt T_{2,A} }\, :=\, \mtt a_N \,+\, \mtt b^{(1)}_N \,+\, \mtt b^{(2)}_N \,,
\end{gathered}
\end{equation}
where $\Gamma_{N} $ has been introduced just before the statement of
Theorem \ref{s29}, and 
\begin{equation*}
\cb{ \mtt a_N}  \, :=\,
\frac{1}{4}\, \mtt r_{4}\, \mtt t_{2,A} 
\,-\, \frac{1}{8}\, \mtt t_{4, A} 
\,+\, \frac{5 }{ 12}\, \mtt r_{3} \, \mtt t_{3,A}
\,-\, \frac{5 }{8}\,  \mtt r^2_{3} \, \mtt t_{2,A}  
\, - \, \, \mtt t_{2,A} \,
\mtt r_{3} \, \mtt s_A
\,+\, \frac{1}{2}\,  \mtt t_{3,A} \, \mtt s_A  \,,
\end{equation*}
\begin{equation*}
\cb{\mtt b^{(1)}_N }  \, :=\,
\frac{3}{8} \, \mtt t_{2,A}^2
\,+\, \frac{1}{8}\, \mtt s^4_{A}
\,-\, \frac{3}{4}\, \mtt s^2_{A} \, \mtt t_{2,A}
\,-\,  \frac{1}{4}\, \mtt t_{2,A} \, \mtt r_{3}\, \mtt s_A \,,
\end{equation*}
\begin{equation*}
\cb{\mtt b^{(2)}_N}  \, :=\, 
\frac{5}{12} \, \mtt r_{3}\, \mtt s^3_A
\,-\, \frac{1}{4} \,  \mtt r_{4}\, \mtt s^2_A
\,+\, \frac{5}{8} \,  \mtt r^2_{3}\, \mtt s^2_A
\,+\,  \frac{1}{8} \,  \mtt r_{5}\, \mtt s_A
\,-\,  \frac{2}{3} \,  \mtt r_{3}\, \mtt r_{4}\,
\mtt s_A \,+\,  \frac{5}{8} \,  \mtt r^3_{3}\, \mtt s_A \,.
\end{equation*}
Notice that we omitted the empirical density $\mss m$ from the
notation.

For $0\le \mss M \le |\Sigma^+_{r}|$, denote by
$\cb{\nu^{\rm gc}_{r,\mss M} }$ the Bernoulli product measure on
$\Omega^+_r$ whose marginal at $x\in \Sigma^+_{r}$ has mean
$\varrho_{\varphi} (x/n )$, where $\varphi$ is selected according to
\eqref{50c}.  In particular,
$E_{\nu^{\rm gc}_{r,\mss M} } [\, \sum_{x\in \Sigma^+_{r}} \eta_x \,]
= \mss M$, and $\mss m = \mss M/|\Sigma^+_{r} |$.

\begin{proposition}
\label{s12b}
For each $q\ge 1$, $0<\delta<1/2$, there exists a finite constant
$C_1 = C_1(\delta, \mf r, q)$ such that
\begin{align*}
\sup_{\delta \le \mss m \le 1-\delta}\, 
\Big|\, E_{\nu^{\rm c}_{r,\mss M}  } [\, \eta_A \,]
\,-\,
E_{\nu^{\rm gc}_{r,\mss M}  } [\, \eta_A \,]
\, \Big( 1 \,+\, \frac{1}{\Gamma_N} \, \mtt T_{1,A}
\,+\, \frac{1}{\Gamma^2_N}\, \mtt T_{2,A} \,\Big)
\,\Big| \,\le\, \frac{C_1}{N^3}
\end{align*}
for all finite subsets $A$ of $\bb Z^d$ and $r \ge 1 $ such that
$A\subset \Sigma^+_r$, $|A|\le q$.
\end{proposition}

\begin{proof}
Fix $ \delta \le \mss m \le 1-\delta$, and recall that $\varphi$ was
selected to get \eqref{50c}.  By definition of the canonical measure,
\begin{equation}
\label{22c}
E_{\nu^{\rm c}_{r,\mss M} } \big[\, \eta_A\,\big]
\,=\,  
\frac{E_{\nu^{\rm gc}_{r,\mss M}  } 
\big[\,  \eta_A \,  \mtt 1\{
\sum_{y\in\Sigma^+_r}  \eta_y = \mss M\}\, \big]}
{\nu^{\rm gc}_{r,\mss M} \big[\, \sum_{y\in\Sigma^+_r}   \eta_y =
\mss M \,\big]} \,\cdot 
\end{equation}
As $\nu^{\rm gc}_{r,\mss M}$ is a product measure, and $\eta_A =1$
means $\eta_y=1$ for $y\in A$, the previous expression is equal to
\begin{equation}
\label{23b}
E_{\nu^{\rm gc}_{r,\mss M} } \big[\, \eta_A \, \big]\;
\frac{\nu^{\rm gc}_{r,\mss M}  
\big[\, \sum_{y\in\Sigma^+_r \setminus A}  \eta_y = \mss M - |A|\big]}
{\nu^{\rm gc}_{r,\mss M}  \big[\, \sum_{y\in\Sigma^+_r}   \eta_y =
\mss M \,\big]}  \,. 
\end{equation}

Consider the denominator of the previous expression.  In Theorem
\ref{s29}, set $p=7$ so that the error in the expansion is of order
$N^{-3}$. By the choice of $\varphi$, the variable $z$ introduced in
\eqref{74} is equal to $0$.  Thus, as $q_{{2j+1}, N} (0) = 0$ for
$j\ge 0$, by Theorem \ref{s29} with $p=7$,
\begin{equation}
\label{24b}
\sqrt{2\pi \, \Gamma_{N}} \,
\nu^{\rm gc}_{r,\mss M} \big[\,
\sum_{y\in\Sigma^+_r}   \eta_y = \mss M \,\big]
\,=\, 1 \, +\, \frac{\sqrt{2\pi} }{N} \, q_{2, N} (0)
\, +\, \frac{\sqrt{2\pi}}{N^2} \, q_{4, N} (0)
\,+\,   R^{(1)}_N(\mss m)  \,,
\end{equation}
where 
\begin{equation}
\label{75}
\sup_{\delta \le \mss m \le 1-\delta}\,
\big|\, R^{(1)}_N(\mss m) \, \big| \,\le\, 
\mf a_1 \, \frac{1}{N^{3}}\,\cdot
\end{equation}
Here and below, we adopt the convention established in Section
\ref{sec4}: $\mf a_1 = \mf a_1(\delta)$ represents a finite constant,
depending only on $\mf r$ and $\delta$, whose value may change from
line to line. The same convention is used for the remainder
$R^{(1)}_N(\mss m)$, which stands for a generic term satisfying
\eqref{75}.

We turn to the numerator of \eqref{23b}. Let
\begin{equation*}
\cb{\mf q_{j, N,A} (z) }\, :=\,  \sqrt{2\pi}\, e^{z^2/2}\,
q_{j, N,A} (z)\,,\quad
\cb{\mf q_{j, N} (z) }\, :=\,  \sqrt{2\pi}\, e^{z^2/2}\,
q_{j, N,\varnothing}  (z)\,,\quad
j \ge 1\,,\, z\in \bb R\,,
\end{equation*}
which cancels the factor $(2\pi)^{-1/2}e^{-z^2/2}$ in the expression
\eqref{78}.  Note that the functions $\mf q_{j, N, A} (\cdot)$ are
polynomials. Moreover, $\mf q_{2j ,N, A}(\cdot)$, $j\ge 1$, are sums
of monomials of even degree, while $\mf q_{2j -1,N, A}(\cdot)$,
$j\ge 1$, are sums of monomials of odd degree.

Since by \eqref{50c}
$\mss M = \sum_{x\in \Sigma^+_r} \varrho_\varphi(x/n)$, by Theorem
\ref{s29} with $p=7$,
\begin{equation}
\label{25b}
\sqrt{ \Gamma_{N,A}}\; \nu^{\rm gc}_{r,\mss M}  
\Big\{\, \sum_{y\in\Sigma^+_r \setminus A}  \eta_y = \mss M -
|A|\, \Big\}
\,=\,
\frac{1}{\sqrt{2\pi}} \, e^{-\mtt x^2/2} \Big\{ \, 1
\,+\, \sum_{j=1}^4 \frac{1}{N^{j/2}_A} \, \mf q_{j, N,A} (\mtt x)\, 
\,+\,   R^{(1)}_N(\mss m)   \, \Big\} \,,
\end{equation}
where $R^{(1)}_N(\mss m) $ satisfies \eqref{75}, $\Gamma_{N,A}$ is
defined in \eqref{79}, 
\begin{equation*}
\cb{ \mtt x } \, :=\, \frac{\mtt s_A}{\sqrt{\Gamma_{N,A}}}\, \,,
\end{equation*}
and $\mtt s_A$ has been introduced in \eqref{77}.  In \eqref{25b}, we
estimated $\exp\{ \mtt x^2/2\} \, R^{(1)}_N(\mss m)$ by
$R^{(1)}_N(\mss m) $.  This is possible because,
$\varrho_\varphi(\cdot)$ being bounded away from $0$ and $1$ uniformly
for $\delta \le \mss m\le 1-\delta$ (see \eqref{56b} for a proof of
this claim in a similar context), $\mtt x$ is uniformly bounded in
$A$, with $|A|\le q$, $N$, and $\delta \le \mss m\le 1-\delta$.
Actually, $|\mtt x| \le \mf a_1 (\delta, q)\, N^{-1/2}$. We also
incorporated in $R^{(1)}_N(\mss m)$ the term $j=5$ of the expansion
provided by Theorem \ref{s29}: as $\mf q_{5, N, A}(\cdot)$ is a sum of
monomials of odd degree, this term is bounded by
$\mf a_1 N^{-5/2}\, |\mtt x| \le \mf a_1 N^{-3}$.

Combining identities \eqref{22c}, \eqref{23b}, \eqref{24b}, \eqref{25b}
yields that
\begin{equation}
\label{33b}
\begin{aligned}
E_{\nu^{\rm c}_{r,\mss M} } \big[\, \eta_A\,\big]
\,=\,  & E_{\nu^{\rm gc}_{r,\mss M} } \big[\, \eta_A\,\big]  \; 
\frac{\sqrt{\Gamma_{N}}} {\sqrt{\Gamma_{N,A}}} \, e^{-\mtt x^2/2}\, 
\Big\{ \,  1 \,+\, \sum_{j=1}^4 \frac{1} {N^{j/2}_A} \,
\mf q_{j, N,A} (\mtt x)\, 
\,+\,  R^{(1)}_N(\mss m)  \, \Big\}
\\
&\times\, 
\Big\{ \, 1
\,-\,  \frac{1}{N} \, \mf q_{2, N} (0)\,
\,-\,  \frac{1 }{N^2} \, \mf q_{4, N} (0)\,
\,+\,  \frac{1}{N^2} \, \mf q_{2, N} (0)^2\, 
\,+\,  R^{(1)}_N(\mss m)  \, \Big\}
\,.
\end{aligned}
\end{equation}
We may further expand this expression by noting that
\begin{equation}
\label{35b}
\begin{gathered}
\frac{\sqrt{\Gamma_{N}}} {\sqrt{\Gamma_{N,A}}}  \,=\,
1\,+\, \frac{1}{2}\, \frac{\mtt t_{2,A}}{\Gamma_{N}}
\,+\, \frac{3}{8}\, \frac{\mtt t^2_{2,A}}{\Gamma^2_{N}}
\,+\,  O \Big( \frac{1}{N^{3}}\Big)\,, 
\\
e^{-\mtt x^2/2}
\,=\, 1\,-\, \frac{1}{2}\, \mtt x^2 \,+\, \frac{1}{8}\, \mtt x^4
\,+\,  O \Big( \frac{1}{N^{3}}\Big)\,,
\end{gathered}
\end{equation}
where $\mtt t_{2,A}$ has been introduced in \eqref{77}.  Rewrite
$\mtt x$ as $(\mtt s_A /\sqrt{\Gamma_{N}}) \sqrt{\Gamma_{N} /\Gamma_{N,A}}$ and
recall the first expansion in \eqref{35b} to conclude that
\begin{equation}
\label{44b}
e^{-\mtt x ^2/2}
\,=\, 1 \,-\, \frac{1}{2}\, \frac{\mtt s_A^2} {\Gamma_{N}}
\,+\, \frac{1}{2}\, \Big(\, \frac{\mtt s_A^4} {4}
- \mtt s_A^2\, \mtt t_{2,A}\,\Big) \, \frac{1} {\Gamma^2_{N}}
\,+\,  O \Big( \frac{1}{N^{3}}\Big)\,.
\end{equation}

To estimate the right-hand side in \eqref{33b}, we collect all terms
of the same order in $ N$ in the previous equation.

\smallskip\noindent{\bf (A0) Zeroth order term.}  We start with the
terms of zero-order in \eqref{33b}. By \eqref{35b}, \eqref{44b}, and
the definition of the polynomials $\mf q_{j,N,A} $, the zero-order
term is given by $E_{\nu^{\rm gc}_{r,\mss M} } \big[\, 
\eta_A\,\big]  $. 

\smallskip\noindent{\bf (A1) First order term.}  Denote by
$\cb{\mf q^{(k)}_{j, N, A}} (\cdot)$ the degree $k$ component of the
polynomial $\mf q_{j, N, A} (\cdot)$. We adopt the same notation for
the polynomials $\mf q^{(k)}_{j,N} (\cdot)$.

By \eqref{35b}, \eqref{44b}, and the definition of the polynomials
$\mf q_{j,N, A} $, the term of order $N^{-1}$ in \eqref{33b} is given by
\begin{equation}
\label{45b}
E_{\nu^{\rm gc}_{r,\mss M} } \big[\, \eta_A\,\big] \,  \Big\{ 
\frac{1}{2\, \Gamma_N} \, \mtt t_{2,A}\, -\, \frac{1}{2\, \Gamma_N}\, \mtt s_A^2
\,+\, \frac{1} {N^{1/2}_A} \,
\mf q^{(1)}_{1, N, A} (\mtt x)\,
\,+\, \frac{1} {N_A} \, \mf q^{(0)}_{2, N, A} (\mtt x)\, 
\,-\, \frac{1} {N} \, \mf q ^{(0)}_{2, N} (0) \Big\} \,.
\end{equation}
The first term comes from the expansion of
$\sqrt{\Gamma_N}/\sqrt{\Gamma_{N,A}}$, the second one from
\eqref{44b}. The next two from the sum inside braces in \eqref{25b}
and the last one from \eqref{24b}.

The difference of the last two terms may be seen to be of order
$N^{-2}$ and should be examined in the next step where we
consider the terms of this order. Therefore, by definition of the
polynomial $\mf q_{1, N, A} (\cdot)$, the expression inside braces
(omitting the difference of the last two terms) in the previous
displayed equation is equal to
\begin{equation*}
\frac{1}{2\, \Gamma_N}\,  (\mtt t_{2,A} - \mtt s_A^2) 
-\, 
\frac{\mtt T_3 - \mtt t_{3,A}  }{2\, \Gamma^2_{N,A}}\, \mtt s_A\,,
\end{equation*}
where $\mtt T_k$, $\mtt t_{k,A} $ have been introduced in equations
\eqref{43b} and \eqref{77}, respectively.

A Taylor
expansion similar to the previous ones yields that
\begin{equation}
\label{40b}
-\, \frac{1}{2}\, 
\frac{\mtt T_3 - \mtt t_{3,A} } {\Gamma^2_{N,A}}\, \mtt s_A
\,=\, -\, \frac{1}{2}\, \frac{\mtt r_{3} \, \mtt s_A} {\Gamma_N}
\,-\,  \frac{1}{2}\,  \frac{1}{\Gamma^2_{N}} \, \big\{ 2 \, \mtt t_{2,A} \,
\mtt r_{3} \,-\, \mtt t_{3,A}\,\big\} \, \mtt s_A
\,+\, O \Big( \frac{1}{N^{3}}\Big)\,,
\end{equation}
where $\mtt r_k$, $k\ge 3$,  has been introduced in \eqref{43b}.

The terms of order $N^{-2}$ in \eqref{40b} are kept for the next step
in the expansion, while the terms of order $N^{-1}$ are combined with
the previous terms of this order. After this operation one concludes
that the term of order $N^{-1}$ in \eqref{33b} is given by
\begin{equation*}
E_{\nu^{\rm gc}_{r,\mss M} } \big[\,  \eta_A\,\big] \,   
\frac{1}{2\, \Gamma_{N}} \, \big\{ \, \mtt t_{2,A} - \mtt s_A^2
\,-\,  \mtt r_{3} \, \mtt s_A\,\big\} \,.
\end{equation*}

\smallskip\noindent{\bf (A2) Second order term.} We first collect all
terms of order $N^{-2}$ which appeared in (A1). By \eqref{45b} and
\eqref{40b}, they are
\begin{equation}
\label{76}
E_{\nu^{\rm gc}_{r,\mss M} } \big[\,  \eta_A\,\big] \,
\Big\{ \, \frac{1} {N_A}
\, \mf q^{(0)}_{2, N,A} (\mtt x)\, 
\,-\, \frac{1} {N} \, \mf q^{(0)}_{2, N} (0) \,-\,
\frac{1}{2}\, \frac{1}{\Gamma^2_{N} } \, \big\{ 2 \, \mtt t_{2,A} \,
\mtt r_{3} \,-\, \mtt t_{3,A}\,\big\} \, \mtt s_A  \,\Big\} \,.
\end{equation}
Since the zero-order term of the polynomial corresponds to the
constants, $\mf q^{(0)}_{2, N,A} (\mtt x)
= \mf q^{(0)}_{2, N, A} (0)$. 
By the explicit formula for $N^{-1}_A\mf q_{2, N, A}$,
\begin{equation*}
\begin{aligned}
&  \frac{1} {N_A}
\, \mf q^{(0)}_{2, N, A} (0)\, 
\,-\, \frac{1} {N} \, \mf q^{(0)}_{2, N} (0)
\\
& \quad =\,
\frac{1}{\Gamma^2_{N}}\, \Big\{\, 
\frac{H_4(0)}{4!}\, \big[ \, 2\, \mtt r_{4}\, \mtt t_{2,A} 
\,-\, \mtt t_{4,A}\, \big] 
\,+\, \frac{H_6(0)  }{ 2\, (3!)^2}\, \, \big[\, 3\, \mtt r^2_{3}
\, \mtt t_{2,A}\, \,-\, 2 \, \, \mtt r_{3} \,
\mtt t_{3,A}\,\big] \Big\}
\, +\, O\Big( \frac{1}{N^3}\Big) 
\\
& \quad =\,
\frac{1}{\Gamma^2_{N}}\, \Big\{\, 
\frac{1}{4}\, \mtt r_{4}\, \mtt t_{2,A} 
\,-\, \frac{1}{8}\, \mtt t_{4, A} 
\,+\, \frac{5 }{ 12}\, \mtt r_{3} \, \mtt t_{3,A}
\,-\, \frac{5 }{8}\,  \mtt r^2_{3} \, \mtt t_{2,A}  \Big\}
\, +\, O\Big( \frac{1}{N^3}\Big) \,.
\end{aligned}
\end{equation*}
In this formula, $H_p(\cdot)$ is the Hermite polynomial of degree $p$,
and $\mtt t_{p,A}$, $\mtt r_{p}$ were introduced in \eqref{77},
\eqref{43b}. Therefore, the expression in \eqref{76} is equal to
\begin{equation*}
E_{\nu^{\rm gc}_{r,\mss M} } \big[\,  \eta_A\,\big] \,
\frac{1}{\Gamma^2_{N}}\, \mtt a_N
\, +\, O\Big( \frac{1}{N^3}\Big) \,,
\end{equation*}
where $\mtt a_N$ is defined just before the statement of the
proposition. 

On the other hand, by \eqref{33b}, \eqref{35b}, \eqref{44b}, and the
definition of the polynomials $\mf q_{j,N, A} $, the remaining terms
of order $N^{-2}$ in \eqref{33b} are given by
$\mtt B_N\, E_{\nu^{\rm gc}_{r,\mss M} } \big[\, \eta_A\,\big] $,
where
\begin{equation*}
\begin{aligned}
\mtt B_N   \, & =\,
\frac{1}{8\, \Gamma^2_{N}}\, \Big\{
3\, \mtt t_{2,A}^2
\,+\, \mtt s^4_{A}
\,-\, 6\, \mtt s^2_{A} \, \mtt t_{2,A}\, \Big\}
\,+\, \frac{1}{2\, \Gamma_{N}}\, \big\{\, \mtt t_{2,A} - \mtt
s^2_{A}\, \big\}\, \frac{1}{\sqrt{N}}\, \mf q_{1,N}^{(1)}(\mtt x) 
\\
& \,+\, \frac{1}{\sqrt{N}}\, \mf q_{1,N}^{(3)}(\mtt x)
\,+\, \frac{1}{N}\, \mf q_{2,N}^{(2)}(\mtt x) 
\,+\, \frac{1}{N^{3/2}}\, \mf q_{3,N}^{(1)}(\mtt x) 
\,-\, \frac{1}{N^{3/2}}\, \mf q_{2,N}^{(0)}(0) \,
\mf q_{1,N}^{(1)}(\mtt x)
\, +\, O\Big( \frac{1}{N^3}\Big)\,. 
\end{aligned}
\end{equation*}
To obtain this identity we observed that
\begin{gather*}
\frac{1} {N_A} \, \mf q_{2,N,A}^{(0)}(\mtt x)
\,-\, \frac{1} {N}  \, \mf q_{2,N}^{(0)}(0)
\,=\, O\Big( \frac{1}{N^2}\Big)\,,
\quad
\frac{1} {N^{2}_A} \, \mf q_{4,N,A}^{(0)}(\mtt x)
\,-\, \frac{1} {N^{2}} \, \mf q_{4,N}^{(0)}(0)
\,=\, O\Big( \frac{1}{N^3}\Big)\,,
\\
\frac{1}{\sqrt{ N_A} }\,  \mf q_{1,N,A}^{(1)}(\mtt x)
\,-\, \frac{1}{\sqrt{N}  }\,  \mf q_{1,N}^{(1)}(\mtt x)
\,=\, O\Big( \frac{1}{N^2}\Big)\,,
\quad
\frac{1}{N^{j/2}_A}\,  \mf q_{j,N,A}^{(4-j)}(\mtt x)
\,-\, \frac{1}{N^{j/2} }\,  \mf q_{j,N}^{(4-j)}(\mtt x)
\,=\, O\Big( \frac{1}{N^3}\Big)\,.
\end{gather*}
By the explicit formulae for the polynomials
$\mf q_{j,N}(\cdot)$,
\begin{equation*}
\begin{aligned}
\mtt B_N  \, & =\,
\frac{1}{8\, \Gamma^2_{N}}\, \Big\{
3\, \mtt t_{2,A}^2
\,+\, \mtt s^4_{A}
\,-\, 6\, \mtt s^2_{A} \, \mtt t_{2,A}\, \Big\}
\,-\, \frac{1}{4\, \Gamma^2_{N}}\, \big\{\, \mtt t_{2,A} - \mtt
s^2_{A}\, \big\}\, \mtt r_{3}\, \mtt s_A
\\
& \,+\, \frac{1}{\Gamma^2_{N}} \,\Big\{\,
\frac{1}{6} \, \mtt r_{3}\, \mtt s^3_A
\,-\, \frac{1}{4} \,  \mtt r_{4}\, \mtt s^2_A
\,+\, \frac{5}{8} \,  \mtt r^2_{3}\, \mtt s^2_A
\,+\,  \frac{1}{8} \,  \mtt r_{5}\, \mtt s_A
\,-\,  \frac{35}{48} \,  \mtt r_{3}\, \mtt r_{4}\,
\mtt s_A
\\
& \qquad\qquad 
\,+\,  \frac{35}{48} \,  \mtt r^3_{3}\, \mtt s_A
\,+\,  \frac{1}{16} \,  \mtt r_{4}\, \mtt r_{3}\, \mtt s_A
\,-\,  \frac{5}{48} \,  \mtt r^3_{3}\, \mtt s_A \,\Big\}
\, +\, O\Big( \frac{1}{N^3}\Big)\,.
\end{aligned}
\end{equation*}
That is, up to an error of order $N^{-3}$, 
\begin{equation*}
\begin{aligned}
\mtt B_N  \, & =\,
\frac{1}{8\, \Gamma^2_{N}}\, \Big\{
3\, \mtt t_{2,A}^2
\,+\, \mtt s^4_{A}
\,-\, 6\, \mtt s^2_{A} \, \mtt t_{2,A}
\,-\,  2 \,  \mtt t_{2,A} \, \mtt r_{3}\, \mtt s_A \, \Big\}
\\
& \,+\, \frac{1}{\Gamma^2_{N}} \,\Big\{\,
\frac{5}{12} \, \mtt r_{3}\, \mtt s^3_A
\,-\, \frac{1}{4} \,  \mtt r_{4}\, \mtt s^2_A
\,+\, \frac{5}{8} \,  \mtt r^2_{3}\, \mtt s^2_A
\,+\,  \frac{1}{8} \,  \mtt r_{5}\, \mtt s_A
\,-\,  \frac{2}{3} \,  \mtt r_{3}\, \mtt r_{4}\,
\mtt s_A \,+\,  \frac{5}{8} \,  \mtt r^3_{3}\, \mtt s_A\,\Big\} \,.
\end{aligned}
\end{equation*}
With the notation introduced before the statement of the proposition,
$\mtt B_N = (1/\Gamma^2_{N})\, (\mtt b^{(1)}_N + \mtt b^{(2)}_N)$. To
complete the proof of the proposition for
$r\ge r_0 (\delta, \mf r,q)$, it remains to collect the previous
estimates. The result holds for $1\le r\le r_0$ by adjusting the
constant $C_1$.
\end{proof}

\begin{remark}
\label{s70}
If the density profile $\varrho(\cdot)$ is constant equal to $\rho$,
then a simple calculation yields that the first order term is equal to
$-\, (1/2N) \, \mtt p_A '' (\rho)\, \chi(\rho) = (1/N) \, \mtt q_A(\rho)$:
\begin{equation*}
E_{\nu^{\rm gc}_{r, \mss M}  } [\, \eta_A \,]
\frac{1}{\Gamma_N} \, \mtt T_{1,A} \,=\,
\frac{1}{N}\, \mtt q_A(\rho) \,.
\end{equation*}
The first order term can also be represented as $(1/N)\, \mtt p_1
(\rho)$, where 
\begin{equation*}
\mtt p_1(\rho)  \,:=\,  \frac{1}{2}\, \Big\{\,
\frac{\gamma_3}{\gamma_2^{2}}\, 
E_{\nu_\rho} \Big[\, \eta_A\,;\, \mc S_A \,\Big] 
\,-\,  \frac{1}{\gamma_2}\, 
E_{\nu_\rho} \Big[\, \eta_A\,;\, \mc S_A^2\,\Big]
\, \Big\} \,,
\end{equation*}
where $\cb{\mc S_A}: = \sum_{x\in A} (\eta_x - \rho)$, and we recover the
formula for the first order correction \eqref{37b} obtained in Lemma
\ref{s05}.

A similar formula can be derived for the second term in the expansion,
which is equal to $(1/N)^2\, \mtt p_2 (\rho)$, where
\begin{equation*}
\begin{aligned}
\mtt p_2(\rho) \;=& \; \frac{3|A|}{2}\, \mtt p_1(\rho)
\;+\;\Big[\frac{2\gamma_3\gamma_4}{3\gamma_2^{4}}
-\frac{5\gamma_3^{3}}{8\gamma_2^{5}}
-\frac{\gamma_5}{8\gamma_2^{3}}\Big]\, E_{\nu_\rho}\big[\eta_A \,;\, \mc
S_A\big]
\\
+&
\;\Big[\frac{5\gamma_3^{2}}{8\gamma_2^{4}}
-\frac{\gamma_4}{4\gamma_2^{3}}\Big]\,
E_{\nu_\rho}\big[\eta_A \,;\, \mc S_A^{2}\big]
\;-\;\frac{5\gamma_3}{12\gamma_2^{3}}\, E_{\nu_\rho}\big[\eta_A \,;\, \mc
S_A^{3}\big] \;+\;\frac{1}{8\gamma_2^{2}}\,
E_{\nu_\rho}\big[\eta_A \,;\, \mc S_A^{4}\big]\,.
\end{aligned}
\end{equation*}
These representations have the advantage that they easily extend to
cylinder functions $h$.

In the inhomogeneous case, the first order term is equal to
\begin{equation*}
\frac{1}{2}\, \frac{1}{\Gamma_{N}} \,  \Big\{\,
\frac{\mtt T_{3} }{\Gamma_{N}} \, 
E_{\nu^{\rm gc}_{r, \mss M}  } \Big[\, \eta_A\,;\, \sum_{x\in A}
[\, \eta_x - \varrho_\varphi (x/n) \,]\, \,\Big] 
\,-\,  
E_{\nu^{\rm gc}_{r, \mss M}  }  \Big[\, \eta_A\,;\, \Big(\sum_{x\in A}
[\, \eta_x - \varrho_\varphi (x/n) \,]\, \Big)^2\,\Big]
\, \Big\} \,,
\end{equation*}
an expression which extends to product inhomogeneous measures the
formula for the first order correction in the equivalence of
ensembles.  We obtained this formula inserting the constant $\mtt s_A$
inside the expectation and then replacing $\varrho_\varphi (y/n) - 1$
by $\varrho_\varphi (y/n) - \eta_y$, which is allowed due to the
presence of the function $\eta_A$ inside the expectation.

\end{remark}

The next consequence of Proposition \ref{s12b} is needed in Section
\ref{sec4}. Assume that $n = r_1\, r_2$ for some positive integers
$r_1\ge 2$, $r_2\ge 2$. Recall the definition of the cubes
$B_{\mtt k}$, $\mtt k \in \Sigma^+_{r_2}$, introduced above equation
\eqref{156}. Fix a cylinder function $h$ and recall from equation
\eqref{95} that $B^o_{\mtt k}$ represents the $h$-interior of the set
$B_{\mtt k}$.

\begin{corollary}
\label{s42}
Let $h$ be a local function, and $J\colon \bb T^d \to \bb R$ a
continuous function, and fix $0<\delta < 1/2$. Then, there exists a
finite constant $\mf a_2=\mf a_2(\delta)$, depending only on $\mf r$,
$h$ and $\delta$, such that
\begin{equation*}
E_{\nu^{\rm c}_{\mtt k,\mss M}  }  \Big [\,
\Big( \sum_{x\in B^o_{\mtt k}} J(x/n)\,  \big\{ \tau_x h
-  E_{\nu^{\rm c}_{\mtt k, \mss M}} [\tau_x h] \,\big\}
\Big)^2  \, \Big]
\,\le\,
\mf a_2 \, \Vert J\Vert^2_\infty \, \, r_1^d
\end{equation*}
for all $\mtt k\in \Sigma^+_{r_2}$, $\delta \le \mss M/|B_{\mtt k}| \le
1-\delta$, $r_1\ge 1$, $r_2\ge 1$.
\end{corollary}

\begin{proof}
The expectation appearing in the statement of the corollary is equal to
\begin{equation*}
\sum_{x, y \in B^o_{\mtt k}} J(x/n)\, J(y/n)\, 
E_{\nu^{\rm c}_{\mtt k,\mss M}  }  \Big [\,
\big\{ \tau_x h -  E_{\nu^{\rm c}_{\mtt k, \mss M}} [\tau_x h] \,\big\}
\,
\big\{ \tau_y h -  E_{\nu^{\rm c}_{\mtt k, \mss M}} [\tau_y h] \,\big\}
\, \Big]\,.
\end{equation*}
By Proposition \ref{s12b}, in the previous formula, we may replace all
expectations with respect to the canonical measure
$\nu^{\rm c}_{\mtt k, \mss M}$ by ones with respect to the
grand-canonical measures $\nu^{\rm gc}_{\mtt k, \mss M}$ at a cost
bounded by $\mf a_2(\delta) \, r_1^{-d}$. After this replacement, all
non-diagonal terms vanish. As there are only $O(r_1^d)$ diagonal
terms, the corollary is proved.
\end{proof}

In the next results we adopt the notation introduced at the beginning
of Section \ref{sec9}. Recall, in particular, the definition of the
measures $\nu^{\rm c} _{2,\mtt j, \mss M} $ and
$\nu^{\rm gc} _{2,\mtt j, \mss M} $.  The assertion of Corollary
\ref{s18} below is an immediate consequence of Proposition \ref{s12b}
applied to $\Sigma^+_r = B_{2,\mtt j}$ for $\mtt t_1 \mtt t_2$
sufficiently large. For the remaining values it is a matter of
adjusting the constant $\mf a_1$ since the left-hand side is bounded.
We just need some notation. Fix $\mtt j \in \Sigma^+_{\mtt t_3}$, and
let
\begin{equation*}
\cb{\Gamma_{2,\mtt j} } \,:=\,
\sum_{y\in  B_{2, \mtt j}}
\chi (\varrho^{(2)}_{\mtt j, \mss m} (y/n)) \,, \quad
\cb{\mtt r_{3, \mtt j}} \,:=\,
\frac{1}{\Gamma_{2,\mtt j}}
\sum_{y\in B_{2,\mtt j}}
\gamma_3(\varrho^{(2)}_{\mtt j, \mss m}(y/n)) \,.
\end{equation*}
For a finite subset $A$ of  $B_{2, \mtt j}$, let
\begin{equation*}
\cb{ \mtt s_{\mtt j, A}  } \,:=\, 
\sum_{y\in A} [\varrho^{(2)}_{\mtt j, \mss m}  (y/n) - 1] \,,
\quad 
\cb{\mtt t_{2,\mtt j,A} } \,:=\,
\sum_{y\in  A} 
\chi (\varrho^{(2)}_{\mtt j, \mss m}   (y/n))\,,
\quad \cb{\mtt T_{1, \mtt j , A} }\,:=\,
\frac{1}{2} \, \big\{ \, \mtt t_{2, \mtt j ,A} - \mtt s_{\mtt j, A}^2
\,-\,  \mtt r_{3, \mtt j} \, \mtt s_{\mtt j, A} \,\big\} \,.
\end{equation*}

\begin{corollary}
\label{s18}
For each $q\ge 1$, $\delta>0$, there exists a finite constant $\mf a_1
= \mf a_1(q, \delta)$ such that
\begin{equation*}
\max_{\mtt j\in \Sigma^+_{\mtt t_3}}
\sup _{\delta \le \mss m \le 1-\delta}
\max_{\mtt k\in \Lambda_{\mtt j} }
\max_{A}\,
\Big | \, E_{\nu^{\rm c} _{2, \mtt j, \mss M} }
\big [ \,   \eta_{A}  \,\big] 
\; -\; E_{\nu^{\rm gc} _{2, \mtt j, \mss M} }
\big [ \, \eta_{A}  \,\big]
\, \Big( \, 1 \,+\, \frac{1}{\Gamma_{2,\mtt j}} \,
\mtt T_{1, \mtt j , A}  \,\Big)
\, \Big|
\; \le \; \frac {\mf a_1}{(\mtt t_1 \mtt t_2)^{2d}}
\end{equation*}
for all $\mtt t_1\ge 2$, $\mtt t_2\ge 2$. In this formula,
$\mss m = \mss M/|\Sigma^+_{\mtt t_1 \mtt t_2}|$, and the last maximum
is carried over all subsets $A$ of $B_{\mtt k}$ with $q$ elements.
\end{corollary}

As $|\, \mtt T_{1, \mtt j , A}/\Gamma_{2,\mtt j}\,|$ is uniformly
bounded by $\mf a_1(q,\delta)/(\mtt t_1 \mtt t_2)^{d}$, we also obtain
from Proposition \ref{s12b} the classical equivalence of ensembles:
\begin{equation}
\label{178}
\max_{\mtt j\in \Sigma^+_{\mtt t_3}}
\sup _{\delta \le \mss m \le 1-\delta}
\max_{A}\,
\Big | \, E_{\nu^{\rm c} _{2, \mtt j, \mss M} }
\big [ \,   \eta_{A}  \,\big] 
\; -\; E_{\nu^{\rm gc} _{2, \mtt j, \mss M} }
\big [ \, \eta_{A}  \,\big]
\, \Big|
\; \le \; \frac {\mf a_1}{(\mtt t_1 \mtt t_2)^{d}}\,,
\end{equation}
where the last maximum is carried over all subsets $A$ of
$B_{2,\mtt j}$ with $q$ elements.

Recall from \eqref{82} the definition of the expectation $\mss m_{\mtt
j, \mtt k}$, $\mtt j\in \Sigma^+_{\mtt t_3}$, $\mtt k\in \Lambda_{\mtt
j}$. 

\begin{corollary}
\label{s71}
For each $\delta>0$, there exists a finite constant
$\mf a_1 =\mf a_1(\delta)$ such that
\begin{gather*}
\max_{\mtt j\in \Sigma^+_{\mtt t_3}}
\sup _{\delta \le \mss m \le 1-\delta}
\max_{\mtt k' \neq \mtt k\in \Lambda_{\mtt j} }
\Big | \, E_{\nu^{\rm c} _{2, \mtt j, \mss M} }
\big [ \,   (\mss m_{\mtt k} - \mss m_{\mtt j, \mtt k}) \,
(\mss m_{\mtt k'}  - \mss m_{\mtt j, \mtt k'})  \,\big] 
\, \Big|
\; \le \; \frac {\mf a_1}{(\mtt t_1 \mtt t_2)^{d}}\,,
\\
\max_{\mtt j\in \Sigma^+_{\mtt t_3}}
\sup _{\delta \le \mss m \le 1-\delta}
\max_{\mtt k\in \Lambda_{\mtt j} }
\Big | \, E_{\nu^{\rm c} _{2, \mtt j, \mss M} }
\big [ \,   (\mss m_{\mtt k} - \mss m_{\mtt j, \mtt k})^2  \,\big] 
\, \Big|
\; \le \; \frac {\mf a_1}{\mtt t_1^{d}}
\end{gather*}
for all $\mtt t_1\ge 2$, $\mtt t_2\ge 2$. 
\end{corollary}

\begin{proof}
Rewrite the expectation appearing in the statement of the corollary as
\begin{gather*}
\frac{1}{|B_{\mtt k}|^2}\, \sum_{x\in B_{\mtt k}}\,
\sum_{y\in B_{\mtt k'}}
\big\{\, E_{\nu^{\rm c} _{2, \mtt j, \mss M} }
\big [ \, \eta_x\, \eta_y   \,\big]
\,-\,
E_{\nu^{\rm c} _{2, \mtt j, \mss M} }  [ \, \eta_x\,]
\,
E_{\nu^{\rm c} _{2, \mtt j, \mss M} }  [ \, \eta_y   \,]\,\big\}\,.
\end{gather*}
Replace the canonical measures by the grand canonical ones, at a cost,
according to \eqref{178}, of order $\mf a_1(\delta)/(\mtt t_1 \mtt
t_2)^{d}$. With the grand canonical measures replacing the canonical
ones, the previous sum vanishes. This completes the proof of the first
assertion. The same proof applies to the second assertion. The
difference is that we have diagonal terms which are roughly estimated
by $1$, providing a worse estimate of order $\mtt t_1^{-d}$ only.
\end{proof}

\section{Inhomogeneous product measures}
\label{sec12}

In this section, we adopt the setting introduced at the beginning of
Section \ref{sec9}, changing the notation. Here, $s_1$, $s_2$, $s_3$
are positive integers such that $s_1s_2s_3=n$. We represent the cubes
$B_{\mtt k}$, $B_{2, \mtt j}$, by $D_{\mtt k}$, $D_{2, \mtt j}$,
respectively. In contrast with Section \ref{sec9}, we admit the
possibility that $s_1=1$. With this choice of $s_1$ and setting
$r_1=s_2$, $r_2=s_3$, the set $D_{2, \mtt j}$,
$\mtt j \in \Sigma^+_{s_3}$, of this section corresponds to the cube
$B_{\mtt j}$ of Section \ref{sec4}.  In particular, we may apply the
results stated in this section for the level-2 cubes $D_{2, \mtt j}$ to
the level-1 cubes $B_{\mtt j}$ of Section \ref{sec4}.

We also adopt the notation introduced in Section \ref{sec4}:
$\cb{\mf a_1}$ stands for a finite and positive constant which depends
only on $\mf r$, $\cb{\mf a_2}$ for one which depends only on $\mf r$
and a cylinder function $h$ (which sometimes will be $\eta_A$),
$\cb{\mf a_{3,p}}$, $p\ge 1$, for one which depends only on $\mf r$, a
cylinder function $h$, and $\Vert \rho\Vert_{C^p}$. All these
constants are allowed to change from line to line. It may happen that
these constants also depend on an extra variable (say some $\delta>0$
which prevents the empirical density to be close to $0$ and $1$). In
this case this dependence appears explicitly in the notation, and we
write, for instance, $\mf a_2(\delta)$ instead of $\mf a_2$.

Fix a density profile $\rho\colon \bb T^d \to (0,1)$ of class
$C^3(\bb T^d)$ satisfying \eqref{83}.  Although all quantities in this
section depend on the density profile $\rho(\cdot)$, it is omitted
from the notation.

Recall from \eqref{97} the definition of the density profile
$\varrho^{(2)}_{\mtt j, \vartheta} (\cdot)$, $0< \vartheta < 1$,
$\mtt j\in \Sigma^+_{s_3}$.  We derive some properties of
$\varrho^{(2)}_{\mtt j, \vartheta} (\cdot)$ needed later in this
section. We first claim that
$\varrho^{(2)}_{\mtt j, \vartheta} (\cdot)$ is bounded away from $0$
and $1$. More precisely, that for each $0<\delta<1/2$, there exists a
finite constant $\mf a_1= \mf a_1(\delta, \mf r) \in (0,1/2)$ such
that
\begin{equation}
\label{56b}
\mf a_1 \,\le\, \varrho^{(2)}_{\mtt j, \vartheta} (x/n) \,\le \, 1\,-\,
\mf a_1
\end{equation}
for all $x\in D_{2, \mtt j}$, $\mtt j\in \Sigma^+_{s_3}$,
$\delta\le \vartheta \le 1-\delta$. We adopt here the same convention
as in the previous sections, and $\mf a_1$ represents a finite
constant which depends only on $\delta$ and $\mf r$, and whose value
may change from line to line.

To prove this assertion, recall from \eqref{42} the definition of the
function $R^{(2)}_{\mtt j}(\cdot)$, and observe that the maps
$\rho(\cdot) \mapsto R^{(2)}_{\mtt j}(\varphi)$,
$\varphi \mapsto R^{(2)}_{\mtt j}(\varphi)$ are both increasing.  To
avoid superscripts we represent below $R^{(2)}_{\mtt j}(\cdot)$,
$\Phi^{(2)}_{\mtt j}(\cdot)$ by $\cb{R_{2,\mtt j}(\cdot)}$,
$\cb{\Phi_{2, \mtt j}(\cdot)}$, respectively.  Thus, by \eqref{50},
\eqref{97}, \eqref{83},
\begin{align*}
\vartheta \,=\, \frac{1}{|D_{2, \mtt j}|}\, 
\sum_{x\in D_{2, \mtt j}}  \varrho^{(2)}_{\mtt j, \vartheta} (x/n )
\,\le\, \frac{e^{\Phi_{2, \mtt j}(\vartheta) }  (1-\mf r) }
{e^{\Phi_{2, \mtt j}(\vartheta) }  (1-\mf r)  + \mf r  }\,\cdot
\end{align*}
A similar lower bound can be obtained for  $\vartheta$. It follows
from these bounds that
\begin{align}
\label{98}
\frac{ \vartheta } { 1-\vartheta } \,
\frac{ \mf r  } { 1-\mf r  }
\,\le\,
e^{\Phi_{2, \mtt j} (\vartheta)}
\,\le\,
\frac{ \vartheta } { 1-\vartheta } \,
\frac { 1-\mf r }  { \mf r } 
\,\cdot
\end{align}
Reporting these bounds to \eqref{97} yields that
\begin{equation*}
\frac{\delta\, \mf r^2}
{\delta\, \mf r^2 + (1-\delta)\, (1-\mf r)^2}
\,\le\, 
\varrho^{(2)}_{\mtt j, \vartheta} (x/n) \,\le\,
\frac{(1-\delta) \, (1-\mf r)^2}
{(1-\delta) \, (1-\mf r)^2  + \delta\, \mf r^2}
\end{equation*}
for all $\delta \le \vartheta \le 1-\delta$, $x\in D_{2, \mtt j}$,
$\mtt j\in \Sigma^+_{s_3}$, as claimed.

The next assertion provides a formula for the gradient of
$\varrho^{(2)}_{\mtt j, \vartheta}$, represented by
$\cb{\nabla \varrho^{(2)}_{\mtt j, \vartheta}}$. We claim that
\begin{equation}
\label{53}
\frac{ (\nabla \varrho^{(2)}_{\mtt j, \vartheta})(x/n) }
{\varrho^{(2)}_{\mtt j, \vartheta} (x/n)^2}
\,=\,
e^{-\Phi_{2, \mtt j}(\vartheta)} \, \frac{ (\nabla \rho) (x/n) }
{\rho (x/n)^2} 
\end{equation}
for all $x\in D_{2, \mtt j}$. Indeed, differentiating both sides of
\eqref{97} yields that
\begin{equation*}
(\nabla \varrho^{(2)}_{\mtt j, \vartheta})(\mtt x) \,=\,
\frac{e^{\Phi_{2, \mtt j}(\vartheta) }  (\nabla \rho) (\mtt x)  }
{\{ e^{\Phi_{2, \mtt j}(\vartheta) }  \rho(\mtt x)  + [1- \rho(\mtt x)]
\}^2}\,
\,=\, \varrho^{(2)}_{\mtt j, \vartheta}(\mtt x)^2 \,
e^{- \Phi_{2, \mtt j}(\vartheta) } 
\frac{ (\nabla \rho) (\mtt x)  }  {\rho (\mtt x) ^2}\,,
\end{equation*}
as claimed. In particular, by \eqref{83}, \eqref{56b} and \eqref{98},
there exists a finite constant $\mf a_1= \mf a_1(\delta, \mf r)$ such
that
\begin{equation}
\label{57}
\big|\, (\nabla \varrho^{(2)}_{\mtt j, \vartheta})(x/n) \,\big|\,\le\,
\mf a_1 \,  \Vert \nabla \rho\Vert_\infty
\end{equation}
for all $\delta \le \vartheta \le 1-\delta$, $x\in D_{2, \mtt j}$,
$\mtt j\in \Sigma^+_{s_3}$.

We will need below to compare the chemical potential
$\Phi_{2, \mtt j}(\cdot)$ to the one associated with a constant
density profile. Let $\bar\rho_{2, \mtt j}$ be the average of
$\rho(\cdot)$ on $D_{2, \mtt j}$ and
$\overline{R}_{2, \mtt j}\colon \bb R \to (0,1)$ be the associated
strictly increasing function given by
\begin{equation}
\label{64}
\cb{\bar\rho_{2, \mtt j}} \,:=\, \frac{1}{|D_{2, \mtt j}|}\, 
\sum_{x\in D_{2, \mtt j}}   \rho(x/n)  \,, \quad
\cb{\overline{R}_{2, \mtt j} (\varphi)}
\,:=\, \frac{e^\varphi \bar\rho_{2, \mtt j}  }
{e^\varphi \bar\rho_{2, \mtt j}  + [1- \bar\rho_{2, \mtt j}] } \,\cdot
\end{equation}
Denote by $\cb{\overline{\Phi}_{2, \mtt j}(\cdot)}$ the inverse of
$\overline{R}_{2, \mtt j}(\cdot)$. For each density
$\vartheta\in (0,1)$,
\begin{equation}
\label{63}
\vartheta \, =\, 
\frac{e^{\overline{\Phi}_{2, \mtt j} (\vartheta) }  \bar\rho_{2, \mtt j}  }
{e^{\overline{\Phi}_{2, \mtt j} (\vartheta) }  \bar\rho_{2, \mtt j} + [1- \bar\rho_{2, \mtt j} ] }\,,
\quad \text{so that}\quad
e^{\overline{\Phi}_{2, \mtt j} (\vartheta) }  \,=\, \frac{\vartheta\, (1-\bar\rho_{2, \mtt j})}
{\bar\rho_{2, \mtt j} (1- \vartheta)}\,
\cdot
\end{equation}

We claim that there exists a finite constant
$\mf a_1= \mf a_1(\mf r, \delta)$ such that
\begin{equation}
\label{59}
\big|\,  e^{\Phi_{2, \mtt j}(\vartheta) -
\overline{\Phi}_{2, \mtt j} (\vartheta) } \,-\, 1 \,\big|\,\le\,
\mf a_1 \,
\Big( \, \Vert \nabla \rho\Vert^2_\infty\,
+ \Vert \nabla^2 \rho\Vert_\infty\, \Big)\, 
\Big(\frac{s_1s_2}{n}\Big)^2
\end{equation}
for all $\delta \le \vartheta \le 1- \delta$,
$\mtt j\in \Sigma^+_{s_3}$, $s_1 \ge 1$, $s_2\ge 2$. Here,
$\cb{\nabla^2\rho}$ represents the Hessian of $\rho(\cdot)$.

To prove \eqref{59}, we first obtain a formula for the expression on
the left-hand side. By \eqref{50} and the second equation in
\eqref{63},
\begin{equation*}
1 \,-\, \vartheta \,=\,
\frac{1}{|D_{2, \mtt j}|}\,  \sum_{x\in D_{2, \mtt j}}
\varrho^{(2)}_{\mtt j, \vartheta}(x/n)\, 
e^{- \Phi_{2, \mtt j} (\vartheta) }\,
\frac{1 -   \rho (x/n)} {\rho(x/n) }
\quad\text{and}\quad
1- \vartheta \,=\,
\vartheta\, e^{- \overline{\Phi}_{2, \mtt j} (\vartheta) }  
\frac {1-\bar\rho_{2, \mtt j}}{\bar\rho_{2, \mtt j} } \cdot
\end{equation*}
Thus, writing that $0 = (1 \,-\, \vartheta) - (1 \,-\, \vartheta)$,
and replacing $\vartheta$ in the last expression by \eqref{50} yields that
\begin{equation*}
0 \,=\, 
\frac{1}{|D_{2, \mtt j}|}\,  \sum_{x\in D_{2, \mtt j}}
\varrho^{(2)}_{\mtt j, \vartheta}(x/n)\, 
\Big\{\, 
e^{- \Phi_{2, \mtt j} (\vartheta) }\,
\frac{1 -   \rho (x/n)} {\rho(x/n) }
- e^{- \overline{\Phi}_{2, \mtt j}(\vartheta) } \,
\frac{1 -   \bar\rho_{2, \mtt j}}
{\bar\rho_{2, \mtt j} } \, \Big\} \,.
\end{equation*}
Adding and subtracting the cross terms yields that
\begin{equation}
\label{62}
e^{\Phi_{2, \mtt j}(\vartheta) - \overline{\Phi}_{2, \mtt j} (\vartheta) } \,-\, 1
\,=\, \frac{1}{(1-\bar\rho_{2, \mtt j}) \, \vartheta}\,
\frac{1}{|D_{2, \mtt j}|}\,
\sum_{x\in D_{2, \mtt j}} \varrho^{(2)}_{\mtt j, \vartheta}(x/n)\,
\frac{ \bar\rho_{2, \mtt j} - \rho(x/n)}{\rho(x/n)} 
\,.
\end{equation}

We turn to the proof of the bound \eqref{59}.  Add and subtract
$\vartheta$ to $\varrho^{(2)}_{\mtt j, \vartheta}(x/n)$.  Since
$\rho(\cdot)$ is Lipschitz continuous, by \eqref{83}, \eqref{50},
\eqref{57}, the right-hand side of \eqref{62} is equal to
\begin{align*}
\frac{1}{1-\bar\rho_{2, \mtt j} }\,
\frac{1}{ |D_{2, \mtt j}| }
\sum_{x\in D_{2, \mtt j}} 
\frac{ \bar\rho_{2, \mtt j} - \rho(x/n)}{\rho(x/n)} \, +\, \mf
R^{(1)}_n \,,
\quad
\text{where}\quad
|\mf R^{(1)}_n| \,\le\, \mf a_1 \, \Vert \nabla \rho\Vert_\infty^2\,
\Big( \frac{s_1s_2}{n}\Big)^2\,.
\end{align*} 
Denote by $\cb{x_{2,\mtt j}}$ the center of the cube $D_{2, \mtt
j}$. By a second-order Taylor expansion,
\begin{align*}
\bar\rho_{2, \mtt j} - \rho(x/n) \,=\,
\frac{1}{ |D_{2, \mtt j}| } \, 
\sum_{y\in D_{2, \mtt j}}
\frac{y-x}{n} \cdot \nabla \rho (x/n) \, +\, \mf R^{(2)}_n (x)
\,=\,
\frac{x_{2,\mtt j} - x}{n} \cdot \nabla \rho (x/n) \, +\, \mf R^{(2)}_n(x)
\end{align*}
because, by symmetry, the sum over $y$ of $y- x_{2,\mtt j} $
vanishes. In this formula, $\mf R^{(2)}_n(x)$ represents a remainder
such that
\begin{equation*}
\max_{x\in D_{2, \mtt j}} |\mf R^{(2)}_n (x) | \,\le\,
\Vert \nabla^2 \rho\Vert_\infty\,
\Big( \frac{s_1s_2}{n}\Big)^2\,.
\end{equation*}
Therefore, the right-hand side of \eqref{62} is equal to
\begin{align*}
& \frac{1}{1-\bar\rho_{2, \mtt j} }\, \frac{1}{ |D_{2, \mtt j}| }
\sum_{x\in D_{2, \mtt j}} 
\frac{x_{2,\mtt j} - x}{n} \cdot \frac{\nabla \rho (x/n) } {\rho(x/n)}   
\, +\, \mf R^{(1)}_n \, +\, \mf R^{(2)}_n
\\
& \quad
\,=\,
\frac{1}{1-\bar\rho_{2, \mtt j} }\,
\frac{1}{ |D_{2, \mtt j}| }
\sum_{x\in D_{2, \mtt j}} 
\frac{x_{2,\mtt j} - x}{n} \cdot \Big\{\,
\frac{\nabla \rho (x/n) } {\rho(x/n)}
\,-\,
\frac{\nabla \rho (x_{2,\mtt j} /n) } {\rho(x_{2,\mtt j} /n)}\,\Big\}
\, +\, \mf R^{(1)}_n \, +\, \mf R^{(2)}_n \,.
\end{align*}
We were allowed to introduce the term
$\nabla \rho (x_{2,\mtt j} /n) /\rho(x_{2,\mtt j} /n)$ because the sum
over $x\in D_{2, \mtt j}$ vanishes by symmetry. In this formula,
$\mf R^{(2)}_n$ is a remainder such that
$|\mf R^{(2)}_n | \,\le\, \mf a_1\, \Vert \nabla^2 \rho\Vert_\infty\,
( s_1s_2/n)^2$. As $\rho(\cdot)$ is of class $C^2(\bb T^d)$, the first
term on the right-hand side is bounded by an expression which
satisfies the same bound as $\mf R^{(2)}_n$. This completes the proof
of \eqref{59}. \smallskip

The next result provides an expansion of
$\varrho^{(2)}_{\mtt j, \vartheta} (x/n) - \vartheta $,
$x\in D_{2, \mtt j}$, $\mtt j\in \Sigma^+_{s_3}$. The point here is
that this difference can be represented as a sum of polynomials in
$\vartheta$ multiplied by functions of the density profile
$\rho(\cdot)$. To state the result, let
\begin{gather*}
\cb{H^{(1)}_{2, \mtt j}  (x/n) } \,:=\,
\frac{\bar\rho_{2, \mtt j}}{1 - \bar\rho_{2, \mtt j}} \,
\frac{x - x_{2, \mtt j}}{n}  
\cdot     \frac{ \nabla \rho (x/n) } { \rho (x/n)^2 }\,,
\\
\cb{H^{(2)}_{2, \mtt j} (x/n)} \,:=\, 
\frac{2} {\chi( \bar\rho_{2, \mtt j})^2} \,
\Big( \frac{x - x_{2, \mtt j}}{n}  
\cdot  \nabla \rho (x/n)  \Big)^2\,,
\end{gather*}
\vskip -.3 truecm
\begin{gather*}
\cb{H^{(3)}_{2, \mtt j} (x/n) } \,:=\,
\frac { 1 }  { \chi (\bar\rho_{2, \mtt j})^2}\,
\frac{1}{|D_{2, \mtt j}|} 
\sum_{y\in D_{2, \mtt j}}
\Big( \frac{y-x}{n} \cdot  \nabla \rho (x/n)  \Big)^2 \,,
\end{gather*}
\vskip -.3 truecm
\begin{align*}
\cb{H^{(4)}_{2, \mtt j} (x/n)} \,:=\,
-\, \frac{1}{ 2}\, \frac{1} { \chi(\bar\rho_{2, \mtt j})} \,
\frac{1}{|D_{2, \mtt j}|} 
\sum_{y\in D_{2, \mtt j}}
\frac{y-x}{n} \cdot \nabla^2\rho (x/n)\, \frac{y-x}{n} \,,
\end{align*}
where, recall, $\nabla^2\rho$ represents the Hessian of
$\rho(\cdot)$. Mind that $H^{(1)}_{2, \mtt j} (\cdot)$ is of order
$(s_1s_2/n)$, while the last three terms are of order $(s_1s_2/n)^2$.

\begin{lemma}
\label{s15}
Fix $0<\delta<1/2$.  There exists a finite constant
$\mf a_1 = \mf a_1(\delta, \mf r)$, depending only on $\delta$ and
$\mf r$, such that
\begin{align*}
\varrho^{(2)}_{\mtt j, \vartheta} (x/n) \,-\, \vartheta \, & =\,
\chi(\vartheta)\, H^{(1)}_{2, \mtt j} (x/n)
\,+\, \chi(\vartheta)\, (1- \vartheta)\, H^{(2)}_{2, \mtt j} (x/n)
\\
& +\, \chi(\vartheta) \, (\vartheta - \bar\rho_{2, \mtt j})\,
H^{(3)}_{2, \mtt j} (x/n) \,+\, \chi(\vartheta) \, H^{(4)}_{2, \mtt j} (x/n)
\,+\, \mf R^{(3)}_n (x) \,,
\end{align*}
for all $x \in D_{2, \mtt j}$, $\mtt j\in \Sigma^+_{s_3}$,
$\delta \le \vartheta \le 1-\delta$,
where
\begin{equation*}
\max_{\mtt j\in \Sigma^+_{s_3}}\, \max_{x \in D_{2, \mtt j}}
\big |\, \mf R^{(3)}_n (x) \,\big|
\,\le\, \mf a_1 \,
\, ( \, 1+ \Vert \rho\Vert_{C^3}\,)^6
\, \Big( \frac{s_1s_2}{n}\Big)^3\,.
\end{equation*}
\end{lemma}

\begin{proof}
Recall that $\mf a_1$ represents a finite constant depending only on
$\delta$ and $\mf r$ whose value may change from line to line.  For
three non-negative integers $i$, $p$, $q$, $\mf R^{p,q}_{n,i}(x)$
represents a remainder such that
\begin{equation*}
\max_{ \mtt j \in \Sigma^+_{s_3} }
\max_{x\in D_{2, \mtt j}} |\, \cb{\mf R^{p,q}_{n,i}(x)} \, | \,\le\,
\mf a_1\, ( \, 1+ \Vert \rho\Vert_{C^i}\,)^p\, \epsilon^q_n \,,
\quad
\text{where}\quad \cb{\epsilon_n} := \frac{ s_1s_2}{n} \, \cdot
\end{equation*}

Fix $\delta \le \vartheta \le 1-\delta$.  Write
$\varrho^{(2)}_{\mtt j,\vartheta} (y/n)$ as $F (\rho (y/n))$, where
$F(a) = a / \{ a + e^{- \Phi_{2, \mtt j}(\vartheta)} (1-a) \}$.  By
\eqref{50} and a Taylor expansion,
\begin{align}
\label{61}
\varrho^{(2)}_{\mtt j,\vartheta} 
(x/n) - \vartheta \, =\,
& -\,  \frac{1}{ |D_{2, \mtt j}| } \sum_{y\in D_{2, \mtt j}}
F' (\rho (x/n))\,
\big\{\, \rho  (y/n)
- \rho  (x/n) \,\big \}
\\
& -\, \frac{1}{2}
\,  \frac{1} { |D_{2, \mtt j}| } \sum_{y\in D_{2, \mtt j}}
F'' (\rho   (x/n))\,
\big\{\, \rho (y/n)
- \rho  (x/n) \,\big \}^2
\,+\, \mf R^{3,3}_{n,1} (x) \,.
\nonumber
\end{align}
Note that
\begin{equation*}
F' (\rho (x/n))\,=\,
e^{-\Phi_{2, \mtt j}(\vartheta)}\, 
\Big( \frac{\varrho^{(2)}_{\mtt j,\vartheta} (x/n)}
{\rho (x/n)} \Big) ^2 \,, \;\;
F'' (\rho (x/n))\,=\, -\, 2\, 
e^{-\Phi_{2, \mtt j}(\vartheta)}\,
[\, 1 - e^{-\Phi_{2, \mtt j}(\vartheta)} \,]\, 
\Big( \frac{\varrho^{(2)}_{\mtt j,\vartheta} (x/n)}
{\rho (x/n)} \Big) ^3 \,.
\end{equation*}
Hence, by \eqref{56b}, \eqref{98}, 
\begin{equation}
\label{165}
\max_{ \mtt j \in \Sigma^+_{s_3} }
\max_{x\in D_{2, \mtt j}} |\, F' (\rho (x/n))\, | \,\le\, \mf a_1\,,
\quad
\max_{ \mtt j \in \Sigma^+_{s_3} }
\max_{x\in D_{2, \mtt j}} |\, F'' (\rho (x/n))\, | \,\le\, \mf a_1\,.
\end{equation}
The bound on the remainder in \eqref{61} follows from \eqref{98},
the explicit formula for $F'''(\cdot)$, and from the fact that
$\mf r \le \rho(\cdot) \le 1-\mf r$.

We examine the first term on the right-hand side of \eqref{61}.  By
the formula for $F' (\rho (x/n))$, and since
$\sum_{y\in D_{2, \mtt j}} (y-x_{2, \mtt j})=0$, a Taylor expansion
yields that this expression can be written as
\begin{gather*}
\Big(\, \frac{\varrho^{(2)}_{\mtt j,\vartheta} (x/n) } { \rho (x/n) }
\Big)^2
e^{- \Phi_{2, \mtt j}(\vartheta)}
\, \Big\{ \frac{x-x_{2, \mtt j}}{n}\cdot \nabla \rho (x/n) \,-\,
\frac{1}{ 2\, |D_{2, \mtt j}| } \sum_{y\in D_{2, \mtt j}}
\frac{y-x}{n} \cdot \nabla^2\rho (x/n)\, \frac{y-x}{n} \,\Big\}
\,+\, \mf R^{1,3}_{n,3} (x) \,.
\end{gather*}
The same reasons as above provide the bound for the remainder.  By
\eqref{98} and \eqref{59}, in the main term of the previous displayed
equation we may replace the chemical potential
$\Phi_{2, \mtt j}(\vartheta)$ by
$\overline{\Phi}_{2, \mtt j}(\vartheta)$ at a cost
$\mf R^{3,3}_{n,2} (x)$.  Therefore, the first term on the right-hand
side of \eqref{61} is equal to
\begin{gather}
\label{164}
\Big(\, \frac{\varrho^{(2)}_{\mtt j,\vartheta} (x/n) } { \rho (x/n) }
\Big)^2
e^{-\overline{ \Phi}_{2, \mtt j}(\vartheta)}
\, \Big\{ \frac{x-x_{2, \mtt j}}{n}\cdot \nabla \rho (x/n) \,-\,
\frac{1}{ 2\, |D_{2, \mtt j}| } \sum_{y\in D_{2, \mtt j}}
\frac{y-x}{n} \cdot \nabla^2\rho (x/n)\, \frac{y-x}{n} \,\Big\}
\,+\, \mf R^{3,3}_{n,3} (x)  \,.
\end{gather}

We now examine the square of the ratio
$\varrho^{(2)}_{\mtt j,\vartheta} (x/n) /\rho (x/n)$ in the previous
formula. Add and subtract $\vartheta$ to
$\varrho^{(2)}_{\mtt j,\vartheta} (x/n) $, expand the square, use
\eqref{61} and \eqref{165} to get that
\begin{gather*}
\Big(\, \frac{\varrho^{(2)}_{\mtt j,\vartheta} (x/n) } { \rho (x/n) }
\Big)^2 \,=\,
\frac{\vartheta^2}{\rho (x/n)^2} 
\,+\, 2\, \vartheta\,  \frac{\varrho^{(2)}_{\mtt j,\vartheta} (x/n)^2}{\rho (x/n)^4}\,
e^{- \Phi_{2, \mtt j}(\vartheta)}\,
\frac{x-x_{2, \mtt j}}{n}\cdot \nabla\rho(x/n) \,
+\, \mf R^{3,2}_{n,2} (x) \,.
\end{gather*}
By \eqref{57}, we may further replace
$\varrho^{(2)}_{\mtt j,\vartheta} (x/n)$ by its average over
$D_{2, \mtt j}$, which is equal to $\vartheta$, at a cost
$\mf R^{2,2}_{n,1} (x) $. Finally, replace $\Phi_{2, \mtt j}$ by
$\overline{\Phi}_{2, \mtt j}$ at an extra cost
$\mf R^{3,3}_{n,2} (x)$.  In conclusion,
\begin{gather*}
\Big(\, \frac{\varrho^{(2)}_{\mtt j,\vartheta} (x/n) } { \rho (x/n) }
\Big)^2 \,=\,
\frac{\vartheta^2}{\rho (x/n)^2} 
\,+\, 2\, \frac{\vartheta^3}{\rho (x/n)^4}\,
e^{- \overline{\Phi}_{2, \mtt j}(\vartheta)}\,
\frac{x-x_{2, \mtt j}}{n}\cdot \nabla\rho(x/n) \,
+\, \mf R^{3,2}_{n,2} (x) \,.
\end{gather*}

Hence, by \eqref{63}, \eqref{164}, which is equal to the first
term on the right-hand side of \eqref{61}, is given by
\begin{equation}
\label{166}
\begin{aligned}
& \chi(\vartheta)\,
\frac{\bar\rho_{2, \mtt j} }{1 - \bar\rho_{2, \mtt j} }\,
\frac{x-x_{2, \mtt j}}{n} \cdot     \frac{ \nabla \rho (x/n) } { \rho (x/n)^2 }
\,+\,
2\, \chi(\vartheta)\, (1- \vartheta)\,
\frac{1} { \chi(\bar\rho_{2, \mtt j} )^2} \,
\Big( \frac{x-x_{2, \mtt j}}{n} \cdot  \nabla \rho (x/n)  \Big)^2 
\\
&\quad \,-\, \frac{1}{ 2}\, \frac{\chi(\vartheta)}
{\chi (\bar\rho_{2, \mtt j})} \,
\frac{1}{  |D_{2, \mtt j}| } \sum_{y\in D_{2, \mtt j}}
\frac{y-x}{n} \cdot \nabla^2\rho (x/n)\, \frac{y-x}{n} 
\,+\, \mf R^{6,3}_{n,3} (x)  \,.
\end{aligned}
\end{equation}
In the second term of the first line we replaced $\rho (x/n)^{-4}$ by
$\bar\rho_{2, \mtt j}^{-4}$ at a cost $\mf R^{6,3}_{n,1} (x)$, and in
the first term of the second line, we replaced $\rho (x/n)^{-2}$ by
$\bar\rho_{2, \mtt j}^{-2}$ at a cost $\mf R^{3,3}_{n,2} (x)$.

We turn to the second term in \eqref{61}.  Since $F''(\rho (x/n))$, in
the formula above \eqref{165}, can be written as
$- \, 2 \, [\varrho^{(2)}_{\mtt j,\vartheta} (x/n)/\rho (x/n) ]^3$
$e^{-\Phi_{2, \mtt j}(\vartheta)}\, [1-e^{-\Phi_{2, \mtt j}
(\vartheta)}]$, the second term in \eqref{61} is equal to
\begin{align*}
& e^{-\overline{\Phi}_{2, \mtt j}(\vartheta)}\,
[\, 1- e^{-\overline{\Phi}_{2, \mtt j}(\vartheta)}\, ] \,
\frac{\vartheta^3} {\rho (x/n)^3}\,
\,  \frac{1} { |D_{2, \mtt j}| } \sum_{y\in D_{2, \mtt j}} \,
\Big( \frac{y-x}{n} \cdot  \nabla \rho (x/n)  \Big)^2
\,+\, \mf R^{4,3}_{n,2} (x) 
\\
&\quad
=\,
\frac { \chi(\vartheta) }  { \chi (\bar\rho_{2, \mtt j})^2}\,
(\vartheta - \bar\rho_{2, \mtt j})\, 
\,  \frac{1} { |D_{2, \mtt j}| } \sum_{y\in D_{2, \mtt j}} \,
\Big( \frac{y-x}{n} \cdot  \nabla \rho (x/n)  \Big)^2
\,+\, \mf R^{5,3}_{n,2} (x)  
\,.
\end{align*}
To obtain the term on the left-hand side of the equality, we performed
many changes.  We replaced $\varrho^{(2)}_{\mtt j,\vartheta} (x/n)$ by
its spatial average, which is $\vartheta$, at a cost
$\mf R^{3,3}_{n,1} (x)$, in view of \eqref{57}. Recalling \eqref{59},
we also changed $\Phi_{2, \mtt j}(\vartheta)$ by
$\overline{\Phi}_{2, \mtt j}(\vartheta)$ at a cost
$\mf R^{4,4}_{n,2} (x)$. Finally we applied a Taylor expansion to the
squared term and recalled the bounds \eqref{98} to obtain a remainder
of type $\mf R^{2,3}_{n,2} (x)$. To obtain the identity, we replaced
$\rho (x/n)^{-3}$ by $\bar\rho_{2, \mtt j}^{-3}$ at a cost
$\mf R^{5,3}_{n,1} (x)$, and we used equation \eqref{63}.

To complete the proof of the lemma, it remains to recall \eqref{61},
\eqref{166}, and the previous formula which yields the second term in
\eqref{61}.  
\end{proof}

\begin{remark}
\label{rm5}
Setting $s_1$ to be $1$, turns the set $D_{2,\mtt j}$ into a set of width
$s_2$. Hence, taking $s_2=\ell_1$, $s_3=L_2$, where $\ell_1$, $L_2$ are the
quantities introduced at the beginning of Section \ref{sec4}, Lemma
\ref{s15} becomes a statement on $\varrho_{\mtt k,\vartheta} (\cdot)$,
$\mtt k\in \Sigma^+_{L_2}$, the density profiles considered in the
second part of Section \ref{sec4}.  The next two corollaries are used
in that section.
\end{remark}

In the next two results $s_1=1$, as noted in their assumptions,
so that $s_2s_3=n$. Moreover, $\cb{B_{\mtt k} = D_{2,\mtt k}}$,
$\mtt k\in \Sigma^+_{s_3}$.  Denote by $\cb{x_{\mtt k }}$ the center
of the cube $B_{\mtt k}$. Fix a finite subset $A$ of $\bb Z^d$, and
recall from \eqref{90} that
$\cb{J_{\mtt k}} := |B^o_{\mtt k}|^{-1} \sum_{x\in B^o_{\mtt k}}
J(x/n)$, where $B^o_{\mtt k}$, introduced in \eqref{95}, represents
the interior of the set $B_{\mtt k}$.  Recall that $\mf a_{3, p}$,
$p\ge 1$, stands for a finite constant which depends on $\delta$,
$\mf r$, a cylinder function $h$ (which sometimes will be $\eta_A$),
and $\Vert \rho\Vert_{C^p}$, whose value may change from line to line.

\begin{corollary}
\label{s35}
Let $s_1=1$.  For each $0<\delta<1/2$, finite subset $A$ of $\bb Z^d$,
function $J$ in $C^1(\bb T^d)$, there exists a finite constant
$\mf a_{3, 3}$, also  depending on $\delta$, such that
\begin{equation*}
\begin{aligned}
\sum_{x\in B^o_{\mtt k}} J(x/n)\,
E_{\nu^{\rm gc} _{\mtt k, \mss M }} [\, \tau_x \eta_A\,] \,
=\, |B_{\mtt k}^o| \, J_{\mtt k}\, \Big\{\, \mss m^{|A|}
\,+\, \mss m^{|A|-1}\,
\frac{\chi(\mss m)} {\chi(\bar\rho_{\mtt k})}\,
\sum_{y\in A} \frac{y}{n}\cdot (\nabla \rho)(x_{\mtt k}/n)\,\Big\} 
\,+\, \mf R_n (\mtt k)
\end{aligned}
\end{equation*}
for all $\delta \le \mss m\le 1-\delta$.  Here
$\mss m = \mss M/|B_{\mtt k}|$, and
\begin{equation*}
\max_{\mtt k \in \Sigma^+_{s_3}}
\max_{\delta \le \mss m \le 1- \delta } 
|\, \mf R_n (\mtt k)\, | \,\le\,  \mf a_{3, 3}\, \Vert J\Vert_{C^1} \,
\, (s_2/n)^2\, s^d_2 \,.
\end{equation*}
\end{corollary}

\begin{proof}
Rewrite the expectation appearing on the left-hand side of the
corollary's statement as
$\prod_{y\in A} \varrho_{\mtt k, \mss m} ((y+x)/n)$ and apply
Lemma \ref{s15} with $s_1=1$ to get that the left-hand side is equal
to
\begin{equation*}
|B^o_{\mtt k}| \, J_{\mtt k}\, \mss m^{|A|}
\,+\, \chi(\mss m)\, \mss m^{|A|-1}
\frac{\bar\rho_{\mtt k}} {1-\bar\rho_{\mtt k}}\,
\sum_{y\in A} \sum_{x\in B^o_{\mtt k}} J(x/n)
\frac{x+y-x_{\mtt k}}{n}\cdot
\frac{(\nabla \rho)([x+y]/n)}{\rho([x+y]/n)^2} 
\,+\, \mf R_n 
\end{equation*}
where
$|\, \mf R_n\, | \,\le\, \mf a_{3, 3}\, \Vert J\Vert_\infty \,
(s_2/n)^2\, s^d_2$.  Note that $|x+y-x_{\mtt k}|/n \le s_2/n$. We may
therefore replace $J(x/n)$ by $J_{\mtt k}$ at a cost bounded by
$\mf a_{3, 1} \, \Vert J\Vert_{C^1}\, (s_2/n)^2 \, s_2^d$, and
$(\nabla \rho)([x+y]/n) /\rho([x+y]/n)^{2}$ by
$(\nabla \rho)(x_{\mtt k}/n) / \bar\rho_{\mtt k}^{2}$ at a cost
$\mf a_{3, 2} \, \Vert J\Vert_{\infty}\, (s_2/n)^2 \, s_2^d$.  After
these replacements, to complete the proof it remains to observe that
$\sum_{x\in B^o_{\mtt k}} (x-x_{\mtt k})=0$ by symmetry.
\end{proof}

Let $G_{\mtt k} \colon \bb T^d \to \bb R$,
$\mtt k \in \Sigma^+_{s_3}$, be given by
\begin{equation*}
\cb{G_{\mtt k}(z/n) } \,:=\, \frac{\bar\rho_{\mtt k}}{1 - \bar\rho_{\mtt k}} \,
\frac{z-x_{\mtt k}}{n} \cdot
\frac{ \nabla \rho (z/n) } { \rho (z/n)^2 }\, \cdot
\end{equation*}
This is the function $H^{(1)}_{2, \mtt j}(\cdot)$ introduced in Lemma
\ref{s15} with $s_1=1$.  Note that $G_k(\cdot)$ restricted to
$B_{\mtt k}$ is bounded by $\mf a_{3, 1} s_2/n$.  Recall near
\eqref{79} the definition of $\Gamma_{\mtt k}:=\Gamma_N$ with respect
to $B_{\mtt k}$, from \eqref{43b} the definition of
$\mtt T^{(1)}_{\mtt k, A} := \mtt T_{1,A}$ with respect to
$B_{\mtt k}$, and from \eqref{90} and \eqref{168} the one of the
polynomial $\mtt q_{A} \colon [0,1]\to \bb R$.

\begin{corollary}
\label{s31}
Let $s_1=1$.  Fix $0<\delta < 1/2$, a finite subset $A$ of $\bb Z^d$,
and a function $J$ in $C^1(\bb T^d)$. There exists a finite constant
$\mf a_{3,3}$, also depending on $\delta$, such that
\begin{equation*}
\begin{aligned}
\sum_{x\in B^o_{\mtt k}} J(x/n)\,
E_{\nu^{\rm gc} _{\mtt k, \mss M }}  [\, \tau_x \eta_A\,] \,
\frac{1} {\Gamma_{\mtt k}}\, \mtt T^{(1)}_{\mtt k, x+A} \, 
=\,  \frac{|B^o_{\mtt k}|}{|B_{\mtt k}|} \,
J_{\mtt k}\, \mtt q_{A} (\mss m)  
\,+\, \mf R_n (\mtt k)
\end{aligned}
\end{equation*}
where $\mss m = \mss M/|B_{\mtt k}|$, and 
\begin{equation*}
\max_{\mtt k \in \Sigma^+_{s_3}} 
\max_{\delta \le \mss m \le 1- \delta } 
\big|\, \mf R_n (\mtt k)\, \big |
\,\le\,  \mf a_{3,3} \, \Vert J\Vert_{C^1}\, \Big\{ \frac{1}{n}
\,+\,  \Big( \frac {s_2}{n} \Big)^2\,\Big\} \,.
\end{equation*}
\end{corollary}

\begin{proof} 
The proof is similar to the one of Corollary \ref{s35}.  If the
expectation was carried out with respect to a Bernoulli product
measure with density $\mss m$, the first order term corresponds to
that in Remark \ref{s70}. To assess the error, expand 
$E_{\nu^{\rm gc} _{\mtt k, \mss M }} [\, \tau_x \eta_A\,] $ and each
term in the definition of $\mtt T^{(1)}_{\mtt k, A}$ by using Lemma
\ref{s15}. A long but simple computation, similar to the one presented
in the previous proof, yields the result. The main observation is that
the sums $\sum_{x\in B^o_{\mtt k}} (x-x_{\mtt k})$,
$\sum_{x\in B_{\mtt k}} (x-x_{\mtt k})$ vanish by symmetry.  For
example, it follows from this remark and the expansion established in
Lemma \ref{s15} that
\begin{equation*}
\frac{1}{ | B_{\mtt k}|} \sum_{x\in B_{\mtt k}}
F(\varrho_{\mtt k, \mss m} (x/n)) \,=\,
F(\mss m) \,+\, \mf R_n 
\end{equation*}
for some error $\mf R_n$ bounded by $\mf a_{3,3} (s_2/n)^2$. It also
follows from this remark that
$|B^o_{\mtt k}|^{-1} \sum_{x\in B^o_{\mtt k}} (x+y-x_{\mtt k})/n =
y/n$, what explains the term $1/n$ in the error; for instance, the
first order term in Corollary \ref{s35} divided by $|B_{\mtt k}|$ is
of order $1/n$.
\end{proof}

\subsection*{Recursive step}

We turn to estimates needed in the induction step of the proof.  Here,
$s_1$ may be large.  In particular, when $s_1$ is large, bounds of
order $(s_1/n)^2$ are not sharp enough and we need to include the next
order term in the expansion to obtain errors bounded by $(s_1/n)^3$.

The next lemma follows from Lemma \ref{s15} by averaging over
$x\in D_{\mtt k}$, for a fixed $\mtt k\in\Lambda_{\mtt j}$.

\begin{lemma}
\label{s20b}
Fix $0<\delta<1/2$.  There exists a finite constant
$\mf a_{3,3}$, also  depending on $\delta$,
such that
\begin{align*}
E_{\nu^{\rm gc} _{2,\mtt j, \mss M} } [\mss m_{\mtt k} ]
\,-\, \mss m_{2, \mtt j}
\, & =\,
\chi(\mss m_{2, \mtt j})\, \mtt Y_{\mtt j, \mtt k}
\,+\, \mf R_n
\,,
\end{align*}
where
\begin{equation*}
\cb{\mtt Y_{\mtt j, \mtt k}} \,=\, H^{(1)}_{\mtt j, \mtt k}  \,+\,
\mtt Y^{(2)}_{\mtt j, \mtt k}\,, \quad
\cb{\mtt Y^{(2)}_{\mtt j, \mtt k}} \,=\, 
(1- \mss m_{2, \mtt j})\, H^{(2)}_{\mtt j, \mtt k} 
\, +\, 
(\mss m_{2, \mtt j} - \bar\rho_{2,\mtt j})\,
H^{(3)}_{\mtt j, \mtt k}
+\, H^{(4)}_{\mtt j, \mtt k}  \,,
\end{equation*}
\vskip -.3 truecm
\begin{gather*}
\cb{H^{(1)}_{\mtt j, \mtt k}  } \,:=\,
\frac{\bar\rho_{2,\mtt j}}{1 - \bar\rho_{2,\mtt j}} \,
\frac{1}{|D_{\mtt k}|} \sum_{x\in D_{\mtt k}}
\frac{x - x_{2,\mtt j}}{n} 
\cdot     \frac{ \nabla \rho (x/n) } { \rho (x/n)^2 }\,,
\\
\cb{H^{(2)}_{\mtt j, \mtt k}  } \,:=\, 
\frac{2} {\chi( \bar\rho_{2,\mtt j})^2} \,
\frac{1}{|D_{\mtt k}|} \sum_{x\in D_{\mtt k}}
\Big( \frac{x - x_{2,\mtt j}}{n} 
\cdot  \nabla \rho (x/n)  \Big)^2\,,
\end{gather*}
\vskip -.3 truecm
\begin{gather*}
\cb{H^{(3)}_{\mtt j, \mtt k}  } \,:=\,
\frac { 1 }  { \chi (\bar\rho_{2,\mtt j})^2}\,
\frac{1} { |D_{2, \mtt j}| } \frac{1}{|D_{\mtt k}|} 
\sum_{y\in D_{2, \mtt j}} \,
\sum_{x\in D_{\mtt k}}
\Big( \frac{y-x}{n} \cdot  \nabla \rho (x/n)  \Big)^2 \,,
\end{gather*}
\vskip -.3 truecm
\begin{align*}
\cb{H^{(4)}_{\mtt j, \mtt k}  } \,:=\,
-\, \frac{1}{ 2}\, \frac{1} { \chi(\bar\rho_{2,\mtt j})} \,
\frac{1} { |D_{2, \mtt j}| } \frac{1}{|D_{\mtt k}|} 
\sum_{y\in D_{2, \mtt j}} \,
\sum_{x\in D_{\mtt k}}
\frac{y-x}{n} \cdot \nabla^2\rho (x/n)\, \frac{y-x}{n} \,,
\end{align*}
and $|\, \mf R_n \,| \,\le\, \mf a_{3,3} \, (s_1s_2 /n)^3$ for all
$\mtt j\in \Sigma^+_{s_3}$, $\mtt k\in \Lambda_{\mtt j}$,
$\delta \le \mss m_{2, \mtt j} \le 1-\delta$.
\end{lemma}

Recall from \eqref{82} that
$\cb{\widetilde{\mss m}_{\mtt j, \mtt k} } = E_{\nu^{\rm gc} _{2,\mtt
j, \mss M} } [\mss m_{\mtt k} ]$. The next result follows from the
previous lemma and a Taylor expansion.

\begin{corollary}
\label{s38}
Fix a smooth function $F\colon \bb R \to \bb R$, $0 < \delta <
1/2$. Then, there exists a constant $\mf a_1 = \mf a_1(\delta, F)$ such
that
\begin{equation*}
\Big | \, F(\widetilde{\mss m}_{\mtt j, \mtt k})\,-\,
F (\mss m_{2, \mtt j}) \, - \,  F'(\mss m_{2, \mtt j})
\, \chi(\mss m_{2, \mtt j})\,  \mtt Y_{\mtt j, \mtt k}
\,-\, \frac{1}{2}\,  F''(\mss m_{2, \mtt j})\,
\big[\, \chi(\mss m_{2, \mtt j})\, H^{(1)}_{\mtt j, \mtt k}
\big]^2 
\,\Big| \, \le \, \mf a_1  \, \Big(\, 
\frac{s_1s_2}{n} \,\Big)^3
\end{equation*}
for all $\mtt j\in \Sigma^+_{s_3}$, $\mtt k \in \Lambda_{\mtt j}$,
$\delta \le \mss M/|D_{2,\mtt j}|\le 1-\delta$, $s_1$, $s_2\ge 1$.
\end{corollary}

The next result provides an estimate of the averages of the density
inside the cube $D_{\mtt k}$.  We claim that under the hypotheses of
Lemma \ref{s15}, for each function $f\colon [0,1]\to \bb R$ of class
$C^2$, there exists a finite constant $\mf a_{3,3}(\delta, f)$ such
that
\begin{equation}
\label{85}
\Big|\,  \frac{1}{|D_{\mtt k}|} \sum_{x\in D_{\mtt k}}
\big\{ \, f(\varrho^{(2)}_{\mtt j, \mss m_{2, \mtt j}} (x/n) )
- f(\widetilde{\mss m}_{\mtt j, \mtt k}  )\,\big\}\, \Big|
\,\le\, \mf a_{3,3}(\delta, f) \, \Big\{ \,\Big( \frac{s_1}{n}\Big)^2
\,+\, \Big( \frac{ s_1s_2}{n}\Big)^6\Big\} \,,
\end{equation}
for all $\mtt j\in \Sigma^+_{s_3}$, $\mtt k\in \Lambda_{\mtt j}$,
$\delta \le \mss m_{2, \mtt j} \le 1-\delta$.

To prove this claim, we first estimate the difference
$\varrho^{(2)}_{\mtt j, \mss m_{2, \mtt j}} (x/n) - \widetilde{\mss
m}_{\mtt j, \mtt k}$. Subtract the left-hand side of Lemma \ref{s20b}
from the one of Lemma \ref{s15}. The resultant on the left-hand side
equals
$\varrho^{(2)}_{\mtt j, \mss m_{2, \mtt j}} - \widetilde{\mss m}_{\mtt
j, \mtt k}$. On the right-hand side we obtain four terms and the
remainder. All terms, with the exception of the first one, are bounded
by $\mf a_{3,2}(\delta) \, (s_1s_2/n)^2$, uniformly in
$x\in D_{\mtt k}$. The difference of the first terms is equal to
\begin{equation*}
\chi(\mss m_{2, \mtt j})\,\frac{\bar\rho_{2, \mtt j}}{1 - \bar\rho_{2, \mtt j}} \,
\frac{1}{|D_{\mtt k}|} \sum_{y\in D_{\mtt k}}
\Big\{ \frac{x - x_{2, \mtt j}}{n}  
\cdot     \frac{ \nabla \rho (x/n) } { \rho (x/n)^2 }
\,-\,
\frac{y - x_{2, \mtt j}}{n}  
\cdot     \frac{ \nabla \rho (y/n) } { \rho (y/n)^2 }\,\Big\}\,.
\end{equation*}
Rewriting $x - x_{2, \mtt j}$ as $(x - y) + (y-x_{2, \mtt j})$, the
term inside braces become
\begin{equation*}
\frac{x - y} {n}  
\cdot     \frac{ \nabla \rho (x/n) } { \rho (x/n)^2 }
\,+\,
\frac{y - x_{2, \mtt j}}{n}  
\cdot
\Big(\, \frac{ \nabla \rho (x/n) } { \rho (x/n)^2 }
\,-\, \frac{ \nabla \rho (y/n) } { \rho (y/n)^2 }
\Big)
\end{equation*}
This expression, and thus, the previous displayed equation, is bounded
by $\mf a_{3,2} (\delta) \, (s_1/n)$.  Thus, as there exists a finite
constant $C_0$ such that $a + (ab)^2 \le C\, (a+ (ab)^3)$ for all
$0\le a, b\le 1$, taking into account the remainders,
$|\varrho^{(2)}_{\mtt j, \mss m_{2, \mtt j}} (x/n) - \widetilde{\mss
m}_{\mtt j, \mtt k} |\le \mf a_{3,3} (\delta) \, \{ (s_1/n) +
(s_1s_2/n)^3\} $ for all $x\in D_{\mtt k}$.

By Taylor's expansion, the left-hand side of \eqref{85} is equal to
\begin{equation*}
f'(\widetilde{\mss m}_{\mtt j, \mtt k}  )
\, \frac{1}{|D_{\mtt k}|} \sum_{x\in D_{\mtt k}}
\big\{ \, \varrho^{(2)}_{\mtt j, \mss m} (x/n) 
- \widetilde{\mss m}_{\mtt j, \mtt k} \,\big\}\, +\, \mf R_n\,,
\end{equation*}
where
$|\mf R_n|\le \mf a_{3,3}(\delta, f)\, \{ (s_1/n) + (s_1s_2/n)^3\}^2$.
By definition of $\widetilde{\mss m}_{\mtt j, \mtt k} $, the first
term vanishes, which completes the proof of the claim \eqref{85}.

Recall that $\mss m_{\mtt k}$ represents the particle's density on
$D_{\mtt k}$: 
\begin{equation*}
\cb{\mss m_{\mtt k}} \,:=\, \frac{1}{|D_{\mtt k}|}
\sum_{x\in   D_{\mtt k}}  \eta_x \,, \quad \mtt k \in
\Sigma^+_{s_2s_3}\;. 
\end{equation*}

\begin{lemma}
\label{s17}
Fix $0<\delta <1/2$, and $p\ge 1$. Then, there exists a constant
$\mf a_1 = \mf a_1(\delta, p)$ such that
\begin{equation*}
\max_{\mtt j\in \Sigma^+_{s_3}}
\sup _{\delta \le \mss m \le 1-\delta}\,
\max_{\mtt k\in \Lambda_{\mtt j} }
\Big | \, 
E_{\nu^{\rm c} _{2, \mtt j, \mss M} }
\big [ \, \mss m_{\mtt k}^p\,\big]  \; -\;
E_{\nu^{\rm gc} _{2, \mtt j, \mss M} }
\big [ \, \mss m_{\mtt k}^p \,\big] \,
\Big( \, 1 \,-\, \frac{ p\, (p-1) }{ (s_1 s_2)^d}
\, \frac{(1-\mss m_{2, \mtt j})^2}
{2\, \chi(\mss m_{2, \mtt j})} \,\Big)
\, \Big|
\; \le \; \frac {\mf a_1}{(s_1s_2)^d}
\,\Big(\, 
\frac {s_1s_2} {n} \,+\,
\frac {1}{s_1^d} \,\Big)
\end{equation*}
for all $s_1\ge 2$, $s_2\ge 2$. In this formula, $\mss m = \mss
M/|D_{2, \mtt j}|$.
\end{lemma}

\begin{proof}
By definition, 
\begin{equation*}
\mss m_{\mtt k}^p \,=\, \frac{1}{s_1^{dp} }
\, \sum_{x_1, \dots, x_p} \eta_{x_1} \cdots \eta_{x_p} \,,
\end{equation*}
where the sum is carried over all sites $x_1 \dots, x_p$ in
$D_{\mtt k}$. Apply Corollary \ref{s18} to $A = \{x_1, \dots, x_p\}$.
Mind that the number of elements of $A$ might be strictly less than
$p$. The first order term in Corollary \ref{s18} is
$E_{\nu^{\rm gc} _{2, \mtt j, \mss M} } [ \, \eta_A \,] $. Summing over
$x_1, \dots, x_p$ yields
$E_{\nu^{\rm gc} _{2, \mtt j, \mss M} } [ \, \mss m_{\mtt k}^p \,] $.

We turn to the second order term.  By Lemma \ref{s15} and a
straightforward calculation (cf. Remark \ref{s70} and
$\mtt p''_A(\mss m) = |A|\big(|A|-1\big){\mss m}^{|A|-2}$), there
exists a finite constant $\mf a_1=  \mf a_1(\delta)$ such that
\begin{equation*}
\Big |\, \frac{1}{\Gamma_{2,\mtt j}} \, \mtt T^{(1)}_{2, \mtt j , A} \,+\,
\frac{ |A|\, (|A|-1) }{ (s_1 s_2)^d}
\, \frac{(1-\mss m_{2, \mtt j})^2}
{2\, \chi(\mss m_{2, \mtt j})} 
\, \Big| \,\le\,  \mf a_1 \, \frac{1}{(s_1 s_2)^d}\,
\frac{s_1s_2}{n}\,\cdot 
\end{equation*}
Hence, in Corollary \ref{s18}, the second order term in the expansion
of $E_{\nu^{\rm c} _{2, \mtt j, \mss M} } [ \, \eta_A \,]$ is equal to
\begin{equation*}
E_{\nu^{\rm gc} _{2, \mtt j, \mss M} } [ \, \eta_A \,] \,
\frac{1}{\Gamma_{2,\mtt j}} \, \mtt T^{(1)}_{2, \mtt j , A}
\,=\,
-\, \frac{1}{s_1^{dp} }
\, \sum_{x_1, \dots, x_p}
E_{\nu^{\rm gc} _{2, \mtt j, \mss M} } [ \,  \eta_{x_1} \cdots \eta_{x_p} \,]
\frac{ |A|\, (|A|-1) }{ (s_1 s_2)^d}
\, \frac{(1-\mss m_{2, \mtt j})^2}
{2\, \chi(\mss m_{2, \mtt j})}  \, +\,
R_n\,,
\end{equation*}
where $A = \{x_1, \dots, x_p\}$, and
$|R_n|\le \mf a_1(\delta, p) \, (s_1 s_2)^{1-d} \,
n^{-1}$.

The number of sets $A$ with non-distinct sites is bounded above by
$C_0 s_1^{d(p-1)}$. Then, considering only the case where all sites
are distinct we get an extra error bounded by
$s_1^{-d} \, (s_1 s_2)^{-d}$. If all sites are distinct $|A|=p$. After
replacing $|A|$ by $p$ we may reconsider all terms in the sum by a
further error of order $s_1^{-d} \, (s_1 s_2)^{-d}$.  This completes
the proof of the lemma.
\end{proof}

Let $\mss Q_{1}(\cdot)$ be the identity, and for $p\ge 2$,
let $\mss Q_p$, $\mss q_{p,i} \colon \bb R \to \bb R$, $0\le i\le 2$,
be the polynomials given by
\begin{equation}
\label{99}
\begin{gathered}
\cb{\mss Q_{p} (m)} \,=\, \sum_{i=0}^2
\frac{1}{|D_{\mtt k}|^i} \, \mss q_{p,i} (m) \,, \quad 
\cb{\mss q_{p,0} (m)} \,=\, m^{p}\,, \quad 
\cb{\mss q_{p,1} (m)} \,=\, {p \choose 2} \,
\chi (m) \,  m^{p-2}\,,
\\
\cb{ \mss q_{p,2} (m)} \,=\ -\,  \Big[\,
2\, {p \choose 3} \,+\, 3\, {p \choose 4} \,\Big]
\, \chi (m)\, m^{p-2}
\,+\, 
\Big[ \,  {p \choose 3}  \,+\, 3\, {p \choose 4}
\,\Big]  \, \chi(m)\, \, m^{p-3} 
\,.
\end{gathered}
\end{equation}
The next result provides an expansion of
$E_{\nu^{\rm gc} _{2,\mtt j, \mss M} } [\mss m_{\mtt k}^p ]$,
$p\ge 2$.  Note that the left-hand side vanishes for $p=1$.

\begin{lemma}
\label{s22}
Fix $0 < \delta < 1/2$, $p\ge 1$. Then, there exists a constant
$\mf a_{3,3} = \mf a_{3,3} (\delta, p)$ such that
\begin{equation*}
\Big|\, E_{\nu^{\rm gc} _{2,\mtt j, \mss M} } [\mss m_{\mtt k}^p ] - 
\mss Q_p (\widetilde{\mss m}_{\mtt j, \mtt k})
\,\Big| \, \le \, \frac{\mf a_{3,3} }{s_1^{d}}   \, \Big\{\, 
\frac{1 }{s_1^{2d}} \,+\, \Big(\frac{s_1}{n}\Big)^2 \,\Big\}
\end{equation*}
for all $\mtt j\in \Sigma^+_{s_3}$, $\mtt k \in \Lambda_{\mtt j}$,
$\delta \le \mss M/|D_{2,\mtt j}|\le 1-\delta$, $s_1$, $s_2\ge 1$.
\end{lemma}

\begin{proof}
Fix $p\ge 2$, and write
\begin{equation*}
\mss m_{\mtt k}^p \,=\, \frac{1}{ |D_{\mtt k} |^{p} }
\,\Big\{
\, \sum_{x_1, \dots, x_p} \eta_{x_1} \cdots \eta_{x_p} \,+\,
\, c_1\, \sum_{x_1, \dots, x_{p-1}} \eta_{x_1} \cdots \eta_{x_{p-1}}
\,+\,
\, c_2\, \sum_{x_1, \dots, x_{p-2}} \eta_{x_1} \cdots \eta_{x_{p-2}}
\, \Big\}  \,+\, \mf R_r \,.
\end{equation*}
Here and below, $\mf R_r$ represents a remainder which satisfies the
bound $|\mf R_r | \le C_0 s_1^{-3d}$ for some finite constant $C_0$
depending only on $p$. In this formula, the sums are carried out over
all distinct sites $x_1 \dots, x_q$ in $D_{\mtt k}$, and
$c_{1} = p(p-1)/2$, $c_2 = c_{2,1}+ c_{2,2}$,
$c_{2,1} = p(p-1)(p-2)/ 6$, $c_{2,2} = p(p-1)(p-2)(p-3)/8$.  Note
that $c_1$ is $p$ choose $2$, $c_{2,1}$ $p$ choose $3$, and
$c_{2,2}$ three times $p$ choose $4$. By definition of
$\varrho^{(2)}_{\mtt j,\vartheta } (\cdot) $, and since
$\nu^{\rm gc} _{2,\mtt j, \mss M}$ is a product measure, taking
expectation in the previous formula yields that
\begin{equation*}
\begin{aligned}
E_{\nu^{\rm gc} _{2,\mtt j, \mss M}}
\big[ \mss m_{\mtt k}^p\,\big]  \,=\, \frac{1}{ |D_{\mtt k} |^{p} }
\,\Big\{ &
\, \sum_{x_1, \dots, x_p} \varrho (x_1/n) \cdots \varrho (x_p/n)   \,+\,
\, c_1\, \sum_{x_1, \dots, x_{p-1}} \varrho (x_1/n)  \cdots
\varrho (x_{p-1}/n)
\\
&\quad 
\,+\,
\, c_2\, \sum_{x_1, \dots, x_{p-2}} \varrho (x_1/n)  \cdots
\varrho (x_{p-2}/n )
\, \Big\}  \,+\, \mf R_r \,,
\end{aligned}
\end{equation*}
where
$\varrho (\cdot) = \varrho^{(2)}_{\mtt j, \mss m_{2, \mtt j} }
(\cdot)$.

Similarly, by definition of $\widetilde{\mss m}_{\mtt j, \mtt k}$,
\begin{align*}
\widetilde{\mss m}_{\mtt j, \mtt k}^p \,=\, 
\frac{1}{ |D_{\mtt k} |^{p} }
\,\Big\{ & 
\sum_{x_1, \dots, x_p} \varrho (x_1/n) \cdots \varrho (x_p/n)   \,+\,
\, c_1\, \sum_{x_1, \dots, x_{p-1}} \varrho (x_1/n) ^2 \varrho (x_2/n)\cdots
\varrho (x_{p-1}/n)
\\
& +\,
\, c_{2,1}\, \sum_{x_1, \dots, x_{p-2}} \varrho (x_1/n)^3  \varrho
(x_2/n) \cdots \varrho (x_{p-2}/n )
\\
& +\, c_{2,2}\, \sum_{x_1, \dots, x_{p-2}} \varrho (x_1/n)^2  \varrho
(x_2/n)^2 \varrho (x_3/n) \cdots \varrho (x_{p-2}/n )
\, \Big\}  \,+\, \mf R_r \,.
\end{align*}

Subtract the right-hand side of the two previous equations.
The first terms cancel. The difference of the second ones is equal to
\begin{align*}
\frac{c_1}{ |D_{\mtt k} |^{p} }
\sum_{x_1, \dots, x_{p-1}} \chi(\varrho (x_1/n))
\, \varrho (x_2/n)\cdots
\varrho (x_{p-1}/n) \,,
\end{align*}
where the sum is carried out over all distinct sites $x_1 \dots, x_{p-1}$
in $D_{\mtt k}$. This expression is equal to
\begin{align*}
& \frac{c_1}{ |D_{\mtt k} |} \,
\Big(\, \frac{1}{ |D_{\mtt k} |}
\sum_{x\in D_{\mtt k}} \chi(\varrho (x/n))  \,\Big)\,
\Big(\, \frac{1}{ |D_{\mtt k} |}
\sum_{x\in D_{\mtt k}} \varrho (x/n)  \,\Big)^{p-2}
\\
&\quad
-\, \frac{c_1 (p-2) }{ |D_{\mtt k} |^2} \,
\Big(\, \frac{1}{ |D_{\mtt k} |}
\sum_{x\in D_{\mtt k}} \chi(\varrho (x/n))\,
\varrho (x/n) \,\Big)\,
\Big(\, \frac{1}{ |D_{\mtt k} |}
\sum_{x\in D_{\mtt k}} \varrho (x/n)  \,\Big)^{p-3}
\\
&\quad
-\, \frac{c_3 }{ |D_{\mtt k} |^2} \,
\Big(\, \frac{1}{ |D_{\mtt k} |}
\sum_{x\in D_{\mtt k}} \chi(\varrho (x/n))\,\,\Big)\,
\Big(\, \frac{1}{ |D_{\mtt k} |}
\sum_{x\in D_{\mtt k}} \varrho (x/n)^2 \,\,\Big)\,
\Big(\, \frac{1}{ |D_{\mtt k} |}
\sum_{x\in D_{\mtt k}} \varrho (x/n)  \,\Big)^{p-4} \, +\, \mf R_r\,,
\end{align*}
where, $c_3 = c_1 (p-2)(p-3)/2 = 2\, c_{2,2}$, and
$c_1(p-2)= 3\, c_{2,1}$. By \eqref{85}, we may replace in the previous
formula $f (\varrho (x/n))$ by
$f(\widetilde{\mss m}_{\mtt j, \mtt k})$ at a cost bounded by
$\mf a_{3,3}(\delta) \, (s_1/n)^2$. Hence, by definition of
$\widetilde{\mss m}_{\mtt j, \mtt k}$, the previous expression is
equal to
\begin{align*}
& \frac{c_1}{ |D_{\mtt k} |} \,
\chi(\widetilde{\mss m}_{\mtt j, \mtt k})\, 
\widetilde{\mss m}_{\mtt j, \mtt k}^{p-2}\,
-\, \frac{c_1 (p-2) }{ |D_{\mtt k} |^2} \,
\chi(\widetilde{\mss m}_{\mtt j, \mtt k})\,
\widetilde{\mss m}_{\mtt j, \mtt k}^{p-2}
-\, \frac{c_3 }{ |D_{\mtt k} |^2} \,
\chi(\widetilde{\mss m}_{\mtt j, \mtt k})\,
\widetilde{\mss m}_{\mtt j, \mtt k}^{p-2}\,
\, + \, \mf R_{r}
\, + \, \mf R_{n,r}\,,
\end{align*}
where $|\mf R_{n,r}| \le \mf a_{3,3}(\delta) \,  s_1^{-d} (s_1/n)^2$.

An analogous argument yields that the difference of the third terms in
the expansion of
$E_{\nu^{\rm gc} _{2,\mtt j, \mss M}} [ \mss m_{\mtt k}^p\,] $ and
$\widetilde{\mss m}_{\mtt j, \mtt k}^p$ is equal to
\begin{align*}
& \frac{c_{2,1}}{ |D_{\mtt k} |^2} \,
\Big( \, \frac{1}{ |D_{\mtt k} |} \,
\sum_{x\in D_{\mtt k}} \chi(\varrho (x/n))\, [ 1 + \varrho (x/n))]   \,\Big)
\, \Big( \, \frac{1}{ |D_{\mtt k} |} \,
\sum_{x\in D_{\mtt k}} \varrho (x/n)  \,\Big)^{p-3}
\\
& \quad +\, \frac{c_{2,2}}{ |D_{\mtt k} |^2} \,
\Big( \, \frac{1}{ |D_{\mtt k} |} \,
\sum_{x\in D_{\mtt k}} \chi(\varrho (x/n))  \,\Big)
\, \Big( \, \frac{1}{ |D_{\mtt k} |} \,
\sum_{x\in D_{\mtt k}} \varrho (x/n)  \,\Big)^{p-3} 
\\
& \quad +\, \frac{c_{2,2}}{ |D_{\mtt k} |^2} \,
\Big( \, \frac{1}{ |D_{\mtt k} |} \,
\sum_{x\in D_{\mtt k}} \chi(\varrho (x/n))  \,\Big)
\Big( \, \frac{1}{ |D_{\mtt k} |} \,
\sum_{x\in D_{\mtt k}} \varrho (x/n)^2  \,\Big)
\, \Big( \, \frac{1}{ |D_{\mtt k} |} \,
\sum_{x\in D_{\mtt k}} \varrho (x/n)  \,\Big)^{p-4} \,+\, \mf R_r \,.
\end{align*}
Apply \eqref{85} once more to obtain that this expression is equal to
\begin{align*}
\frac{c_{2,1}}{ |D_{\mtt k} |^2} \,
\chi(\widetilde{\mss m}_{\mtt j, \mtt k})\, [ 1 +
\widetilde{\mss m}_{\mtt j, \mtt k} ] \, \widetilde{\mss m}_{\mtt j,
\mtt k}^{p-3} \,+\,
\frac{c_{2,2}}{ |D_{\mtt k} |^2} \,
\chi(\widetilde{\mss m}_{\mtt j, \mtt k}) \,
[\, 1 + \widetilde{\mss m}_{\mtt j, \mtt k} \, ]\, 
\widetilde{\mss m}_{\mtt j, \mtt k}^{p-3}   \,+\, \mf R_r  \, + \, \mf R_{n,r}\,.
\end{align*}
To complete the proof, it remains to recollect all previous estimates.
\end{proof}

The next statement is a consequence of Corollary \ref{s38} and the
previous result. 

\begin{corollary}
\label{s39}
Fix $0 < \delta < 1/2$, $p\ge 1$. Then, there exists a constant
$\mf a_{3,3} = \mf a_{3,3}(p, \delta)$ such that
\begin{equation*}
\begin{aligned}
\sum_{\mtt k\in\Lambda_{\mtt j}} J_{\mtt k} \,
E_{\nu^{\rm gc} _{2,\mtt j, \mss M} } [\mss m_{\mtt k}^p ]
\, &=\, s_2^d\, J_{2,\mtt j}\, \mss Q_p(\mss m_{2,\mtt j})
\,+\, \mss Q'_p(\mss m_{2,\mtt j}) \,
\chi (\mss m_{2,\mtt j})\,
\sum_{\mtt k\in\Lambda_{\mtt j}}  [\, J_{\mtt k} - J_{2, \mtt j}\,] \,
\mtt Y_{\mtt j, \mtt k}  \\
& +\, \frac{1}{2}\, \mss Q''_p(\mss m_{2,\mtt j}) \,
\chi(\mss m_{2,\mtt j})^2\,  \sum_{\mtt k\in\Lambda_{\mtt j}}  J_{\mtt k} \,
H^{(1)}_{\mtt j, \mtt k} \,
(\mtt Y_{\mtt j, \mtt k}  + \mtt Y^{(2)}_{\mtt j, \mtt k} )
\\
& +\, \frac{1}{6}\, \mss Q'''_p(\mss m_{2,\mtt j}) \,
\chi(\mss m_{2,\mtt j})^3 \sum_{\mtt k\in\Lambda_{\mtt j}}  J_{\mtt k} \,
[\,H^{(1)}_{\mtt j, \mtt k} \,]^3\,
\,+\, \mf R_n
\end{aligned}
\end{equation*}
where
\begin{equation*}
|\, \mf R_n | \, \le \, \mf a_{3,3} \,
\Vert J\Vert_\infty\,  s_2^d\, \Big\{ \, \Big( \frac{s_1s_2 }{n} \Big)^4
\,+\,  \frac{1 }{s_1^{d}}   \, \Big[ \, 
\frac{1 }{s_1^{2d}} \,+\, \Big(\frac{s_1 }{n} \Big)^2\,\Big] \Big\}
\end{equation*}
for all $\mtt j\in \Sigma^+_{s_3}$, 
$\delta \le \mss M/|D_{2,\mtt j}|\le 1-\delta$, $s_1$, $s_2\ge 1$.
\end{corollary}

\begin{proof}
According to Lemma \ref{s22},
\begin{equation*}
\sum_{\mtt k\in\Lambda_{\mtt j}} J_{\mtt k} \, E_{\nu^{\rm gc} _{2,\mtt
j, \mss M} } [\mss m_{\mtt k}^p ]
\,=\, 
\sum_{\mtt k\in\Lambda_{\mtt j}} J_{\mtt k}\,
\mss Q_{p} (\widetilde{\mss m}_{\mtt j, \mtt k}) \,+\, \mf R^{(1)}_n
\,,
\end{equation*}
where
$|\mf R^{(1)}_n|\le \mf a_{3,3}(\delta)\, (s_2/s_1)^{d} \, \{\,
s_1^{-2d} \,+\, (s_1/n)^{2} \}$.

Fix a smooth function $F\colon [0,1]\to \bb R$. By Lemma \ref{s20b},
$\widetilde{\mss m}_{\mtt j, \mtt k} - \mss m_{2,\mtt j}$ is of order
$s_1s_2/n$. Therefore, by a Taylor expansion,
\begin{equation*}
\begin{aligned}
\sum_{\mtt k\in\Lambda_{\mtt j}} J_{\mtt k}\, F(\widetilde{\mss m}_{\mtt j, \mtt k})
\, & =\, F(\mss m_{2,\mtt j}) \, \sum_{\mtt k\in\Lambda_{\mtt j}} J_{\mtt k}\, 
\,+\, F'(\mss m_{2,\mtt j}) \sum_{\mtt k\in\Lambda_{\mtt j}}
[\, J_{\mtt k} \, - J_{2,\mtt j}\,]\,
[\, \widetilde{\mss m}_{\mtt j, \mtt k} - \mss m_{2,\mtt j} \,]
\\
& +\, \sum_{i=2}^3 \frac{1}{i!}\, 
F^{(i)} (\mss m_{2,\mtt j}) \, \sum_{\mtt k\in\Lambda_{\mtt j}} J_{\mtt k}\, 
[\, \widetilde{\mss m}_{\mtt j, \mtt k} - \mss m_{2,\mtt j} \,]^i
\,+\, \mf R^{(2)}_n \,,
\end{aligned}
\end{equation*}
where $F^{(i)}$ represents the $i$-th derivative of $F$, and
$|\mf R^{(2)}_n|\le \mf a_{3,3}\, \Vert J\Vert_\infty\, s_2^d\,
(s_1s_2/n)^4$. Moreover, in the second term on the right-hand side we
were allowed to introduce $J_{2,\mtt j}$ because the sum
$\sum_{\mtt k\in\Lambda_{\mtt j}} [\, \widetilde{\mss m}_{\mtt j, \mtt
k} - \mss m_{2,\mtt j} \,]$ vanishes. To complete the proof, it
remains to recall the statement of Lemma \ref{s20b} to replace in the
previous formula
$\widetilde{\mss m}_{\mtt j, \mtt k} - \mss m_{2,\mtt j} $ by
$\chi (\mss m_{2,\mtt j})\, \mtt Y_{\mtt j, \mtt k} + \mf R_n$.
\end{proof}

The next result, for the canonical expectation, is a consequence of
Lemmata \ref{s17}, \ref{s22} and Corollary \ref{s39}.

\begin{corollary}
\label{s30b}
Fix $0<\delta <1/2$, and $p\ge 1$. Then, there exists a constant
$\mf a_{3,3} = \mf a_{3,3} (p, \delta)$ such that
\begin{equation}
\label{110}
\begin{aligned}
\sum_{\mtt k\in\Lambda_{\mtt j}} J_{\mtt k} \,
E_{\nu^{\rm c} _{2,\mtt j, \mss M} } [\mss m_{\mtt k}^p ]
\, &=\, s_2^d\, J_{2,\mtt j}\, \mss Q_p(\mss m_{2,\mtt j})
\,+\, \mss Q'_p(\mss m_{2,\mtt j}) \,
\chi (\mss m_{2,\mtt j})\,
\sum_{\mtt k\in\Lambda_{\mtt j}}  [\, J_{\mtt k} - J_{2, \mtt j}\,] \,
\mtt Y_{\mtt j, \mtt k}  \\
& +\, \frac{1}{2}\, \mss Q''_p(\mss m_{2,\mtt j}) \,
\chi(\mss m_{2,\mtt j})^2\,  \sum_{\mtt k\in\Lambda_{\mtt j}}  J_{\mtt k} \,
H^{(1)}_{\mtt j, \mtt k} \,
(\mtt Y_{\mtt j, \mtt k}  + \mtt Y^{(2)}_{\mtt j, \mtt k} )
\\
& +\, \frac{1}{6}\, \mss Q'''_p(\mss m_{2,\mtt j}) \,
\chi(\mss m_{2,\mtt j})^3 \sum_{\mtt k\in\Lambda_{\mtt j}}  J_{\mtt k} \,
[\,H^{(1)}_{\mtt j, \mtt k} \,]^3
\\
&- \, \frac {1}{2}\, s_2^d\, J_{2,\mtt j}\, \frac{ p\, (p-1) }{ (s_1 s_2)^d}
\, \chi(\mss m_{2, \mtt j})  \, \mss m_{2,\mtt j}^{p-2}
\,+\, \mf R_n
\end{aligned}
\end{equation}
where
\begin{equation*}
|\, \mf R_n | \, \le \, \mf a_{3,3} \, \Vert J\Vert_\infty
\, s_2^d\, \Big\{ \, \Big(
\frac{s_1s_2 }{n} \Big)^4
\,+\,  \frac{s_1s_2 }{n} \, \frac{1} {(s_1s_2)^d } \, 
\,+\,  \frac{1 }{s_1^{d}}   \, \Big[ \, 
\frac{1 }{s_1^{2d}} \,+\, \frac{1 }{(s_1s_2)^{d}} \,+\,
\Big(\frac{s_1 }{n} \Big)^2\,\Big] \Big\}
\end{equation*}
for all $\mtt j\in \Sigma^+_{s_3}$, 
$\delta \le \mss M/|D_{2,\mtt j}|\le 1-\delta$, $s_1$, $s_2\ge 1$.

\end{corollary}

\subsection*{Radon-Nikodym derivatives}

In this subsection, we present estimates on the Radon-Nikodym
derivative of the canonical measure projected on a small cube with
respect to the grand canonical measure. 

\begin{lemma}
\label{s10}
For each $0<\delta<1/2$, there exists a finite constant
$\mf a_1 = \mf a_1(\delta)$ such that
\begin{equation*}
\nu^{\rm c} _{2,  \mtt j, \mss M}
\big[\, \{\eta : \eta_x = \xi_x
\;\forall x\in D_{\mtt k}\,\}\,\big]
\,\le\,
\mf a_1 \,
\nu^{\rm gc} _{2, \mtt j,  \mss M} 
\big[\, \{\eta : \eta_x = \xi_x
\;\forall x\in D_{\mtt k} \,\}\,\big]
\end{equation*}
for all $\mtt j \in \Sigma^+_{s_3} $, $\mtt k \in \Lambda_{\mtt j}$,
$\delta \le \mss m \le 1-\delta$, $\xi \in \{0,1\}^{D_{\mtt k} }$, and
$s_1\ge 2$, $s_2\ge 2$. Here,
$\mss m = \mss M/|D_{2,\mtt j} |$.
\end{lemma}

\begin{proof}
By definition of the canonical measure, the left-hand side of the
equation appearing in the statement of the lemma is equal to
\begin{equation*}
\nu^{\rm gc} _{2,\mtt j, \mss M} 
\big[\, \{\eta : \eta_x = \xi_x
\;\forall x\in D_{\mtt k}\,\}\,\big] \;
\frac{\nu^{\rm gc} _{2,\mtt j, \mss M} 
\big[\, \sum_{x\in D_{2,\mtt j, \mtt k}}  \eta_x = \mss M -
|\xi|\, \big]}
{\nu^{\rm gc} _{2,\mtt j, \mss M}  \big[\, \sum_{x\in D_{2,\mtt j}}   \eta_x =
\mss M \,\big]}\,,
\end{equation*}
where
$\cb{D_{2, \mtt j, \mtt k} } = D_{2,\mtt j} \setminus D_{\mtt
k}$, $\cb{|\xi| } = \sum_{x\in D_{\mtt k}} \xi_x$. It remains to bound
the previous ratio. This is done with the help of the local central
limit theorem.

Since the chemical potential of the product measure has been chosen
for the density of particles on $D_{2, \mtt j}$ to be
$\mss m = \mss M/|D_{2, \mtt j} |$, see \eqref{50}, by Theorem
\ref{s29} with $\Sigma^+_r = D_{2, \mtt j}$,
\begin{equation*}
\sqrt{2\pi \, \Gamma_{2, \mtt j}} \, \nu^{\rm gc} _{2,\mtt j, \mss M}  \big[\,
\sum_{x\in D_{2, \mtt j}}   \eta_x = \mss M \,\big]
\,=\, 1 \,+\, \frac{\sqrt{2\pi} }{|D_{2, \mtt j} |^{1/2}}\,
q_{1,L, \varnothing}(0) \,+\, O\Big(\, \frac{1}{ |D_{2, \mtt j} |}\,\Big)\,,
\end{equation*}
where
$\cb{\Gamma_{2, \mtt j}} := \sum_{x\in D_{2, \mtt j}} \chi (\,
\varrho^{(2)}_{\mtt j, \mss m} (x/n)\,)$,
$\varrho^{(2)}_{\mtt j, \mss m} (x/n)$, $x\in D_{2, \mtt j}$, is the
density profile defined in \eqref{97}, $L = |D_{2, \mtt j} |$, and
$q_{1,L, A}(\cdot) $ has been introduced in \eqref{78}.  As
$q_{1,L, \varnothing}(0)=0$, by Theorem \ref{s29}, the right-hand side
is bounded below by $1/2$ for $s_1s_2$ sufficiently large, uniformly
in $\delta \le \mss m \le 1 -\delta$, $\mtt j \in \Sigma^+_{s_3} $.

We turn to the numerator. By the same Theorem \ref{s29},
\begin{equation*}
\sqrt{2\pi \, \Gamma_{2, \mtt j, \mtt k}} \, \nu^{\rm gc} _{2,\mtt j, \mss M}  \Big[\,
\sum_{x\in D_{2, \mtt j, \mtt k}}   \eta_x = \mss M - |\xi| \,\Big]
\,=\, e^{-z^2/2} \,+\, \frac{\sqrt{2\pi} }{|D_{2, \mtt j, \mtt k} |^{1/2}}\,
q_{1, L, D_{\mtt k}}(z) \,+\, O\Big(\, \frac{1}{ |D_{2, \mtt j, \mtt k} |}\,\Big)\,,
\end{equation*}
where
$\cb{\Gamma_{2, \mtt j, \mtt k}} := \sum_{x\in D_{2, \mtt j, \mtt k}}
\chi\big(\, \varrho^{(2)}_{\mtt j, \mss m} (x/n)\,\big)$,
\begin{equation*}
z\,=\,
\frac{1}{ \sqrt{\Gamma_{2, \mtt j, \mtt k} } } \,
\Big\{\, \sum_{x\in D_{\mtt k}}
\varrho^{(2)}_{\mtt j, \mss m} (x/n) \,
-\, |\xi |\, \Big\}\,.
\end{equation*}
The first term on the right-hand side in the penultimate displayed
equation is bounded by $1$. By the definition \eqref{78} of
$q_{1,L, D_{\mtt k}}(\cdot)$, the second term is bounded by
$C_0\, |D_{2, \mtt j, \mtt k}|^{-1/2}\, H_3(z)\, \lambda_{3, L,
D_{\mtt k}} (\varphi)$. Here, $\varphi$ is the chemical potential
selected in \eqref{21b} to get \eqref{50} with $\mss m$ instead of
$\vartheta$.  By \eqref{79}, the formula for the third Hermite
polynomial, and the one above for $z$, this expression is less than or
equal to
$C_0\, \Gamma_{2,\mtt j, \mtt k}^{-2} \, |D_{\mtt k}|\, \sum_{x\in
D_{2, \mtt j, \mtt k}} \, \gamma_{3, x}(\varphi)$, where, recall,
$\gamma_{3,x} (\varphi)$ is the third cumulant of a Bernoulli
distribution with density $\varrho^{(2)}_{\mtt j, \mss m} (x/n)$.  By
\eqref{56b}, and the explicit expressions for the cumulants of a
Bernoulli distribution,
$C_0\, \Gamma_{2,\mtt j, \mtt k}^{-2} \, |D_{\mtt k}|\, \sum_{x\in
D_{2, \mtt j, \mtt k}} \, \gamma_{3, x}(\varphi) \leq \mf
a_1(\delta)$. Actually, it is bounded by
$\mf a_1(\delta) |D_{\mtt k}|/ \Gamma_{2,\mtt j, \mtt k}$, but this
better bound does not improve the estimate.  To complete the proof for
$s_1s_2$ sufficiently large, it remains to recollect all previous
estimates.

For the other cases, note that, for each
fixed $s_1$, $s_2$, the
right-hand side of the equation appearing in the statement of the lemma
is bounded away from $0$, uniformly in $\delta \le \mss m\le
1-\delta$. Thus,  to extend the result for small $s_1s_2$, it
is enough to adjust the constant $\mf a_1$.
\end{proof}

Recall from Theorem \ref{s26}, and \eqref{82}, the
definitions of $\mss m_{\mtt k}$, $\mss m_{\mtt j, \mtt k} $.

\begin{corollary} 
\label{s25}
For each $0<\delta<1/2$, there exists a finite constant
$\mf a_1=\mf a_1(\delta)$ such that
\begin{equation*}
\nu^{\rm c} _{2,\mtt j, \mss M}
\big \{ \, \big | \,\mss m_{\mtt k} - \mss m_{\mtt j, \mtt k} \big | \; \ge 
\;   a \big\}
\; \le \;  \mf a_1  \exp \big\{ - (a^2/2) \,s_1^d \, \big\}
\end{equation*}
for all $\mf a_1\,  (s_1 s_2)^{-d} \le a\le 1$,
$\mtt j \in \Sigma^+_{s_3}$, $\mtt k \in \Lambda_{\mtt j}$,
$\delta \le \mss m_{2, \mtt j} \le 1-\delta$, and
$s_1\ge 2$, $s_2\ge 2$.  
\end{corollary}

\begin{proof}
By Lemma \ref{s10}, 
\begin{equation*}
\nu^{\rm c} _{2,\mtt j, \mss M}
\big \{ \, \big | \,\mss m_{\mtt k} - \mss m_{\mtt j, \mtt k} \big | \; \ge 
\;   a \big\}
\,\le\, \mf a_1 \, 
\nu^{\rm gc} _{2,\mtt j, \mss M}
\big \{ \, \big | \,\mss m_{\mtt k} - \mss m_{\mtt j, \mtt k} \big |
\; \ge  \;   a \big\}
\end{equation*}
for some finite constant $\mf a_1 = \mf a_1(\delta)$. By Lemma
\ref{s17} with $p=1$,
$| \, \widetilde{\mss m}_{\mtt j, \mtt k} - \mss m_{\mtt j, \mtt k} \,
| \,\le\, \mf a_1\,  (s_1s_2)^{-d}$. Thus, the previous expression is less
than or equal to
\begin{equation*}
\mf a_1 \, \nu^{\rm gc} _{2,\mtt j, \mss M}
\big \{ \, \,\mss m_{\mtt k} - \widetilde{\mss m}_{\mtt j, \mtt k} 
\; \ge  \;   a/2 \,  \big\}
\,+ \,
\mf a_1 \, \nu^{\rm gc} _{2,\mtt j, \mss M}
\big \{ \, \,\mss m_{\mtt k} - \widetilde{\mss m}_{\mtt j, \mtt k} 
\; \le  \;   -\, a/2 \,  \big\}\,
\end{equation*}
provided $a \ge 2 \,\mf a_1 \, (s_1s_2)^{-d}$. These two expressions
are estimated in the same way. We present the details for the first
one.

We proceed as in the classical proof of the upper bound large
deviations for i.i.d. random variables. There is just an extra
difficulty coming from the fact that the variables are not identically
distributed.  Recall the notation introduced at the beginning of
Sections \ref{sec9}, \ref{sec12}.  By Chebyshev exponential
inequality, the first term of the previous displayed equation is
bounded by
\begin{equation}
\label{81}
\begin{aligned}
& \mf a_1 \, 
e^{ -\theta \, s_1^d \, (a/2)}
\, E_{\nu^{\rm gc} _{2,\mtt j, \mss M} }
\big [ \,
e^{  \theta\, s_1^d\,  (\mss m_{\mtt k}
- \widetilde{\mss m}_{\mtt j, \mtt k})  } \big]
\\
&\quad
\,=\, \mf a_1\, 
\exp -\, s_1^d \,
\Big\{  \theta \, [ (a/2)  + \bar\varrho_{2, \mtt j, \mtt k}]  \,-\,   
\frac{1}{|D_{\mtt k}|} \sum_{x\in D_{\mtt k} } \ln \big[\,
\varrho^{(2)}_{\mtt j, \mss m}(x/n) \, e^\theta + 1 -
\varrho^{(2)}_{\mtt j, \mss m}(x/n) \, \big] \, \Big\} 
\end{aligned}
\end{equation}
for every $\theta>0$. In this formula, we adopted the notation
introduced at the beginning of this section, and replaced
$\widetilde{\mss m}_{\mtt j, \mtt k}$ by
$\cb{\bar\varrho_{2, \mtt j, \mtt k}} := |D_{\mtt k}|^{-1} \sum_{x\in
D_{\mtt k} } \varrho^{(2)}_{\mtt j, \mss m} (x/n)$, where
$\mss m = \mss m_{2, \mtt j}$.

It is easily seen that the expression inside braces is strictly
concave in $\theta$. Let $\theta (a)$ be the parameter which maximizes
the expression, so that
\begin{equation}
\label{80}
(a/2)  + \bar\varrho_{2, \mtt j, \mtt k}  \,=\,   
\frac{1}{|D_{\mtt k}|} \sum_{x\in D_{\mtt k} }
\frac{\varrho^{(2)}_{\mtt j, \mss m}(x/n) \, e^{\theta(a)} }
{\varrho^{(2)}_{\mtt j, \mss m}(x/n) \, e^{\theta (a)} + 1 -
\varrho^{(2)}_{\mtt j, \mss m}(x/n) } \;\cdot
\end{equation}
Note that $\theta$ depends on $\mss m$ and $\rho(\cdot)$.  By
definition of $\bar\varrho_{2, \mtt j, \mtt k}$, $\theta(0)=0$.  With
this notation, the expression inside braces on the right-hand side of
\eqref{81} becomes
\begin{equation*}
R(a)\,:=\, 
\theta(a)  \, [ (a/2)  + \bar\varrho_{2, \mtt j, \mtt k}]  \,-\,   
\frac{1}{|D_{\mtt k}|} \sum_{x\in D_{\mtt k} } \ln \big[\,
\varrho^{(2)}_{\mtt j, \mss m}(x/n) \, e^{\theta(a)} + 1 -
\varrho^{(2)}_{\mtt j, \mss m}(x/n) \, \big] \,.
\end{equation*}
Since $\theta(0) =0$, $R(0)=0$. A computation yields that
$R'(a)= \theta(a)/2$. In particular, $R'(0)=0$, and $R''(a) = (1/2)
\, \theta'(a)$. By \eqref{80},
\begin{equation*}
\frac{1}{2}   \,=\,  \theta'(a)\, 
\frac{1}{|D_{\mtt k}|} \sum_{x\in D_{\mtt k} }
\chi ( \rho_x(a))   \, \quad \text{where}\;\;
\rho_x(a)  \,=\, 
\frac{\varrho^{(2)}_{\mtt j, \mss m}(x/n) \, e^{\theta(a)} }
{\varrho^{(2)}_{\mtt j, \mss m}(x/n) \, e^{\theta (a)} + 1 -
\varrho^{(2)}_{\mtt j, \mss m}(x/n) } \;\cdot
\end{equation*}
Thus, $1/2 \le (1/4) \, \theta'(a)$, so that  $\theta'(a) \ge 2$, and 
$R''(a) \ge 1$, which completes the proof of the corollary.
\end{proof}

\section{The space $L^1(0,T; \mc H_{- \mtt r})$}
\label{sec-c3}

In this section, we present some results on the tightness of a family
of probability measures on $L^1(0,T; \mc H_{- \mtt r})$, $T>0$,
$\mtt r\in \bb R$, and the characterization of measures on this
space. Here, the time interval $[0,T]$ is endowed with the Borel
$\sigma$-algebra $\ms B = \ms B([0,T])$ to turn the measure space
$([0,T], \ms B, \lambda)$, where $\lambda$ is the Lebesgue measure, 
into a countably generated measure space.

\subsection*{Compact sets in $L^1(0,T; \mc B)$}

Let $\mc B$ be a separable Banach space.  By \cite[Example 1.2.28 and
Proposition 1.2.29]{hnvw16}, $L^1(0,T; \mc B)$ is a separable Banach
space, and therefore a Polish space.  By Simon \cite[Theorem 1]{s86},
a subset $\mc F$ of $L^1(0, T; \mc B)$ is relatively compact if, and
only if,
\begin{itemize}
\item[(a)] For all $0<t_1<t_2<T$, the set
$\{ \int_{t_1}^{t_2} X(t)\, dt : X \in \mc F\} $ is a relatively compact
subset of $\mc B$.

\item[(b)] $\lim_{h\to 0} \sup_{X\in \mc F} \int_0^{T-h} \| X(t+h)
-X(t)\|_{\mc B} \; dt = 0$.
\end{itemize}

For $0<\sigma<1$, let
\begin{align*}
\cb{W^{\sigma, 1} (0, T ; \mc B ) } \, :=\,
\Big\{\, X \in L^1 (0, T ; \mc B) : \int_0^T \int_0^T
\frac{ \Vert X (t) - X(s)\Vert_{\mc B}} {|t-s|^{1+\sigma}}
\;ds\,dt < \infty\,\Big\}\,.
\end{align*}
For $X\in W^{\sigma, 1} (0, T ; \mc B)$, define the seminorm
\begin{align*}
\cb{ \Vert X \Vert_{\dot W^{\sigma, 1} (0, T ; \mc B)} } \, :=\,
\int_0^T \int_0^T
\frac{ \Vert X (t) - X(s)\Vert_{\mc B}} {|t-s|^{1+\sigma}}
\;ds\,dt\,.
\end{align*}
By \cite[Lemma 5]{s86}, there exists a finite constant
$C_0=C_0(\sigma, T)$ such that
\begin{align}
\label{222}
\int_0^{T-h} \| X(t+h) -X(t)\|_{\mc B} \; dt \,\le\, C_0\, h^\sigma\,
\Vert X \Vert_{\dot W^{\sigma, 1} (0, T ; \mc B)} 
\end{align}
for all $0<h<T$ and $X\in L^1 (0, T ; \mc B) $. In particular, a
subset $\mc F$ of $L^1(0, T; \mc B)$ is relatively compact if it
satisfies (a) and
\begin{itemize}
\item[(b')] $\sup_{X\in \mc F}
\Vert X \Vert_{\dot W^{\sigma, 1} (0, T ; \mc B)}  < \infty$.
\end{itemize}

\subsection*{Tightness}

In this subsection, we establish sufficient conditions for the
tightness of a sequence of probability measures $(\bb Q_n:n\ge 1)$ on
$L^1(0, T; \mc H_{- \mtt r})$.

As $L^1(0, T; \mc H_{- \mtt r})$ is complete and separable, each
probability measure is tight.  As $\mc H_{- \mtt r'}$ is compactly
embedded in $\mc H_{- \mtt r}$ if $\mtt r'< \mtt r$, we have the
following sufficient condition for the tightness of sequences of
probability measures on $L^1(0, T; \mc H_{- \mtt r})$.

\begin{lemma}
\label{s89}
Fix $\mtt r\in \bb R$.  A sequence of measures $(\bb Q_n:n\ge 1)$ on
$L^1(0, T; \mc H_{- \mtt r})$ is tight if there exists
$\mtt r'< \mtt r$ such that for all $\epsilon>0$,
\begin{gather}
\tag{A}
\lim_{A\to\infty} 
\limsup_{n\to\infty} \bb Q_{n} \Big[ \, \sup_{0<t_1<t_2<T} \Big\|\int_{t_1}^{t_2}
\mtt x (t)  \, dt \, \Big\|_{\mc H_{-\mtt r'}} >A \,\Big ] \,=\,  0\,,
\\
\lim_{h\to 0} 
\limsup_{n\to\infty}  \bb Q_{n}
\Big[ \,  \int_0^{T-h}
\Vert \mtt x(t+h) - \mtt x(t) \Vert_{\mc H_{- \mtt r}}dt
> \epsilon  \, \Big] = 0\;.
\tag{B}
\end{gather}
\end{lemma}

The next result follows from the previous lemma and \eqref{222}.

\begin{lemma}
\label{s63}
Fix $\mtt r\in \bb R$.  A sequence of measures $(\bb Q_n:n\ge 1)$ on
$L^1(0, T; \mc H_{- \mtt r})$ is tight if there exist
$\mtt r'< \mtt r$ and $0<\sigma<1$ such that
\begin{gather}
\tag{A}
\lim_{A\to\infty} 
\limsup_{n\to\infty} \bb Q_{n} \Big[ \, \sup_{0<t_1<t_2<T} \Big\|\int_{t_1}^{t_2}
\mtt x (t)  \, dt \, \Big\|_{\mc H_{-\mtt r'}} >A \,\Big ] \,=\,  0\,,
\\
\lim_{A\to\infty} 
\limsup_{n\to\infty}  \bb Q_{n}
\Big[ \,  \Vert \mtt x \Vert_{\dot W^{\sigma, 1} (0, T ; \mc H_{- \mtt r})}  > A  \,
\Big] = 0\;.
\tag{B}
\end{gather}
\end{lemma}

The following result should be compared with the
Kolmogorov--\v{C}entsov theorem. In that theorem, the constant $b$
appearing in Lemma \ref{s57} is required to satisfy $b>1$, whereas
here it suffices that $b>0$. This is the main advantage of working
with the weaker topology of $L^1(0,T ; \mc H_{- \mtt r})$ rather than
$C([0,T],\mc H_{- \mtt r})$.

\begin{lemma}
\label{s57}
Suppose that there exist a finite constant $C_0$ and $b>0$ such that
\begin{equation*}
\sup_{n\ge 1}  \,
E_{\bb Q_{n}} \big[ \,  
\Vert \mtt x (t) - \mtt x(s)\Vert_{\mc H_{- \mtt r}}\, \big] \, \le \,
C_0\, |t-s|^b
\end{equation*}
for all $0\le s<t\le T$. Then, condition {\rm (B)} of the previous lemma is
fulfilled for $0<\sigma<b$.
\end{lemma}

\begin{proof}
Fix $0<\sigma<b$.  By Chebyshev's inequality and the definition of the
norm $\Vert \,\cdot\, \Vert_{\dot W^{\sigma, 1} (0, T ; \mc H_{- \mtt r})} $,
\begin{equation*}
\bb Q_{n} \Big[ \,
\Vert \mtt x \Vert_{\dot W^{\sigma, 1} (0, T ; \mc H_{- \mtt r})}  >A \,\Big ]
\,\le\,
\frac{1}{A}\,
\int_0^T \int_0^T
\frac{1} {|t-s|^{1+\sigma}}\,
E_{\bb Q_{n}} \big[ \,  
\Vert \mtt x (t) - \mtt x(s)\Vert_{\mc H_{- \mtt r}}\, \big]
\;ds\,dt\,.
\end{equation*}
By hypothesis, the right-hand side is bounded by $(2C_0/A) [a
(1+a)]^{-1} T^{1+a}$, where $a=b-\sigma$. As this bound does not
depend on $n$ and vanishes as $A\to\infty$, condition {\rm (B)} of
Lemma \ref{s63} holds.
\end{proof}

The condition of the previous lemma can be expressed in a simpler
form. Let $\{\omega_m: m\in \bb Z^d\}$ be the sequence given by
$\cb{\omega_m = \omega\, (1+\|m\|)^{-(d+1)}}$, where $\omega$ is set
by $\sum_{m\in \bb Z^d}\omega_m=1$. We claim that for every
$\mtt r\in \bb R$ and every $Y \in \mc H_{-\mtt r}$,
\begin{equation}
\label{236}
\Vert Y \Vert_{\mc H_{- \mtt r}} \;\le\;
\sum_{m\in \bb Z^d} \gamma_m^{-\mtt r/2}\,  \omega^{-1/2}_m \,
\big|\, \< Y ,\, \phi_m\>\, \big| \;.
\end{equation}

To prove this assertion, denote by $a$ the right-hand side of
\eqref{236}. Clearly, we may assume that $a<\infty$.  For each
$m\in\bb Z^d$, the $m$-th term is bounded by $a$. Thus,
$\gamma_m^{-\mtt r/2} |\<Y,\phi_m\>| \le \omega_m^{1/2}\, a$.  Taking
the square, and summing over $m$ yields that
$\Vert Y\Vert^2_{\mc H_{-\mtt r}} = \sum_m \gamma^{-\mtt r}_m
|\<Y,\phi_m\>|^2 \le a^2 \sum_m \omega_m = a^2$, as claimed.

\begin{lemma}
\label{s64}
Suppose that there exist a sequence $\{\beta_m : m\in\bb Z^d\}$ and
$b>0$ such that
\begin{align*}
\sup_{n\ge 1}
E_{\bb Q_{n} } \Big[ \,    \, |\, \< \mtt x (t) - \mtt x(s) ,\, \phi_m\>\, | \,
\Big] \,\le\, \beta_m \, |t-s|^b
\end{align*}
for all $m\in\bb Z^d$, $0\le s<t\le T$ and
\begin{align*}
\sum_{m\in \bb Z^d} \beta_m \, \gamma_m^{-\mtt r/2}\,  \omega^{-1/2}_m
\,<\,\infty\,. 
\end{align*}
Then, condition {\rm (B)} of the Lemma \ref{s63} is fulfilled for
$0<\sigma<b$.
\end{lemma}

\begin{proof}
Taking expectations in \eqref{236} and applying Fubini yields that
\begin{equation}
\label{223}
E_{\bb Q_{n}} \big[ \,
\Vert \mtt x (t) - \mtt x(s)\Vert_{\mc H_{- \mtt r}}\, \big]
\;\le\;
\sum_{m\in \bb Z^d} \gamma_m^{-\mtt r/2}\,  \omega^{-1/2}_m \,
E_{\bb Q_{n} } \Big[ \, |\, \< \mtt x (t) - \mtt x(s) ,\, \phi_m\>\, |
\, \Big] 
\end{equation}
for all $0\le s < t\le T$ and $n\ge 1$.

By assumption, the right-hand side of \eqref{223} is less than or
equal to $C_0 \, |t-s|^b$ for some finite constant $C_0$. This bound
holds for all $n\ge 1$ and $0\le s\le t \le T$. This proves that the
condition of Lemma \ref{s57} holds.
\end{proof}

A similar argument provides a simple sufficient condition for
requirement (A) in Lemmata \ref{s89} and \ref{s63}.

\begin{lemma}
\label{s66}
Condition (A) in Lemmata \ref{s89} and \ref{s63} holds provided
\begin{equation}
\label{158}
\limsup_{n\to\infty}
E_{\bb Q_{n}}  \Big[ \, \sum_{m\in \bb Z^d}
(\gamma^{-\mtt r'}_m/\omega_m)^{1/2}
\,  \int_{0}^{T}
\big|\,  \< \mtt x (t)  ,\, \phi_m\> \big | \, \, dt \, \Big ]
\,<\, \infty \,.
\end{equation}
\end{lemma}

\begin{proof}
By definition of the
$\mc H_{-\mtt r'}$-norm, 
\begin{align*}
\bb Q_{n} \Big[ \, \sup_{0<t_1<t_2<T} \Big\|\int_{t_1}^{t_2}
\mtt x (t)  \, dt \, \Big\|_{\mc H_{-\mtt r'}} >A \,\Big ]
\, & \le\, 
\bb Q_{n}  \Big[ \, \sup_{0<t_1<t_2<T}  \sum_{m\in \bb Z^d}
\gamma^{-\mtt r'}_m\,  \Big|\,  \int_{t_1}^{t_2}
\< \mtt x (t)  ,\, \phi_m\> \, dt \, \Big |^2 > A^2 \,\Big ]
\\
& \le\, 
\bb Q_{n}  \Big[ \, \sum_{m\in \bb Z^d}
\gamma^{- \mtt r'}_m\,  \sup_{0<t_1<t_2<T}   \Big|\,  \int_{t_1}^{t_2}
\< \mtt x (t)  ,\, \phi_m\> \, dt \, \Big |^2 > A^2 \,\Big ]\,.
\end{align*}
Recall the definition of the sequence $\{\omega_m : m\in \bb Z^d\}$
and the argument employed in the previous lemma to estimate this
expression by
\begin{align*}
& \sum_{m\in \bb Z^d} \bb Q_{n}  \Big[ \, 
(\gamma^{- \mtt r'}_m/\omega_m)
\,  \sup_{0<t_1<t_2<T}   \Big|\,  \int_{t_1}^{t_2}
\< \mtt x (t)  ,\, \phi_m\> \, dt \, \Big |^2 > A^2 \,\Big ]
\\
& =\,
\sum_{m\in \bb Z^d} \bb Q_{n}  \Big[ \, 
(\gamma^{-\mtt r'}_m/\omega_m)^{1/2}
\,  \sup_{0<t_1<t_2<T}   \Big|\,  \int_{t_1}^{t_2}
\< \mtt x (t)  ,\, \phi_m\> \, dt \, \Big | > A \,\Big ]
\,.
\end{align*}
By Chebyshev's inequality this expression is less than or equal to
\begin{align*}
\frac{1}{A} \, E_{\bb Q_{n}}  \Big[ \, \sum_{m\in \bb Z^d}
(\gamma^{-\mtt r'}_m/\omega_m)^{1/2}
\,  \int_{0}^{T}
\big|\,  \< \mtt x (t)  ,\, \phi_m\> \big | \, \, dt \, \Big ]\,.
\end{align*}
By \eqref{158}, the $\limsup_n$ of this expression vanishes as
$A\to\infty$. This completes the proof of the lemma.
\end{proof}

We summarize the previous results in the next corollary, which
provides simple conditions for tightness of a sequence of measures in
$L^1(0,T; \mc H_{- \mtt r})$.

\begin{corollary}
\label{s67}
Fix $\mtt r\in \bb R$, and a sequence of measures $(\bb Q_n:n\ge 1)$ on
$L^1(0,T ; \mc H_{- \mtt r})$.  Suppose that there exists
$\mtt r'< \mtt r$ such that
\begin{equation}
\tag{A}
\limsup_{n\to\infty}
E_{\bb Q_{n}}  \Big[ \, \sum_{m\in \bb Z^d}
(\gamma^{-\mtt r'}_m/\omega_m)^{1/2}
\,  \int_{0}^{T}
\big|\,  \< \mtt x (t)  ,\, \phi_m\> \big | \, \, dt \, \Big ]
\,<\, \infty \,.
\end{equation}
Suppose, furthermore, that there exist a sequence
$\{\beta_m : m\in\bb Z^d\}$ and $b>0$ such that
\begin{align}
\tag{B}
\sup_{n\ge 1}
E_{\bb Q_{n} } \Big[ \,    \, |\, \< \mtt x (t) - \mtt x(s) ,\, \phi_m\>\, | \,
\Big] \,\le\, \beta_m \, |t-s|^b
\end{align}
for all $m\in\bb Z^d$, $0\le s<t\le T$ and
\begin{align}
\label{239}
\sum_{m\in \bb Z^d} \beta_m \, \gamma_m^{-\mtt r/2}\,  \omega^{-1/2}_m
\,<\,\infty\,. 
\end{align}
Then, the sequence $(\bb Q_n:n\ge 1)$ on $L^1(0,T; \mc H_{- \mtt r})$ is
tight.
\end{corollary}

We may also state a sufficient criterion with respect to condition (B)
in Lemma \ref{s89} by similar arguments, following the same route as
above.

\begin{corollary}
\label{s90}
Fix $\mtt r\in \bb R$, and a sequence of measures $(\bb Q_n:n\ge 1)$ on
$L^1(0,T ; \mc H_{- \mtt r})$.  Suppose that there exists
$\mtt r'< \mtt r$ such that condition \eqref{158} holds, and that
there exist sequences $\{\beta_m : m\in\bb Z^d\}$ and
$\{b_n: n\geq 1\}$ such that $\lim_n b_n=0$, and
\begin{align}
\tag{B}
E_{\bb Q_{n} } \Big[ \,    \, |\, \< \mtt x (t+h) - \mtt x(t) ,\, \phi_m\>\, | \,
\Big] \,\le\, \beta_m \, (h + b_n)
\end{align}
for all $m\in\bb Z^d$, $n\ge 1$, $0<h<T$, $0<t< T-h$, and
\begin{align*}
\sum_{m\in \bb Z^d}\beta_m \, \gamma_m^{-\mtt r/2}\,  \omega^{-1/2}_m
\,<\,\infty\,. 
\end{align*}
Then, the sequence $(\bb Q_n:n\ge 1)$ on $L^1(0,T; \mc H_{- \mtt r})$ is
tight.
\end{corollary}

\begin{proof}
Let
\begin{equation*}
\mf b \;:=\; \sum_{m\in \bb Z^d}\beta_m \, \gamma_m^{-\mtt r/2}\,
\omega^{-1/2}_m \;<\; \infty\,,
\end{equation*}
which is finite by assumption.  It is enough to check conditions {\rm
(A)} and {\rm (B)} of Lemma \ref{s89}.  Condition {\rm (A)} holds by
Lemma \ref{s66}.

We turn to condition {\rm (B)} of Lemma \ref{s89}. Recall the estimate
\eqref{223}.  Fix $0<h<T$ and $n\ge 1$.  Applying \eqref{223} with the
pair $(t, t+h)$ in place of $(s,t)$, and then the hypothesis {\rm
(B)}, we obtain, for all $0<t<T-h$,
\begin{equation}
\label{224}
E_{\bb Q_{n}} \big[ \,
\Vert \mtt x (t+h) - \mtt x(t)\Vert_{\mc H_{- \mtt r}}\, \big]
\;\le\;
\sum_{m\in \bb Z^d} \gamma_m^{-\mtt r/2}\,  \omega^{-1/2}_m \, \beta_m \,
( h + b_n) \;=\; \mf b \, ( h + b_n) \,.
\end{equation}

By Chebyshev's inequality and Fubini's theorem, for every
$\epsilon>0$,
\begin{equation*}
\bb Q_{n}
\Big[ \,  \int_0^{T-h}
\Vert \mtt x(t+h) - \mtt x(t) \Vert_{\mc H_{- \mtt r}}\, dt
> \epsilon  \, \Big]
\;\le\;
\frac{1}{\epsilon}\, \int_0^{T-h}
E_{\bb Q_{n}} \big[ \,
\Vert \mtt x (t+h) - \mtt x(t)\Vert_{\mc H_{- \mtt r}}\, \big] \; dt \,.
\end{equation*}
By \eqref{224}, the right-hand side is bounded by
$(T/\epsilon)\, \mf b\, ( h + b_n)$.  As $b_n\to 0$,
\begin{equation*}
\limsup_{n\to\infty}  \bb Q_{n}
\Big[ \,  \int_0^{T-h}
\Vert \mtt x(t+h) - \mtt x(t) \Vert_{\mc H_{- \mtt r}}\, dt
> \epsilon  \, \Big]
\;\le\; \frac{T\, \mf b}{\epsilon}\; h \,,
\end{equation*}
an expression which vanishes as $h\to 0$.  This is condition {\rm (B)}
of Lemma \ref{s89}, and completes the proof of the corollary.
\end{proof}

\begin{remark}
\label{s91}
The interest of Corollary \ref{s90} lies in the extra term $b_n$,
which permits the modulus of continuity in time to deteriorate as
$n\to\infty$, provided it does so at a rate which vanishes in the
limit.  On the other hand, one loses the H\"older exponent: the bound
has to be linear in $h$.

The reason is that Corollary \ref{s90} appeals to Lemma \ref{s89}
directly, whereas Corollary \ref{s67} passes through the fractional
Sobolev norm of Lemma \ref{s63}.  Bounding
$\Vert \mtt x \Vert_{\dot W^{\sigma, 1} (0, T ; \mc H_{- \mtt r})}$ requires
the integrability of $|t-s|^{b-\sigma-1}$, which fails for a bound of
the form $\beta_m (|t-s| + b_n)$, since the term $b_n$ contributes
$b_n \int_0^T\!\!\int_0^T |t-s|^{-1-\sigma} ds\, dt = \infty$.
\end{remark}

We conclude this subsection with a result which avoids estimates of
the expectation. 

\begin{lemma}
\label{s92}
Fix $\mtt r \in \bb R$ and a sequence $(\mtt x^n : n\ge 1)$ of random
elements of $L^1(0,T;\mc H_{-\mtt r})$.  Suppose that there exist
$\mtt r' < \mtt r$, $b>0$, a sequence $\{\beta_m : m \in \bb Z^d\}$
with
\begin{equation*}
\mf b \;:=\; \sum_{m\in \bb Z^d} \beta_m \, \gamma_m^{-\mtt r'/2}\,
\omega^{-1/2}_m \;<\;\infty\;,
\end{equation*}
and non-negative random variables $\Theta_n$ such that
$\lim_{A\to\infty} \limsup_{n\to\infty} \bb P[\Theta_n > A] = 0$, and,
almost surely,
\begin{equation}
\label{237}
\int_0^T \big|\, \< \mtt x^n (t) ,\, \phi_m\>\, \big| \, dt \;\le\;
\Theta_n\, \beta_m \;, \qquad
\big|\, \< \mtt x^n (t) - \mtt x^n(s) ,\, \phi_m\>\, \big| \;\le\;
\Theta_n\, \beta_m\, |t-s|^{b}
\end{equation}
for all $m \in \bb Z^d$ and $0\le s<t\le T$.  Then the sequence of
laws of $\mtt x^n$ is tight in $L^1(0,T;\mc H_{-\mtt r})$.
\end{lemma}

\begin{proof}
It is enough to check the two conditions of Lemma \ref{s89}.  By
\eqref{236} and the first bound in \eqref{237}, for all $0<t_1<t_2<T$,
\begin{equation*}
\Big\| \int_{t_1}^{t_2} \mtt x^n (t)\, dt \Big\|_{\mc H_{-\mtt r'}}
\;\le\;
\sum_{m} \gamma_m^{-\mtt r'/2} \omega_m^{-1/2}
\int_0^T \big|\, \< \mtt x^n (t) ,\, \phi_m\>\, \big| \, dt
\;\le\; \mf b \; \Theta_n \;.
\end{equation*}
Since this bound is uniform in $t_1$, $t_2$,
$\bb Q_n[\sup_{t_1<t_2}\Vert \int_{t_1}^{t_2}\mtt x^n\Vert_{-\mtt r'}
>A] \le \bb P[\Theta_n > A/\mf b]$, and condition (A) of Lemma
\ref{s89} follows from the property of the sequence $\Theta_n$.

Similarly, by \eqref{236}, the second bound in \eqref{237}, as
$\gamma_m \ge 1$ and $\mtt r>\mtt r'$,
\begin{equation*}
\int_0^{T-h} \Vert \mtt x^n(t+h) - \mtt x^n(t)\Vert_{\mc H_{-\mtt r}}
\, dt \;\le\; T\, \mf b\, \Theta_n \, h^{b}
\end{equation*}
for all $0<h<T$.  Therefore, for every $\epsilon>0$, by the property
of the sequence $\Theta_n$,
\begin{equation*}
\lim_{h\to 0} \limsup_{n\to\infty} \bb Q_n\Big[ \int_0^{T-h}
\Vert \mtt x^n(t+h) - \mtt x^n(t)\Vert_{\mc H_{-\mtt r}}\, dt >
\epsilon \Big]
\;\le\;
\lim_{h\to 0}  \limsup_{n\to\infty} \bb P\Big[ \Theta_n
> \frac{\epsilon}{T \mf b\,
h^{b}} \Big] \,=\, 0 \;.
\end{equation*}
This proves condition (B) of Lemma \ref{s89}.
\end{proof}

\subsection*{Characterization of probability measures}

By Theorems 1.3.10 and 1.3.21 in \cite{hnvw16}, the dual of
$L^1(0,T; \mc H_{- \mtt r})$ is $L^\infty(0,T; \mc H_{\mtt r})$: for any bounded
linear functional $\Phi\colon L^1(0,T; \mc H_{- \mtt r}) \to \bb R$ there
exists $\mf g\in L^\infty(0,T; \mc H_{\mtt r})$ such that
\begin{equation*}
\Phi (X) \, =\, \Phi_{\mf g} (X)
\quad
\text{for all} \;\; X \in L^1(0,T; \mc H_{- \mtt r})
\,, \;\;
\text{where} \;\;
\cb{\Phi_{\mf g} (X)} \,:=\,
\int_0^T \<X_s , \mf g(s)\> \, ds \,.
\end{equation*}

Denote by $\cb{\ms A}$ the space of functions
$\mtt F \colon L^1 (0,T; \mc H_{- \mtt r}) \to \bb R$ for which there exist
$p\ge 1$, $\mf g_1, \dots, \mf g_p \in L^\infty (0,T; \mc H_{\mtt r})$, and a
continuous bounded function $F\colon \bb R^p \to \bb R$ such that
\begin{equation}
\label{182}
\mtt F  (X) \,=\, F( \Phi_{\mf g_1} (X), \dots, \Phi_{\mf g_p}
(X))\,. 
\end{equation}
Each element of $\ms A$ is clearly a bounded continuous function.  As
$L^\infty (0,T; \mc H_{\mtt r})$ is the dual of $L^1 (0,T; \mc H_{- \mtt r})$, it
separates the points of $L^1 (0,T; \mc H_{- \mtt r})$, and so does $\ms A$.
The set $\ms A$ is also an algebra.

\begin{lemma}
\label{s73}
Let $\bb P_1$, $\bb P_2$ be two probability measures on
$L^1 (0,T; \mc H_{- \mtt r})$. Suppose that
\begin{equation*}
E_{\bb P_1} \big[\, \mtt F(X) \,\big]
\,=\,
E_{\bb P_2} \big[\, \mtt F(X) \,\big]
\end{equation*}
for all $\mtt F\in \ms A$. Then, $\bb P_1 = \bb P_2$.
\end{lemma}

\begin{proof}
We have seen at the beginning of this section that
$L^1 (0,T; \mc H_{- \mtt r})$ is a separable Banach space. It is thus a
complete and separable metric space. Since $\ms A$ is an algebra of
bounded and continuous functions which separates points, the result
follows from \cite[Proposition 3.4.5]{ek09}.
\end{proof}

\begin{lemma}
\label{s72}
For any $\mf g\in L^\infty(0,T; \mc H_{\mtt r})$ and any $\epsilon>0$,
there exist $q'\ge 1$, disjoint open intervals $I_j$ in $[0,T]$, and
elements $\mf g_j\in \mc H_{\mtt r}$, $1\le j\le q'$, such that
$\mf g_\epsilon \in L^\infty(0,T; \mc H_{\mtt r})$ defined by
\begin{gather*}
\mf g_\epsilon (s) = \sum_{j=1}^{q'} \mf g_j \mtt 1_{I_j} (s) 
\;\;\text{satisfies the bounds}
\\ 
\esssup_{0\le s\le T} \Vert \mf g_\epsilon (s)\Vert_{\mc H_{\mtt r}}
\,\le\;
\esssup_{0\le s\le T} \Vert \mf g (s)\Vert_{\mc H_{\mtt r}} \,,\quad
\esssup_{s\in I} \Vert \mf g_\epsilon (s) - \mf g (s)\Vert_{\mc
H_{\mtt r}} \,\le\; \epsilon\,, \;\;
\lambda ([0,T] \setminus I) \,\le\, \epsilon\,,
\end{gather*}
where $I = \cup_j I_j$, and $\lambda$ stands for the Lebesgue
measure. 
\end{lemma}

\begin{proof}
By \cite[Lemma 1.2.19]{hnvw16}, the simple functions are dense in
$ L^\infty(0,T; \mc H_{\mtt r})$ with respect to the convergence in
measure: for any $\mf g\in L^\infty(0,T; \mc H_{\mtt r})$ and any
$\epsilon>0$, there exist $q\ge 1$, measurable subsets $E_j$ of
$[0,T]$, and elements $\mf g_j\in \mc H_{\mtt r}$, $1\le j\le q$, such that
$\mf g_\epsilon \in L^\infty(0,T; \mc H_{\mtt r})$ defined by
\begin{gather*}
\mf g_\epsilon (s) = \sum_{j=1}^q \mf g_j \mtt 1_{E_j} (s) 
\;\;\text{satisfies the bounds}
\\ 
\esssup_{0\le s\le T} \Vert \mf g_\epsilon (s)\Vert_{\mc H_{\mtt r}}
\,\le\;
\esssup_{0\le s\le T} \Vert \mf g (s)\Vert_{\mc H_{\mtt r}} \,,\quad
\esssup_{s\in E} \Vert \mf g_\epsilon (s) - \mf g (s)\Vert_{\mc
H_{\mtt r}} \,\le\; \epsilon\,, \;\;
\lambda ([0,T] \setminus E) \,\le\, \epsilon\,,
\end{gather*}
where $E = \cup_j E_j$, and $\lambda$ stands for the Lebesgue
measure. 

Since the Lebesgue measure on $[0,T]$ is regular, we may replace in
the previous assertion simple functions by step functions, that is,
replace the measurable sets $E_j$, $1\le j\le q$, by open disjoint
intervals $I_{k}$, $1\le k\le q'$, at the price of enlarging the
exceptional set by at most $\epsilon$.
\end{proof}

Denote by $\cb{\ms A_{\rm fdd}}$ the space of functions
$\mtt F \colon L^1 (0,T; \mc H_{- \mtt r}) \to \bb R$ for which there exist
$p\ge 1$, $g_1, \dots, g_p \in \mc H_{\mtt r}$,
$0\le s_1 < t_1 \le s_2 < \cdots < t_{p-1} \le s_p <t_p\le T$, and a
continuous bounded function $F\colon \bb R^p \to \bb R$ such that
\begin{equation*}
\mtt F  (X) \,=\, F\Big( \int_{s_1}^{t_1} \< X_s , g_1\> \, ds ,
\dots,
\int_{s_p}^{t_p} \< X_s , g_p\> \, ds\, \Big)\,.
\end{equation*}
Each element of $\ms A_{\rm fdd}$ is clearly a bounded continuous
function.  

\begin{corollary}
\label{s74}
Let $\bb P_1$, $\bb P_2$ be two probability measures on
$L^1 (0,T; \mc H_{- \mtt r})$. Suppose that
\begin{equation*}
E_{\bb P_1} \big[\, \mtt F(X) \,\big]
\,=\,
E_{\bb P_2} \big[\, \mtt F(X) \,\big]
\end{equation*}
for all all $\mtt F\in \ms A_{\rm fdd}$. Then, $\bb P_1 = \bb P_2$.
\end{corollary}

\begin{proof}
By Lemma \ref{s73}, it is enough to extend the identity to functions
$\mtt F$ in $\ms A$. Fix such a function. By definition, there exists
$p\ge 1$, $\mf g_1, \dots, \mf g_p \in L^\infty (0,T; \mc H_{\mtt r})$, and
a continuous bounded function $F\colon \bb R^p \to \bb R$ for which
\eqref{182} holds. We assume that $p=1$, the same proof applies to the
general case. Hence,
\begin{equation*}
\mtt F(X) \,=\, F \Big( \int_{0}^{T} \< X_s , \mf g (s)\> \, ds\,\Big)
\end{equation*}
for some continuous function $F$ and $\mf g\in  L^\infty (0,T; \mc
H_{\mtt r})$. For $n\ge 1$, denote by $\mf g_n\in  L^\infty (0,T; \mc
H_{\mtt r})$ the approximation given by Lemma \ref{s72} with
$\epsilon=1/n$. Let
\begin{equation*}
\mtt F_n(X) \,=\, F \Big( \int_{0}^{T} \< X_s , \mf g_n(s) \> \, ds\,\Big)\,.
\end{equation*}
By definition of $\mf g_n$, $\mtt F_n$ belongs to $\ms A_{\rm fdd}$,
and, for each $X\in L^1 (0,T; \mc H_{- \mtt r})$, 
\begin{equation*}
\lim_{n\to\infty} \int_{0}^{T} \< X_s , \mf g_n(s) \> \, ds  \, = \,
\int_{0}^{T} \< X_s , \mf g(s) \> \, ds\,.
\end{equation*}
Thus, as $F$ is continuous, 
$\lim_n \mtt F_n (X) = \mtt F(X)$. By the dominated convergence
theorem, and by hypothesis,
\begin{equation*}
E_{\bb P_1} \big[\, \mtt F(X) \,\big] \,=\,
\lim_{n\to\infty} E_{\bb P_1} \big[\, \mtt F_n (X) \,\big] 
\,=\,
\lim_{n\to\infty} E_{\bb P_2} \big[\, \mtt F_n (X) \,\big]
\,=\, E_{\bb P_2} \big[\, \mtt F(X) \,\big] 
\,,
\end{equation*}
which completes the proof of the corollary.
\end{proof}

\section{Ergodic properties of the exclusion process}
\label{sec10}

Recall that $\Sigma_p\subset \bb Z^d$ represents the cube
$\{-p, \dots, p\}^d$. Fix a function $\rho\colon \bb R^d \to [0,1]$ of
class $C^1$. Assume that
$\sup_{x\in [-1,1]^d} |(\nabla \rho)(x)| < \infty$, and that there
exists $0<\mf a_0<1/2$ such that $\mf a_0 \le \rho(x) \le 1-\mf a_0$
for all $x\in [-1,1]^d$.

Fix a sequence $\ell_n\to\infty$ such that $\ell_n/n\to 0$.  In this
section, $\cb{\nu^n_{\rho(\cdot)}}$ represents the product measure on
$\varOmega_n = \{0,1\}^{\Sigma_{\ell_n}}$ whose marginals are given by
$\nu^n_{\rho(\cdot)} \{ \eta: \eta_x=1\} = \rho(x/n)$,
$x\in \Sigma_{\ell_n}$. 

Denote by $\ms L_{n}$ the generator of the exclusion dynamics on
$\varOmega_n$ reversible with respect to the measure
$\nu^n_{\rho (\cdot)}$:
\begin{gather*}
(\ms L_n g) (\xi) \,=\, \sum_{j=1}^d \sum_{x: \{x,x+e_j\} \subset
\Sigma_{\ell_n}} 
\widehat {c}_{x,x+e_j}(\xi)\, 
[\, g(T_{x,x+e_j} \xi) - g(\xi)\,] 
\\
\widehat {c}_{x,x+e_j}(\xi) \,=\, \exp \big\{ - (1/2) [\, H_{x+e_j} - H_x \,]\,
[\xi_{x+e_j} - \xi_x]\,\big\}\,, \quad H_z \,=\,
\ln \frac{ \rho(z/n)}{1-\rho(z/n)}  \,.
\end{gather*}

Denote by $\varOmega_{n, \mss M}$,
$0\le \mss M \le |\Sigma_{\ell_n} |$, the configurations of
$\varOmega_n$ with $\mss M $ particles:
$\cb{\varOmega_{n,\mss M} } :=\{\xi\in \varOmega_n: \sum_{x\in
\Sigma_{\ell_n} } \xi_x = \mss M\}$, and by $\nu^{\rm c}_{n, \mss M}$
the measure $\nu^n _{\rho (\cdot)}(\cdot)$ conditioned to
$\varOmega_{n,\mss M} $:
\begin{equation}
\label{91b}
\cb{\nu^{\rm c}_{n, \mss M} (\xi') } \,=\,
\frac{\nu^n_{\rho(\cdot)} (\xi')}
{\sum_{\zeta \in \varOmega_{n,\mss M}}
\nu^n_{\rho(\cdot)} (\zeta)} \,, \quad
\xi' \,\in\, \varOmega_{n,\mss M}\,.
\end{equation}

Denote by $L_+^2(\nu^{\rm c}_{n, \mss M} )$ the space of non-negative
elements of $L^2(\nu^{\rm c}_{n, \mss M} )$, and by
$I_{n} (\, \cdot \,;\, \nu^{\rm c}_{n, \mss M} )\colon L_+^2(\nu^{\rm
c}_{n, \mss M} ) \to \bb R_+$ the rate functional given by
\begin{align*}
\cb{I_{n} (g  \,;\, \nu^{\rm c}_{n, \mss M} )}
\, :=\, \< \, (-\ms L_n) \sqrt{g} \,,\,  \sqrt{g}\, \>_{\nu^{\rm c}_{n, \mss M}} \,.
\end{align*}
This functional is the Donsker-Varadhan level $2$ large deviations rate
functional if $g$ has mean one, as well as the integrated carr\'e du
champ. It is explicitly given by
\begin{gather*}
I_{n} (g  \,;\, \nu^{\rm c}_{n, \mss M} )
\,=\, \frac{1}{2}\, \sum_{j=1}^d \sum_{x: \{x,x+e_j\} \subset
\Sigma_{\ell_n}}  \int \widehat{c}_{x,x+e_j}(\xi)\,
\Big[\, \sqrt{ g(T_{x,x+e_j} \xi) } - \sqrt{ \vphantom{T_{x,y} } g(\xi) }
\, \Big]^2
\, d \nu^{\rm c}_{n, \mss M}  \,.
\end{gather*}

The next result is \cite[Theorem 2.3]{cmr02}. As $\nu^n_{\rho(\cdot)}$
is a product measure (in the notation of \cite{cmr02} the potential
$\varPhi$ vanishes) condition ${\rm USMT} (1,m,1)$ holds.

\begin{theorem}
\label{s83}
There exists a positive constant $\mf c_{\rm LS}$ such that
\begin{equation*}
\int f\, \ln f\, \, d \nu^{\rm c}_{n, \mss M}   \,\le\,
\mf c_{\rm LS}\, n^2\, I_{n} (f \,;\, \nu^{\rm c}_{n, \mss M} )
\end{equation*}
for all densities $f$ with respect to $\nu^{\rm c}_{n, \mss M}$, $0\le
\mss M\le |\Sigma_{n}|$, $n\ge 1$.
\end{theorem}

Taking $f= 1+ \epsilon \, g$ for a zero-mean function $g$, and
expanding in $\epsilon$ yields the spectral gap stated below.

\begin{corollary}
\label{s84}
There exists a positive constant $\mf c_{\rm SG}$ (bounded by
$(1/2) \mf c_{\rm LS}$) such that
\begin{equation*}
\int g^2 \, \, d \nu^{\rm c}_{n, \mss M}   \,\le\,
\mf c_{\rm SG}\, n^2\, I_{n} (g \,;\, \nu^{\rm c}_{n, \mss M} )
\end{equation*}
for all functions $g$ such that
$\int g \, d\nu^{\rm c}_{n, \mss M} =0$, $0\le \mss M\le |\Sigma_{n}|$,
$n\ge 1$.
\end{corollary}

\section{Decomposition of cylinder functions}
\label{sec8}

Throughout this section, $\rho\colon \bb T^d \to (0,1)$ is a fixed
smooth density profile. Let
\begin{align}
\label{09}
{\color{blue} \xi_\varnothing } \,=\, 1\,, \quad 
{\color{blue} \xi_D}  \,=\, \prod_{x\in D} [\eta_x-\rho(x/n)]\,, \quad
\cb{w_{y} } \, := \, w_{\rho (\cdot) , y}  \, := \,
\frac{\eta_y - \rho(y/n)}{ \rho (y/n) [1-\rho(y/n)]} \,,
\end{align}
$D$ a
finite subset of $\bb T^d_n$.

We start with an elementary identity. Let $h$ be a cylinder function
and denote by $\cb{\mtt s(h)}$ its support. Thus, $h$ depends on the
variables $\{\eta_x : x\in \mtt s(h)\}$ and only on these variables.
Recall the definition of the polynomial $\mtt p_h$ introduced in
\eqref{05}. Assume that the density profile $\rho$ is constant equal
to $\alpha$.  An elementary computation yields that
\begin{equation}
\label{148}
\frac{d}{d\alpha} E_{\nu^n_\alpha} [ h ] \,=\,
\sum_{z\in \mtt s(h)} E_{\nu^n_\alpha} [ h\, w_z ]\,.
\end{equation}

Consider a cylinder function $f\colon \{0,1\}^{\bb Z^d} \to \bb R$.
Denote by $A\subset \bb Z^d$ its support, and assume that $n$ is large
enough so that $A\subset \bb T^d_n$.  Let $\Pi^1_{\rho(\cdot)}$,
$\Pi^2_{\rho(\cdot)}$ be the operators given by
\begin{equation}
\label{140}
\begin{gathered}
(\Pi^1_{\rho(\cdot)} f) (\eta) \; =\;
\sum_{z\in A} E_{\nu^n_{\rho(\cdot)} }  \big[\, f\, w_z\,]   \, \xi_z 
\;, \\
(\Pi^{2}_{\rho(\cdot)} f) (\eta) \; =\;
\sum_{D\subset A \,,\, |D|\ge 2}
E_{\nu^n_{\rho(\cdot)} }  \big[\, f\, w_D\,]   \, \xi_D \;,
\end{gathered}
\end{equation}
where $w_D = \prod_{x\in D} w_x$.  Note that we can take the first sum
over $\bb T^d_n$ because
$E_{\nu^n_{\rho(\cdot)} } \big[\, f\, w_z\,] =0$ whenever $z$ lies
outside the support of $f$. A similar observation applies to the
second sum.

As $A$ is the support of $f$, there exist constants $c_B$,
$B\subset A$, such that
\begin{equation}
\label{13}
f(\eta) \;=\; \sum_{B \subset A} c_B \, \eta_B\;,
\end{equation}
where $\eta_\varnothing = 1$, $\eta_B = \prod_{x\in B} \eta_x$ and the
sum is performed over all subsets $B$ of $A$.  With this representation,
\begin{equation}
\label{16}
(\Pi^1_{\rho(\cdot)}  f) (\eta) \; =\; 
\sum_{z\in A} \xi_z \,
\sum_{\{z\} \subset B \subset A} c_B \,   \rho_{B\setminus \{z\}}
\,, \quad \text{where}\;\;
\rho_D \,:=\,  \prod_{x\in D} \rho(x/n)\,, \;\; D\subset \bb T^d_n\,, 
\end{equation}
because, by \eqref{13},
\begin{equation*}
E_{\nu^n_{\rho(\cdot)} }  \big[\, f\, w_z\,]
\,=\, \sum_{B \subset A} c_B E_{\nu^n_{\rho(\cdot)} }
\big[\, \eta_B \, w_z\,]
\,=\,  \sum_{\{z\}\subset B \subset A} c_B \, \rho_{B\setminus \{z\}}\,.
\end{equation*}
Similarly, 
\begin{equation}
\label{15}
(\Pi^2_{\rho(\cdot)}  f) (\eta) \; =\;
\sum_{C\subset A, |C|\ge 2} \xi_C \,
\sum_{C\subset B \subset A} c_B \,   \rho_{B\setminus C} \,.
\end{equation}

The cylinder function $f$ is said to have {\it degree $k \ge 0$} in
$L^2(\nu^n_{\rho(\cdot )} ) $ if there exists a finite collection of
subsets $D$ of $\bb Z^d$ with cardinality $k$ and real numbers $c'_D$
such that
\begin{equation*}
f(\eta) \;=\; \sum_{D } c'_D\, \xi_D \;.
\end{equation*}
Note that $\Pi^1_{\rho(\cdot)} f $ corresponds to the terms of degree
$1$ of $f - E_{\nu^n_{\rho(\cdot)} } [f] $ and $\Pi^2_{\rho(\cdot)} f$
to the ones of degree greater than or equal to $2$.

\begin{lemma}
\label{l42}
For each cylinder function $f$, 
\begin{equation*}
f(\eta) \,-\, E_{\nu^n_{\rho(\cdot)} } [f]  \,=\,
(\Pi^1_{\rho(\cdot)}  f)  (\eta)
\;+\; (\Pi^2_{\rho(\cdot)} f) (\eta) \;.
\end{equation*}
\end{lemma}

\begin{proof}
By \eqref{13}, 
\begin{equation*}
f(\eta) \,-\, E_{\nu^n_{\rho(\cdot)} } [f] \,=\,
\sum_{B \subset A} c_B \, [\, \eta_B - \rho_B]\,.
\end{equation*}
Rewrite the previous sum as
\begin{equation*}
\sum_{B \subset A} c_B \, \Big \{\, \prod_{z\in B} [\, \xi_z +
\rho(z/n) \,] - \rho_B\, \Big\} \,=\,
\sum_{B \subset A} c_B \,  \sum_{\varnothing \subsetneq C\subset B} \xi_C \,
\rho_{B\setminus C}
\,=\, \sum_{\varnothing \subsetneq C\subset A} \xi_C \,
\sum_{C\subset B \subset A} c_B \,   \rho_{B\setminus C} \,.
\end{equation*}
Separating the singletons from the rest yields that
\begin{equation*}
f(\eta) \,-\, E_{\nu^n_{\rho(\cdot)} } [f]  \,=\,
\sum_{z\in A} \xi_z \,
\sum_{\{z\} \subset B \subset A} c_B \,   \rho_{B\setminus \{z\}}
\,+\,
\sum_{C\subset A, |C|\ge 2} \xi_C \,
\sum_{C\subset B \subset A} c_B \,   \rho_{B\setminus C}
\end{equation*}
To complete the proof of the lemma, it remains to recall the
identities \eqref{16}, \eqref{15}.
\end{proof}

\begin{remark}
\label{rm2}
The above computation shows that the term
$E_{\nu^n_{\rho(\cdot)} } [f]$ removes the constants from $f(\eta)$,
and that the expansion
\begin{equation*}
f \,-\,  E_{\nu^n_{\rho(\cdot)} } [f] \,-\,  (\Pi^1_{\rho(\cdot)} f) ( \eta)
\end{equation*}
only contains terms of degree $2$ or larger.
\end{remark}

\subsection*{The diffusion coefficient}

We prove in this subsection an identity for the cylinder functions
$h_j$ used in the semi-martingale decomposition of the density field.

\begin{lemma}
\label{s34}
For all $0<\rho<1$, $1\le j\le d$, $x\neq 0$,
$E_{\nu_{\rho} } \, [\, \tau_{x} h_{j} \, [\, \eta_0 -\rho\,] \, ]
\,=\,0$, and
\begin{equation*}
E_{\nu_{\rho} } \, \big[\, h_{j}
\, [\, \eta_0 -\rho\,] \, \big] \,=\, \chi(\rho) \,
E_{\nu_{\rho} } \, \big[\, c_{0,e_j}(\eta) \, \big] \,.
\end{equation*}
\end{lemma}

\begin{proof}
Fix $1\le j\le d$.  The proof is based on the following
observation. Let $\mtt j_{y,y+\bs e_k}$, $y\in \bb Z^d$,
$1\le k\le d$, be the instantaneous current over the bond
$(y,y+\bs e_k)$:
$\cb {\mtt j_{y,y+\bs e_k}} := c_{y,y+\bs e_k} (\eta) \, [\, \eta_y -
\eta_{y+\bs e_k} \,] $.  As $c_{0,e_k} (\cdot)$ does not depend on
$\eta_0$, $\eta_{e_k}$,
\begin{equation}
\label{137}
E_{\nu_{\rho} } \, \Big[\, \mtt j_{x, x+ e_k} (\eta) 
\,  [\, \eta_0 -\rho\,] \, \Big] \,=\, 0 \quad \text{if $x$ and
$x+e_k\neq 0$}\,.
\end{equation}

By \eqref{11b}, 
\begin{equation*}
\begin{aligned}
& \tau_{z(x,p)}  h_{j} \, =\,  \tau_{z(x,0)}  h_j
\,-\, \sum_{r=0}^{p- 1} \mtt j _{z(x,r),
z(x,r+1) }\,,\quad \text{for}\;\; p\ge 1
\\
& \tau_{z(x,p)}  h_{j} \, =\, \tau_{z(x,0)}  h_j
\,+\, \sum_{r=p}^{- 1} \mtt j _{z(x,r),
z(x,r+1) }\,,\quad \text{for}\;\; p\le -\, 1\,,
\end{aligned}
\end{equation*}
where
$\cb{z(x,r)} : = \sum_{i\neq j} x_i\, \bs e_i\, +\, r\, \bs e_j$.
Therefore, if $p\ge 1$, 
\begin{equation}
\label{139}
E_{\nu_{\rho} } \, \Big[\, \tau_{z(x,p)}   h_{j}
\, [\, \eta_0 -\rho\,] \, \Big] \,=\,
E_{\nu_{\rho} } \, \Big[\, \tau_{z(x,0)}  h_{j}
\, [\, \eta_0 -\rho\,] \, \Big] \,-\,
\sum_{r=0}^{p- 1}
E_{\nu_{\rho} } \, \Big[\, 
\, [\, \eta_0 -\rho\,] \, \mtt j _{z(x,r), z(x,r+1)}
\,\Big]\,.
\end{equation}
A similar formula holds for  $p\le -1$.

Fix $k\neq j$. In this case, if $x_k\neq 0$, then $z(x,r) \neq 0$ for
all $r\in \bb Z$. Hence, by \eqref{137}, the second term on the
right-hand side of the previous identity vanishes. Thus, for all
$x\in \bb Z^d$ such that $x_k\neq 0$,
\begin{equation*}
E_{\nu_{\rho} } \, [\, \tau_{z(x,p)}   h_{j} \, [\, \eta_0 -\rho\,] \, ]
\,=\,
E_{\nu_{\rho} } \, [\, \tau_{z(x,0)}    h_{j} \, [\, \eta_0 -\rho\,] \, ] \,.
\end{equation*}
Consider $p$ sufficiently large for the support of $\tau_{z(x,p)}  h_{j}$
not to contain $0$. In this case, the left-hand side of the previous
identity vanishes. This implies that the right-hand side vanishes,
which in turn yields that the left-hand side vanishes for all $p$.
Therefore,
\begin{equation}
\label{138}
E_{\nu_{\rho} } \, [\, \tau_x h_{j} \, [\, \eta_0 -\rho\,] \, ]
\,=\, 0\quad\text{if}\;\; x_k \neq 0 \quad \text{for some}\;\;
k\neq j \,.
\end{equation}

It remains to consider the case $x= p \bs e_j$ for some $p\in\bb Z$.
By \eqref{139} and \eqref{137},
\begin{equation*}
E_{\nu_{\rho} } \, \Big[\, \tau_{p\bs e_j}  h_{j}
\, [\, \eta_0 -\rho\,] \, \Big] \,=\,
E_{\nu_{\rho} } \, \Big[\, h_{j}
\, [\, \eta_0 -\rho\,] \, \Big] \, - \, \chi(\rho) \,
E_{\nu_{\rho} } \, \Big[\, c_{0,e_j}(\eta) \, \Big] 
\end{equation*}
for $p\neq 0$. For $p$ large the left-hand side vanishes. This proves
that the right-hand side vanishes, and, therefore, that the left-hand
side vanishes for all $p\neq 0$.
\end{proof}

\begin{corollary}
\label{s21}
For every $1\le j, k\le d$, $0<\rho<1$, 
\begin{equation*}
E_{\nu_{\rho} } \, \Big[\, h_{j} \sum_{x\in \bb Z^d} x_k\, 
\,  [\, \eta_x -\rho\,] \, \Big] \,=\, 0\,.
\end{equation*}
As $h_j$ is a cylinder function, we may replace the sum over $\bb Z^d$
by a sum over a finite set $A$ which contains the support of $h_j$.
\end{corollary}

\begin{proof}
Fix $1\le j\le d$.  
As $\nu_\rho$ is translation invariant, the expression appearing in
the statement of the corollary is equal to
\begin{equation*}
-\, \sum_{x\in A} x_k\, E_{\nu_{\rho} } \, \Big[\, \tau_x h_{j}
\, [\, \eta_0 -\rho\,] \, \Big] \,.
\end{equation*}
The assertion thus follows from the previous lemma.
\end{proof}

\begin{corollary}
\label{s24}
Let $\rho\colon \bb T^d \to (0,1)$ be a density profile in
$C^2(\bb T^d)$ and such that $\mf r \le \rho(x) \le 1-\mf r$ for some
$0<\mf r < 1/2$. Then, there exists a finite constant $C_0$, depending
only on $\mf r$, $\Vert \nabla \rho \Vert_\infty$, $\Vert \nabla^2 \rho
\Vert_\infty$, such that
\begin{equation}
\label{136}
\Big|\, E_{\nu^n_{\rho( \cdot)}} \, [\, h_{j} \,]
\, -\, E_{\nu_{\rho(0)}} \, [\, h_{j} \,]\, \Big|\, \le \, C_0 \,
\frac{1}{n^2} \,\cdot
\end{equation}
\end{corollary}

\begin{proof}
We may rewrite the first expectation as
\begin{equation*}
E_{\nu^n_{\rho( \cdot)}} \, [\, h_{j} \,] \,=\, 
E_{\nu_{\rho} } \, \Big[\, h_{j}  \,
\exp\Big\{ {\sum_{x\in A} \eta_x \, \ln \frac{\rho(x/n)}{\rho} \,+\,
[1-\eta_x]\, } \ln \frac{1-\rho(x/n)}{1-\rho} \,\Big\}\, \Big]\,,
\end{equation*}
where $A$ is a cube centered at the origin, sufficiently large to
contain the support of $h_j$ and $\rho = \rho(0)$.  As the exponential
has mean $1$ with respect to $\nu_{\rho} $, we may write the
difference \eqref{136} as
\begin{equation*}
E_{\nu_{\rho} } \, \Big[\, \{ \, h_{j}
-  E_{\nu_{\rho} }  [h_{j} ]  \} \,
\exp\Big\{ {\sum_{x\in A} \eta_x \, \ln \frac{\rho(x/n)}{\rho} \,+\,
[1-\eta_x]\, } \ln \frac{1-\rho(x/n)}{1-\rho} \,\Big\}\, \Big]\,.
\end{equation*}
Since $\rho(\cdot)$ is Lipschitz-continuous, the expression inside
braces is of order $1/n$. This expectation is thus equal to
\begin{equation*}
E_{\nu_{\rho} } \, \Big[\, \{ \, h_{j}
-  E_{\nu_{\rho} }  [h_{j} ]  \} \,
\sum_{x\in A} [\eta_x -\rho] \, \ln \frac{\rho(x/n)}{\rho} \,-\,
[\eta_x - \rho]\,  \ln \frac{1-\rho(x/n)}{1-\rho} \, \Big]
\,+\, \mf R_n \,,
\end{equation*}
where $\mf R_n$ is a remainder, whose value may change from line to
line and such that $|\mf R_n|\le C_0 n^{-2}$ for some finite constant
$C_0$ depending only on $\mf r$, $\Vert \nabla \rho \Vert_\infty$,
$\Vert \nabla^2 \rho \Vert_\infty$. We subtracted the expectation of
$\eta_x$ because $h_j - E_{\nu_{\rho} } [h_{j} ]$ has zero-mean. We
may now remove the expectation of $h_j$ and expand $\rho(\cdot)$ to
write the previous equation as
\begin{equation*}
\frac{1}{n}\, \sum_{k=1}^d (\partial_{x_k} \rho)(0)
\sum_{x\in A} E_{\nu_{\rho} } \, \Big[\, h_{j}
\, x_k\, \frac{\eta_x -\rho}{\chi(\rho)} \, \, \Big]
\,+\, \mf R_n \,.
\end{equation*}
By Corollary \ref{s21}, the first term vanishes, which completes the
proof of the corollary.
\end{proof}

\section{Martingales}

We introduce in this section some martingales used throughout the
article. Fix a test function $F\colon \bb T^d \to \bb R$ and recall
from \eqref{30} the semi-martingale decomposition of the process
$t\mapsto X^n_t(F)$:
\begin{equation}
\label{30b}
X^n_t(F)  \, =\;  X^n_0(F) \; + \; \int_0^t
(\partial_s + L_n)  X^n_s(F)\, ds\, +\, M^n_t(F)  \;.
\end{equation}
By \cite[Lemma A1.5.1]{kl}, the process $t \mapsto M^n_t(F)$ is a
martingale, and
\begin{equation}
\label{144}
M^n_t(F)^2 \;-\; \int_0^t \Gamma^n (\eta^n(s), F) \; ds
\end{equation}
is also a martingale, where $\Gamma^n(\eta, F)$ is twice the carr\'e du
champ of $F$:
\begin{equation}
\label{143}
\cb{\Gamma^n(\eta, F)} \; :=\; L_n X^n(F)^2 \,-\, 2 \, X^n(F) \, L_n X^n(F)
\;.
\end{equation}
An elementary computation yields that
\begin{equation*}
\Gamma^n(\eta, F) \;=\;
\frac{1}{n^d} \sum_{j=1}^d \sum_{x\in \bb T^d_n}
c_{0,e_j}  (\tau_x \eta) \, [\eta_{x+e_j} - \eta_x]^2\,
[(\nabla^n_j F)(x/n)]^2  \;,
\end{equation*}
where $\nabla^n_j F$ stands for the discrete partial derivative in the
$j$-th direction given by
\begin{equation}
\label{163}
{\color{blue} (\nabla^n_j F)(x/n)} := n[F((x+e_j)/n) - F(x/n)] \,.
\end{equation}
Thus, $|\,\Gamma^n(\eta, F)\,|$ is bounded by
$C_0 \Vert \nabla F \Vert^2_\infty$ for some finite constant $C_0$,
and
\begin{equation}
\label{134}
\bb E^n_{\eta} \big[\, M^n_t(F)^2\,\big] \;\le\; C_0\, t\, 
\Vert \nabla F \Vert^2_\infty
\end{equation}
for all $t>0$, $\eta\in \Omega_n$, $n\ge 1$.

We may also consider time-dependent functions.  Fix a function
$H \in C^{1,2} ([0,T] \times \bb T^d)$, and denote by
$\mss M^n_t(H_\cdot)$ the martingale defined by
\begin{equation*}
\cb{\mss M^{n}_t(H_\cdot)}  \;:=\; X^n_t(H_{t}) \;-\; X^n_0(H_{0}) \;-\; \int_0^t
(\partial_s + L_n) X^n_s(H_{s})\, ds\;.
\end{equation*}
In this formula, the time-derivative $\partial_s$ in
$\partial_s X^n_s(H_{s})$ acts on the density field $X^n_s$ and on the
test function $H_{s}$.

Fix up to the end of this section $0<t\le T$. Recall from \eqref{159b}
that we denote by $(P_{s,t}: 0\le s \le t)$ the semigroup associated
to the linear operators $(\mc A_r: r\ge 0)$ introduced in \eqref{41}.
Thus, for each $G\in C^\infty (\bb T^d)$,
$\mf u(s,x) = (P_{s,t} G)(x)$, $0\le s\le t$, is the solution of the
linear equation
\begin{equation}
\label{159}
\left\{
\begin{aligned}
& \partial_s \mf u \,+\, \mc A_s \mf u\, =\, 0 \,,
\quad 0\le s\le t\,, 
\\
& \mf u (t) \,=\, G\,.
\end{aligned}
\right.
\end{equation}
Fix $F\in H^{\beta+2}(\bb T^d)$ for some non-integer real number
$\beta>0$.  Denote by $\mf v \colon [0,t] \times \bb T^d \to \bb R$
the solution of equation \eqref{159b}.  Let
$\mf u \colon [0,t] \times \bb T^d \to \bb R$ be given by
$\mf u(s,x) = \mf v(t-s,x)$. Clearly, $\mf u$ is the solution to
\eqref{159} with $G=F$.

Fix $G$ in $C^\infty(\bb T^d)$, and denote by $M^{n,t}_s (G)$,
$0\le s\le t$, the martingale defined by
\begin{equation}
\label{51b}
\cb{M^{n,t}_s(G)}  \;:=\; \mss M^n_s (P_{\cdot, t} G)\,=\, 
X^n_s(P_{s,t}G) \;-\; X^n_0(P_{0,t}G) \;-\; \int_0^s
(\partial_r + L_n) X^n_r(P_{r,t}G)\, dr\;.
\end{equation}
As above, the time-derivative $\partial_r$ in
$\partial_r X^n_r(P_{r,t} G)$ acts on the density field $X^n_r$ and on
the test function $P_{r,t} G$. Taking $s=t$ in \eqref{51b} yields a
representation for $X^n_t(G)$:
\begin{equation*}
X^n_t(G) \;=\; X^n_t(P_{t,t}G) \;=\; X^n_0(P_{0,t}G)
\;+\; M^{n,t}_t(G)
\;+\; \int_0^t (\partial_r + L_n) \, X^n_r(P_{r,t} G)\; dr
\;.
\end{equation*}

\begin{lemma}
\label{s56}
Assume that $d\le 3$ and that
$H_n (\mu_n \,|\, \nu^n_{u_0(\cdot)}) = o(n^{d/2})$.
Fix $0\le t\le T$, $G$ in $C^\infty(\bb T^d)$. Then,
$X^n_t(G) - X^n_0(P_{0,t}G) -M^{n,t}_t(G)$ converges to $0$ in
$L^1(\bb P_{\mu_n})$.
\end{lemma}

\begin{proof}
Recall formula \eqref{145} for $(\partial_r + L_n) \, X^n_r(F)$. Here
as the test function is time dependent, we have an extra term. More
precisely, as $\partial_r P_{r,t} G\,=\, -\, \mc A_r P_{r,t} G$,
\begin{equation*}
(\partial_r + L_n) X^n_r(P_{r,t}G ) \;=\;
R^n_r(P_{r,t} G) \;+\; B^n_r(P_{r,t} G) \, +\, \mss D^{(n)}_r (P_{r,t} G)\;, 
\end{equation*}
where $\mss D^{(n)}_r (\cdot)$ is the deterministic remainder
introduced in \eqref{145} and satisfying the bound
\begin{equation}
\label{189}
|\, \mss D^{(n)}_r (P_{r,t} G) \,| \,\le\,
\mf a_{3,4} (\{h_j : 1\le j\le d\}) \, \Vert
P_{r,t} G\Vert_{C^2} \, n^{(d/2)-2}\,.
\end{equation}
In the penultimate displayed equation, for a function
$H$ in $C^4(\bb T^d)$, $R_r^n(H)$, $B_r^n(H)$ are given by
\begin{equation}
\label{147}
\begin{gathered}
\cb{R^n_r(H)} \; :=\; \frac{1}{n^{d/2}} \sum_{j=1}^d
\sum_{x\in\bb T_n^d}
\big\{ \, (\Delta^n_j H)  -
(\partial^2_{x_j}  H ) \,\big\}(x/n)\,
D_j (u(r,x/n))  \,[\, \eta_x^n(r) - u(r,x/n) \,]  \;,
\\
\cb{B^n_r(H)} \; :=\; \frac{1}{n^{d/2}} \sum_{j=1}^d
\sum_{x\in\bb T_n^d}
(\Delta^n_j H)  (x/n)\, 
\Pi_{u(r, \cdot), h_j}  (x, \eta^n(r))  \,.
\end{gathered}
\end{equation}
By a fourth order Taylor expansion,
$| (\Delta^n_j H) - (\partial^2_{x_j} H )| \le C_0 \, \Vert
H\Vert_{C^4}\, n^{-2}$ for some finite constant $C_0$. Thus,
$R^n_r(P_{r,t} G)$ fulfills the same bound \eqref{189} with the $C^2$
norm replaced by the $C^4$ norm.  Therefore,
\begin{equation}
\label{146}
X^n_t(G) \;=\;  X^n_0(P_{0,t}G)
\;+\; M^{n,t}_t(G) \;+\; \int_0^t
B^n_r(P_{r,t} G)  \; dr \, +\, \int_0^t \mf R^n_r (P_{r,t} G) \, dr \;,
\end{equation}
where the remainder
$\cb{\mf R^n_r (H)} := \mss D^{(n)}_r (H) + R^n_r(H) $ satisfies
\eqref{189} with the $C^2$ norm replaced by the $C^4$ norm.  To
complete the proof, it remains to recall that, by Theorem \ref{t02},
$\int_0^t B^n_r(P_{r,t} G) \; dr$ converges to $0$ in
$L^1(\bb P_{\mu_n})$.
\end{proof}

We turn to a formula for the martingale $s\mapsto M^{n,t}_s (F)$.
Recall from \eqref{30b} the definition of the martingale $M^n_t(F)$,
$F\in C^\infty(\bb T^d)$. Let
$H\colon \bb R_+ \times \bb T^d \to \bb R$ be a smooth function, and
denote by $H_t \colon \bb T^d_n \to \bb R$ the function given by
$H_t(x) = H(t,x/n)$, and by $\cb{M^n_t(H_s)}$ the value of the
martingale $M^n_t(F)$ for $F=H_s$.  Similarly, and in contrast with
equation \eqref{51b}, $[(\partial_s + L_n ) X^n_s] (H_s)$ stands for
the expression $(\partial_s + L_n ) X^n_s(F)$ computed at $F=H_s$. To
distinguish this situation from the case \eqref{51b} in which the
partial derivative acts on the test function, we added a square
bracket.  The proof of \cite[Lemma 8.3]{jl23} yields that
\begin{equation*}
M^n_t(H_t) \;=\; X^n_t(H_t) \,-\, X^n_0(H_0) \;-\;
\int_0^t \big\{ \, \big[ (\partial_s + L_n )  X^n_s\big] (H_s)
\,+\, X^n_s(\partial_s H_s)\,\big\}\; ds
\;+\; \int_0^t M^n_s(\partial_s H_s)\; ds\;.
\end{equation*}
This expression can be rewritten as
\begin{equation}
\label{185}
M^n_t(H_t) \;=\; X^n_t(H_t) \,-\, X^n_0(H_0) \;-\;
\int_0^t  (\partial_s + L_n )  X^n_s (H_s) \; ds
\;+\; \int_0^t M^n_s(\partial_s H_s)\; ds\;,
\end{equation}
where now, as in \eqref{51b}, the partial derivative $\partial_s$
in $(\partial_s + L_n ) X^n_s (H_s)$ acts both on the density field
and on the test function.

The previous identity provides a representation of the martingale
$s\mapsto M^{n,t}_s (F)$, introduced in \eqref{51b}, in terms of the
martingales $s\mapsto M^n_s (G)$. To derive it, fix $0<t\le T$, and a
function $F$ in $C^\infty(\bb T^d)$. Let $H(s,x) = (P_{s,t} F)(x)$,
$0\le s\le t$. In equation \eqref{51b} with $s=t$, replace the
right-hand side by the first three terms on the right-hand side of
\eqref{185}.  This yields that
\begin{equation}
\label{56}
M^{n,t}_t(F) \;=\; M^n_t(F)  \,-\, 
\int_0^t M^n_s(\partial_s P_{s,t} F )\; ds
\,=\, 
M^n_t(F) \;+\;  \int_0^t M^n_s(\mc A_s P_{s,t}  F)\; ds 
\end{equation}
because $\partial_s P_{s,t} F = - \mc A_s \, P_{s,t} F$. Mind that, as
the operators $\mc A_s$ depend on $s$ through $u(s,\cdot)$, they do
not commute with the semigroup, so that $\mc A_s \, P_{s,t}$ cannot be
replaced by $P_{s,t}\, \mc A_s$.

We conclude this section with a formula for the time integral of
$M^{n,t}_t(F)$.  We claim that for all $0\le s<t \le T$, and
$G\in C^\infty(\bb T^d)$,
\begin{equation}
\label{183}
\begin{aligned}
\int_{s}^{t} M^{n,r}_r (G) \, dr \, & =\,
\int_{s}^{t}  M^{n}_r (G) \;dr
\,+\, \int_{0}^{s} 
M^{n}_{u} \Big(\, \int_s^t  \mc A_u \, P_{u,r} \, G \, dr \,\Big)
\, du
\\
\, & +\,
\int_{s}^{t}
M^{n}_{u} \Big(\, \int_u^t  \mc A_u \, P_{u,r} \, G\, dr \Big)
\, du \,.
\end{aligned}
\end{equation}
By \eqref{56}, the left-hand side is equal to
\begin{equation*}
\int_{s}^{t}  M^{n}_r (G) \;dr  \,+\,
\int_{s}^{t}  dr \, \int_0^r M^{n}_{u} (\mc A_u \, P_{u,r} \, G) \, du\,   \,.
\end{equation*}
Consider the second term of this sum.  A change of the order of
integration yields that it is equal to
\begin{align*}
\int_{0}^{s}  du \, \int_s^t M^{n}_{u} (\mc A_u \, P_{u,r} \, G) \, dr
\,+\,
\int_{s}^{t}  du \, \int_u^t  M^{n}_{u} (\mc A_u \, P_{u,r} \, G) \, dr\,.
\end{align*}
The linear dependence of the martingale $M^n(F)$ on its argument $F$
completes the proof of the claim \eqref{183}.

In the case where the operator $\mc A$ is time-independent, we may
simplify further the right-hand side of \eqref{183} noting that
$\mc A  P_{u,r} \, G = P_{r-u} \mc A\, G = (d/dr) P_{r-u} \,G$.

\subsection*{Acknowledgments}

C. L. has been partially supported by FAPERJ CNE E-26/201.117/2021, by
CNPq, Brazil Bolsa de Produtividade em Pesquisa PQ 305779/2022-2, and
by the French National Research Agency (ANR) via Project
ANR-25-CE40-6875-01 (ANR DySLoS).  S. S. was partially supported by a
Simons sabbatical grant.

The A.I. Claude was used to improve the some of the English grammar
and to find typos.

\end{document}